\documentclass[a4paper,12pt,amsart,frenchb]{article}

\usepackage{amsmath,amsbsy,amsfonts,amssymb}
\usepackage[french]{babel}
\newcommand{\ass}[2]{\vskip0.3cm\noindent
{\bf {#1}}. { \sl {#2}}\vskip0.3cm\noindent
}

\begin{document}

\title{  Calcul d'une valeur d'un facteur $\epsilon$  par une formule int\'egrale}
\author{J.-L. Waldspurger}
\date{18 avril 2012}
\maketitle

\bigskip

{\bf Introduction}

\bigskip

Soient $F$ un corps local non archim\'edien de caract\'eristique nulle, $V$ un espace vectoriel de dimension finie sur $F$, muni d'une forme quadratique non d\'eg\'en\'er\'ee $q$, $D_{0}$ une droite de $V$ qui n'est pas isotrope pour $q$, $W$ l'orthogonal de $D_{0}$ dans $V$. Notons $G$ le groupe sp\'ecial orthogonal de $V$ et $H$ celui de $W$, que l'on identifie \`a un sous-groupe de $G$.  Soient $\pi$, resp. $\rho$, une repr\'esentation admissible irr\'eductible de $G(F)$, resp. $H(F)$, que l'on r\'ealise dans un espace $E_{\pi}$, resp. $E_{\rho}$. Notons $Hom_{H(F)}(\pi,\rho)$ l'espace des applications lin\'eaires $\varphi:E_{\pi}\to E_{\rho}$ telles que $\varphi\circ\pi(h)=\rho(h)\circ\varphi$ pour tout $h\in H(F)$. Notons $m(\rho,\pi)$ la dimension de cet espace. D'apr\`es un th\'eor\`eme de Aizenbud, Gourevitch, Rallis et Schiffmann ([AGRS]), on a $m(\rho,\pi)\leq 1$. Supposons $\pi$ et $\rho$ temp\'er\'ees. Dans les articles [W2] et [W3], on a \'etabli une formule qui calcule $m(\rho,\pi)$ comme somme d'int\'egrales sur des sous-tores non n\'ecessairement maximaux de $H$ de fonctions qui se d\'eduisent des caract\`eres de $\pi$ et $\rho$. D'apr\`es les travaux encore en cours d'Arthur, la th\'eorie de l'endoscopie tordue relie les repr\'esentations $\pi$ et $\rho$, ou plus exactement les $L$-paquets qui contiennent ces repr\'esentations, \`a des repr\'esentations autoduales de groupes lin\'eaires. Indiquons plus pr\'ecis\'ement de quel groupe lin\'eaire il s'agit, par exemple pour la repr\'esentation $\pi$. Notons $d$ la dimension de $V$. Si $d$ est pair, le groupe est $GL_{d}$. Si $d$ est impair, $G$ appara\^{\i}t usuellement comme groupe endoscopique de $GL_{d-1}$ tordu. Mais $G$ est aussi un groupe endoscopique de $GL_{d}$ tordu et, pour ce que nous faisons, il semble que cette deuxi\`eme interpr\'etation soit plus pertinente. D'apr\`es la conjecture locale de Gross-Prasad, le nombre $m(\rho,\pi)$ doit \^etre reli\'e \`a un  facteur $\epsilon$ de la paire de repr\'esentations de groupes lin\'eaires correspondant \`a la paire $(\rho,\pi)$. Cela sugg\`ere l'existence d'une formule int\'egrale, parall\`ele \`a celle \'evoqu\'ee ci-dessus, qui calcule ce facteur $\epsilon$ de paire. Inversement, une telle formule devrait permettre, via la th\'eorie de l'endoscopie tordue, de prouver la conjecture locale de Gross-Prasad pour les repr\'esentations temp\'er\'ees. Le but de l'article est d'\'etablir cette formule int\'egrale.

Oublions maintenant les objets introduits ci-dessus, qui n'ont servi que de motivation. On conserve toutefois le corps $F$. Soient $r$ et $m$ deux entiers  positifs ou nuls. Posons $d=m+1+2r$, $G=GL_{d}$, $H=GL_{m}$.  Soit $\pi$ une repr\'esentation admissible irr\'eductible et temp\'er\'ee de $G(F)$. On suppose $\pi$ autoduale, c'est-\`a-dire qu'elle est isomorphe \`a sa contragr\'ediente $\pi^{\vee}$.  Soit $\rho$ une repr\'esentation de $H(F)$ v\'erifiant des conditions similaires. Notons $\theta_{d}$ l'automorphisme de $G$ d\'efini par $\theta_{d}(g)=J_{d}{^tg}^{-1}J_{d}^{-1}$, o\`u $J_{d}$ est la matrice antidiagonale de coefficients $(J_{d})_{i,d+1-i}=(-1)^{i}$. Introduisons le groupe non connexe $G\rtimes \{1,\theta_{d}\}$ et sa composante connexe $\tilde{G}=G\theta_{d}$. Puisque $\pi$ est autoduale, elle se prolonge en une repr\'esentation du groupe non connexe $G(F)\rtimes \{1,\theta_{d}\}$. Elle se prolonge m\^eme de deux fa\c{c}ons. Fixons un caract\`ere $\psi:F\to {\mathbb C}^{\times}$ continu et non trivial.   La th\'eorie des mod\`eles de Whittaker permet de choisir l'un des prolongements. On note $\tilde{\pi}$ la restriction de ce prolongement \`a $\tilde{G}(F)$. On effectue des constructions analogues pour $H$ et $\rho$. Selon Jacquet, Piatetskii-Shapiro et Shalika, on d\'efinit le facteur $\epsilon(s,\pi\times \rho,\psi)$ pour $s\in {\mathbb C}$. Notons $\omega_{\pi}$ et $\omega_{\rho}$ les caract\`eres centraux de $\pi$ et $\rho$. Soit enfin $\nu$ un \'el\'ement de $F^{\times}$. On pose
$$\epsilon_{\nu}(\tilde{\pi},\tilde{\rho})=\omega_{\pi}((-1)^{[m/2]}2\nu)\omega_{\rho}((-1)^{1+[d/2]}2\nu)\epsilon(1/2,\pi\times \rho,\psi).$$
C'est ce terme que nous allons calculer par une formule int\'egrale.

{\bf Remarque.} Pour \'eviter un pi\`ege, signalons que l'\'equation fonctionnelle locale n'entra\^{\i}ne pas que ce terme est un signe $\pm 1$, mais seulement que c'est une racine quatri\`eme de l'unit\'e.

  Pour manier plus ais\'ement les objets $\tilde{G}$ et $\tilde{H}$, on les interpr\`ete comme des groupes tordus. Soient $V$ un espace vectoriel de dimension $d$ sur $F$, $W$ un sous-espace de dimension $m$ et $Z$ un sous-espace de $V$ suppl\'ementaire de $W$. Par le choix d'une base de $V$, $G$ s'identifie au groupe $GL(V)$ des automorphismes lin\'eaires de $V$. Notons $V^*$ l'espace dual de $V$. Alors $\tilde{G}$ s'identifie \`a l'espace $Isom(V,V^*)$ des isomorphismes lin\'eaires de $V$ sur $V^*$ ou, si l'on pr\'ef\`ere, \`a l'espace des formes bilin\'eaires non d\'eg\'en\'er\'ees sur $V$. Le groupe $G$ agit \`a droite et \`a gauche sur $\tilde{G}$ par
$$(g,\tilde{x},g')\mapsto {^tg}^{-1}\circ\tilde{x}\circ g'$$
pour $g,g'\in G$ et $\tilde{x}\in \tilde{G}$. On note simplement $(g,\tilde{x},g')\mapsto g\tilde{x}g'$ ces actions.  Le point base $\theta_{d}$ s'identifie \`a une forme symplectique si $d$ est pair, \`a une forme quadratique si $d$ est impair. On renvoie \`a 2.1 pour plus de pr\'ecision. De m\^eme, $\tilde{H}$ s'identifie \`a $Isom(W,W^*)$. Fixons une forme quadratique non d\'eg\'en\'er\'ee $\tilde{\zeta}$ sur $Z$. On suppose qu'elle est somme orthogonale d'une forme hyperbolique et de la forme $x\mapsto 2\nu x^2$ de dimension $1$. On interpr\`ete $\tilde{\zeta}$ comme un \'el\'ement de $Isom(Z,Z^*)$. On plonge alors $\tilde{H}$ dans $\tilde{G}$: un \'el\'ement $\tilde{y}\in Isom(W,W^*)$ s'identifie \`a l'\'el\'ement de $Isom(V,V^*)$ qui envoie $W$ sur $W^*$, $Z$ sur $Z^*$, et dont les restrictions \`a $W$, resp. $Z$, co\"{\i}ncident avec $\tilde{y}$, resp. $\tilde{\zeta}$.

Tout \'el\'ement $\tilde{x}\in \tilde{G}$ d\'efinit un automorphisme $\theta_{\tilde{x}}$ de $G$ caract\'eris\'e par l'\'egalit\'e $\tilde{x}g=\theta_{\tilde{x}}(g)\tilde{x}$. On d\'efinit la notion de sous-tore maximal de $\tilde{G}$ de la fa\c{c}on suivante. Soient $T$ un sous-tore maximal de $G$ d\'efini sur $F$ et $B$ un sous-groupe de Borel de $G$, contenant $T$ mais pas forc\'ement d\'efini sur $F$. Notons $\tilde{T}$ le sous-ensemble des \'el\'ements $\tilde{x}\in \tilde{G}$ tels que $\theta_{\tilde{x}}$ conserve $T$ et $B$. C'est un espace principal homog\`ene pour  chacune des actions de $T$ \`a droite ou \`a gauche.   Pour  $\tilde{x}\in \tilde{T}$, la restriction \`a $T$ de $\theta_{\tilde{x}}$ ne d\'epend pas de $\tilde{x}$. On note cet automorphisme $\theta_{\tilde{T}}$. Nous dirons que $\tilde{T}$ est un sous-tore maximal de $\tilde{G}$ si $\tilde{T}(F)$ n'est pas vide. 

Consid\'erons une d\'ecomposition $W=W'\oplus W''$. Posons $H'=GL(W')$, $\tilde{H}'=Isom(W',W^{_{'}*})$. Soit $\tilde{T}'$ un sous-tore maximal de $\tilde{H}'$, auquel est associ\'e un sous-tore maximal $T'$ de $H'$, et soit $\tilde{\zeta}_{H,T}\in Isom(W'',W^{_{''}*})$ une forme quadratique. On impose les conditions suivantes:

- la dimension de $W'$ est paire;

- $\tilde{T}'$ est anisotrope, c'est-\`a-dire que le seul sous-tore d\'eploy\'e contenu dans le sous-ensemble des \'el\'ements de $T'$ fixes par $\theta_{\tilde{T}'}$ est \'egal \`a $\{1\}$;

- le groupe sp\'ecial orthogonal de la forme $\tilde{\zeta}_{H,T}$ sur $W''$ est quasi-d\'eploy\'e ainsi que celui de la forme $\tilde{\zeta}_{G,T}= \tilde{\zeta}_{H,T}\oplus \tilde{\zeta}$ sur $W''\oplus Z$.

On note $\tilde{T}$ l'ensemble des \'el\'ements $\tilde{y}\in\tilde{H}$ tels que $\tilde{y}(W')=W^{_{'}*}$, $\tilde{y}(W'')=W^{_{''}*}$, que la restriction de $\tilde{y}$ \`a $W'$ appartienne \`a $\tilde{T}'$ et que la restriction de $\tilde{y}$ \`a $W''$ co\"{\i}ncide avec $\tilde{\zeta}_{H,T}$. On note $\underline{\cal T}$ l'ensemble des sous-ensembles $\tilde{T}$ de $\tilde{H}$ obtenus de cette fa\c{c}on. Le groupe $H(F)$ agit par conjugaison sur $\tilde{H}$. Cette action conserve l'ensemble $\underline{\cal T}$. On fixe un ensemble de repr\'esentants ${\cal T}$ de l'ensemble des orbites.

Soit $\tilde{T}\in \underline{\cal T}$, reprenons pour cet \'el\'ement les objets d\'efinis ci-dessus, en notant simplement $T=T'$. L'action de $T(F)$ par conjugaison conserve $\tilde{T}(F)$. On note $\tilde{T}(F)_{/\theta}$ l'ensemble des orbites. C'est une vari\'et\'e analytique sur $F$ sur laquelle on d\'efinit une certaine mesure. On d\'efinit aussi sur cette vari\'et\'e deux fonctions $\Delta^{\tilde{H}}$ et $\Delta_{r}$, des valeurs absolues de certains d\'eterminants. On d\'efinit \'egalement un groupe de Weyl $W(H,\tilde{T})$. Tous ces termes sont \'el\'ementaires, on renvoie \`a 1.4 et 3.2 pour des d\'efinitions pr\'ecises. Ce qui est plus crucial est d'associer \`a $\tilde{\pi}$ et $\tilde{\rho}$ deux fonctions $c_{\tilde{\pi}}$ et $c_{\tilde{\rho}}$ sur $\tilde{T}(F)_{/\theta}$. Soit $\tilde{t}$ un \'el\'ement de $\tilde{T}(F)$ en position g\'en\'erale, notons $G_{\tilde{t}}$ la composante neutre du sous-groupe des points fixes par $\theta_{\tilde{t}}$ dans $G$. Ce groupe se d\'ecompose en $T_{\theta}\times SO(\tilde{\zeta}_{G,T})$, o\`u $T_{\theta}=T\cap G_{\tilde{t}}$ et $SO(\tilde{\zeta}_{G,T})$ est le groupe sp\'ecial orthogonal de la forme quadratique $\tilde{\zeta}_{G,T}$ introduite ci-dessus.  A $\tilde{\pi}$ est associ\'e un caract\`ere $\Theta_{\tilde{\pi}}$ qui est une fonction localement int\'egrable sur $\tilde{G}(F)$. Plus pr\'ecis\'ement, d'apr\`es un r\'esultat d'Harish-Chandra g\'en\'eralis\'e au cas non connexe par Clozel, ce caract\`ere admet au voisinage de $\tilde{t}$ un d\'eveloppement en combinaison lin\'eaire de transform\'ees de Fourier d'int\'egrales orbitales nilpotentes. C'est-\`a-dire, notons $\mathfrak{g}_{\tilde{t}}$ l'alg\`ebre de Lie de $G_{\tilde{t}}$ et $Nil(\mathfrak{g}_{\tilde{t}})$ l'ensemble des orbites nilpotentes dans $\mathfrak{g}_{\tilde{t}}(F)$. Il existe un voisinage $\omega$ de $0$ dans  $\mathfrak{g}_{\tilde{t}}(F)$ et, pour tout ${\cal O}\in Nil(\mathfrak{g}_{\tilde{t}})$, il existe un nombre complexe $c_{\tilde{\pi},{\cal O}}(\tilde{t})$ de sorte que, pour toute fonction $\varphi\in C_{c}^{\infty}(\mathfrak{g}_{\tilde{t}}(F))$ \`a support dans $\omega$, on ait l'\'egalit\'e
$$\int_{\mathfrak{g}_{\tilde{t}}(F)}\Theta_{\tilde{\pi}}(\tilde{t}exp(X))\varphi(X)dX=\sum_{{\cal O}\in Nil(\mathfrak{g}_{\tilde{t}})}c_{\tilde{\pi},{\cal O}}(\tilde{t})\int_{{\cal O}}\hat{\varphi}(X)dX,$$
o\`u $\hat{\varphi}$ est la transform\'ee de Fourier de $\varphi$. Evidemment, pour que cette formule ait un sens, on doit d\'efinir pr\'ecis\'ement la transformation de Fourier ainsi que les diverses mesures. Remarquons que les orbites nilpotentes dans $\mathfrak{g}_{\tilde{t}}(F)$ sont les m\^emes que celles dans l'alg\`ebre de Lie $\mathfrak{so}(\tilde{\zeta}_{G,T})(F)$. Supposons d'abord $d$ impair. Alors, parce que $dim(W')$ est paire, l'espace $W''\oplus Z$ de la forme $\tilde{\zeta}_{G,T}$ est de dimension impaire. Puisque cette forme est quasi-d\'eploy\'ee, il y a une unique orbite nilpotente r\'eguli\`ere dans $\mathfrak{so}(\tilde{\zeta}_{G,T})(F)$. On la note ${\cal O}_{reg}$ et on pose $c_{\tilde{\pi}}(\tilde{t})=c_{\tilde{\pi},{\cal O}_{reg}}(\tilde{t})$. Supposons maintenant $d$ pair.  Alors $\mathfrak{so}(\tilde{\zeta}_{G,T})(F)$ poss\`ede en g\'en\'eral plusieurs orbites nilpotentes r\'eguli\`eres. Mais on peut en s\'electionner une de la fa\c{c}on suivante. On peut d\'ecomposer l'espace  $W''\oplus Z$  muni de sa forme $\tilde{\zeta}_{G,T}$ en somme orthogonale d'un hyperplan $X$ et d'une droite $D_{0}$ sur laquelle la forme quadratique est \'equivalente \`a $x\mapsto 2\nu x^2$. Nos hypoth\`eses impliquent que le groupe sp\'ecial orthogonal $SO(X)$ est quasi-d\'eploy\'e. Puisque $ dim(X)$ est impaire, $\mathfrak{so}(X)(F)$ poss\`ede une unique orbite nilpotente r\'eguli\`ere.  Fixons un point de cette orbite et notons ${\cal O}_{\nu}$ l'orbite dans $\mathfrak{so}(\tilde{\zeta}_{G,T})(F)$ qui contient ce point. C'est encore une orbite nilpotente r\'eguli\`ere. On pose $c_{\tilde{\pi}}(\tilde{t})=c_{\tilde{\pi},{\cal O}_{\nu}}(\tilde{t})$. On a ainsi d\'efini une fonction $c_{\tilde{\pi}}$ presque partout sur $\tilde{T}(F)$. Cette fonction est invariante par conjugaison par $T(F)$ et peut \^etre consid\'er\'ee comme une fonction sur $\tilde{T}(F)_{/\theta}$. Par une construction similaire, on d\'efinit une fonction $c_{\tilde{\rho}}$ presque partout sur le m\^eme ensemble.

Posons
$$\epsilon_{geom,\nu}(\tilde{\pi},\tilde{\rho})=\sum_{\tilde{T}\in {\cal T}} \vert W(H,\tilde{T})\vert ^{-1}\int_{\tilde{T}(F)_{/\theta}}c_{\tilde{\pi}}(\tilde{t})c_{\tilde{\rho}}(\tilde{t})D^{\tilde{H}}(\tilde{t})\Delta_{r}(\tilde{t})d\tilde{t}.$$
Cette expression est absolument convergente. Notre r\'esultat principal est le th\'eor\`eme 7.1 dont voici l'\'enonc\'e.

\ass{Th\'eor\`eme}{ Soit  $\pi$, resp. $\rho$, une repr\'esentation admissible, irr\'eductible, temp\'er\'ee et autoduale de $G(F)$, resp. $H(F)$.  Alors on a l'\'egalit\'e
$$\epsilon_{geom,\nu}(\tilde{\pi},\tilde{\rho})=\epsilon_{\nu}(\tilde{\pi},\tilde{\rho}).$$}

Comme nous l'avons expliqu\'e, notre motivation est la conjecture locale de Gross-Prasad. Nous ignorons si cette fa\c{c}on, plut\^ot compliqu\'ee, de calculer un facteur $\epsilon$ peut avoir d'autres applications.

La d\'emonstration reprend celle de [W2] et [W3]. Donnons de tr\`es br\`eves indications dans le cas o\`u $r=0$ et $d=m+1$ (en fait, ce cas ne peut pas \^etre trait\'e \`a part car la preuve utilise une r\'ecurrence compliqu\'ee sur le couple $(d,m)$). Consid\'erons une fonction $\tilde{f}\in C_{c}^{\infty}(\tilde{G}(F))$ qui est tr\`es cuspidale, cf. 1.7. On d\'efinit une suite $(\Omega_{N})_{N\geq1}$ de sous-ensembles ouverts compacts de $H(F)\backslash G(F)$ v\'erifiant les propri\'et\'es usuelles
$$\Omega_{N}\subset \Omega_{N+1}\text{ pour tout N et }\bigcup_{N} \Omega_{N}=H(F)\backslash G(F).$$
On note $\kappa_{N}$ la fonction caract\'eristique de l'image r\'eciproque de $\Omega_{N}$ dans $G(F)$. Cela \'etant, on pose
$$J_{N}(\Theta_{\tilde{\rho}},\tilde{f})=\int_{G(F)}\int_{\tilde{H}(F)}\Theta_{\tilde{\rho}}(\tilde{y})\tilde{f}(g^{-1}\tilde{y}g)d\tilde{y}\kappa_{N}(g)dg.$$
Cette int\'egrale est absolument convergente. Comme pour la formule des traces locale d'Arthur, il y a deux fa\c{c}ons de calculer la limite de cette expression quand $N$ tend vers l'infini. L'une, que l'on peut qualifier de "g\'eom\'etrique", et qui est la r\'eplique de celle de [W2]; l'autre, que l'on peut qualifier de "spectrale", qui s'appuie sur la formule de Plancherel pour le groupe $G(F)$ et qui est la r\'eplique de celle de [W3]. Ces deux voies conduisent \`a une \'egalit\'e
$$J_{geom}(\Theta_{\tilde{\rho}},\tilde{f})=lim_{N\to \infty}J_{N}(\Theta_{\tilde{\rho}},\tilde{f})=J_{spec}(\Theta_{\tilde{\rho}},\tilde{f}).$$
Les deux expressions extr\^emes contiennent des distributions (en $\tilde{f}$) qui ne sont pas invariantes: des int\'egrales orbitales pond\'er\'ees et des caract\`eres pond\'er\'es. Le proc\'ed\'e habituel d'Arthur permet par r\'ecurrence d'en d\'eduire d'autres expressions qui ne contiennent plus que des distributions invariantes, et qui continuent d'\^etre \'egales entre elles. Supposons que $\tilde{\pi}$ est  "elliptique", le cas g\'en\'eral s'en d\'eduisant assez facilement. On prend pour $\tilde{f}$ un pseudo-coefficient de $\tilde{\pi}$. Alors les deux expressions "invariantes" ci-dessus deviennent respectivement $\epsilon_{geom,\nu}(\tilde{\pi},\tilde{\rho})$ et $\epsilon_{\nu}(\tilde{\pi},\tilde{\rho})$, ce qui prouve l'\'egalit\'e de ces deux termes.

Expliquons pourquoi appara\^{\i}t le terme $\epsilon(1/2,\pi\times \rho,\psi)$. Notons $E_{\pi}$ et $E_{\rho}$ des espaces dans lesquels se r\'ealisent  $\pi$ et $\rho$. Consid\'erons l'espace $Hom_{H(F)}(\pi,\rho)$ des applications lin\'eaires $\varphi:E_{\pi}\to E_{\rho}$ telles que $\varphi\circ \pi(h)=\rho(h)\circ\varphi$ pour tout $h\in H(F)$. D'apr\`es [AGRS], cet espace est de dimension $1$.   Pour $\varphi\in Hom_{H(F)}(\pi,\rho)$ et $\tilde{y}\in \tilde{H}(F)$, l'application lin\'eaire $\tilde{\rho}(\tilde{y})^{-1}\circ\varphi\circ\tilde{\pi}(\tilde{y})$ appartient encore \`a $Hom_{H(F)}(\pi,\rho)$ et ne d\'epend pas de $\tilde{y}$. Il existe donc un nombre $c\in {\mathbb C}^{\times}$ tel que
$$\tilde{\rho}(\tilde{y})^{-1}\circ\varphi\circ\tilde{\pi}(\tilde{y})=c\varphi$$
pour tous $\varphi\in Hom_{H(F)}(\pi,\rho)$ et $\tilde{y}\in \tilde{H}(F)$. C'est ce nombre $c$ qui appara\^{\i}t naturellement dans nos calculs spectraux. Or le th\'eor\`eme 2.7 de [JPSS] permet de calculer $c$: \`a des termes \'el\'ementaires pr\`es, c'est $\epsilon(1/2,\pi\times \rho,\psi)$.

Voici le contenu de l'article. La premi\`ere section rassemble des d\'efinitions et r\'esultats g\'en\'eraux sur les "groupes tordus", selon la terminologie de Labesse. La deuxi\`eme introduit plus pr\'ecis\'ement les groupes $GL_{d}$ tordus. La partie "g\'eom\'etrique" de notre formule int\'egrale est trait\'ee dans la section 3. La section 4 contient les majorations n\'ecessaires pour prouver les diverses convergences d'int\'egrales  utilis\'ees dans la section 6.  La section 5 \'etablit le r\'esultat \'evoqu\'e ci-dessus, \`a savoir que $\epsilon(1/2,\pi\times \rho,\psi)$ mesure la compatibilit\'e des deux prolongements $\tilde{\pi}$ et $\tilde{\rho}$. La partie "spectrale" de la formule int\'egrale est trait\'ee dans la section 6. Le th\'eor\`eme principal est prouv\'e dans la septi\`eme et derni\`ere section. Ainsi qu'on l'a d\'ej\`a dit, nos preuves sont parall\`eles \`a celles de [W2] et [W3], au point d'\^etre parfois identiques. Pour \'epargner le lecteur, ou par paresse, on s'est souvent content\'e d'indiquer sommairement les modifications \`a apporter ou m\^eme simplement d'affirmer les r\'esultats sans d\'emonstration. 

{\bf Remarque sur la notation.} Dans les articles [W2] et [W3], on avait  utilis\'e la lettre $\theta$ pour noter les caract\`eres ou quasi-caract\`eres ($\theta_{\pi}$, $\theta_{f}$ etc...). Ici, cette lettre sera r\'eserv\'ee aux automorphismes des groupes lin\'eaires. Les caract\`eres ou quasi-caract\`eres seront not\'es par la lettre $\Theta$.

Je remercie R. Beuzart pour m'avoir signal\'e une erreur dans une premi\`ere version de ce texte.

\bigskip

\section{Notations, groupes tordus}

\bigskip

\subsection{Notations g\'en\'erales}

Soit $F$ un corps local non archim\'edien de caract\'eristique nulle. On note $\vert .\vert _{F}$ la valeur absolue usuelle de $F$, $val_{F}$ la valuation, $\mathfrak{o}_{F}$ l'anneau des entiers et $\mathfrak{p}_{F}$ son id\'eal maximal. On fixe une cl\^oture alg\'ebrique $\bar{F}$ de $F$ et une uniformisante $\varpi_{F}$ de $F$.

Toutes les vari\'et\'es alg\'ebriques seront suppos\'ees d\'efinies sur $F$, sauf mention explicite du contraire. De m\^eme pour les actions d'un groupe alg\'ebrique sur une vari\'et\'e. Soit $G$ un groupe alg\'ebrique r\'eductif connexe. On note $A_{G}$ le plus grand tore d\'eploy\'e central dans $G$, $X(G)$ le groupe des caract\`eres de $G$ d\'efinis sur $F$, ${\cal A}_{G}=Hom(X(G),{\mathbb R})$ et ${\cal A}_{G}^*=X(G)\otimes_{{\mathbb Z}}{\mathbb R}$ le dual de ${\cal A}_{G}$. On note $\mathfrak{g}$ l'alg\`ebre de Lie de $G$ et
$$\begin{array}{ccc}G\times \mathfrak{g}&\to&\mathfrak{g}\\ (g,X)&\mapsto& gXg^{-1}\\ \end{array}$$
l'action adjointe.

Quand on d\'efinit un objet relatif au groupe $G$, on peut pr\'eciser si besoin est la notation en introduisant la lettre $G$ en exposant.  Par exemple, les int\'egrales orbitales pond\'er\'ees sont not\'ees $J_{M}(x,f)$ s'il n'y a pas d'ambig\" uit\'e sur le groupe ambiant $G$, ou $J_{M}^G(x,f)$ si cela semble pr\'ef\'erable.

Pour toute bijection $\theta$ d'un ensemble $X$ dans lui-m\^eme, on note $X^{\theta}$ le sous-ensemble des points fixes.

\bigskip
\subsection{Groupes tordus}

On appelle groupe tordu un couple $(G,\tilde{G})$ v\'erifiant les conditions qui suivent. Le terme $G$ est un groupe r\'eductif connexe. Le terme $\tilde{G}$ est une vari\'et\'e alg\'ebrique telle que $\tilde{G}(F)\not= \emptyset$. Il y a deux actions de groupe alg\'ebrique de $G$ sur $\tilde{G}$, \`a droite et \`a gauche, not\'ees
$$\begin{array}{ccccccc}G\times\tilde{G}&\to&\tilde{G}&&\tilde{G}\times G&\to&\tilde{G}\\(g,\tilde{x})&\mapsto&g\tilde{x}&&(\tilde{x},g)&\mapsto&\tilde{x}g.\\ \end{array}$$
Chacune d'elles fait de $\tilde{G}$ un espace principal homog\`ene sur $G$.

Consid\'erons un tel groupe tordu. Notons $Aut(G)$ le groupe des automorphismes de $G$. Il existe une unique application alg\'ebrique
$$\begin{array}{ccc}\tilde{G}&\to&Aut(G)\\ \tilde{x}&\mapsto&\theta_{\tilde{x}}\\ \end{array}$$
de sorte que l'on ait l'\'egalit\'e $\tilde{x}g=\theta_{\tilde{x}}(g)\tilde{x}$ pour tous $\tilde{x}\in \tilde{G}$ et $g\in G$. De $\theta_{\tilde{x}}$ se d\'eduit des automorphismes de $X(G)$, de $A_{G}$, de ${\cal A}_{G}$ etc... qui ne d\'ependent pas de $\tilde{x}\in \tilde{G}(F)$. On les note $\theta_{\tilde{G}}$. On suppose

{\bf Hypoth\`ese:} $\theta_{\tilde{G}}$ est d'ordre fini.

Le groupe $G$ op\`ere sur $\tilde{G}$ par conjugaison: $(g,\tilde{x})\mapsto g\tilde{x}g^{-1}$. Pour tout sous-ensemble $\tilde{X}$ de $\tilde{G}$, on note $Norm_{G}(\tilde{X})$ le normalisateur de $\tilde{X}$ dans $G$ et $Z_{G}(\tilde{X})$ son centralisateur. Si $\tilde{X}$ est r\'eduit \`a un point $\{\tilde{x}\}$, on note simplement ces groupes $Z_{G}(\tilde{x})$ et on note $G_{\tilde{x}}$ la composante neutre de $Z_{G}(\tilde{x})$. La vari\'et\'e $\tilde{G}$ op\`ere sur $G$ par $(\tilde{x},g)\mapsto \theta_{\tilde{x}}(g)$. Si $X$ est un sous-ensemble de $G$, on d\'efinit alors son normalisateur $N_{\tilde{G}}(X)$ et son centralisateur $Z_{\tilde{G}}(X)$, avec la variante $Z_{\tilde{G}}(x)$ quand $X$ est r\'eduit \`a un point $\{x\}$.

On note $A_{\tilde{G}}$ le plus grand sous-tore de $A_{G}$ sur lequel $\theta_{\tilde{G}}$ agit trivialement. On note ${\cal A}_{\tilde{G}}$, resp. ${\cal A}^*_{\tilde{G}}$, le sous-groupe des \'el\'ements de ${\cal A}_{G}$, resp. ${\cal A}^*_{G}$ fix\'es par $\theta_{\tilde{G}}$. On note $a_{\tilde{G}}$ la dimension de ${\cal A}_{\tilde{G}}$. On d\'efinit un homomorphisme
$H_{\tilde{G}}:G(F)\to {\cal A}_{\tilde{G}}$ par $H_{\tilde{G}}(g)(\chi)=log(\vert \chi(g)\vert _{F})$ pour tous $g\in G(F)$ et $\chi\in X(G)^{\theta_{\tilde{G}}}$.
 
 On note $\tilde{G}_{ss}$ le sous-ensemble des \'el\'ements semi-simples de $\tilde{G}$, c'est-\`a-dire des $\tilde{x}\in \tilde{G}$ tels qu'il existe un sous-tore maximal $T$ de $G$, d\'efini sur $\bar{F}$, et un sous-groupe de Borel $B$ de $G$, d\'efini sur $\bar{F}$ et contenant $T$, tels que $\theta_{\tilde{x}}$ conserve $T$ et $B$. Si $\tilde{x}\in \tilde{G}_{ss}(F)$, on pose
$$D^{\tilde{G}}(\tilde{x})=\vert det(1-\theta_{\tilde{x}})_{\vert \mathfrak{g}/\mathfrak{g}_{\tilde{x}}}\vert _{F}.$$
On note $\tilde{G}_{reg}$ l'ensemble des \'el\'ements fortement r\'eguliers de $\tilde{G}$, c'est-\`a-dire l'ensemble des \'el\'ements $\tilde{x}$ tels que $Z_{G}(\tilde{x})$ soit ab\'elien et $G_{\tilde{x}}$ soit un tore.

On appelle sous-groupe parabolique tordu $(P,\tilde{P})$ un couple v\'erifiant les conditions suivantes. Le terme $P$ est un sous-groupe parabolique de $G$. Le terme $\tilde{P}$ est son normalisateur dans $\tilde{G}$ et on suppose $\tilde{P}(F)\not=\emptyset$. Pour un tel couple, on appelle composante de L\'evi tordue un couple $(M,\tilde{M})$ tel que $M$ soit une composante de L\'evi de $P$ et $\tilde{M}$ est l'intersection des normalisateurs de $P$ et $M$ dans $\tilde{G}$. On appelle groupe de L\'evi tordu de $(G,\tilde{G})$ un couple $(M,\tilde{M})$ tel qu'il existe un sous-groupe parabolique tordu $(P,\tilde{P})$ dont $(M,\tilde{M})$ est une composante de L\'evi tordue. On v\'erifie qu'un groupe de L\'evi tordu est un groupe tordu (c'est-\`a-dire que $\tilde{M}(F)\not=\emptyset$, cf. [W1] 1.6). Pour un tel groupe de L\'evi tordu, on note ${\cal P}(\tilde{M})$ l'ensemble des sous-groupes paraboliques tordus $(P,\tilde{P})$ de composante de L\'evi tordue $(M,\tilde{M})$, ${\cal F}(\tilde{M})$ l'ensemble des sous-groupes paraboliques tordus $(Q,\tilde{Q})$ tels que $M\subset Q$ et $\tilde{M}\subset \tilde{Q}$ et ${\cal L}(\tilde{M})$ l'ensemble des groupes de L\'evi tordus $(L,\tilde{L})$ tels que $M\subset L$ et $\tilde{M}\subset \tilde{L}$. Soit $(Q,\tilde{Q}) $ un sous-groupe parabolique tordu. On \'ecrit simplement $\tilde{Q}=\tilde{L}U$ pour signifier que:

- $\tilde{L}$ est le second membre d'une composante de L\'evi tordue $(L,\tilde{L})$ de $(Q,\tilde{Q})$;

- $U$ est le radical unipotent de $Q$.  

Si $(M,\tilde{M})$ est un groupe de L\'evi fix\'e et que $(Q,\tilde{Q})$ appartient \`a ${\cal F}(\tilde{M})$, il sera de plus suppos\'e que $\tilde{L}$ contient $\tilde{M}$.

Comme dans le cas non tordu, les L\'evi tordus se caract\'erisent comme les commutants de tores d\'eploy\'es. Pr\'ecis\'ement, soit $A$ un sous-tore d\'eploy\'e de $G$, notons $M$ son commutant dans $G$ et $Z_{\tilde{G}}(A)$ son commutant dans $\tilde{G}$. Supposons $Z_{\tilde{G}}(A)(F)\not=\emptyset$. Alors $(M,Z_{\tilde{G}}(A))$ est un L\'evi tordu. Prouvons cela. On sait bien que $M$ est un L\'evi de $G$. Choisissons un \'el\'ement $x_{*}\in X_{*}(A)$ en position g\'en\'erale et posons $a=x_{*}(\varpi_{F})$. Alors $M$ est le commutant de $a$. Notons $P$ le sous-groupe parabolique de $G$, de composante de L\'evi $M$, dont le radical unipotent est engendr\'e par les sous-groupes radiciels associ\'es aux racines $\alpha$ de $A_{M}$ telles que $\vert \alpha(a)\vert _{F}>1$. Soit $\tilde{x}\in Z_{\tilde{G}}(A)$. Puisque $\theta_{\tilde{x}}$ fixe $a$, il conserve aussi $M$ et $P$. Donc $\tilde{x}$ appartient \`a l'ensemble $\tilde{M}$ d\'efini comme plus haut, et $Z_{\tilde{G}}(A)\subset \tilde{M}$. Puisque $Z_{\tilde{G}}(A)(F)$ n'est pas vide, $\tilde{M}(F)$ ne l'est pas non plus, ce qui est la condition pour que $(M,\tilde{M})$ soit un groupe de L\'evi tordu. De plus, $Z_{\tilde{G}}(A)$ et $\tilde{M}$ sont \'evidemment tous deux des espaces principaux homog\`enes pour le groupe $M$. Ils sont donc \'egaux, ce qui prouve l'assertion. Inversement, tout groupe de L\'evi tordu $(M,\tilde{M})$ est le commutant au sens ci-dessus du tore $A_{\tilde{M}}$.

Soit $P_{min}$ un sous-groupe parabolique minimal de $G$ et $M_{min}$ une composante de L\'evi de $P_{min}$. On peut compl\'eter ces donn\'ees, de fa\c{c}on unique, en un sous-groupe parabolique tordu $(P_{min},\tilde{P}_{min})$ et une composante de L\'evi tordue $(M_{min},\tilde{M}_{min})$. En effet,  $\tilde{P}_{min}$ est le normalisateur de $P_{min}$ dans $\tilde{G}$ et $\tilde{M}_{min}$ est le normalisateur de $M_{min}$ dans $\tilde{P}_{min}$. On doit voir que $\tilde{M}_{min}(F)\not=\emptyset$. Soit $\tilde{y}\in \tilde{G}(F)$. Alors le couple $(\theta_{\tilde{y}}(P_{min}),\theta_{\tilde{y}}(M_{min}))$ est form\'e d'un  sous-groupe parabolique minimal de $G$ et d'une composante de L\'evi de ce sous-groupe. Deux tels couples sont conjugu\'es par un \'el\'ement de $G(F)$. Choisissons donc $g\in G(F)$ tel que la conjugaison par $g$ envoie $(\theta_{\tilde{y}}(P_{min}),\theta_{\tilde{y}}(M_{min}))$ sur $(P_{min},M_{min})$. Posons $\tilde{x}=g\tilde{y}$. Alors $\tilde{x}$ appartient \`a $\tilde{M}_{min}(F)$.

On appelle sous-tore maximal tordu un couple $(T,\tilde{T})$ v\'erifiant les conditions suivantes. Le terme $T$ est un sous-tore maximal de $G$. Le terme $\tilde{T}$ est une sous-vari\'et\'e de $\tilde{G}$. L'ensemble $\tilde{T}(F)$ n'est pas vide et il existe un sous-groupe de Borel $B$ de $G$, d\'efini sur $\bar{F}$, contenant $T$, tel que $\tilde{T}$ soit l'intersection des normalisateurs de $T$ et $B$ dans $\tilde{G}$. Cela entra\^{\i}ne que l'on a $\tilde{T}=T\tilde{x}=\tilde{x}T$ pour tout $\tilde{x}\in \tilde{T}$. La restriction \`a $T$ de l'automorphisme $\theta_{\tilde{x}}$ ne d\'epend pas de $\tilde{x}$, on la note $\theta_{\tilde{T}}$ ou simplement $\theta$ si cela ne cr\'ee pas d'ambig\" uit\'e. On note $T_{\theta}$ la composante neutre du sous-groupe des points fixes $T^{\theta}$.

Evidemment, pour un couple $(G,\tilde{G})$, ou   $(P,\tilde{P})$  etc... le premier terme $G$ ou $P$ etc... est uniquement d\'etermin\'e par le second $\tilde{G}$ ou $\tilde{P}$ etc... Dans la suite de l'article, on parlera simplement du groupe tordu $\tilde{G}$, ou du sous-groupe parabolique tordu $\tilde{P}$  etc... Les termes $G$ ou $P$ etc... seront utilis\'es si besoin est sans les d\'efinir explicitement. D'autre part, on peut utiliser les d\'efinitions que l'on vient d'introduire dans le cas o\`u $\tilde{G}=G$ muni des multiplications \`a droite et \`a gauche. On supprime alors les $\tilde{}$; par exemple, on d\'efinit les ensembles ${\cal P}(M)$, ${\cal L}(M)$, la fonction $H_{P}$ etc...

{\bf Remarque.} Les espaces tordus ont \'et\'e introduits par Labesse. D'autres auteurs pr\'ef\`erent \'etudier les groupes non connexes. Un espace tordu $(G,\tilde{G})$ appara\^{\i}t comme une composante d'un tel groupe. En effet, fixons $\tilde{x}\in \tilde{G}(F)$.  D'apr\`es l'hypoth\`ese de finitude faite plus haut, on peut choisir un entier $n\geq1$ tel que l'image de $\theta_{\tilde{x}}^n$ dans le groupe des automorphismes ext\'erieurs de $G$ soit \'egale \`a $1$. Donc $\theta_{\tilde{x}}^n$ est l'automorphisme int\'erieur associ\'e \`a un \'el\'ement du groupe adjoint $G_{ad}(F)$. L'image de $G(F)$ dans $G_{ad}(F)$ est d'indice fini. Quitte \`a accro\^{\i}tre $n$, on peut donc supposer que $\theta_{\tilde{x}}^n$ est l'automorphisme int\'erieur associ\'e \`a un \'el\'ement $x\in G(F)$. Fixons un tel $x$. Consid\'erons le groupe ab\'elien libre $\tilde{X}$ engendr\'e par $\tilde{x}$. Il agit sur $G$: $\tilde{x}^m$ agit par $\theta_{\tilde{x}}^m$. Consid\'erons le produit semi-direct $G\rtimes \tilde{X}$, puis le plus petit sous-groupe distingu\'e de ce produit contenant $x^{-1}\tilde{x}^n$. Notons $G^+$ le quotient. Alors $G^+$ est un groupe lin\'eaire alg\'ebrique de composante neutre $G$ et $\tilde{G}$ s'identifie \`a la composante connexe de $G^+$ qui contient $\tilde{x}$. Cette remarque permet d'appliquer les r\'esultats d\'emontr\'es  dans la litt\'erature pour des groupes non connexes.

\bigskip

\subsection{Les espaces ${\cal A}_{\tilde{M}}$}

Soit $\tilde{G}$ un groupe tordu. Fixons un L\'evi minimal $M_{min}$ de $G$. On note ${\cal L}^{\tilde{G}}$ l'ensemble des L\'evis tordus $\tilde{M}$ de $\tilde{G}$ tels que $M_{min}\subset M$. On d\'efinit le groupe de Weyl usuel $W^G=Norm_{G(F)}(M_{min})/M_{min}$. 

On munit ${\cal A}_{M_{min}}$ d'un produit scalaire invariant par l'action du groupe de Weyl $W^G$.  Pour tout L\'evi tordu $\tilde{M}$, on en d\'eduit par conjugaison et restriction un produit scalaire sur ${\cal A}_{\tilde{M}}$.  Si  $\tilde{L}$ est un L\'evi tordu  tel que $\tilde{M}\subset \tilde{L}$, on note ${\cal A}_{\tilde{M}}^{\tilde{L}}$ l'orthogonal de ${\cal A}_{\tilde{L}}$ dans ${\cal A}_{\tilde{M}}$. On note $\zeta\mapsto \zeta_{\tilde{L}}$ et $\zeta\mapsto \zeta^{\tilde{L}}$ les projections orthogonales de ${\cal A}_{\tilde{M}}$ sur ${\cal A}_{\tilde{L}}$, resp. sur ${\cal A}_{\tilde{M}}^{\tilde{L}}$.  Fixons un sous-groupe compact sp\'ecial $K$ de $G(F)$ en bonne position relativement \`a $M_{min}$ et supposons $M_{min}\subset M$.  Soit $\tilde{P}=\tilde{M}U\in {\cal P}(\tilde{M})$. On d\'efinit une fonction $H_{\tilde{P}}:G(F)\to {\cal A}_{\tilde{M}}$ par $H_{\tilde{P}}(g)=H_{\tilde{M}}(m)$ pour $g=muk$, avec $m\in M(F)$, $u\in U(F)$ et $k\in K$.

On fixe une extension de $H_{\tilde{G}}$ \`a $\tilde{G}(F)$, c'est-\`a-dire une fonction que l'on note encore $H_{\tilde{G}}:\tilde{G}(F)\to {\cal A}_{\tilde{G}}$ telle que $H_{\tilde{G}}(g\tilde{x}g')=H_{\tilde{G}}(g)+H_{\tilde{G}}(\tilde{x})+H_{\tilde{G}}(g')$ pour tous $g,g'\in G(F)$ et $\tilde{x}\in \tilde{G}(F)$. Pour tout L\'evi tordu $\tilde{M}$, on en d\'eduit une fonction analogue $H_{\tilde{M}}:\tilde{M}(F)\to {\cal A}_{\tilde{M}}$ de la fa\c{c}on suivante. Posons $W^{\tilde{M}}=Norm_{G(F)}(\tilde{M})/M(F)$. Fixons, ainsi qu'il est loisible, une extension $H'_{\tilde{M}}$ de $H_{\tilde{M}}$ \`a $\tilde{M}(F)$ de sorte que la compos\'ee de cette fonction et de la projection orthogonale de ${\cal A}_{\tilde{M}}$ sur ${\cal A}_{\tilde{G}}$ co\"{\i}ncide avec la restriction de $H_{\tilde{G}}$ \`a $\tilde{M}(F)$.  Pour $g\in Norm_{G(F)}(\tilde{M})$ et $\tilde{x}\in \tilde{M}(F)$,  le terme $H'_{\tilde{M}}(g\tilde{x}g^{-1})$ ne d\'epend que de l'image $w$ de $g$ dans $W^{\tilde{M}}$. Notons-le $H'_{\tilde{M}}(w\tilde{x}w^{-1})$. D'autre part, le groupe $W^{\tilde{M}}$ agit naturellement dans ${\cal A}_{\tilde{M}}$. On pose
$$H_{\tilde{M}}(\tilde{x})=\vert W^{\tilde{M}}\vert ^{-1}\sum_{w\in W^{\tilde{M}}}w^{-1}H'_{\tilde{M}}(w\tilde{x}w^{-1}).$$
L'application $H_{\tilde{M}}$ ainsi d\'efinie sur $\tilde{M}(F)$ est un prolongement de l'application not\'ee de la m\^eme fa\c{c}on sur $M(F)$. Elle ne d\'epend pas de l'application auxiliaire $H'_{\tilde{M}}$. Si $\tilde{L}\in {\cal L}(\tilde{M})$,  la restriction de $H_{\tilde{L}}$ \`a $\tilde{M}(F)$ co\"{\i}ncide avec la compos\'ee de $H_{\tilde{M}}$ et de la projection orthogonale de ${\cal A}_{\tilde{M}}$ sur ${\cal A}_{\tilde{L}}$. Pour $g\in G(F)$, la conjugaison par $g$ d\'efinit un isomorphisme de ${\cal A}_{\tilde{M}}$ sur ${\cal A}_{g\tilde{M}g^{-1}}$ et on a l'\'egalit\'e $H_{g\tilde{M}g^{-1}}(g\tilde{x}g^{-1})=gH_{\tilde{M}}(\tilde{x})g^{-1}$ pour tout $\tilde{x}\in \tilde{M}(F)$.

 \bigskip

\subsection{Mesures}

Soient $G$ un groupe  r\'eductif connexe, $M_{min}$ un groupe de L\'evi minimal de $G$ et $K$ un sous-groupe compact sp\'ecial de $G(F)$ en bonne position relativement \`a $M_{min}$. Fixons une mesure de Haar sur $G(F)$ et munissons $K$ de la mesure de masse totale $1$ (remarquons que l'on ne suppose pas que cette mesure sur $K$ soit \'egale \`a la restriction de la mesure sur $G(F)$). Alors, pour tout $P=MU\in {\cal F}(M_{min})$, Arthur a d\'efini une mesure de Haar sur les groupes $M(F)$ et $U(F)$, cf. [A1] paragraphe 1. Soit $T$ un tore. Si $T$ est d\'eploy\'e, on munit $T(F)$ de la mesure telle que le plus grand sous-groupe compact de $T(F)$ soit de mesure $1$. En g\'en\'eral, on munit $A_{T}(F)$ de cette mesure puis $T(F)$ de celle pour laquelle la mesure quotient sur $T(F)/A_{T}(F)$ soit de masse totale $1$.

Fixons un caract\`ere continu non trivial $\psi$ de $F$. On munit $F$ de la mesure autoduale pour $\psi$. Fixons une forme bilin\'eaire sym\'etrique non d\'eg\'en\'er\'ee $<.,.>$ sur $\mathfrak{g}(F)\times \mathfrak{g}(F)$, invariante par l'action adjointe de $G(F)$.  Pour toute sous-alg\`ebre $\mathfrak{h}$ de $\mathfrak{g}$ telle que la restriction de $<.,.>$ \`a $\mathfrak{h}(F)$ soit non d\'eg\'en\'er\'ee, on munit $\mathfrak{h}(F)$ de la mesure de Haar autoduale pour le bicaract\`ere $(X,Y)\mapsto \psi(<X,Y>)$. Soit $H$ un sous-groupe de $G$, supposons que l'on a muni $H(F)$ d'une mesure de Haar  $dh$ et que la restriction de $<.,.>$ soit non d\'eg\'en\'er\'ee. La mesure sur $\mathfrak{h}(F)$ se rel\`eve par l'exponentielle en une mesure de Haar $d'h$ sur $H(F)$ qui n'a aucune raison d'\^etre celle que l'on a fix\'ee. On note $\nu(H)$ la constante telle $dh=\nu(H)d'h$ (en fait, nous n'utiliserons cette notation que dans le cas o\`u $H$ est un sous-tore de $G$).

Supposons que $G$ soit la premi\`ere composante d'un groupe tordu $(G,\tilde{G})$. La mesure sur $G(F)$ en d\'etermine une sur l'espace principal homog\`ene $\tilde{G}(F)$.   Consid\'erons un sous-tore maximal tordu $(T,\tilde{T})$.  Le groupe $T(F)$ agit par conjugaison sur $\tilde{T}(F)$. Notons $\tilde{T}(F)_{/\theta}$ l'ensemble des orbites. Il est naturellement muni d'une structure de groupe de Lie sur $F$ et d'une action de $T(F)$ par multiplication \`a gauche. Pour tout $\tilde{t}\in \tilde{T}(F)_{/\theta}$, l'application
$$\begin{array}{ccc}T_{\theta}(F)&\to&\tilde{T}(F)_{/\theta}\\ t&\mapsto&t\tilde{t}\\ \end{array}$$
est un isomorphisme local. Il existe une mesure sur $\tilde{T}(F)_{/\theta}$, invariante par l'action de $T(F)$ et ind\'ependante du choix de $\tilde{t}$, telle que cette application conserve localement les mesures. On munit $\tilde{T}(F)_{/\theta}$ de cette mesure. On pose
$$W(G,\tilde{T})=Norm_{G(F)}(\tilde{T})/T(F).$$
On dit que $\tilde{T}$ est elliptique dans $\tilde{G}$ si $A_{\tilde{T}}=A_{\tilde{G}}$.
Fixons un ensemble ${\cal T}(\tilde{G})$, resp. ${\cal T}_{ell}(\tilde{G})$, de repr\'esentants des classes de conjugaison par $G(F)$ dans l'ensemble des sous-tores maximaux tordus, resp. et elliptiques, de $\tilde{G}$. On fixe des ensembles analogues pour tout L\'evi tordu. La formule de Weyl prend l'une ou l'autre des formes suivantes
$$\int_{\tilde{G}(F)}\tilde{f}(\tilde{x})d\tilde{x}=\sum_{\tilde{T}\in {\cal T}(\tilde{G})}\vert W(G,\tilde{T})\vert ^{-1}[T(F)^{\theta}:T_{\theta}(F)]^{-1}\int_{\tilde{T}(F)_{/\theta}}\int_{T_{\theta}(F)\backslash G(F)}\tilde{f}(g^{-1}\tilde{t}g)dgD^{\tilde{G}}(\tilde{t})d\tilde{t}$$
$$\qquad =\sum_{\tilde{M}\in {\cal L}^{\tilde{G}}}\vert W^{M}\vert \vert W^{G}\vert ^{-1}\sum_{\tilde{T}\in {\cal T}_{ell}(\tilde{M})}\vert W(M,\tilde{T})\vert ^{-1}[T(F)^{\theta}:T_{\theta}(F)]^{-1}$$
$$\qquad \int_{\tilde{T}(F)_{/\theta}}\int_{T_{\theta}(F)\backslash G(F)}\tilde{f}(g^{-1}\tilde{t}g)dgD^{\tilde{G}}(\tilde{t})d\tilde{t}$$
pour toute fonction $\tilde{f}\in C_{c}^{\infty}(\tilde{G}(F))$.

On note ${\cal A}_{\tilde{G},F}$, resp. ${\cal A}_{A_{\tilde{G}},F}$, l'image par l'application $H_{\tilde{G}}$ de $G(F)$, resp. $A_{\tilde{G}}(F)$. On note ${\cal A}_{\tilde{G},F}^{\vee}$, resp. ${\cal A}_{A_{\tilde{G}},F}^{\vee}$, le r\'eseau dans ${\cal A}_{\tilde{G}}^*$ form\'e des $\lambda$ tels que $\lambda(\zeta)\in 2\pi{\mathbb Z}$ pour tout $\zeta\in {\cal A}_{\tilde{G},F}$, resp. $\zeta\in {\cal A}_{A_{\tilde{G}},F}$. On munit le quotient $i{\cal A}_{\tilde{G}}^*/i{\cal A}_{A_{\tilde{G}},F}^{\vee}$ de la mesure de Haar de masse totale $1$. On pose $i{\cal A}_{\tilde{G},F}^*=i{\cal A}_{\tilde{G}}^*/i{\cal A}_{\tilde{G},F}^{\vee}$ et on le munit de la mesure telle que l'application naturelle $i{\cal A}_{\tilde{G},F}^*\to i{\cal A}_{\tilde{G}}^*/i{\cal A}_{A_{\tilde{G}},F}^{\vee}$ pr\'eserve localement les mesures.

\bigskip

\subsection{Int\'egrales orbitales pond\'er\'ees}

Dans la suite de cette section, on fixe un groupe tordu $(G,\tilde{G})$. On  suppose fix\'es un L\'evi minimal $M_{min}$ de $G$ et un sous-groupe compact sp\'ecial $K$ de $G(F)$ en bonne position relativement \`a $M_{min}$. Soit $\tilde{M}\in {\cal L}^{\tilde{G}}$. La th\'eorie des $(G,M)$-familles se g\'en\'eralise au cas tordu en une th\'eorie des $(\tilde{G},\tilde{M})$-familles. En particulier, soit $g\in G(F)$. De la famille de points $(H_{\tilde{P}}(g))_{\tilde{P}\in {\cal P}(\tilde{M})}$ se d\'eduit une $(\tilde{G},\tilde{M})$-famille $(v_{\tilde{P}}(g))_{\tilde{P}\in {\cal P}(\tilde{M})}$: on pose $v_{\tilde{P}}(g,\lambda)=e^{-\lambda(H_{\tilde{P}}(g))}$ pour $\lambda\in i{\cal A}_{\tilde{M}}^*$. De cette $(\tilde{G},\tilde{M})$-famille se d\'eduit un nombre $v_{\tilde{M}}(g)$, cf. [A2] p.37. Soient $\tilde{f}\in C_{c}^{\infty}(\tilde{G}(F))$ et $\tilde{x}\in \tilde{M}(F)\cap \tilde{G}_{reg}(F)$. On d\'efinit l'int\'egrale orbitale pond\'er\'ee
$$J_{\tilde{M}}(\tilde{x},\tilde{f})=D^{\tilde{G}}(\tilde{x})^{1/2}\int_{G_{\tilde{x}}(F)\backslash G(F)}\tilde{f}(g^{-1}\tilde{x}g)v_{\tilde{M}}(g)dg.$$
Ces int\'egrales v\'erifient les m\^emes conditions de r\'egularit\'e et de croissance que dans le cas non tordu. 

\bigskip

\subsection{Quasi-caract\`eres}

Soit $\tilde{x}\in \tilde{G}_{ss}(F)$.  On d\'efinit la notion de bon voisinage $\omega\subset \mathfrak{g}_{\tilde{x}}(F)$: la d\'efinition est la m\^eme qu'en [W2] 3.1, en ajoutant quelques $\tilde{}$. On note $Nil(\mathfrak{g}_{\tilde{x}})$ l'ensemble des orbites nilpotentes dans $\mathfrak{g}_{\tilde{x}}(F)$. Pour tout ${\cal O}\in Nil(\mathfrak{g}_{\tilde{x}})$, l'int\'egrale orbitale sur ${\cal O}$ a pour transform\'ee de Fourier une distribution localement int\'egrable, que l'on note $X\mapsto \hat{j}({\cal O},X)$. Evidemment, on doit d\'efinir correctement les mesures, on renvoie pour cela \`a [W2] 1.2.

Soit $\Theta$ une fonction d\'efinie presque partout sur $\tilde{G}(F)$ et invariante par conjugaison par $G(F)$. On dit que c'est un quasi-caract\`ere si, pour tout $\tilde{x}\in \tilde{G}_{ss}(F)$, il existe un bon voisinage $\omega$ de $0$ dans $\mathfrak{g}_{\tilde{x}}(F)$ et, pour tout ${\cal O}\in Nil(\mathfrak{g}_{\tilde{x}})$, il existe $c_{\Theta,{\cal O}}(\tilde{x})\in {\mathbb C}$ de sorte que l'on ait l'\'egalit\'e
$$\Theta(\tilde{x}exp(X))=\sum_{{\cal O}\in Nil(\mathfrak{g}_{\tilde{x}})}c_{\Theta,{\cal O}}(\tilde{x})\hat{j}({\cal O},X)$$
pour presque tout $X\in \omega$. Cette d\'efinition est la m\^eme qu'en [W2] 4.1. Dans cette r\'ef\'erence, on avait not\'e $\theta$ les quasi-caract\`eres. Pour une raison \'evidente, on croit bon ici de modifier cette notation.

\bigskip

\subsection{Fonctions tr\`es cuspidales}

Soit $\tilde{f}\in C_{c}^{\infty}(\tilde{G}(F))$. On dit que $\tilde{f}$ est tr\`es cuspidale si et seulement si, pour tout sous-groupe parabolique tordu propre $\tilde{P}=\tilde{M}U$ de $\tilde{G}$ et pour tout $\tilde{x}\in \tilde{M}(F)$, on a l'\'egalit\'e
$$\int_{U(F)}\tilde{f}(\tilde{x}u)du=0.$$

Les int\'egrales orbitales pond\'er\'ees des fonctions tr\`es cuspidales poss\`edent les m\^emes propri\'et\'es que dans le cas non tordu. En particulier, soient $\tilde{f}$ une fonction tr\`es cuspidale et $\tilde{x}\in \tilde{G}_{reg}(F)$. Notons $\tilde{M}(\tilde{x})$ le commutant de $A_{G_{\tilde{x}}}$ dans $\tilde{G}$, au sens de 1.2. C'est un  L\'evi tordu de $\tilde{G}$. Quitte \`a conjuguer $M_{min}$ et $K$, on peut supposer $ \tilde{M}(\tilde{x})\in {\cal L}^{\tilde{G}}$. On pose
$$\Theta^J_{\tilde{f}}(\tilde{x})=(-1)^{a_{\tilde{M}(\tilde{x})}-a_{\tilde{G}}}D^{\tilde{G}}(\tilde{x})^{-1/2}J_{\tilde{M}(\tilde{x})}(\tilde{x},\tilde{f}).$$
Cela ne d\'epend pas de la conjugaison effectu\'ee. La fonction $\Theta^J_{\tilde{f}}$ ainsi d\'efinie est invariante par conjugaison par $G(F)$.

\ass{Proposition}{Soit $\tilde{f}\in C_{c}^{\infty}(\tilde{G}(F))$ une fonction tr\`es cuspidale. Alors $\Theta^J_{\tilde{f}}$ est un quasi-caract\`ere.}

La preuve est exactement la m\^eme que dans le cas non tordu, cf. [W2] corollaire 5.9.

\bigskip

\subsection{Distributions locales associ\'ees \`a une fonction tr\`es cuspidale}

Soit $\tilde{x}\in \tilde{G}_{ss}(F)$. Supposons que $G_{\tilde{x}}$ soit le produit de deux groupes r\'eductifs connexes $G_{\tilde{x}}=G'\times G''$, chacun conserv\'e par $Z_{G}(\tilde{x})(F)$. Tout \'el\'ement $X\in \mathfrak{g}_{\tilde{x}}(F)$ se d\'ecompose en la somme d'un \'el\'ement de $\mathfrak{g}'(F)$ et d'un \'el\'ement de $\mathfrak{g}''(F)$. On note $X=X'+X''$ cette d\'ecomposition. On note $f\mapsto f^{\sharp}$ la transformation de Fourier partielle dans $C_{c}^{\infty}(\mathfrak{g}_{\tilde{x}}(F))$ relative \`a la deuxi\`eme variable, c'est-\`a-dire
$$f^{\sharp}(X)=\int_{\mathfrak{g''}(F)}f(X'+Y'')\psi(<Y'',X''>)dY''.$$
Soit $\omega\subset \mathfrak{g}_{\tilde{x}}(F)$ un bon voisinage de $0$. On suppose que $\omega=\omega'+\omega''$, avec $\omega'\subset \mathfrak{g}'(F)$ et $\omega''\subset \mathfrak{g}''(F)$.

Soit $f\in C_{c}^{\infty}(\tilde{G}(F))$ une fonction tr\`es cuspidale. On d\'efinit une fonction $\Theta^J_{\tilde{f},\tilde{x},\omega}$ sur $\mathfrak{g}_{\tilde{x},reg}(F)$ par
$$\Theta^J_{\tilde{f},\tilde{x},\omega}(X)=\left\lbrace\begin{array}{cc}0,&\text{ si }X\not\in \omega,\\ \Theta^J_{\tilde{f}}(\tilde{x}exp(X)),&\text{ si }X\in \omega.\\ \end{array}\right.$$
Pour $g\in G(F)$, on d\'efinit $^g\tilde{f}_{\tilde{x},\omega}\in C_{c}^{\infty}(\mathfrak{g}_{\tilde{x}}(F))$ par
$$^g\tilde{f}_{\tilde{x},\omega}(X)=\left\lbrace\begin{array}{cc}0,&\text{ si }X\not\in \omega,\\  \tilde{f}(g^{-1}\tilde{x}exp(X)g),&\text{ si }X\in \omega.\\ \end{array}\right.$$
On pose $^g\tilde{f}_{\tilde{x},\omega}^{\sharp}=(^g\tilde{f}_{\tilde{x},\omega})^{\sharp}$. Cette fonction v\'erifie la propri\'et\'e suivante:

(1) soit $\tilde{P}=\tilde{M}U$ un sous-groupe parabolique tordu tel que $\tilde{P}\not=\tilde{G}$ et $\tilde{x}\in \tilde{M}(F)$; alors
$$\int_{U(F)}{^{u}\tilde{f}}_{\tilde{x},\omega}^{\sharp}(X)du=0$$
pour tout $X\in \mathfrak{m}_{\tilde{x}}(F)$.

La preuve est la m\^eme que celle de [W2] lemme 5.5(i).

Soit $\tilde{M}$ un L\'evi de $\tilde{G}$  tel que $\tilde{x}\in \tilde{M}(F)$. Quitte \`a conjuguer $M_{min}$ et $K$, on peut supposer $  \tilde{M}\in {\cal L}^{\tilde{G}}$. On d\'efinit une fonction $J^{\sharp}_{\tilde{M},\tilde{x},\omega}(.,\tilde{f})$ sur $\mathfrak{m}_{\tilde{x}}(F)\cap \mathfrak{g}_{\tilde{x},reg}(F)$ par
$$J^{\sharp}_{\tilde{M},\tilde{x},\omega}(X,\tilde{f})=D^{G_{\tilde{x}}}(X)^{1/2}\int_{G_{\tilde{x}exp(X)}(F)\backslash G(F)}{^g\tilde{f}}^{\sharp}_{\tilde{x},\omega}(X)v_{\tilde{M}}(g)dg.$$
Cela ne d\'epend pas de la conjugaison effectu\'ee. On a

(2) Si $A_{\tilde{M}}\subsetneq A_{G_{\tilde{x}exp(X)}}$, on a $J^{\sharp}_{\tilde{M},\tilde{x},\omega}(X,\tilde{f})=0$.

La preuve est la m\^eme que celle de [W2] lemme 5.5(ii).

 Soit $X\in \mathfrak{g}_{\tilde{x},reg}(F)$. Notons $\tilde{{\bf M}}(X)$ le commutant de $A_{G_{\tilde{x}exp(X)}}$ dans $\tilde{G}$. C'est un L\'evi tordu de $\tilde{G}$ qui contient $\tilde{x}$. On pose
$$\Theta^{J,\sharp}_{\tilde{f},\tilde{x},\omega}(X)=(-1)^{a_{\tilde{{\bf M}}(X)}-a_{\tilde{G}}}D^{G_{\tilde{x}}}(X)^{-1/2}J^{\sharp}_{\tilde{{\bf M}}(X),\tilde{x},\omega}(X,\tilde{f}).$$
La fonction $\Theta^{J,\sharp}_{\tilde{f},\tilde{x},\omega}$ ainsi d\'efinie est invariante par conjugaison par $G_{\tilde{x}}(F)$.

On d\'efinit la notion de quasi-caract\`ere sur une alg\`ebre de Lie de m\^eme que l'on a d\'efini cette notion sur un groupe ou un groupe tordu.

\ass{Proposition}{Les fonctions $\Theta^J_{\tilde{f},\tilde{x},\omega}$ et $\Theta^{J,\sharp}_{\tilde{f},\tilde{x},\omega}$ sont des quasi-caract\`eres sur $\mathfrak{g}_{\tilde{x}}(F)$. La seconde est la transform\'ee de Fourier partielle de la premi\`ere, c'est-\`a-dire que, pour toute $\varphi\in C_{c}^{\infty}(\mathfrak{g}_{\tilde{x}}(F))$,on a l'\'egalit\'e
$$\int_{\mathfrak{g}_{\tilde{x}}(F)}\Theta^{J,\sharp}_{\tilde{f},\tilde{x},\omega}(X)\varphi(X)dX=\int_{\mathfrak{g}_{\tilde{x}}(F)}\Theta^J_{\tilde{f},\tilde{x},\omega}(X)\varphi^{\sharp}(X)dX.$$}

 La premi\`ere assertion concernant la fonction $\Theta^J_{\tilde{f},\tilde{x},\omega}$ r\'esulte de la proposition du paragraphe pr\'ec\'edent. La deuxi\`eme assertion se d\'emontre comme dans le cas non tordu, cf. [W2] proposition 5.8. Puisque $\Theta^J_{\tilde{f},\tilde{x},\omega}$ est un quasi-caract\`ere \`a support compact modulo conjugaison, il existe une fonction $\varphi\in C_{c}^{\infty}(\mathfrak{g}_{\tilde{x}}(F))$, tr\`es cuspidale, telle que $\Theta^J_{\tilde{f},\tilde{x},\omega}$ soit \'egal au quasi-caract\`ere $\Theta^J_{\varphi}$ associ\'e \`a $\varphi$, cf. [W2] proposition 6.4. Le lemme 6.1 de [W2] se g\'en\'eralise imm\'ediatement aux transform\'ees de Fourier partielles: la transform\'ee de Fourier partielle de $\Theta^J_{\varphi}$ est $\Theta^J_{\varphi^{\sharp}}$, et c'est un quasi-caract\`ere. Puisque cette transform\'ee de Fourier est $\Theta^{J,\sharp}_{\tilde{f},\tilde{x},\omega}$, cela entra\^{\i}ne la premi\`ere assertion concernant cette fonction. $\square$
 
 \bigskip
 
 \subsection{Repr\'esentations de groupes tordus}
 
 On appelle repr\'esentation de $\tilde{G}(F)$ un triplet $(\pi,\tilde{\pi},E)$, o\`u $E$ est un espace vectoriel complexe, $\pi$ une repr\'esentation de $G(F)$ dans $E$ et $\tilde{\pi}$ une application de $\tilde{G}(F)$ dans le groupe $Aut(E)$ des automorphismes ${\mathbb C}$-lin\'eaires de $E$ telle que $\pi(g)\tilde{\pi}(\tilde{x})\pi(g')=\tilde{\pi}(gxg')$ pour tous $g,g'\in G(F)$ et $\tilde{x}\in \tilde{G}(F)$. Deux repr\'esentations $(\pi_{1},\tilde{\pi}_{1},E_{1})$ et $(\pi_{2},\tilde{\pi}_{2},E_{2})$ sont dites \'equivalentes s'il existe un isomorphisme ${\mathbb C}$-lin\'eaire $A:E_{1}\to E_{2}$ et un \'el\'ement $B\in Aut(E_{1})$, commutant \`a la repr\'esentation $\pi_{1}$, de sorte que $\pi_{1}(g)=A^{-1}\pi_{2}(g)A$ pour tout $g\in G(F)$ et $\tilde{\pi}_{1}(\tilde{x})=BA^{-1}\tilde{\pi}_{2}(\tilde{x})A$ pour tout $\tilde{x}\in \tilde{G}(F)$.  
 
 Soit $(\pi,\tilde{\pi},E)$ une repr\'esentation de $\tilde{G}(F)$. On dit qu'elle est lisse, ou admissible, si $\pi$ l'est. On dit  qu'elle est unitaire s'il existe un produit hermitien d\'efini positif sur $E$ tel que $\tilde{\pi}$ prenne ses valeurs dans le groupe unitaire pour ce produit. Nous dirons qu'elle est temp\'er\'ee si elle est unitaire, si $\pi$ est de longueur finie et si toutes les composantes irr\'eductibles de $\pi$ sont temp\'er\'ees. On d\'efinit la contragr\'ediente $(\pi^{\vee},\tilde{\pi}^{\vee},E^{\vee})$: $(\pi^{\vee},E^{\vee})$ est la contragr\'ediente de $(\pi,E)$ et $\tilde{\pi}^{\vee}$ est d\'efinie par $<\tilde{\pi}^{\vee}(\tilde{x})\check{e},e>=<\check{e},\tilde{\pi}(\tilde{x})^{-1}e>$. Pour $\lambda\in {\cal A}_{\tilde{G}}^*\otimes_{{\mathbb R}}{\mathbb C}$, on d\'efinit $(\pi_{\lambda},\tilde{\pi}_{\lambda},E)$ par $\tilde{\pi}_{\lambda}(\tilde{x})=e^{\lambda(H_{\tilde{G}}(\tilde{x}))}\tilde{\pi}(\tilde{x})$. On dit que la repr\'esentation est $G(F)$-irr\'eductible si $\pi$ est irr\'eductible. Fixons un point quelconque $\tilde{x}\in \tilde{G}(F)$. La donn\'ee d'une repr\'esentation $G(F)$-irr\'eductible est  simplement  la donn\'ee de deux objets
 
 - une repr\'esentation irr\'eductible $(\pi,E)$ de $G(F)$ telle que la repr\'esentation $g\mapsto \pi(\theta_{\tilde{x}}(g))$ soit \'equivalente \`a $\pi$;
 
 - un op\'erateur $\tilde{\pi}(\tilde{x})$ de $E$ tel que $\pi(\theta_{\tilde{x}}(g))=\tilde{\pi}(\tilde{x})\pi(g)\tilde{\pi}(\tilde{x})^{-1}$ pour tout $g\in G(F)$.
 
  En pratique, nous noterons simplement $\tilde{\pi}$ la repr\'esentation $(\pi,\tilde{\pi},E)$. Nous noterons sans plus de commentaire $\pi$ la repr\'esentation de $G(F)$ associ\'ee et $E_{\tilde{\pi}}$ l'espace $E$. On note $Temp(\tilde{G})$ l'ensemble des classes d'isomorphie de repr\'esentations admissibles irr\'eductibles et temp\'er\'ees de $\tilde{G}(F)$, que nous identifions \`a un ensemble de repr\'esentants de ces classes.  
  
  Soient $\tilde{P}=\tilde{M}U$ un sous-groupe parabolique tordu de $\tilde{G}$ et $\tilde{\tau}$ une repr\'esentation admissible de $\tilde{M}(F)$. On d\'efinit la repr\'esentation induite $Ind_{P}^{G}(\tilde{\tau})$ de $\tilde{G}(F)$. Elle se r\'ealise dans l'espace habituel de la repr\'esentation $Ind_{P}^G(\tau)$, que l'on note $E_{P,\tau}^G$. Soient $e\in E_{P,\tau}^G$ et $\tilde{x}\in \tilde{G}(F)$. Fixons $\tilde{y}\in \tilde{M}(F)$, soit $\gamma$ l'\'el\'ement de $G(F)$ tel que $\tilde{x}=\gamma\tilde{y}$. Alors on a l'\'egalit\'e
  $$(Ind_{P}^G(\tilde{\tau},\tilde{x})e)(g)=\delta_{P}(\tilde{y})^{1/2}\tilde{\tau}(\tilde{y})e(\theta_{\tilde{y}}^{-1}(g\gamma))$$
  pour tout $g\in G(F)$, o\`u $\delta_{P}$ est \'etendu de fa\c{c}on naturelle \`a $\tilde{M}(F)$. Supposons que $M$ contienne $M_{min}$. On peut aussi  r\'ealiser la repr\'esentation induite dans l'espace ${\cal K}_{P,\tau}^G$ des restrictions \`a $K$ des \'el\'ements de $E_{P,\tau}^G$. Supposons $\tilde{\tau}$ temp\'er\'ee. On d\'efinit une $(\tilde{G},\tilde{M})$-famille $({\cal R}_{\tilde{P}'}(\tau))_{\tilde{P}'\in {\cal P}(\tilde{M})}$, qui prend ses valeurs dans l'espace des op\'erateurs de ${\cal K}_{P,\tau}^G$, par
  $${\cal R}_{\tilde{P}'}(\tau,\lambda)=R_{P'\vert P}(\tau)^{-1}R_{P'\vert P}(\tau_{\lambda})$$
  pour $\lambda\in i{\cal A}_{\tilde{M}}^*$, o\`u $R_{P'\vert P}(\tau)$ est l'op\'erateur d'entrelacement normalis\'e, cf. [A3] th\'eor\`eme 2.1. On associe \`a cette $(\tilde{G},\tilde{M})$-famille un op\'erateur ${\cal R}_{\tilde{M}}(\tilde{\tau})$. On d\'efinit le caract\`ere pond\'er\'e de $\tilde{\tau}$: c'est la distribution
  $$\tilde{f}\mapsto J_{\tilde{M}}^{\tilde{G}}(\tilde{\tau},\tilde{f})=trace({\cal R}_{\tilde{M}}(\tau)Ind_{P}^{G}(\tilde{f}))$$
  sur $C_{c}^{\infty}(\tilde{G}(F))$. Dans le cas o\`u $\tilde{M}=\tilde{G}$, c'est le caract\`ere habituel de $\tilde{\tau}$ et on le note plut\^ot $\tilde{f}\mapsto \Theta_{\tilde{\tau}}(\tilde{f})$. D'apr\`es [C] theorem 2, ce caract\`ere est une distribution localement int\'egrable.

  \bigskip
  \subsection{Int\'egrales orbitales pond\'er\'ees invariantes}
  
  En utilisant les caract\`eres pond\'er\'es de la section pr\'ec\'edente, le proc\'ed\'e habituel d'Arthur permet de  construire des int\'egrales orbitales pond\'er\'ees invariantes \`a l'aide des int\'egrales $\tilde{f}\mapsto J_{\tilde{M}}(\tilde{x},\tilde{f})$. Rappelons la construction. Pour $\zeta\in {\cal A}_{\tilde{G}}$, notons ${\bf 1}_{H_{\tilde{G}}=\zeta}$ la fonction caract\'eristique de l'ensemble des $ \tilde{x}\in \tilde{G}(F)$ tels que $H_{\tilde{G}}(\tilde{x})=\zeta$. Notons ${\cal H}_{ac}(\tilde{G}(F))$ l'espace des fonctions $\tilde{f}:\tilde{G}(F)\to {\mathbb C}$ telles que
  
  (1) $\tilde{f}$ est biinvariante par un sous-groupe ouvert compact de $G(F)$;
  
  (2) pour tout   tout $\zeta\in {\cal A}_{\tilde{G}}$, la fonction ${\bf 1}_{H_{\tilde{G}}=\zeta}\tilde{f} $ est \`a support compact sur $\tilde{G}(F)$.
  
  Pour $\tilde{f}\in C_{c}^{\infty}(\tilde{G}(F))$ et $\tilde{M}\in {\cal L}^{\tilde{G}}$, Arthur montre qu'il existe une fonction $\phi_{\tilde{M}}(\tilde{f})\in {\cal H}_{ac}(\tilde{M}(F))$ telle que, pour toute repr\'esentation $\tilde{\pi}\in Temp(\tilde{M})$ et tout $\zeta\in H_{\tilde{M}}(\tilde{M}(F))\subset {\cal A}_{\tilde{M}}$, on ait l'\'egalit\'e
  $$\int_{i{\cal A}_{\tilde{M},F}}J_{\tilde{L}}(\tilde{\pi}_{\lambda},\tilde{f})exp(-\lambda(\zeta))d\lambda=\int_{i{\cal A}_{\tilde{M},F}}\Theta_{\tilde{\pi}_{\lambda}}(\phi_{\tilde{M}}(\tilde{f}){\bf 1}_{H_{\tilde{M}}=\zeta})exp(-\lambda(\zeta))d\lambda$$
  $$=mes(i{\cal A}_{\tilde{M},F})\Theta_{\tilde{\pi}}(\phi_{\tilde{M}}(\tilde{f}){\bf 1}_{H_{\tilde{M}}=\zeta}).$$
  On fixe de telles fonctions $\phi_{\tilde{M}}(\tilde{f})$. Pour $\tilde{x}\in \tilde{M}(F)\cap \tilde{G}_{reg}(F)$, on d\'efinit l'int\'egrale orbitale invariante $I_{\tilde{M}}(\tilde{x},\tilde{f})$ par r\'ecurrence sur $a_{\tilde{M}}-a_{\tilde{G}}$ par la formule
  $$I_{\tilde{M}}(\tilde{x},\tilde{f})=J_{\tilde{M}}(\tilde{x},\tilde{f})-\sum_{\tilde{L}\in {\cal L}(\tilde{M}), \tilde{L}\not=\tilde{G}}I_{\tilde{M}}^{\tilde{L}}(\tilde{x},\phi_{\tilde{L}}(\tilde{f}){\bf 1}_{H_{\tilde{L}}=H_{\tilde{L}}(\tilde{x})}).$$

    Soit $\tilde{f}\in C_{c}^{\infty}(\tilde{G}(F))$. On lui associe une fonction $\Theta_{\tilde{f}}$ sur $\tilde{G}_{reg}(F)$ de la fa\c{c}on suivante. Soit $\tilde{x}\in \tilde{G}_{reg}(F)$. Notons $\tilde{M}(\tilde{x})$ le commutant de $A_{G_{\tilde{x}}}$ dans $\tilde{G}$. Choisissons $g\in G(F)$ tel que $g\tilde{M}(\tilde{x})g^{-1}\in {\cal L}^{\tilde{G}}$. On pose
  $$\Theta_{\tilde{f}}(\tilde{x})=(-1)^{a_{\tilde{M}(\tilde{x})}-a_{\tilde{G}}}D^{\tilde{G}}(\tilde{x})^{-1/2}I_{g\tilde{M}(\tilde{x})g^{-1}}(g\tilde{x}g^{-1},\tilde{f}).$$
 Cela ne d\'epend pas du choix de $g$.
  
Pour $\tilde{f}\in C_{c}^{\infty}(\tilde{G}(F))$, on dit que $\tilde{f}$ est cuspidale si et seulement si, pour tout groupe de L\'evi tordu $\tilde{M}\subsetneq \tilde{G}$ et pour tout $\tilde{x}\in \tilde{G}_{reg}(F)\cap \tilde{M}(F)$, on a $J_{\tilde{G}}(\tilde{x},\tilde{f})=0$.  On a la propri\'et\'e suivante, que nous \'enon\c{c}ons avant de l'expliquer:

(3) pour toute fonction cuspidale $\tilde{f}\in C_{c}^{\infty}(\tilde{G}(F))$, pour tout $\tilde{M}\in{\cal L}^{\tilde{G}}$ et pour tout \'el\'ement $\tilde{x}\in \tilde{M}(F)\cap \tilde{G}_{reg}(F)$ qui est elliptique dans $\tilde{M}(F)$, on a l'\'egalit\'e
$$D^{\tilde{G}}(\tilde{x})^{-1/2}(-1)^{a_{\tilde{M}}-a_{\tilde{G}}}I_{\tilde{M}}(\tilde{x},\tilde{f})=\sum_{{\cal O}\in \{\Pi_{ell}(\tilde{G})\}}c({\cal O})\int_{i{\cal A}_{\tilde{G},F}^*}\Theta_{\tilde{\pi}_{\lambda}}(\tilde{x})\Theta_{(\tilde{\pi}_{\lambda})^{\vee}}(\tilde{f})d\lambda.$$

Que $\tilde{x}$ soit elliptique dans $\tilde{M}(F)$ signifie que $A_{G_{\tilde{x}}}=A_{\tilde{M}}$. L'ensemble $\Pi_{ell}(\tilde{G})$ est un certain sous-ensemble de celui des repr\'esentations temp\'er\'ees virtuelles de $\tilde{G}(F)$. Il est stable par l'action $\tilde{\pi}\mapsto \tilde{\pi}_{\lambda}$ de $i{\cal A}_{\tilde{G}}^*$. L'ensemble $\{\Pi_{ell}(\tilde{G})\}$ est l'ensemble des orbites de cette action. Pour chaque orbite ${\cal O}$, on choisit un point base $\tilde{\pi}\in {\cal O}$. Enfin, $c({\cal O})$ est un certain coefficient. La somme est en fait finie car, pour tout sous-groupe ouvert compact $K'$ de $G(F)$, il n'y a qu'un nombre fini d'orbites pour lesquelles une repr\'esentation de l'orbite admet des invariants non nuls par $K'$.  

Cf. [W6] th\'eor\`eme 7.1 et, dans le cas non tordu,  [A4] th\'eor\`eme 5.1.

\ass{Lemme}{Soit $\tilde{f}\in C_{c}^{\infty}(\tilde{G}(F))$ une fonction cuspidale. Alors $\Theta_{\tilde{f}}$ est un quasi-caract\`ere de $\tilde{G}(F)$.}

La preuve est la m\^eme qu'en [W3] 2.5, en utilisant [C] theorem 3.

\bigskip

\subsection{Un lemme d'annulation}

Soient $\pi$ une repr\'esentation admissible de $G(F)$ et $B$ une forme bilin\'eaire sur $E_{\pi^{\vee}}\times E_{\pi}$. Soit $A\in End_{{\mathbb C}}(E_{\pi})$ un op\'erateur lisse, c'est-\`a-dire qu'il existe un sous-groupe ouvert compact $K'$ de $G(F)$ tel que $\pi(k)A\pi(k')=A$ pour tous $k,k'\in K'$. Pour un tel sous-groupe $K'$, fixons une base ${\cal B}^{K'}$ du sous-espace $E_{\pi}^{K'}$ et introduisons la base duale $\{\check{e}; e\in {\cal B}^{K'}\}$ de $E_{\pi^{\vee}}^{K'}$ (pour l'accouplement usuel, not\'e $<.,.>$, entre ces deux espaces). Posons
$$trace_{B}(A)=\sum_{e\in {\cal B}^{K'}}B(\check{e},A(e)).$$
Ce terme ne d\'epend ni du choix de $K'$, ni de celui de la base.
  
  Soient $\tilde{P}=\tilde{M}U$ un sous-groupe parabolique tordu et $\tau$ une repr\'esentation admissible de $M(F)$. On suppose
  
  (1) $\tilde{P}\not=\tilde{G}$.
  
  Supposons $\pi=Ind_{P}^G(\tau)$, $E_{\pi}=E_{P,\tau}^G$ et $E_{\pi^{\vee}}=E_{P,\tau^{\vee}}^G$. Pour un \'el\'ement $\tilde{y}\in \tilde{M}(F)$, consid\'erons la condition suivante:
  
  $(H)_{\tilde{y}}$ soient $e\in E_{P,\tau}^G$ et $e'\in E_{P,\tau^{\vee}}^G$ tels que $e'(\theta_{\tilde{y}}(g))\otimes e(g)=0$ pour tout $g\in G(F)$; alors $B(e',e)=0$.
  
  Consid\'erons aussi la condition:
  
  $(H)$ la repr\'esentation $\tau$ se prolonge en une repr\'esentation $\tilde{\tau}$ de $\tilde{M}(F)$; soient $e\in E_{P,\tau}^G$ et $e'\in E_{P,\tau^{\vee}}^G$ tels que $e'(g)\otimes e(g)=0$ pour tous $g\in G(F)$; alors $B(e',e)=0$.
  
  Soit $\tilde{f}\in C_{c}^{\infty}(\tilde{G}(F))$. Supposons
  
  (2) $\tilde{f}$ est tr\`es cuspidale.
  
   Pour un \'el\'ement $\tilde{x}\in \tilde{G}(F)$, on note $\tilde{f}_{\tilde{x}}$ la fonction sur $G(F)$ d\'efinie par $\tilde{f}_{\tilde{x}}(g)=\tilde{f}(g\tilde{x})$.
  
  \ass{Lemme}{On suppose v\'erifi\'ees les conditions (1) et (2).
  
  (i) Soit $\tilde{y}\in \tilde{M}(F)$, supposons  v\'erifi\'ee la condition $(H)_{\tilde{y}}$. Alors $trace_{B}(Ind_{P}^G(\tau,\tilde{f}_{\tilde{y}}))=0$.
  
  (ii) Supposons v\'erifi\'ee la condition $(H)$. Alors $trace_{B}(Ind_{P}^G(\tilde{\tau},\tilde{f}))=0$.}
  
  Preuve. On ne perd rien \`a supposer que $K$ est en bonne position relativement \`a $M$. Consid\'erons la situation de (i). Fixons un sous-groupe ouvert compact $K'$ de $K$ tel que $\tilde{f}$ soit biinvariante par $K'$ et par $K''=\theta_{\tilde{y}}(K')$.   On fixe un ensemble de repr\'esentants $\Gamma'$ de l'ensemble des doubles classes $P(F)\backslash G(F)/K'$. L'ensemble $\Gamma''=\theta_{\tilde{y}}(\Gamma')$ est un ensemble de repr\'esentants de l'ensemble des doubles classes $P(F)\backslash G(F)/K''$. Choisissons une base ${\cal B}^{K'}$ de $(E_{P,\tau}^G)^{K'}$ telle que, pour tout $e'\in {\cal B}^{K'}$, il existe $\gamma'\in \Gamma'$ de sorte que le support de $e$ soit inclus dans $P(F)\gamma'K'$. L'\'el\'ement correspondant $\check{e}'$ de la base duale v\'erifie la m\^eme propri\'et\'e, avec le m\^eme $\gamma'$. Choisissons de m\^eme une base ${\cal B}^{K''}$ de $(E_{P,\tau}^G)^{K''}$. Pour tout $e'\in {\cal B}^{K''}$, on a l'\'egalit\'e
  $$Ind_{P}^G(\tau,\tilde{f}_{\tilde{y}})e''=\sum_{e'\in {\cal B}^{K'}}<\check{e}',Ind_{P}^G(\tau,\tilde{f}_{\tilde{y}})e''>e',$$
  d'o\`u
  $$(3) \qquad trace_{B}(Ind_{P}^G(\tau,\tilde{f}_{\tilde{y}}))=\sum_{e'\in {\cal B}^{K'}, e''\in {\cal B}^{K''}}B(\check{e}'',e')<\check{e}',Ind_{P}^G(\tau,\tilde{f}_{\tilde{y}})e''>.$$
  Fixons $e'\in {\cal B}^{K'}$ et $e''\in {\cal B}^{K''}$. Soient $\gamma'\in \Gamma'$ et $\gamma''\in \Gamma''$ tels que le support de $e'$ soit contenu dans $P(F)\gamma'K'$ et celui de $e''$ soit contenu dans $P(F)\gamma''K''$. Si $\gamma''\not=\theta_{\tilde{y}}(\gamma')$, on a $\check{e}''(\theta_{\tilde{y}}(g))\otimes e'(g)=0$ pour tout $g\in G(F)$. Donc $B(\check{e}'',e')=0$ d'apr\`es $(H)_{\tilde{y}}$. Supposons $\gamma''=\theta_{\tilde{y}}(\gamma')$. On calcule
  $$<\check{e}',Ind_{P}^G(\tau,\tilde{f}_{\tilde{y}})e''>=\int_{K}\int_{G(F)}\tilde{f}_{\tilde{y}}(g)<\check{e}'(h),e''(hg)>dg\,dh$$
  $$\qquad =\int_{K}\int_{G(F)}\tilde{f}_{\tilde{y}}(h^{-1}g)<\check{e}'(h),e''(g)>dg\,dh$$
  $$\qquad =\int_{K}\int_{K}\int_{M(F)}\int_{U(F)}\tilde{f}(h^{-1}muk\tilde{y})<\check{e'}(h),\tau(m)e''(k)>\delta_{P}(m)^{1/2}du\,dm\,dk\,dh.$$
  Par les changements de variables $k\mapsto \theta_{\tilde{y}}(k)$, $u\mapsto \theta_{\tilde{y}}(u)$, on obtient
  $$<\check{e}',Ind_{P}^G(\tau,\tilde{f}_{\tilde{y}})e''> =\int_{K}\int_{\theta_{\tilde{y}}^{-1}(K)}\int_{M(F)}\int_{U(F)}\tilde{f}(h^{-1}m\tilde{y}uk)$$
  $$<\check{e'}(h),\tau(m)e''(\theta_{\tilde{y}}(k))>\delta_{P}(m\tilde{y})^{1/2}du\,dm\,dk\,dh.$$
  Fixons $h,k,m$ et supposons $<\check{e'}(h),\tau(m)e''(\theta_{\tilde{y}}(k))>\not=0$. Alors $h\in P(F)\gamma'K'$ et $\theta_{\tilde{y}}(k)\in P(F)\gamma''K''$. Cette derni\`ere condition entra\^{\i}ne $k\in P(F)\gamma'K'$. Mais alors $k\in P(F)hK'$. Ecrivons $k=m'u'hk'$, avec $m'\in M(F)$, $u'\in U(F)$, $k'\in K'$. Consid\'erons l'int\'egrale int\'erieure de la formule ci-dessus. Puisque $\tilde{f}$ est invariante \`a droite par $K'$, le $k'$ dispara\^{\i}t. Par le changement de variables $u\mapsto m'uu^{_{'}-1}m^{_{'}-1}$, cette int\'egrale devient
  $$\delta_{P'}(m')^{1/2}\int_{U(F)}\tilde{f}(h^{-1}m\tilde{y}m'uh)du.$$
  Elle est nulle puisque $\tilde{f}$ est tr\`es cuspidale. Donc $<\check{e}',Ind_{P}^G(\tau,\tilde{f}_{\tilde{y}})e''>=0$. Tous les termes de (3) sont donc nuls, ce qui prouve l'assertion (i) de l'\'enonc\'e.
  
  Consid\'erons la situation de (ii). Fixons $\tilde{y}\in \tilde{M}(F)$. On a l'\'egalit\'e $Ind_{P}^G(\tilde{\tau},\tilde{f})=Ind_{P}^G(\tau,\tilde{f}_{\tilde{y}})Ind_{P}^G(\tilde{\tau},\tilde{y})$. D\'efinissons une forme bilin\'eaire $B'$ sur $E_{P,\tau^{\vee}}^G\times E_{P,\tau}^G$ par
  $$B'(e',e)=B(Ind_{P}^G(\tilde{\tau},\tilde{y})^{-1}e',e).$$
  On v\'erifie formellement que $trace_{B}(Ind_{P}^G(\tilde{\tau},\tilde{f}))=trace_{B'}(Ind_{P}^G(\tau,\tilde{f}_{\tilde{y}}))$ et que la condition $(H)$ pour $B$ entra\^{\i}ne la condition $(H)_{\tilde{y}}$ pour $B'$. Il ne reste plus qu'\`a appliquer le (i) \`a la forme $B'$ pour conclure. $\square$
  
  {\bf Remarque.} Quand $\pi$ est unitaire, on peut consid\'erer des formes sesquilin\'eaires sur $E_{\pi}\times E_{\pi}$ plut\^ot que des formes bilin\'eaires sur $E_{\pi^{\vee}}\times E_{\pi}$.
  
  \bigskip
  
  \subsection{Induction de quasi-caract\`eres}
  
  Soit $\tilde{M}$ un L\'evi tordu de $\tilde{G}$. Soit $\Delta^{\tilde{M}}$ une distribution sur $\tilde{M}(F)$ invariante par conjugaison par $M(F)$. On peut d\'efinir une distribution induite $\Delta=Ind_{M}^G(\Delta^{\tilde{M}})$ sur $\tilde{G}(F)$, qui est invariante par conjugaison. La d\'efinition est similaire \`a celle du cas non tordu. Quitte \`a conjuguer $\tilde{M}$, on suppose $\tilde{M}\in {\cal L}^{\tilde{G}}$. On fixe $\tilde{P}=\tilde{M}U\in {\cal P}(\tilde{M})$. Pour $\tilde{f}\in C_{c}^{\infty}(\tilde{G}(F))$, on d\'efinit $\tilde{f}_{P}\in C_{c}^{\infty}(\tilde{M}(F))$ par
  $$\tilde{f}_{P}(\tilde{m})=\delta_{P}(\tilde{m})^{1/2}\int_{K}\int_{U(F)}\tilde{f}(k^{-1}\tilde{m}uk)du\,dk.$$
  On pose $\Delta(\tilde{f})=\Delta^{\tilde{M}}(\tilde{f}_{P})$. Cela ne d\'epend pas des choix effectu\'es. Dans le cas o\`u $\Delta^{\tilde{M}}$ est un quasi-caract\`ere sur $\tilde{M}(F)$, $\Delta$ est aussi un quasi-caract\`ere sur $\tilde{G}(F)$. Le d\'eveloppement de $\Delta$ au voisinage d'un point semi-simple de $\tilde{G}(F)$ se calcule en fonction des d\'eveloppements de $\Delta^{\tilde{M}}$. Enon\c{c}ons le r\'esultat qui nous int\'eresse, qui est le m\^eme que dans le cas non tordu ([W3] lemme 2.3).
  
  \ass{Lemme}{Soient $\Theta^{\tilde{M}}$ un quasi-caract\`ere de $\tilde{M}(F)$ et $\Theta=Ind_{M}^G(\Theta^{\tilde{M}})$. Alors
  
  (i) $\Theta$ est un quasi-caract\`ere sur $\tilde{G}(F)$;
  
  (ii) soient $\tilde{x}$ un \'el\'ement semi-simple de $\tilde{G}(F)$ et ${\cal O}\in Nil(\mathfrak{g}_{\tilde{x}})$ une orbite r\'eguli\`ere; on a l'\'egalit\'e
  $$c_{\Theta,{\cal O}}(\tilde{x})=\sum_{\tilde{x}'\in {\cal X}^{\tilde{M}}(\tilde{x})}\sum_{g\in \Gamma_{\tilde{x}'}/G_{\tilde{x}}(F)}\sum_{{\cal O}'\in Nil(\mathfrak{m}_{\tilde{x}'})}D^{\tilde{G}}(\tilde{x})^{-1/2}D^{\tilde{M}}(\tilde{x}')^{1/2}$$
  $$[Z_{M}(\tilde{x}')(F):M_{\tilde{x}'}(F)]^{-1}[g{\cal O}:{\cal O}']c_{\Theta^{\tilde{M}},{\cal O}'}(\tilde{x}').$$}
  
  Expliquons les notations. On a fix\'e un ensemble ${\cal X}^{\tilde{M}}(\tilde{x})$ de repr\'esentants des classes de conjugaison par $M(F)$ dans l'intersection de $\tilde{M}(F)$ avec la classe de conjugaison de $\tilde{x}$ par $G(F)$. Pour $\tilde{x}'\in {\cal X}^{\tilde{M}}(\tilde{x})$, $\Gamma_{\tilde{x}'}$ est l'ensemble des $g\in G(F)$ tels que $g\tilde{x}g^{-1}=\tilde{x}'$. Pour un tel $g$, la conjugaison par $g$ transporte ${\cal O}$ en une orbite $g{\cal O}$ dans $\mathfrak{g}_{\tilde{x}'}(F)$. Pour une orbite nilpotente ${\cal O}'$ de $\mathfrak{m}_{\tilde{x}'}(F)$, on pose $[g{\cal O}:{\cal O}']=1$ si $g{\cal O}$ est incluse dans l'orbite induite de ${\cal O}'$, $[g{\cal O}:{\cal O}']=0$ sinon.
  
  \bigskip
  
  \subsection{ Propri\'et\'es des fonctions tr\`es cuspidales}

  En utilisant le lemme du paragraphe 1.11, on peut adapter les preuves des lemmes 2.2, 2.6 et 2.7 de [W3].   On obtient les r\'esultats suivants. 
  
  \ass{Lemme}{(i) Soit $\tilde{f}\in C_{c}^{\infty}(\tilde{G}(F))$ une fonction tr\`es cuspidale, soient $\tilde{M}\in {\cal L}^{\tilde{G}}$, $\tilde{\tau}$ une repr\'esentation temp\'er\'ee de $\tilde{M}(F)$, $\tilde{L}\in {\cal L}(\tilde{M})$ et $\tilde{Q}\in {\cal F}(\tilde{L})$. Supposons $\tilde{L}\not=\tilde{M}$ ou $\tilde{Q}\not=\tilde{G}$. Alors $J_{\tilde{L}}^{\tilde{Q}}(\tilde{\tau},\tilde{f})=0$.
  
  (ii) Soit $\tilde{f}\in C_{c}^{\infty}(\tilde{G}(F))$ une fonction tr\`es cuspidale. Alors, pour tout $\tilde{L}\in {\cal L}^{\tilde{G}}$, la fonction $\phi_{\tilde{L}}(\tilde{f})$ est cuspidale. Admettons l'hypoth\`ese de 1.10. Alors on a l'\'egalit\'e
  $$\Theta^J_{\tilde{f}}=\sum_{\tilde{L}\in {\cal L}^{\tilde{G}}}\vert W^L\vert \vert W^G\vert ^{-1}(-1)^{a_{\tilde{L}}-a_{\tilde{G}}}Ind_{L}^G(\Theta_{\phi_{\tilde{L}}(\tilde{f})}).$$
  
(iii)  Soit $\tilde{f}\in C_{c}^{\infty}(\tilde{G}(F))$ une fonction cuspidale. Alors il existe une fonction tr\`es cuspidale $\tilde{f}'\in C_{c}^{\infty}(\tilde{G}(F))$ telle que $J_{\tilde{G}}(\tilde{x},\tilde{f}')=J_{\tilde{G}}(\tilde{x},\tilde{f})$ pour tout $\tilde{x}\in \tilde{G}_{reg}(F)$.}
  
  \bigskip

  \section{Le groupe $GL_{d}$ tordu}
  
  \bigskip
  
  \subsection{Description du groupe tordu}
  
  Soient $d\geq1$ un entier et $V$ un espace vectoriel sur $F$ de dimension $d$. Notons $V^*$ le dual de $V$, $G=GL(V)$ le groupe des automorphismes $F$-lin\'eaires de $V$ et $\tilde{G}=Isom(V,V^*)$ l'ensemble des isomorphismes $F$-lin\'eaires de $V$ sur $V^*$. Le groupe $G$ agit \`a droite et \`a gauche sur $\tilde{G}$ par
  $$(g,\tilde{x},g')\mapsto {^tg}^{-1}\circ \tilde{x}\circ g',$$
  o\`u $^tg$ est le transpos\'e de $g$. Le couple $(G,\tilde{G})$ est un groupe tordu.  Remarquons que $\tilde{G}(F)$ s'identifie \`a l'ensemble des formes bilin\'eaires non d\'eg\'en\'er\'ees sur $V$: \`a $\tilde{x}\in \tilde{G}(F)$, on associe la forme $(v',v)\mapsto <v',\tilde{x}v>$.
  
  Soit $\tilde{x}\in \tilde{G}(F)$ un \'el\'ement semi-simple. Son commutant dans $G$ est produit d'un groupe symplectique, d'un groupe orthogonal et de groupes qui sont des restrictions \`a la Weil de groupes lin\'eaires ou unitaires. La composante  orthogonale  se construit de la fa\c{c}on suivante. Consid\'erons l'\'el\'ement $x={^t\tilde{x}}^{-1}\tilde{x}$ de $G(F)$.  Il est semi-simple. Notons $V''_{\tilde{x}}$ le sous-espace de $V$ propre pour l'action de $x$, associ\'e \`a la valeur propre $1$. Notons $V'_{\tilde{x} }$ l'unique suppl\'em\'entaire de $ V''_{\tilde{x}}$ qui soit invariant par $x$. Alors $V=  V'_{\tilde{x}}\oplus V''_{\tilde{x}} $ et $V^*=V_{\tilde{x}}^{_{'}*}\oplus V_{\tilde{x}}^{_{''}*}$. L'\'el\'ement $\tilde{x}$ envoie $V'_{\tilde{x}}$ dans $V_{\tilde{x}}^{_{'}*}$ et $V''_{\tilde{x}}$ dans $V_{\tilde{x}}^{_{''}*}$. Sa restriction \`a $V''_{\tilde{x}}$ est une forme  quadratique. La composante orthogonale de $Z_{G}(\tilde{x})$ est le groupe  orthogonal de cette restriction.  
    
  Fixons une base $(v_{i})_{i=1,...,d}$ de $V$. On peut alors identifier $G$ au groupe $GL_{d}$ des matrices $d\times d$ inversibles. Pour $g\in GL_{d}$, on note ses coefficients $g_{i,j}$, pour $i,j=1,...,d$. On introduit les sous-groupes suivants: $B_{d}$ le  sous-groupe de Borel triangulaire sup\'erieur, $A_{d}$ le sous-tore diagonal, $U_{d}$ le radical unipotent de $B_{d}$, $K_{d}$ le sous-groupe compact sp\'ecial de $GL_{d}(F)$ form\'e des matrices \`a coefficients entiers et de d\'eterminant  de valuation nulle. On note $\xi$ le caract\`ere du groupe $U_{d}(F)$ d\'efini par
  $$\xi(u)=\psi(\sum_{i=1,...,d-1}u_{i,i+1}).$$
  Introduisons la base duale $(v_{i}^*)_{i=1,...,d}$ de $V^*$. Notons $\boldsymbol{\theta}_{d}$ l'\'el\'ement de $\tilde{G}(F)$ d\'efini par $\boldsymbol{\theta}_{d}(v_{i})=(-1)^{i+[(d+1)/2]}v^*_{d+1-i}$ pour tout $i$. Notons simplement $\theta_{d}$ l'automorphisme de $G$ associ\'e \`a $\boldsymbol{\theta}_{d}$. On a $\theta_{d}(g)=J_{d}{^tg}^{-1}J_{d}^{-1}$, o\`u $J_{d}$ est la matrice antidiagonale de coefficients $(J_{d})_{i,d+1-i}=(-1)^{i}$. Cet automorphisme $\theta_{d}$ conserve les sous-groupes que l'on vient d'introduire. Il conserve aussi le caract\`ere $\xi$. Le normalisateur  de $B_{d}$ dans $\tilde{G}$ est $\tilde{B}_{d}=B_{d}\boldsymbol{\theta}_{d}$. Le normalisateur commun de $B_{d}$ et $A_{d}$ est $\tilde{A}_{d}=A_{d}\boldsymbol{\theta}_{d}$.
  
  Pour $g\in GL_{d}(F)$, on pose
   $$\sigma(g)=sup(\{1\}\cup\{log(\vert g_{ij}\vert _{F}); i,j=1,...,d\}\cup\{log(\vert (g^{-1})_{ij}\vert _{F});i,j=1,...,d\}).$$
    Pour $g\in G(F)$, on d\'efinit $\sigma(g)$ en identifiant $G$ \`a $GL_{d}$ par le choix d'une base. La fonction $\sigma$ d\'epend \'evidemment du choix de la base, mais d'une fa\c{c}on inessentielle. Pour tout r\'eel $b$, on note ${\bf 1}_{\sigma< b}$, resp. ${\bf 1}_{\sigma\geq b}$, la fonction caract\'eristique du sous-ensemble des $g\in G(F)$ tels que $\sigma(g)< b$, resp. $\sigma(g)\geq b$
    
    Remarquons que ${\cal A}_{\tilde{G}}=\{0\}$. Pour tout L\'evi tordu $\tilde{M}$ de $\tilde{G}$, le proc\'ed\'e de 1.3 munit $\tilde{M}(F)$ d'une extension canonique de la fonction $H_{\tilde{M}}$ \`a $\tilde{M}(F)$.  On v\'erifie que $H_{\tilde{A}_{d}}(\boldsymbol{\theta}_{d})=0$.

    \bigskip
    
    \subsection{Mod\`ele de Whittaker}
    
    Fixons une base $(v_{i})_{i=1,...,d}$ de $V$. Soit $(\pi,E)$ une repr\'esentation admissible irr\'eductible de $G(F)$. On appelle fonctionnelle de Whittaker une application lin\'eaire $\phi:E\to {\mathbb C}$ telle que $\phi(\pi(u)e)=\xi(u)\phi(e)$ pour tous $u\in U_{d}(F)$ et $e\in E$. Comme on le sait, l'espace de ces fonctionnelles est de dimension au plus $1$. Supposons que cette dimension soit $1$ et fixons une fonctionnelle de Whittaker non nulle $\phi$. Pour $e\in E$, on d\'efinit une fonction $W_{e}$ sur $G(F)$ par $W_{e}(g)=\phi(\pi(g)e)$. L'espace des fonctions $W_{e}$, pour $e\in E$, est le mod\`ele de Whittaker de $\pi$.
    
    On note $\theta_{d}(\pi)$ la repr\'esentation $g\mapsto \pi(\theta_{d}(g))$. Cette repr\'esentation est isomorphe \`a la contragr\'ediente $\pi^{\vee}$. La condition pour que $\pi$ se prolonge en une repr\'esentation $\tilde{\pi}$ de $\tilde{G}(F)$ est que $\theta_{d}(\pi)$ soit isomorphe \`a $\pi$. Supposons qu'il en soit ainsi et soit $\tilde{\pi}$ un prolongement. Supposons aussi que $\pi$ admette un mod\`ele de Whittaker. L'application $\phi\mapsto \phi\circ\tilde{\pi}(\boldsymbol{\theta}_{d})$ conserve la droite des fonctionnelles de Whittaker. Il existe donc un nombre complexe $w(\tilde{\pi},\psi)\in {\mathbb C}^{\times}$ tel que $\phi\circ\tilde{\pi}(\boldsymbol{\theta}_{d})=w(\tilde{\pi},\psi)\phi$ pour toute fonctionnelle de Whittaker $\phi$. Un calcul de changement de base montre que ce nombre ne d\'epend pas de la base $(v_{i})_{i=1,...,d}$ choisie. Par contre, il d\'epend de $\psi$. Soit $b\in F^{\times}$, notons $\psi^b$ le caract\`ere $z\mapsto \psi(bz)$ de $F$, dont on d\'eduit un caract\`ere $\xi^b$ de $U_{d}(F)$. Notons $D^b$ la matrice diagonale de coefficients diagonaux $D^b_{i,i}=b^{d-i}$. Si $\phi$ est une fonctionnelle de Whittaker relative \`a $\xi$, alors $\phi\circ\pi(D^b)$ est une telle fonctionnelle relative \`a $\xi^b$.  On a 
    $$\phi\circ\pi(D^b)\circ\tilde{\pi}(\boldsymbol{\theta}_{d}) =\phi\circ\tilde{\pi}(\boldsymbol{\theta}_{d})\circ\pi(\theta_{d}(D^b)).$$
    Mais $\theta_{d}(D^b)$ est \'egal \`a $D^b$ multipli\'e par la matrice centrale de coefficients diagonaux \'egaux \`a $b^{d-1}$. On en d\'eduit l'\'egalit\'e 
    
    (1) $w(\tilde{\pi},\psi^b)=\omega_{\pi}(b)^{d-1}w(\tilde{\pi},\psi)$, 
    
   \noindent o\`u $\omega_{\pi}$ est le caract\`ere central de $\pi$ (c'est un caract\`ere d'ordre au plus $2$).
    On a d\'efini  en 1.9 la contragr\'ediente $\tilde{\pi}^{\vee}$ de $\tilde{\pi}$. Supposons $\pi$ temp\'er\'ee. On a la relation
    
    (2)  $w(\tilde{\pi}^{\vee},\psi^-)w(\tilde{\pi},\psi)=1,$

 \noindent   o\`u $\psi^-=\psi^b$ pour $b=-1$. En effet, la condition $\pi^{\vee}\simeq \pi$ signifie qu'il existe une forme bilin\'eaire invariante sur $E\times E$. On sait construire cette forme bilin\'eaire.  Introduisons comme ci-dessus l'espace des fonctions de Whittaker $W_{e}$ pour $e\in E$ et introduisons l'espace similaire des fonctions $W_{e}^-$ relatif au caract\`ere $\xi^-$.  On peut prendre pour forme bilin\'eaire la forme
    $$(e',e)\mapsto <e',e>=\int_{GL_{d-1}(F)}W_{e'}^-(\gamma)W_{e}(\gamma)d\gamma,$$
    o\`u $GL_{d-1}$ est identifi\'e au groupe des automorphismes lin\'eaires de l'hyperplan de $V$ engendr\'e par $v_{1},...,v_{d-1}$. Par d\'efinition, on a $<\tilde{\pi}^{\vee}(\boldsymbol{\theta}_{d})e',\tilde{\pi}(\boldsymbol{\theta}_{d})e>=<e',e>$.  Il suffit d'expliciter ces deux termes \`a l'aide de la formule int\'egrale ci-dessus pour obtenir (2).      
    
    On montre de m\^eme en utilisant l'unitarit\'e de $\tilde{\pi}$ que $\vert w(\tilde{\pi},\psi)\vert =1$.

    \bigskip
    
    \subsection{Repr\'esentations induites}
    
    Soit $\tilde{L}$ un L\'evi tordu de $\tilde{G}$. Il existe une d\'ecomposition
    $$V=V_{u}\oplus...\oplus V_{1}\oplus V_{0}\oplus V_{-1}\oplus ...\oplus V_{-u}$$
    de sorte que $\tilde{L}(F)$ soit l'ensemble des $\tilde{x}\in \tilde{G}(F)$ tels que $\tilde{x}(V_{j})=V_{-j}^*$ pour tout $j=-u,...,u$. On a
     $$L=GL_{d_{u}}\times...\times GL_{d_{1}}\times GL_{d_{0}}\times GL_{d_{1}}\times...\times GL_{d_{u}},$$
     o\`u $d_{j}=dim(V_{j})=dim(V_{-j})$.  On peut choisir une base $(v_{i})_{i=1,...,d}$ de $V$ de sorte que $V_{j}$ ait pour base les vecteurs $v_{i}$ pour $i=d_{u}+...+d_{j+1}+1,...,d_{u}+...+d_{j}$. Alors $\boldsymbol{\theta}_{d}$ appartient \`a $\tilde{L}(F)$. Soit $\tilde{\sigma}\in Temp(\tilde{L})$. On a
     $$\sigma=\sigma_{u}\otimes...\otimes \sigma_{1}\otimes\sigma_{0}\otimes \sigma_{-1}\otimes...\otimes \sigma_{-u},$$
     o\`u $\sigma_{j}$ est une repr\'esentation irr\'eductible et temp\'er\'ee de $GL_{d_{j}}(F)$ et $\theta_{d_{j}}(\sigma_{j})\simeq \sigma_{-j}$. Pour tout $j=1,...,u$, choisissons un op\'erateur unitaire $A_{j}:E_{\sigma_{j}}\to E_{\sigma_{-j}}$ tel que $A_{j}\sigma_{j}(\theta_{d_{j}}(x_{j}))=\sigma_{-j}(x_{j})A_{j}$ pour tout $x_{j}\in GL_{d_{j}}(F)$. Notons $\tilde{G}_{0}$ l'analogue de $\tilde{G}$ quand $d$ est remplac\'e par $d_{0}$. Alors il existe un unique prolongement $\tilde{\sigma }_{0}$ de $ \sigma_{0}$  \`a $\tilde{G}_{0}(F)$ tel que, pour
    $$  e=e_{u}\otimes...\otimes e_{1}\otimes e_{0}\otimes e_{-1}\otimes...\otimes e_{-u}\in E_{\sigma},$$
    on ait
    $$\tilde{\sigma}(\boldsymbol{\theta}_{d})(e)=A_{u}^{-1}e_{-u}\otimes...\otimes A_{1}^{-1}e_{-1}\otimes \tilde{\sigma}_{0}(\boldsymbol{\theta}_{d_{0}})(e_{0})\otimes A_{1}e_{1}\otimes...\otimes A_{u}e_{u}.$$
   Cela ne d\'epend pas du choix des $A_{j}$.  Introduisons la  repr\'esentation induite $\tilde{\pi}=Ind_{L}^G(\tilde{\sigma})$. On v\'erifie que $w(\tilde{\pi},\psi)=w(\tilde{\sigma}_{0},\psi)$.
       
    Il est utile de calculer le caract\`ere de la repr\'esentation $\tilde{\sigma}$.  Pour $\tilde{x}\in \tilde{L}(F)$ et $j=-u,...,u$, on note $\tilde{x}_{j}:V_{j}\to V_{-j}^*$ la restriction de $\tilde{x}$ \`a $V_{j}$.    
    
\ass{Lemme}{Pour tout $\tilde{x}\in \tilde{L}(F)$ assez r\'egulier, on a l'\'egalit\'e
$$\Theta_{\tilde{\sigma}}(\tilde{x})=\Theta_{\tilde{\sigma}_{0}}(\tilde{x}_{0})\prod_{j=1,...,u}\Theta_{\sigma_{j}}((-1)^{d+1}{^t\tilde{x}}_{-j}^{-1}\tilde{x}_{j}).$$}

Preuve. Pour $j=-u,...,u$ soit $f_{j}\in C_{c}^{\infty}(GL_{d_{j}}(F))$. D\'efinissons $\tilde{f}$ sur $\tilde{L}(F)$ par $\tilde{f}(x\boldsymbol{\theta}_{d})=\prod_{j=-u,...,u}f_{j}(x_{j})$ pour $x\in L(F)$, o\`u les $x_{j}$ sont les diff\'erentes composantes de $x$. Pour 
$$ e=e_{u}\otimes...\otimes e_{1}\otimes e_{0}\otimes e_{-1}\otimes...\otimes e_{-u}\in E_{\sigma},$$ 
 on a l'\'egalit\'e
 $$\tilde{\sigma}(\tilde{f})e=\sigma_{u}(f_{u})A_{u}^{-1}e_{-u}\otimes...\otimes \sigma_{1}(f_{1})A_{1}^{-1}e_{-1}\otimes \sigma_{0}(f_{0})\tilde{\sigma}_{0}(\boldsymbol{\theta}_{d_{0}})(e_{0})\otimes \sigma_{-1}(f_{-1})A_{1}e_{1}\otimes...\otimes \sigma_{-u}(f_{-u})A_{u}e_{u}.$$
 Pour $j=1,...,u$, l'op\'erateur 
 $$e_{j}\otimes e_{-j}\mapsto\sigma_{j}(f_{j})A_{j}^{-1}e_{-j}\otimes \sigma_{-j}(f_{-j})A_{j}e_{j}$$
 a m\^eme trace que l'op\'erateur
 $$e_{j}\mapsto\sigma_{j}(f_{j})A_{j}^{-1}\sigma_{-j}(f_{-j})A_{j}e_{j}=\sigma_{j}(f_{j}\star {^{\theta}f}_{-j})e_{j}$$
 o\`u $^{\theta}f_{-j}(y)=f_{-j}(\theta_{d_{j}}(y))$ et $\star$ est la convolution. Sa trace est
 $$ \int_{GL_{d_{j}}(F)^2}\Theta_{\sigma_{j}}(x_{j}\theta_{d_{j}}(x_{-j}))f_{j}(x_{j})f_{-j}(x_{-j})dx_{j}\,dx_{-j}.$$
 D'autre part, la trace de l'op\'erateur $\sigma_{0}(f_{0})\tilde{\sigma}(\boldsymbol{\theta}_{d_{0}})$ est
 $$ \int_{GL_{d_{0}}(F)}\Theta_{\tilde{\sigma}_{0}}(x_{0}\boldsymbol{\theta}_{d_{0}})f_{0}(x_{0})dx_{0}.$$
 La trace de l'op\'erateur $\tilde{\sigma}(\tilde{f})$ est  donc \'egale  \`a
 $$\int_{\tilde{L}(F)} \Theta_{\tilde{\sigma}_{0}}(x_{0}\boldsymbol{\theta}_{d_{0}})\big(\prod_{j=1,...,u}\Theta_{\sigma_{j}}(x_{j}\theta_{d_{j}}(x_{-j}))\big)\tilde{f}(\tilde{x})d\tilde{x},$$
 o\`u on a \'ecrit $\tilde{x}=x\boldsymbol{\theta}_{d}$ avec $x\in L(F)$. On v\'erifie que $x_{0}\boldsymbol{\theta}_{d_{0}}=\tilde{x}_{0}$ et que $x_{j}\theta_{d_{j}}(x_{-j})=(-1)^{d+1}{^t\tilde{x}}_{-j}^{-1}\tilde{x}_{j}$ pour tout $j=1,...,u$. Cela entra\^{\i}ne l'\'egalit\'e de l'\'enonc\'e. $\square$
 
 \bigskip

\subsection{Repr\'esentations elliptiques des groupes lin\'eaires tordus}

 Soit $\tilde{L}$ un L\'evi tordu de $\tilde{G}$ que l'on \'ecrit comme dans le paragraphe pr\'ec\'edent.
 D\'ecrivons l'ensemble $\Pi_{ell}(\tilde{L})$ qui intervient en 1.10. Pour $j=1,...,u$, soit $\sigma_{j}$ une repr\'esentation admissible irr\'eductible de $GL_{d_{j}}(F)$, de la s\'erie discr\`ete. Soit $Q_{0}$ un sous-groupe parabolique standard de $GL_{d_{0}}$, de L\'evi
 $$L_{0}=GL_{d'_{1}}\times...\times GL_{d'_{s}}.$$
Pour $j=1,...,s$,  soit $\tau_{j}$ une repr\'esentation admissible irr\'eductible de $GL_{d'_{j}}(F)$, de la s\'erie discr\`ete. On suppose $\theta_{d'_{j}}(\tau_{j})\simeq \tau_{j}$ et $\tau_{k}\not\simeq \tau_{j}$ si $j\not=k$. Posons $\sigma_{0}=Ind_{Q_{0}}^{GL_{d_{0}}}(\tau_{1}\otimes...\otimes \tau_{s})$ et 
 $$\sigma=\sigma_{u}\otimes...\otimes \sigma_{1}\otimes \sigma_{0}\otimes \theta_{d_{1}}(\sigma_{1})\otimes...\otimes \theta_{d_{u}}(\sigma_{u}).$$
Alors $\sigma$ se prolonge en une repr\'esentation $\tilde{\sigma}$ de $\tilde{L}(F)$ qui appartient \`a $\Pi_{ell}(\tilde{L})$. L'ensemble $\Pi_{ell}(\tilde{L})$ est l'ensemble des repr\'esentations qui se construisent par ce proc\'ed\'e. 

Poursuivons avec la repr\'esentation $\tilde{\sigma}$ que l'on vient de construire. Notons ${\cal O}$ son orbite  pour l'action $\tilde{\sigma}\mapsto \tilde{\sigma}_{\lambda}$ de $i{\cal A}_{\tilde{L}}^*$. On note $i{\cal A}_{{\cal O}}^{\vee}$ le stabilisateur de $\tilde{\sigma}$ dans $i{\cal A}_{\tilde{L}}^*$. On note $s(\tilde{\sigma})=s({\cal O})=s$ et
 $$c({\cal O})=2^{-s({\cal O})-a_{\tilde{L}}}[i{\cal A}_{{\cal O}}^{\vee}:i{\cal A}_{\tilde{L},F}^{\vee}]^{-1} .$$
 Ce terme est le coefficient qui intervient dans la relation 1.10(3), cf. [W6] th\'eor\`eme 7.1.

 \bigskip
 
 \subsection{Facteurs $\epsilon$}
 
 Soit $m\in {\mathbb N}$ un entier tel que $m<d$ et $m$ et $d$ soient de parit\'es distinctes. Posons $H=GL_{m}$. Soient $\tilde{\pi}\in Temp(\tilde{G})$ et $\tilde{\rho}\in Temp(\tilde{H})$. On introduit le facteur $\epsilon(s,\pi\times \rho,\psi)$ de [JPSS] th\'eor\`eme 2.7. Soit $\nu\in F^{\times}$. On pose
 $$\epsilon_{\nu}(\tilde{\pi},\tilde{\rho})=\epsilon_{-\nu}(\tilde{\rho},\tilde{\pi})=w(\tilde{\pi},\psi)w(\tilde{\rho},\psi)\omega_{\pi}((-1)^{[m/2]}2\nu)\omega_{\rho}((-1)^{1+[d/2]}2\nu)\epsilon(1/2,\pi\times \rho,\psi).$$
 Ce terme ne d\'epend pas de $\psi$. Cela r\'esulte de la relation 2.2(1) et de l'\'egalit\'e bien connue 
 $$(1) \qquad \epsilon(1/2,\pi\times \rho,\psi^b)=\omega_{\pi}(b)^m\omega_{\rho}(b)^d\epsilon(1/2,\pi\times \rho,\psi).$$
 
 \ass{ Remarques}{ (a) Cette d\'efinition est dissym\'etrique en $\pi$ et $\rho$.
 
 (b) Dans le cas particulier o\`u $m=0$,  on consid\`ere par convention que $\omega_{\rho}=1$ et $\epsilon(1/2,\pi\times \rho,\psi)=1$. }
 
 Soit $\tilde{L}$ un L\'evi tordu de $\tilde{G}$ que l'on \'ecrit comme en 2.3. Soit $\tilde{\sigma}\in Temp(\tilde{L})$. On \'ecrit
 $$\sigma=\sigma_{u}\otimes...\otimes\sigma_{1}\otimes \sigma_{0}\otimes \theta_{d_{1}}(\sigma_{1})\otimes...\otimes \theta_{d_{u}}(\sigma_{u}).$$
 Le prolongement $\tilde{\sigma}$ de $\sigma$ d\'etermine  un prolongement $\tilde{\sigma}_{0}$ de $\sigma_{0}$ comme en 2.3. On pose
$$\epsilon_{\nu}(\tilde{\sigma},\tilde{\rho})=w(\tilde{\sigma}_{0},\psi)w(\tilde{\rho},\psi)\omega_{\sigma_{0}}((-1)^{[m/2]}2\nu)\omega_{\rho}((-1)^{[d_{0}/2]+1}2\nu)(\prod_{j=1,...,u}\omega_{\sigma_{j}}(-1)^m)\epsilon(1/2,\sigma_{0}\times \rho,\psi).$$ 
On v\'erifie que
$$\epsilon_{\nu}(\tilde{\sigma},\tilde{\rho})= \epsilon_{\nu}(\tilde{\sigma}_{0},\tilde{\rho})\prod_{j=1,...,u}\omega_{\sigma_{j}}(-1)^m .$$
Remarquons que le terme $\epsilon_{\nu}(\tilde{\sigma},\tilde{\rho})$ ne d\'epend que de l'orbite ${\cal O}$ de $\tilde{\sigma}$ sous l'action de $i{\cal A}_{\tilde{L}}^*$.   Soit $\tilde{Q}\in {\cal P}(\tilde{L})$, posons $\tilde{\pi}=Ind_{Q}^G(\tilde{\sigma})$. On a l'\'egalit\'e

(2) $\epsilon_{\nu}(\tilde{\pi},\tilde{\rho})=\epsilon_{\nu}(\tilde{\sigma},\tilde{\rho})$.
 
 Cela r\'esulte ais\'ement des \'egalit\'es bien connues suivantes:
$$\epsilon(s,\pi\times \rho,\psi)=\epsilon(s,\sigma_{0}\times \rho,\psi)\prod_{j=1,...,u}\epsilon(s,\sigma_{j}\times \rho,\psi)\epsilon(s,\theta_{d_{j}}(\sigma_{j})\times \rho,\psi),$$
$$ \epsilon(1/2,\sigma_{j}\times \rho,\psi)\epsilon(1/2,\sigma_{j}^{\vee}\times \rho^{\vee},\psi)=\omega_{\sigma_{j}}(-1)^m\omega_{\rho}(-1)^{d_{j}},$$
et du fait que $\rho^{\vee}=\rho$.

 \bigskip
 
 \section{La partie g\'eom\'etrique de la formule int\'egrale}
 
 \bigskip
 
 \subsection{Plongements de groupes lin\'eaires tordus}
 
 Soient $V$ un espace vectoriel sur $F$ de dimension $d\geq1$, $r$ un entier naturel tel que $2r+1\leq d$, $W$ un sous-espace de $V$ de dimension $m=d-2r-1$, $(z_{i})_{i=-r,...,r}$ une base d'un suppl\'ementaire $Z$ de $W$ dans $V$, et enfin $\nu$ un \'el\'ement de $F^{\times}$. On note $G=GL(V)$, $H=GL(W)$  et on introduit les groupes tordus $\tilde{G}=Isom(V,V^*)$, $\tilde{H}=Isom(W,W^*)$. L'espace dual $V^*$ se d\'ecompose en $V^*=W^*\oplus Z^*$. Notons $(z^*_{i})_{i=-r,...,r}$ la base de $Z^*$ duale de $(z_{i})_{i=-r,...,r}$ et $\tilde{\zeta}$ l'\'el\'ement de $Isom(Z,Z^*)$ d\'efini par $\tilde{\zeta}(z_{i})=(-1)^{i}2\nu z_{-i}^*$ pour tout $i=-r,...,r$. Le groupe $H$ est un sous-groupe de $G$. On introduit un plongement $\iota:\tilde{H}\to \tilde{G}$. Pour $\tilde{y}\in \tilde{H}$, $\iota(\tilde{y})$ envoie $W$ dans $W^*$ et $Z$ dans $Z^*$. La restriction de $\iota(\tilde{y})$ \`a $W$, resp. $Z$, co\"{\i}ncide avec $\tilde{y}$, resp. $\tilde{\zeta}$. Le plongement $\iota$ ainsi d\'efini est \'equivariant pour les actions de $H$. Pour simplifier les notations, on identifie tout \'el\'ement de $\tilde{H}$ \`a son image par $\iota$. On pose $V_{0}=W\oplus Fz_{0}$ et on note $G_{0}$ le groupe des automorphismes $F$-lin\'eaires de $V_{0}$.
 
 Notons $P$ le sous-groupe parabolique de $G$ form\'e des \'el\'ements qui conservent le drapeau
 $$Fz_{r}\subset Fz_{r}\oplus Fz_{r-1}\subset...\subset Fz_{r}\oplus...\oplus Fz_{1}\subset Fz_{r}\oplus...\oplus Fz_{1}\oplus V_{0}$$
 $$\subset Fz_{r}\oplus...\oplus Fz_{1}\oplus V_{0}\oplus Fz_{-1}\subset...\subset Fz_{r}\oplus...\oplus Fz_{1}\oplus V_{0}\oplus Fz_{-1}\oplus...\oplus Fz_{-r} .$$
 Notons $A$ le tore isomorphe \`a $GL_{1}^{2r}$ form\'e des \'el\'ements qui conservent chaque droite $Fz_{\pm i}$ pour $i=1,...,r$ et agissent trivialement sur $V_{0}$. Pour $a\in A$ et $i=\pm 1,...,\pm r$, on note $a_{i}$ la valeur propre de $a$ sur $z_{i}$.
 Notons $U$ le radical unipotent de $P$ et $M$ sa composante de L\'evi  qui contient $A$. On a $M=AG_{0}$, en particulier $M$ contient $H$. Notons $\tilde{P}$ le normalisateur de $P$ dans $\tilde{G}$ et $\tilde{M}$ le normalisateur commun de $P$ et $M$. Le groupe $\tilde{M}$ contient $\tilde{H}$, en particulier $\tilde{M}(F)\not=\emptyset$, donc $\tilde{M}$ est un L\'evi tordu de $\tilde{G}$. On d\'efinit un caract\`ere $\xi$ de $U(F)$ par
 $$\xi(u)=\psi(\sum_{i=-r,...,r-1}u_{i+1,i}),$$
 o\`u $u_{i+1,i}=<z_{-i-1}^*,uz_{i}>$. On v\'erifie que $\xi$ est invariant par $\theta_{\tilde{y}}$ pour tout $\tilde{y}\in \tilde{H}(F)$.
 
 Fixons un $\mathfrak{o}_{F}$-r\'eseau $R$ dans $V$ qui est somme d'un $\mathfrak{o}_{F}$-r\'eseau de $V_{0}$ et d'un $\mathfrak{o}_{F}$-r\'eseau de base sur $\mathfrak{o}_{F}$ form\'ee de vecteurs proportionnels aux $z_{\pm i}$ pour $i=1,...,r$. Notons $K$ le stabilisateur de $R$ dans $G(F)$ et $R^{\vee}$ le r\'eseau dual de $R$ dans $V^*$. On a l'\'egalit\'e $G(F)=P(F)K$. Pour un entier $N\geq 0$, introduisons la fonction $\kappa_{N}$ sur $G(F)$ d\'efinie de la fa\c{c}on suivante. Elle est invariante \`a gauche par $U(F)$ et \`a droite par $K$.  Sur $M(F)$, c'est la fonction caract\'eristique de l'ensemble des $m\in M(F)$ qui s'\'ecrivent $m=ag_{0}$, avec $a\in A(F)$, $g_{0}\in G_{0}(F)$, tels que:
 
 - pour tout $i=\pm 1,...,\pm r$, $\vert val_{F}(a_{i})\vert\leq N$;  
 
 - $g_{0}^{-1}z_{0}\in \mathfrak{p}_{F}^{-N}R$ et $^tg_{0}z_{0}^*\in \mathfrak{p}_{F}^{-N}R^{\vee}$.
 
 \bigskip
 
 \subsection{Les ingr\'edients de la formule g\'eom\'etrique}
 
 Consid\'erons une d\'ecomposition $W=W'\oplus W''$ telle que la dimension de $W'$ soit paire. Posons $H'=GL(W')$ et $\tilde{H}'=Isom(W',W^{_{'}*})$. Soit $\tilde{T}'$ un sous-tore maximal de $\tilde{H}'$ qui est anisotrope, c'est-\`a-dire que $A_{\tilde{T}'}=\{1\}$. Soit $\tilde{\zeta}_{H,T}\in Isom(W'',W^{_{''}*})$ une forme bilin\'eaire sym\'etrique.  Notons $V''=W''\oplus Z$ et $\tilde{\zeta}_{G,T}\in Isom(V'',V^{_{''}*})$ la somme directe de $\tilde{\zeta}_{H,T}$ et $\tilde{\zeta}$. On suppose que les groupes sp\'eciaux orthogonaux des formes quadratiques $\tilde{\zeta}_{H,T}$ et $\tilde{\zeta}_{G,T}$ sont quasi-d\'eploy\'es. Notons $\tilde{T}$ l'ensemble des \'el\'ements $\tilde{x}\in\tilde{H}$ tels que $\tilde{x}(W')=W^{_{'}*}$ et $\tilde{x}(W'')=W^{_{''}*}$, que la restriction de $\tilde{x}$ \`a $W'$ appartienne \`a $\tilde{T}'$ et que la restriction de $\tilde{x}$ \`a $W''$ soit \'egale \`a $\tilde{\zeta}_{H,T}$. On note $\underline{{\cal T}}$ l'ensemble des sous-ensembles $\tilde{T}$ de $\tilde{H}$ obtenus de cette fa\c{c}on.

 Pour un tel objet $\tilde{T}$, que l'on peut consid\'erer comme un sous-tore tordu de $\tilde{H}$, en g\'en\'eral non maximal, on note $T$ le tore not\'e ci-dessus $T'$, que l'on peut consid\'erer comme un sous-groupe de $H$ (un \'el\'ement de $T$ fixe tout point de $W''$). On d\'efinit comme dans le cas d'un sous-tore tordu maximal les ensembles $T_{\theta}$ et $\tilde{T}(F)_{/\theta}$.  Il est utile de remarquer que, parce que $dim(W')$ est paire, on a l'\'egalit\'e $T(F)^{\theta}=T_{\theta}(F)$. On munit $T_{\theta}(F)$ de la mesure de Haar de masse totale $1$. Comme en 1.4, on en d\'eduit une mesure sur $\tilde{T}(F)_{/\theta}$. Le normalisateur  $Norm_{H}(\tilde{T})$  de $\tilde{T}$ dans $H$  contient $T\times SO(\tilde{\zeta}_{H,T})$, ce dernier groupe  \'etant le groupe sp\'ecial orthogonal de la forme quadratique $\tilde{\zeta}_{H,T}$ sur $W''$. On pose
 $$W(H,\tilde{T})=Norm_{H}(\tilde{T})(F)/(T(F)\times SO(\tilde{\zeta}_{H,T})(F)).$$
 Soit $\tilde{t}\in \tilde{T}(F)$. L'\'el\'ement $t={^t\tilde{t}}^{-1}\tilde{t}$ appartient \`a $GL(W')$. On pose
 $$\Delta_{r}(\tilde{t})=\vert 2\vert _{F}^{r^2+r+rdim(W'')}\vert det((1-t)_{\vert W'}\vert^r _{F}.$$

 Soit $\Gamma$ un quasi-caract\`ere de $\tilde{G}(F)$. On en d\'eduit comme en [W2] 7.3 une fonction $c_{\Gamma}$ d\'efinie presque partout sur $\tilde{T}(F)$. Rappelons la d\'efinition. Soit $\tilde{t}\in \tilde{T}(F)$ en position g\'en\'erale. Son commutant connexe $G_{\tilde{t}}$ dans $G$ est $T_{\theta}\times SO(\tilde{\zeta}_{G,T})$. Si $d$ est impair, la dimension de $V''$ est elle-aussi impaire et l'alg\`ebre de Lie $\mathfrak{g}_{\tilde{t}}(F)$ poss\`ede une unique orbite nilpotente r\'eguli\`ere que l'on note ${\cal O}_{reg}$. On pose $c_{\Gamma}(\tilde{t})=c_{\Gamma,{\cal O}_{reg}}(\tilde{t})$, cf. 1.6. Si $d$ est pair et si la dimension de $V''$ est $\leq2$, la m\^eme construction s'applique. Si $d$ est pair et si la dimension de $V''$ est $\geq4$, il y a plusieurs orbites nilpotentes r\'eguli\`eres dans $\mathfrak{g}_{\tilde{t}}(F)$. Mais on sait les param\'etrer, cf. [W2] 7.1. En particulier, on peut associer \`a $\nu$ une orbite ${\cal O}_{\nu}$. On pose $c_{\Gamma}(\tilde{t})=c_{\Gamma,{\cal O}_{\nu}}(\tilde{t})$. Soit maintenant $\Theta$ un quasi-caract\`ere de $\tilde{H}(F)$. On lui associe une fonction $c_{\Theta}$ d\'efinie presque partout sur $\tilde{T}(F)$. La construction est la m\^eme que ci-dessus, \`a ceci pr\`es que l'on remplace $V''$ par $W''$ et, dans la derni\`ere \'eventualit\'e ci-dessus, le param\`etre $\nu$ par $-\nu$. Les fonctions $c_{\Gamma}$ et $c_{\Theta}$ sont invariantes par conjugaison par $T(F)$. On a
 
 \ass{Lemme}{Les fonctions $c_{\Gamma}$ et $c_{\Theta}$ sont localement constantes sur un ouvert de Zariski non vide de $\tilde{T}(F)$. La fonction
 $$\tilde{t}\mapsto c_{\Theta}(\tilde{t})c_{\Gamma}(\tilde{t})D^{\tilde{H}}(\tilde{t})\Delta_{r}(\tilde{t})$$
 est localement int\'egrable sur $\tilde{T}(F)_{/\theta}$.}
  
  La preuve est la m\^eme que celle de la proposition 7.3 de [W2].
  
 Le groupe $H(F)$ agit par conjugaison sur $\underline{{\cal T}}$. On fixe un ensemble ${\cal T}$ de repr\'esentants des orbites.
 
 \bigskip
 
 \subsection{La formule g\'eom\'etrique}
 
 Soient $\Theta$ un quasi-caract\`ere sur $\tilde{H}(F)$ et $\tilde{f}\in C_{c}^{\infty}(\tilde{G}(F))$ une fonction tr\`es cuspidale.   Pour $g\in G(F)$, posons
 $$J(\Theta,\tilde{f},g)=\int_{\tilde{H}(F)}\int_{U(F)}\Theta(\tilde{y})\tilde{f}(g^{-1}\tilde{y}ug)\xi(u)\,du\,d\tilde{y}.$$
 Pour tout entier $N\geq1$, posons
 $$J_{N}(\Theta,\tilde{f})=\int_{H(F)U(F)\backslash G(F)}J(\Theta,\tilde{f},g)\kappa_{N}(g)dg.$$
 Pour $\tilde{T}\in \underline{{\cal T}}$, on d\'efinit la fonction $c_{\Theta}$ sur $\tilde{T}(F)$.  A $\tilde{f}$, on associe le quasi-caract\`ere $\Theta^J_{\tilde{f}}$ sur $\tilde{G}(F)$, cf. 1.7, puis une fonction $c_{\Theta^J_{\tilde{f}}}$ sur $\tilde{T}(F)$. On note simplement $c_{\tilde{f}}$ cette fonction. Posons
 $$J_{geom}(\Theta,\tilde{f})=\sum_{\tilde{T}\in {\cal T}}\vert W(H,\tilde{T})\vert ^{-1} \int_{\tilde{T}(F)_{/\theta}}c_{\Theta}(\tilde{t})c_{\tilde{f}}(\tilde{t})D^{\tilde{H}}(\tilde{t})\Delta_{r}(\tilde{t})d\tilde{t}.$$

 \ass{Th\'eor\`eme}{On a l'\'egalit\'e
 $$lim_{N\to \infty}J_{N}(\Theta,\tilde{f})=J_{geom}(\Theta,\tilde{f}).$$}
 
 La d\'emonstration occupe la fin de la section. Remarquons que l'\'egalit\'e du th\'eor\`eme ne d\'epend pas des mesures choisies, \`a l'exception de celles sur les tores compacts que l'on a suppos\'ees de masse totale $1$ (la mesure sur $G(F)$ intervient directement dans le membre de gauche, mais aussi dans la d\'efinition du quasi-caract\`ere $\Theta^J_{\tilde{f}}$ associ\'e \`a $\tilde{f}$). Il est plus commode dans cette section de  normaliser les mesures de la fa\c{c}on suivante:
 
 - on munit $\mathfrak{g}(F)$ de la forme bilin\'eaire $(X,X')\mapsto \frac{1}{2}trace(XX')$ et tout sous-espace non d\'eg\'en\'er\'e de la mesure autoduale;
 
 - on munit arbitrairement d'une mesure de Haar tout autre sous-espace de $\mathfrak{g}(F)$ (par exemple le radical nilpotent d'un sous-groupe parabolique);
 
 - on rel\`eve ces mesures aux groupes via l'exponentielle.
 
 Dans le cas d'un tore compact, la mesure de masse totale $1$ n'est plus $dt$, mais $\nu(T)dt$.
 
 \bigskip
 
 \subsection{Localisation}
 
 Par partition de l'unit\'e, on peut localiser le probl\`eme de la fa\c{c}on suivante. On fixe un point $\tilde{x}\in \tilde{G}_{ss}(F)$ et un bon voisinage $\omega$ de $0$ dans $\mathfrak{g}_{\tilde{x}}(F)$, aussi petit que l'on veut. On note $\Omega=(\tilde{x}exp(\omega))^G=\{g^{-1}\tilde{x}exp(X)g; g\in G(F), X\in \omega\}$. On suppose le support de $\tilde{f}$ inclus dans $\Omega$.
 
 Supposons d'abord que la classe de conjugaison de $\tilde{x}$  par $G(F)$ ne coupe pas $\tilde{H}(F)$. Alors, pour tout $\tilde{T}\in {\cal T}$, le quasi-caract\`ere $\Theta^J_{\tilde{f}}$ est nul au voisinage de $\tilde{T}(F)$, donc la fonction $c_{\tilde{f}}$ est nulle. Donc $J_{geom}(\Theta,\tilde{f})=0$. Pour d\'emontrer l'\'egalit\'e du th\'eor\`eme, il suffit de prouver que $J(\Theta,\tilde{f},g)=0$ pour tout $g\in G(F)$. Il suffit de prouver que $\tilde{H}(F)U(F)$ ne coupe pas $\Omega$. Puisque la partie semi-simple d'un \'el\'ement de $\tilde{H}(F)U(F)$ est conjugu\'ee \`a un \'el\'ement de $\tilde{H}(F)$ et que $\Omega$ est stable par l'op\'eration consistant \`a prendre les parties semi-simples, il suffit de prouver qu'aucun \'el\'ement semi-simple de $\Omega$ n'est conjugu\'e \`a un \'el\'ement de $\tilde{H}(F)$. On se rappelle qu'en 2.1, on a associ\'e \`a $\tilde{x}$ une forme quadratique sur le sous-espace propre $V''_{\tilde{x}}$ de $V$   associ\'e \`a la valeur propre $1$ de $^t\tilde{x}^{-1}\tilde{x}$. Notons $(V''_{\tilde{x}},\tilde{x})$  l'espace quadratique form\'e de ce sous-espace propre et de la forme quadratique d\'efinie par la restriction de $\tilde{x}$. Dire que la classe de conjugaison de $\tilde{x}$ par $G(F)$ ne coupe pas $\tilde{H}(F)$ revient \`a dire qu'aucun sous-espace de $(V''_{\tilde{x}},\tilde{x})$ n'est isomorphe \`a l'espace quadratique $(Z,\tilde{\zeta})$. Soit $X\in \omega$ un \'el\'ement semi-simple, posons $\tilde{x}'=\tilde{x}exp(X)$. Si $\omega$ est assez petit, $V''_{\tilde{x}'}$ est le noyau de $X$ dans $V''_{\tilde{x}}$ et  la forme quadratique $\tilde{x}'$ sur ce noyau est la restriction de $\tilde{x}$. Donc $(V''_{\tilde{x}'},\tilde{x}')$ ne contient aucun sous-espace isomorphe \`a $(Z,\tilde{\zeta})$. Un tel \'el\'ement ne peut pas \^etre  conjugu\'e \`a un \'el\'ement de $\tilde{H}(F)$. Cela d\'emontre l'assertion.
 
 On suppose d\'esormais que la classe de conjugaison de $\tilde{x}$ coupe $\tilde{H}(F)$. On ne perd rien \`a supposer plus simplement $\tilde{x}\in \tilde{H}(F)$. On note $W''$, resp. $V''$, le sous-espace propre de $W$, resp. $V$, associ\'e \`a la valeur propre $1$ de $x=^t\tilde{x}^{-1}\tilde{x}$. On munit ces espaces des formes quadratiques d\'efinies par $\tilde{x}$. On a la d\'ecomposition orthogonale $V''=W''\oplus Z$. On note $O(V'')$, $SO(V'')$ etc... les groupes orthogonaux et sp\'eciaux orthogonaux de ces espaces. On note $W'$ l'unique suppl\'ementaire de $W''$ dans $W$ stable par $x$ et $H'=GL(W')$. On a les \'egalit\'es $H_{\tilde{x}}=H'_{\tilde{x}}SO(W'')$, $H_{\tilde{x}}=H'_{\tilde{x}}SO(V'')$ (par abus de notations, on note $H'_{\tilde{x}}$ le commutant connexe dans $H'$ de la restriction de $\tilde{x}$ \`a $W'$). On suppose $\omega=\omega'\times \omega''$, o\`u $\omega'\subset \mathfrak{h}'_{\tilde{x}}(F)$ et $\omega''\subset \mathfrak{so}(V'')(F)$.
 
 On d\'efinit le quasi-caract\`ere $\Theta_{\tilde{x},\omega}$ de $\mathfrak{h}_{\tilde{x}}(F)$ par
 $$\Theta_{\tilde{x},\omega}(X)=\left\lbrace\begin{array}{cc}\Theta(\tilde{x}exp(X)),&\text{ si }X\in \omega,\\ 0,&\text{ sinon. }\\ \end{array}\right.$$
 Pour $g\in G(F)$, on d\'efinit une fonction $^g\tilde{f}_{\tilde{x},\omega}\in C_{c}^{\infty}(\mathfrak{g}_{\tilde{x}}(F))$ comme en 1.8, puis une fonction $^g\tilde{f}_{\tilde{x},\omega}^{\xi}\in C_{c}^{\infty}(\mathfrak{h}_{\tilde{x}}(F))$ par
 $$^g\tilde{f}_{\tilde{x},\omega}^{\xi}(X)=\int_{\mathfrak{u}_{\tilde{x}}(F)}{^g\tilde{f}}_{\tilde{x},\omega}(X+N)\xi(N)dN.$$
 On a not\'e ici $\xi$ le caract\`ere $N\mapsto \xi(exp(N))$ de $\mathfrak{u}_{\tilde{x}}(F)$. Posons
 $$J_{\tilde{x},\omega}(\Theta,\tilde{f},g)=\int_{\mathfrak{h}_{\tilde{x}}(F)}\Theta_{\tilde{x},\omega}(X) \,{^g\tilde{f}}_{\tilde{x},\omega}^{\xi}(X)dX,$$
 puis
 $$J_{\tilde{x},\omega,N}(\Theta,\tilde{f})=\int_{U_{\tilde{x}}(F)H_{\tilde{x}}(F)\backslash G(F)}J_{\tilde{x},\omega}(\Theta,\tilde{f},g)\kappa_{N}(g)dg.$$
 
 Notons $\underline{{\cal T}}'$ l'ensemble des sous-tores maximaux de $H'_{\tilde{x}}$ qui sont anisotropes, c'est-\`a-dire ne contiennent pas de sous-tore d\'eploy\'e autre que $\{1\}$. Le couple d'espaces quadratiques $(W'',V'')$ v\'erifient les conditions de [W2] 7.1. On d\'efinit un ensemble $\underline{{\cal T}}''$ de sous-tores (en g\'en\'eral non maximaux) de $SO(W'')$ comme en [W2] 7.3. On v\'erifie que l'application $\tilde{T}\mapsto T_{\theta}$ est une bijection du sous-ensemble des \'el\'ements $\tilde{T}\in \underline{{\cal T}}$ qui contiennent $\tilde{x}$ sur l'ensemble $\underline{{\cal T}}'\times \underline{{\cal T}}''$. Pour $\tilde{T}$ dans l'ensemble de d\'epart, les fonctions $c_{\Theta}$ et $c_{\tilde{f}}$ sont d\'efinies sur $\tilde{T}(F)$. On d\'efinit des fonctions $c_{\Theta,\tilde{x},\omega}$ et $c_{\tilde{f},\tilde{x},\omega}$ sur $\mathfrak{t}_{\theta}(F)$: elles sont nulles hors de $\mathfrak{t}_{\theta}(F)\cap \omega$; pour $X\in \mathfrak{t}_{\theta}(F)\cap \omega$,
 $$c_{\Theta,\tilde{x},\omega}(X)=c_{\Theta}(\tilde{x}exp(X)),\,\,c_{\tilde{f},\tilde{x},\omega}(X)=c_{\tilde{f}}(\tilde{x}exp(X)).$$

 Fixons un ensemble de repr\'esentants ${\cal T}_{\tilde{x}}$ des classes de conjugaison par $H_{\tilde{x}}(F)$ dans  $\underline{{\cal T}}'\times \underline{{\cal T}}''$.  On d\'efinit une fonction $\Delta''$ sur $\mathfrak{h}_{\tilde{x}}(F)$ par
 $$\Delta''(X)=\vert det(X_{\vert W''/W''(X)})\vert _{F},$$
 o\`u $W''(X)$ est le noyau de $X$ agissant dans $W''$. Posons
 $$J_{\tilde{x},\omega,geom}(\Theta,\tilde{f})=\sum_{I\in{\cal T}_{\tilde{x}}}\vert W(H_{\tilde{x}},I)\vert ^{-1}\nu(I)\int_{\mathfrak{i}(F)}c_{\Theta,\tilde{x},\omega}(X)c_{\tilde{f},\tilde{x},\omega}(X)D^{H_{\tilde{x}}}(X)\Delta''(X)^rdX.$$
 
 Posons enfin
 $$C(\tilde{x})=\vert 2\vert _{F}^{r^2+r+rdim(W'')}[Z_{H}(\tilde{x})(F):H_{\tilde{x}}(F)]^{-1}D^{\tilde{H}}(\tilde{x})\vert det((1-x)_{\vert W''})\vert _{F}^r.$$
 
 \ass{Lemme}{On a les \'egalit\'es
 $$J_{N}(\Theta,\tilde{f})=C(\tilde{x})J_{\tilde{x},\omega,N}(\Theta,\tilde{f}),\,\,J_{geom}(\Theta,\tilde{f})=C(\tilde{x})J_{\tilde{x},\omega,geom}(\Theta,\tilde{f}).$$}
 
 Preuve. Soit $g\in G(F)$. En utilisant la formule de Weyl, on a
 $$(1) \qquad J(\Theta,\tilde{f},g)=\sum_{\tilde{T}\in {\cal T}(\tilde{H})}\vert W(H,\tilde{T})\vert ^{-1}[T(F)^{\theta}:T_{\theta}(F)]^{-1}$$
 $$\qquad \int_{\tilde{T}(F)_{/\theta}}\int_{T_{\theta}(F)\backslash H(F)} {^g\tilde{f}}^{\xi}(h^{-1}\tilde{t}h)dh\Theta(\tilde{t})D^{\tilde{H}}(\tilde{t})d\tilde{t},$$
 o\`u, pour $\tilde{y}\in \tilde{H}(F)$, on a pos\'e
 $$^g\tilde{f}^{\xi}(\tilde{y})=\int_{U(F)}\tilde{f}(g^{-1}\tilde{y}ug)\xi(u)du.$$
 Pour tout sous-tore maximal $I$ de $H_{\tilde{x}}$, notons $T_{I}$ son commutant dans $H$ et $\tilde{T}_{I}=T_{I}\tilde{x}$, qui est un sous-tore maximal tordu de $\tilde{H}$.
 Introduisons un ensemble ${\cal T}(H_{\tilde{x}})$ de repr\'esentants des classes de conjugaison par $H_{\tilde{x}}(F)$ dans l'ensemble des sous-tores maximaux de $H_{\tilde{x}}$. Pour deux sous-tores maximaux tordus $\tilde{T}$ et $\tilde{T}'$ de $\tilde{H}$, notons
 $$W(\tilde{T},\tilde{T}')=\{h\in H(F); h\tilde{T}h^{-1}=\tilde{T}'\}/T(F).$$
 Tout \'el\'ement $w\in W(\tilde{T},\tilde{T}')$ induit une bijection de $\tilde{T}(F)_{/\theta}$ sur $\tilde{T}'(F)_{/\theta}$. On va prouver les propri\'et\'es suivantes:
 
 (2) soient $\tilde{T}\in {\cal T}(\tilde{H})$ et $\tilde{t}\in \tilde{T}(F)_{/\theta}$, en position g\'en\'erale; supposons
 $$\int_{T_{\theta}(F)\backslash H(F)} {^g\tilde{f}}^{\xi}(h^{-1}\tilde{t}h)dh\not=0;$$
 alors il existe $I\in {\cal T}(H_{\tilde{x}})$ et $w\in W(\tilde{T}_{I},\tilde{T})$ tel que 
 $$\tilde{t}\in w(\tilde{x}exp(\mathfrak{i}(F)\cap \omega));$$
 
 (3) soit $\tilde{T}\in {\cal T}(\tilde{H})$ et, pour $i=1,2$, soit $I_{i}\in {\cal T}(H_{\tilde{x}})$ et $w_{i}\in W(\tilde{T}_{I_{i}},\tilde{T})$; alors les sous-ensembles $w_{1}(\tilde{x}exp(\mathfrak{i}_{1}(F)\cap \omega))$ et $w_{2}(\tilde{x}exp(\mathfrak{i}_{2}(F)\cap \omega))$ de $\tilde{T}(F)_{/\theta}$ sont disjoints ou confondus;
 
 (4) soient $\tilde{T}\in {\cal T}(\tilde{H})$, $I_{1}\in {\cal T}(H_{\tilde{x}})$ et $w_{1}\in W(\tilde{T}_{I_{1}},\tilde{T})$; alors le nombre de couples $(I_{2},w_{2})$ tels que $I_{2}\in {\cal T}(H_{\tilde{x}})$, $w_{2}\in W(\tilde{T}_{I_{2}},\tilde{T})$ et $w_{2}(\tilde{x}exp(\mathfrak{i}_{2}(F)\cap \omega))=w_{1}(\tilde{x}exp(\mathfrak{i}_{1}(F)\cap \omega))$ est \'egal \`a
 $$\vert W(H_{\tilde{x}},I_{1})\vert [Z_{H}(\tilde{x})(F):H_{\tilde{x}}(F)][T(F)^{\theta}:T_{\theta}(F)]^{-1}.$$
 
 Sous les hypoth\`eses de (2), on voit comme au d\'ebut du paragraphe que la classe de conjugaison par $G(F)$ de $\tilde{t}$ coupe l'ensemble $\tilde{x}exp(\omega)$. Soient $X\in \omega$ et $g\in G(F)$ tels que $g\tilde{t}g^{-1}=\tilde{x}exp(X)$. L'\'el\'ement $g$ se restreint en une isom\'etrie entre les espaces quadratiques $(V''_{\tilde{t}},\tilde{t})$ et $(V''_{\tilde{x}exp(X)},\tilde{x}exp(X))$. Le premier de ces espaces contient $(Z,\tilde{\zeta})$ tandis que le second est inclus dans $(V''_{\tilde{x}},\tilde{x})=(V'',\tilde{x})$, qui contient lui-aussi $(Z,\tilde{\zeta})$. D'apr\`es le th\'eor\`eme de Witt, on peut trouver un \'el\'ement $g''\in O(V'')(F)$ tel que $g''g$ conserve $Z$ et y agisse par l'identit\'e. Rappelons que $O(V'')\subset Z_{G}(\tilde{x})$. Quitte \`a remplacer $g$ par $g''g$ et $X$ par $g''Xg^{_{''}-1}$, on se ram\`ene au cas o\`u $g$ conserve $Z$ et agit par l'identit\'e sur cet espace. Consid\'erons l'\'egalit\'e $g{^t\tilde{t}}^{-1}\tilde{t}g^{-1}={^t(\tilde{x}}exp(X))^{-1}\tilde{x}exp(X)$. Sur $Z$, le premier terme agit par l'identit\'e. Le second agit par la restriction \`a $Z$ de $exp(2X)$. Cette restriction est donc l'identit\'e, donc la restriction de $X$ \`a $Z$ est nulle. Un \'el\'ement de $\mathfrak{g}_{\tilde{x}}$ dont la restriction \`a $Z$ est nulle appartient \`a $\mathfrak{h}_{\tilde{x}}$. Donc $X\in \mathfrak{h}_{\tilde{x}}(F)$.  L'\'el\'ement $g$ \'etant une isom\'etrie de $(V''_{\tilde{t}},\tilde{t})$ sur $(V''_{\tilde{x}exp(X)},\tilde{x}exp(X))$ et conservant $Z$, il envoie l'orthogonal $W''_{\tilde{t}}$ de $Z$ dans $V''_{\tilde{t}}$ sur l'orthogonal $W''_{\tilde{x}exp(X)}$ de $Z$ dans $V''_{\tilde{x}exp(X)}$. Il envoie aussi $V'_{\tilde{t}}=W'_{\tilde{t}}$ sur $V'_{\tilde{x}exp(X)}=W'_{\tilde{x}exp(X)}$. Puisque $W=W'_{\tilde{t}}\oplus W''_{\tilde{t}}=W'_{\tilde{x}exp(X)}\oplus W''_{\tilde{x}exp(X)}$, $g$ conserve $W$, donc $g\in H(F)$. Quitte \`a multiplier $g$ par un \'el\'ement de $H_{\tilde{x}}(F)$, on peut supposer que $X\in \mathfrak{i}(F)$ pour un \'el\'ement $I\in {\cal T}(H_{\tilde{x}})$. Alors l'\'el\'ement $g^{-1}$ d\'efinit un \'el\'ement $w$ de $W(\tilde{T}_{I},\tilde{T})$ pour lequel $\tilde{t}$ appartient \`a $w(\tilde{x}exp(\mathfrak{i}(F)\cap \omega))$. Cela prouve (2).
 
 Sous les hypoth\`eses de (3), supposons que les ensembles $w_{1}(\tilde{x}exp(\mathfrak{i}_{1}(F)\cap \omega))$ et $w_{2}(\tilde{x}exp(\mathfrak{i}_{2}(F)\cap \omega))$  ne sont pas disjoints. On peut relever $w_{1}$ et $w_{2}$ en des \'el\'ements $h_{1}$ et $h_{2}$ de $H(F)$ tels que 
 $$h_{1}(\tilde{x}exp(\mathfrak{i}_{1}(F)\cap \omega))h_{1}^{-1}\cap h_{2}(\tilde{x}exp(\mathfrak{i}_{2}(F)\cap \omega))h_{2}^{-1}\not=\emptyset.$$
 Posons $h=h_{2}^{-1}h_{1}$. D'apr\`es les propri\'et\'es des bons voisinages, $h$ appartient \`a $Z_{G}(\tilde{x})(F)$. Donc la conjugaison par $h$ conserve $\omega$. D'autre part $y$ conjugue $\tilde{T}_{I_{1}}$ en $\tilde{T}_{I_{2}}$, donc aussi $I_{1}$ en $I_{2}$. Alors
 $$h(\tilde{x}exp(\mathfrak{i}_{1}(F)\cap \omega))h^{-1}=\tilde{x}exp(\mathfrak{i}_{2}(F)\cap \omega),$$
 puis
 $$h_{1}(\tilde{x}exp(\mathfrak{i}_{1}(F)\cap \omega))h_{1}^{-1}=h_{2}h(\tilde{x}exp(\mathfrak{i}_{1}(F)\cap \omega))h^{-1}h_{2}^{-1}=h_{2}(\tilde{x}exp(\mathfrak{i}_{2}(F)\cap \omega))h_{2}^{-1}.$$
 Cela prouve (3).
 
 Sous les hypoth\`eses de (4), posons
 $${\cal Y}=\{h\in Z_{H}(\tilde{x})(F); hI_{1}h^{-1}\in {\cal T}(H_{\tilde{x}})\}/(Z_{H}(\tilde{x})(F)\cap T_{I_{1}}(F)).$$
 La preuve de (2) montre que l'application $h\mapsto (I_{2}=hI_{1}h^{-1}, w_{2}=w_{1}h^{-1})$ est une surjection de ${\cal Y}$ sur l'ensemble des paires dont on veut calculer le nombre. On voit que l'application est aussi injective. Le nombre cherch\'e est donc $\vert {\cal Y}\vert $.  Remarquons que $Z_{H}(\tilde{x})(F)\cap T_{I_{1}}(F)=T_{I_{1}}(F)^{\theta}$. Il y a une application naturelle
 $${\cal Y}\to H_{\tilde{x}}(F)\backslash Z_{H}(\tilde{x})(F)/T_{I_{1}}(F)^{\theta}.$$
 Elle est surjective et ses fibres sont en bijection avec $W(H_{\tilde{x}},I_{1})$. D'autre part, $H_{\tilde{x}}$ est distingu\'e dans $Z_{H}(\tilde{x})$ et son intersection avec $T_{I_{1}}(F)$ est $I_{1}(F)$. Le quotient ci-dessus a donc pour nombre d'\'el\'ements
 $$[Z_{H}(\tilde{x})(F):H_{\tilde{x}}(F)][T_{I_{1}}(F)^{\theta}:I_{1}(F)]^{-1}.$$
 Enfin, $[T_{I_{1}}(F)^{\theta}:I_{1}(F)]=[T(F)^{\theta}:T_{\theta}(F)]$. L'assertion (4) r\'esulte de ces calculs.
 
On utilise ces trois assertions pour transformer la formule (1).  On doit rappeler que, d'apr\`es le choix de nos mesures, pour $I\in {\cal T}(H_{\tilde{x}})$, l'application  
$$\begin{array}{ccc}\mathfrak{i}(F)\cap \omega&\to&\tilde{T}_{I}(F)_{/\theta}\\ X&\mapsto&\tilde{x}exp(X)\\ \end{array}$$
pr\'eserve les mesures. On obtient
$$J(\Theta,\tilde{f},g)=[Z_{H}(\tilde{x})(F):H_{\tilde{x}}(F)]^{-1}\sum_{I\in {\cal T}(H_{\tilde{x}})}\vert W(H_{\tilde{x}},I)\vert ^{-1}\sum_{\tilde{T}\in {\cal T}(\tilde{H})}\vert W(H,\tilde{T})\vert ^{-1}\sum_{w\in W(\tilde{T}_{I},\tilde{T})}$$
$$\qquad \int_{ \mathfrak{i}(F)\cap \omega}\int_{T_{\theta}(F)\backslash H(F)} {^g\tilde{f}}^{\xi}(h^{-1}w\tilde{x}exp(X)w^{-1}h)dh\Theta(w\tilde{x}exp(X)w^{-1})D^{\tilde{H}}(w\tilde{x}exp(X)w^{-1})dX.$$
On a ici identifi\'e les \'el\'ements $w$ \`a des repr\'esentants dans $H(F)$. Les \'el\'ements $w$ disparaissent de cette formule par le changement de variables $h\mapsto wh$. Evidemment, pour tout $I$, il y a un seul $\tilde{T}$ tel que $W(\tilde{T}_{I},\tilde{T})$ soit non vide et, pour ce $\tilde{T}$, le nombre d'\'el\'ements de $W(\tilde{T}_{I},\tilde{T})$ est le m\^eme que celui de $W(H,\tilde{T})$. On obtient
$$J(\Theta,\tilde{f},g)=[Z_{H}(\tilde{x})(F):H_{\tilde{x}}(F)]^{-1}\sum_{I\in {\cal T}(H_{\tilde{x}})}\vert W(H_{\tilde{x}},I)\vert ^{-1} \int_{ \mathfrak{i}(F)\cap \omega}\int_{I(F)\backslash H(F)} $$
$${^g\tilde{f}}^{\xi}(h^{-1}\tilde{x}exp(X)h)dh\Theta(\tilde{x}exp(X))D^{\tilde{H}}(\tilde{x}exp(X))dX.$$
On a $D^{\tilde{H}}(\tilde{x}exp(X))=D^{\tilde{H}}(\tilde{x})D^{H_{\tilde{x}}}(X)$, d'o\`u
$$ J(\Theta,\tilde{f},g)=[Z_{H}(\tilde{x})(F):H_{\tilde{x}}(F)]^{-1}D^{\tilde{H}}(\tilde{x})\int_{H_{\tilde{x}}(F)\backslash H(F)}\sum_{I\in {\cal T}(H_{\tilde{x}})}\vert W(H_{\tilde{x}},I)\vert ^{-1}$$
$$ \int_{ \mathfrak{i}(F)\cap \omega}\int_{I(F)\backslash H_{\tilde{x}}(F)} {^{hg}\tilde{f}}^{\xi}(y^{-1}\tilde{x}exp(X)y)dy\Theta(\tilde{x}exp(X))D^{H_{\tilde{x}}}(X)dX\,dh.$$
 On peut remplacer les int\'egrations sur $\mathfrak{i}(F)\cap \omega$ par des int\'egrations sur $\mathfrak{i}(F)$, \`a condition de remplacer la fonction $\Theta(\tilde{x}exp(X))$ par $\Theta_{\tilde{x},\omega}(X)$.
Par la formule de Weyl appliqu\'ee dans $\mathfrak{h}_{\tilde{x}}(F)$,  on obtient
$$ J(\Theta,\tilde{f},g)=[Z_{H}(\tilde{x})(F):H_{\tilde{x}}(F)]^{-1}D^{\tilde{H}}(\tilde{x})\int_{H_{\tilde{x}}(F)\backslash H(F)}\int_{\mathfrak{h}_{\tilde{x}}(F) }\Theta_{\tilde{x},\omega}(X)\,{^{hg}\tilde{f}}^{\xi}(\tilde{x}exp(X))dX.$$
Pour $X\in \omega$, on a
$${^{hg}\tilde{f}}^{\xi}(\tilde{x}exp(X))=\int_{U(F)}{^{hg}\tilde{f}}(\tilde{x}exp(X)u)\xi(u)du$$
$$\qquad =\int_{U_{\tilde{x}}(F)\backslash U(F)}\int_{U_{\tilde{x}}(F)}{^{hg}\tilde{f}}(\tilde{x}exp(X)uv)\xi(uv)du\,dv.$$
Pour $u\in U_{\tilde{x}}(F)$, l'application $v\mapsto \theta_{\tilde{x}exp(X)u}^{-1}(v^{-1})v$ est un isomorphisme de $U_{\tilde{x}}(F)\backslash U(F)$ sur lui-m\^eme. Son jacobien est la valeur absolue du d\'eterminant de $1-\theta_{\tilde{x}}^{-1}$ agissant sur $\mathfrak{u}(F)/\mathfrak{u}_{\tilde{x}}(F)$. Notons $d(\tilde{x})$ ce d\'eterminant.  Par changement de variables, et en remarquant que
$$\xi(\theta_{\tilde{x}exp(X)u}^{-1}(v^{-1})v)=1,$$
on obtient
$${^{hg}\tilde{f}}^{\xi}(\tilde{x}exp(X))=\vert d(\tilde{x})\vert _{F}\int_{U_{\tilde{x}}(F)\backslash U(F)}\int_{U_{\tilde{x}}(F)}{^{hg}\tilde{f}}(v^{-1}\tilde{x}exp(X)uv)\xi(u)du\,dv$$
$$\qquad =\vert d(\tilde{x})\vert _{F}\int_{U_{\tilde{x}}(F)\backslash U(F)}\int_{U_{\tilde{x}}(F)}{^{vhg}\tilde{f}}(\tilde{x}exp(X)u)\xi(u)du\,dv.$$
On remplace $u$ par $exp(N)$ avec $N\in \mathfrak{u}_{\tilde{x}}(F)$. D'apr\`es les propri\'et\'es de $\omega$, la condition $X\in \omega$ entra\^{\i}ne $X+N\in \omega$ pour tout $N$. On obtient alors
$${^{hg}\tilde{f}}^{\xi}(\tilde{x}exp(X))=\vert d(\tilde{x})\vert _{F}\int_{U_{\tilde{x}}(F)\backslash U(F)}{^{vhg}\tilde{f}}_{\tilde{x},\omega}^{\xi}(X)dv.$$
Reportons cette expression dans (5), multiplions l'expression obtenue par $\kappa_{N}(g)$, puis int\'egrons sur $g\in H(F)U(F)\backslash G(F)$. On obtient
$$J_{N}(\Theta,\tilde{f})=[Z_{H}(\tilde{x})(F):H_{\tilde{x}}(F)]^{-1}D^{\tilde{H}}(\tilde{x})\vert d(\tilde{x})\vert _{F}J_{\tilde{x},\omega,N}(\Theta,\tilde{f}).$$
Pour obtenir la premi\`ere \'egalit\'e de l'\'enonc\'e, il reste \`a calculer $d(\tilde{x})$. Pour cela, on peut \'etendre le corps des scalaires  et se placer sur $\bar{F}$. On peut fixer une base $(v_{i})_{i=1,...,d}$ de $V$  et  une famille $(\lambda_{i})_{i=1,...,d}$ d'\'el\'ements de $\bar{F}^{\times}$ de sorte que:

- $Z$ ait pour base  la famille des $v_{i}$ pour $i\in\{1,...,r\}\cup\{d-r,...,d\}$;

- le sous-groupe $P$ est le stabilisateur du drapeau
$$Fv_{1}\subset Fv_{1}\oplus Fv_{2}\subset...\subset Fv_{1}\oplus ...\oplus Fv_{r}\subset Fv_{1}\oplus ...\oplus Fv_{d-r}$$
$$\subset Fv_{1}\oplus ...\oplus Fv_{d-r+1}\subset...\subset Fv_{1}\oplus ...\oplus Fv_{d} ;$$

- pour tout $i$, $\tilde{x}v_{i}=\lambda_{i}v_{d+1-i}^*$.

L'\'el\'ement $x=^t\tilde{x}^{-1}\tilde{x}$ est l'\'el\'ement diagonal qui envoie $v_{i}$ sur $\mu_{i}=\frac{\lambda_{i}}{\lambda_{d+1-i}}v_{i}$. Le sous-espace $W'$ est engendr\'e par les $v_{i}$ tels que $\mu_{i}\not=1$. Puisque $W'\cap Z=\{0\}$, on a $\mu_{i}=1$ si $v_{i}\in Z$. Introduisons la base  $(E_{i,j})_{i,j=1,...,d}$ de $\mathfrak{g}(F)$, form\'ee des matrices $E_{i,j}$ n'ayant qu'un coefficient non nul, le $(i,j)$-i\`eme, lequel vaut $1$. On calcule $\theta_{\tilde{x}}^{-1}(E_{i,j})=-\frac{\lambda_{d+1-i}}{\lambda_{d+1-j}}E_{d+1-j,d+1-i}$. L'espace $\mathfrak{u}(F)$ se d\'ecompose en sommes des sous-espaces suivants, qui sont stables par $\theta_{\tilde{x}}^{-1}$:

- les espaces $FE_{i,j}\oplus FE_{d+1-j,d+1-i}$ pour $i=1,...,r$, $i<j\leq r$ ou $d+1-r\leq j< d+1-i$;

- les droites $FE_{i,d+1-i}$ pour $i=1,...,r$;

- les espaces $FE_{i,j}\oplus FE_{d+1-j,d+1-i}$, pour $i=1,...,r$, $j=r+1,...,d-r$.

Consid\'erons un sous-espace du troisi\`eme type. La matrice de $1-\theta_{\tilde{x}}^{-1}$ s'y \'ecrit
$$\left(\begin{array}{cc}1&\frac{\lambda_{j}}{\lambda_{i}}\\ \frac{\lambda_{d+1-i}}{\lambda_{d+1-j}}&1\\ \end{array}\right).$$
Si $\mu_{j}\not=1$, c'est-\`a-dire si $v_{j}\in V'$, cette matrice est inversible  donc le sous-espace ne coupe pas $\mathfrak{u}_{\tilde{x}}$. Le d\'eterminant de la matrice est $1-\mu_{d+1-j}$. La contribution de ces sous-espaces \`a $d(\tilde{x})$ est donc $det((1-x)_{\vert W'})^r$. Si $\mu_{j}=1$, la matrice a un noyau de dimension $1$, qui est l'intersection du sous-espace avec $\mathfrak{u}_{\tilde{x}}$. L'autre valeur propre est $2$ et la contribution du sous-espace \`a $d(\tilde{x})$ est $2$. On voit de m\^eme que la contribution de chaque sous-espace du premier ou du deuxi\`eme type est $2$. Il reste \`a compter le nombre de sous-espaces qui contribuent par $2$. On trouve $r^2+r+rdim(W'')$. D'o\`u
$$d(\tilde{x})=2^{r^2+r+rdim(W'')}det((1-x)_{\vert W'})^r,$$
ce qui ach\`eve la preuve de la premi\`ere \'egalit\'e de l'\'enonc\'e.

Pour  $I\in {\cal T}_{\tilde{x}}$, on note $\tilde{T}_{I}$ l'unique \'el\'ement de $\underline{\cal T}$ qui contient $\tilde{x}$ et est tel que $\tilde{T}_{I,\theta}=I$. Pour deux \'el\'ements $\tilde{T}_{1}$ et $\tilde{T}_{2}$ de $\underline{\cal T}$, on pose
$$W(\tilde{T}_{1},\tilde{T}_{2})=\{h\in H(F); h\tilde{T}_{1}h^{-1}=\tilde{T}_{2}\}/(T_{1}(F)\times SO(\tilde{\zeta}_{H,T_{1}})(F)).$$
On a les propri\'et\'es suivantes:

(6) soient $\tilde{T}\in {\cal T}$ et $\tilde{t}\in \tilde{T}(F)_{/\theta}$ en position g\'en\'erale; supposons $c_{\tilde{f}}(\tilde{t})\not=0$; alors il existe $I\in {\cal T}_{\tilde{x}}$ et $w\in W(\tilde{T}_{I},\tilde{T})$ tel que $\tilde{t}\in w(\tilde{x}exp(\mathfrak{i}(F)\cap \omega))$;

(7) soit $\tilde{T}\in {\cal T}$ et, pour $i=1,2$, soient $I_{i}\in {\cal T}_{\tilde{x}}$ et $w_{i}\in W(\tilde{T}_{I_{i}},\tilde{T})$; alors les sous-ensembles $w_{1}(\tilde{x}exp(\mathfrak{i}_{1}(F)\cap \omega))$ et $w_{2}(\tilde{x}exp(\mathfrak{i}_{2}(F)\cap \omega))$ de $\tilde{T}(F)_{/\theta}$ sont disjoints ou confondus;
 
 (8) soient $\tilde{T}\in {\cal T}$, $I_{1}\in {\cal T}_{\tilde{x}}$ et $w_{1}\in W(\tilde{T}_{I_{1}},\tilde{T})$; alors le nombre de couples $(I_{2},w_{2})$ tels que $I_{2}\in {\cal T}_{\tilde{x}}$, $w_{2}\in W(\tilde{T}_{I_{2}},\tilde{T})$ et $w_{2}(\tilde{x}exp(\mathfrak{i}_{2}(F)\cap \omega))=w_{1}(\tilde{x}exp(\mathfrak{i}_{1}(F)\cap \omega))$ est \'egal \`a
 $$\vert W(H_{\tilde{x}},I_{1})\vert [Z_{H}(\tilde{x})(F):H_{\tilde{x}}(F)] .$$
 
 Ces propri\'et\'es se prouvent comme (2), (3) et (4) (en se rappelant qu'ici $T(F)^{\theta}=T_{\theta}(F)$). On laisse les preuves au lecteur. Comme ci-dessus, elles permettent de transformer $J_{geom}(\Theta,\tilde{f})$  en l'expression
 $$J_{geom}(\Theta,\tilde{f})= [Z_{H}(\tilde{x})(F):H_{\tilde{x}}(F)]D^{\tilde{H}}(\tilde{x})\sum_{I\in {\cal T}_{\tilde{x}}}\vert W(H_{\tilde{x}},I)\vert ^{-1}\nu(I)$$
 $$\int_{\mathfrak{i}(F)}c_{\Theta,\tilde{x},\omega}(X)c_{\tilde{f},\tilde{x},\omega}(X)D^{H_{\tilde{x}}}(X)\Delta_{r}(\tilde{x}exp(X))dX.$$
 Soit $I\in {\cal T}_{\tilde{x}}$ et $X\in \mathfrak{i}(F)\cap \omega$, en position g\'en\'erale. Notons $W''=\,'W''\oplus \,''W''$ la d\'ecomposition de $W''$ associ\'ee \`a $I$, posons $\tilde{t}=\tilde{x}exp(X)$ et $t={^t\tilde{t}}^{-1}\tilde{t}$. On a
 $$\Delta_{r}(\tilde{t})=\vert 2\vert _{F}^{r^2+r+rdim(''W'')}\vert det((1-t)_{\vert W'\oplus 'W''})\vert _{F}^r.$$
 L'\'el\'ement $1-t$ agit sur $W'$ comme $1-x$. Il agit sur $'W''$ comme $1-exp(2X)$ et la valeur absolue du d\'eterminant de cette application est la m\^eme que celle du d\'eterminant de $2X$. C'est donc $\vert 2\vert _{F}^{dim('W'')}\Delta''(X)$. On obtient
 $$\Delta_{r}(\tilde{t})=\vert 2\vert _{F}^{r^2+r+rdim(W)}\vert det((1-x)_{\vert W'})\vert_{F}^r\Delta''(X)^r.$$
 Alors l'\'egalit\'e ci-dessus devient
 $$J_{geom}(\Theta,\tilde{f})=C(\tilde{x})J_{\tilde{x},\omega,geom}(\Theta,\tilde{f}).$$
 Cela ach\`eve la preuve. $\square$
 
 \bigskip
 
 \subsection{Une premi\`ere expression de $J_{\tilde{x},\omega,N}(\Theta,\tilde{f})$}
 
 Le lemme pr\'ec\'edent nous ram\`ene \`a prouver l'\'egalit\'e
 $$(1) \qquad lim_{N\to \infty}J_{\tilde{x},\omega,N}(\Theta,\tilde{f})=J_{\tilde{x},\omega,geom}(\Theta,\tilde{f}).$$
  Posons $H''_{\tilde{x}}=SO(W'')$, $G''_{\tilde{x}}=SO(V'')$. On a $H_{\tilde{x}}=H'_{\tilde{x}}H''_{\tilde{x}}$ et  $G_{\tilde{x}}=H'_{\tilde{x}}G''_{\tilde{x}}$. Un certain nombre d'objets relatifs \`a ces groupes se d\'ecomposent conform\'ement \`a ces d\'ecompositions en produits. On affectera leurs composantes d'un $'$ ou d'un $''$. Par exemple, on d\'ecompose tout \'el\'ement $X\in \mathfrak{h}_{\tilde{x}}(F)$ en $X=X'+X''$, avec $X'\in \mathfrak{h}'_{\tilde{x}}(F)$ et $X''\in \mathfrak{h}''_{\tilde{x}}(F)$. Puisque $\Theta$ est un quasi-caract\`ere, on peut supposer, en prenant $\omega$ assez petit, que $\Theta_{\tilde{x},\omega}$ est la restriction \`a $\mathfrak{h}_{\tilde{x}}(F)\cap \omega$ d'une combinaison lin\'eaire de transform\'ees de Fourier d'int\'egrales orbitales nilpotentes. Gr\^ace \`a la "conjecture de Howe", on peut remplacer les int\'egrales orbitales nilpotentes par des int\'egrales orbitales associ\'ees \`a des \'el\'ements semi-simples r\'eguliers et m\^eme "en position g\'en\'erale". L'\'egalit\'e que l'on veut prouver \'etant lin\'eaire en $\Theta$, on peut fixer un \'el\'ement $S\in \mathfrak{h}''_{\tilde{x}}(F)$, en position g\'en\'erale, et supposer que, pour tout $X\in \mathfrak{h}_{\tilde{x}}(F)\cap \omega$, on a l'\'egalit\'e
  $$\Theta_{\tilde{x},\omega}(X)=\hat{j}_{S}(X')\hat{j}^{H''_{\tilde{x}}}(S,X''),$$
  o\`u la seconde fonction est la transform\'ee de Fourier de l'int\'egrale orbitale associ\'ee \`a $S$ et la premi\`ere est une telle transform\'ee de Fourier associ\'ee \`a un \'el\'ement de $\mathfrak{h}'_{\tilde{x}}(F)$ que l'on n'a pas besoin de pr\'eciser. Pour $g\in G(F)$, on a l'\'egalit\'e
  $$J_{\tilde{x},\omega}(\Theta,\tilde{f},g)=\int_{\mathfrak{h}'_{\tilde{x}}(F)\times \mathfrak{h}''_{\tilde{x}}(F)}\Theta_{\tilde{x},\omega}(X){^g\tilde{f}}_{\tilde{x},\omega}^{\xi}(X)dX.$$
  En utilisant la formule de Weyl pour $\mathfrak{h}_{\tilde{x}}'(F)$, on obtient
  $$J_{\tilde{x},\omega}(\Theta,\tilde{f},g)=\sum_{I'\in {\cal T}(H'_{\tilde{x}})}\vert W(H'_{\tilde{x}},I')\vert ^{-1}\int_{\mathfrak{i}'(F)}\hat{j}_{S}(X')D^{H'_{\tilde{x}}}(X')$$
  $$\int_{I'(F)\backslash H'_{\tilde{x}}(F)}\int_{\mathfrak{h}''_{\tilde{x}}(F)}\hat{j}^{H''_{\tilde{x}}}(S,X''){^g\tilde{f}}_{\tilde{x},\omega}^{\xi}(h^{_{'}-1}X'h'+X'')dX''\,dh'\,dX'.$$
  D'o\`u
  $$J_{\tilde{x},\omega,N}(\Theta,\tilde{f})=\sum_{I'\in {\cal T}(H'_{\tilde{x}})}\vert W(H'_{\tilde{x}},I')\vert ^{-1}\int_{\mathfrak{i}'(F)}\hat{j}_{S}(X')D^{H'_{\tilde{x}}}(X')$$
  $$\int_{I'(F)H_{\tilde{x}}''(F)U_{\tilde{x}}(F)\backslash G(F)}\int_{\mathfrak{h}''_{\tilde{x}}(F)}\hat{j}^{H''_{\tilde{x}}}(S,X''){^g\tilde{f}}_{\tilde{x},\omega}^{\xi}(X'+X'')dX''\kappa_{N}(g)dg\,dX'.$$
  Les deux derni\`eres int\'egrales peuvent s'\'ecrire
  $$\int_{I'(F)G''_{\tilde{x}}(F)\backslash G(F)}\int_{H''_{\tilde{x}}(F)U_{\tilde{x}}(F)\backslash G''_{\tilde{x}}(F)}\int_{\mathfrak{h}''_{\tilde{x}}(F)} \hat{j}^{H''_{\tilde{x}}}(S,X''){^{g''g}\tilde{f}}_{\tilde{x},\omega}^{\xi}(X'+X'')dX''\kappa_{N}(g''g)dg''\,dg.$$
  On peut  utiliser le paragraphe 9 de [W2] pour calculer les deux int\'egrales int\'erieures ci-dessus. Pour retrouver les notations de cette r\'ef\'erence, on doit poser $v_{i}=z_{i}$ pour $i=0,...,r$, $v_{-i}=z_{-i}/(2\nu(-1)^{i})$ pour $i=1,...,r$ et d\'efinir des constantes $\xi_{i}$, $i=0,...,r-1$, de sorte que, pour $u\in U_{\tilde{x}}(F)$, on ait
  $$\sum_{i=-r,...,r-1}<z_{-i-1}^*,uz_{i}>=\sum_{i=0,...,r-1}\xi_{i}<v_{-i-1},\tilde{x}uv_{i}>.$$
  On a d\'efini en [W2] 9.2 et 9.3 un \'el\'ement $\Xi\in \mathfrak{g}''_{\tilde{x}}(F)$ et un sous-espace $\Sigma$ de $\mathfrak{g}''_{\tilde{x}}(F)$. Avec les notations ci-dessus, $\Xi$ est l'\'el\'ement de $\mathfrak{g}''_{\tilde{x}}(F)$ qui annule $W''$ et v\'erifie $\Xi v_{i+1}=\xi_{i}v_{i}$ pour tout $i=0,...,r-1$. L'espace $\Sigma$ est la somme directe de $\mathfrak{a}(F)\cap \mathfrak{g}''_{\tilde{x}}(F)$, de l'orthogonal de $\mathfrak{h}''_{\tilde{x}}(F)$ dans $\mathfrak{g}''_{\tilde{x}}(F)\cap \mathfrak{g}_{0}(F)$ et de $\mathfrak{u}_{\tilde{x}}(F)$. Pour tout sous-tore maximal $I''$ de $G''_{\tilde{x}}$, on note $\mathfrak{i}''(F)^S$ le sous-ensemble des \'el\'ements de $\mathfrak{i}''(F)$ qui sont conjugu\'es \`a un \'el\'ement de $\Xi+S+\Sigma$ par un \'el\'ement de $G''_{\tilde{x}}(F)$. Pour tout $X''\in \mathfrak{i}''(F)^S$, on fixe un \'el\'ement $\gamma_{X''}\in G''_{\tilde{x}}(F)$ tel que $\gamma_{X''}^{-1}X''\gamma_{X''}\in \Xi+S+\Sigma$. On peut imposer \`a $\gamma_{X''}$ diverses contraintes de croissance et de r\'egularit\'e. On introduit la fonction ${^g\tilde{f}}_{\tilde{x},\omega}^{\sharp}$ comme en 1.8 et la formule 9.8(2) de [W2] appliqu\'ee \`a nos deux derni\`eres int\'egrales ci-dessus conduit \`a l'\'egalit\'e:
      $$(2) \qquad J_{\tilde{x},\omega,N}(\Theta,\tilde{f})=\sum_{I\in {\cal T}(G_{\tilde{x}})}\nu(A_{I''})^{-1}\vert W(G_{\tilde{x}},I)\vert ^{-1}\int_{\mathfrak{i}'(F)\times \mathfrak{i}''(F)^S}\hat{j}_{S}(X')D^{G'_{\tilde{x}}}(X')D^{G''_{\tilde{x}}}(X'')^{1/2}$$
      $$\int_{I'(F)A_{I''}(F)\backslash G(F)}{^g\tilde{f}}_{\tilde{x},\omega}^{\sharp}(X'+X'')\kappa_{N,X''}(g)dg\,dX''\,dX',$$
    o\`u
      $$\kappa_{N,X''}(g)=\nu(A_{I''})\int_{A_{I''}(F)}\kappa_{N}(\gamma_{X''}^{-1}ag)da.$$
      On voit comme en [W2] 10.1 que cette expression est absolument convergente et born\'ee par une puissance de $N$.
      
      \bigskip
      \subsection{Changement de fonction de troncature}

      Fixons $I\in {\cal T}(G_{\tilde{x}})$.  On a d\'efini en [W2] 10.1 trois polyn\^omes non nuls sur $\mathfrak{i}(F)$. Pour $\epsilon>0$, on note $\mathfrak{i}(F)[\leq\epsilon]$ l'ensemble des $X\in \mathfrak{i}(F)$ pour lesquels on a $\vert Q(X)\vert _{F}\leq \epsilon$ pour l'un au moins de ces polyn\^omes $Q$. On note $\mathfrak{i}(F)[>\epsilon]$ le compl\'ementaire. On a
      
      (1) il existe un entier $b\geq1$ tel que 
      $$\int_{(\mathfrak{i}'(F)\times \mathfrak{i}''(F)^S)\cap \mathfrak{i}(F)[\leq N^{-b}]}\vert \hat{j}_{S}(X')\vert D^{G'_{\tilde{x}}}(X')D^{G''_{\tilde{x}}}(X'')^{1/2}$$
      $$\int_{I'(F)A_{I''}(F)\backslash G(F)}\vert {^g\tilde{f}}_{\tilde{x},\omega}^{\sharp}(X'+X'')\vert \kappa_{N,X''}(g)dg\,dX''\,dX'\,<<\,N^{-1}$$
      pour tout $N\geq1$.
      
  La notation $<<$ signifie qu'il existe une constante $c>0$ telle que le membre de gauche soit  inf\'erieure ou \'egale \`a $c$ fois le membre de droite. La relation ci-dessus se prouve comme le (ii) du lemme 10.1 de [W2]. On fixe $b$ v\'erifiant cette relation.

      Notons $\tilde{M}_{\natural}$ le commutant de $A_{I''}$ dans $\tilde{G}$. C'est un L\'evi tordu de $\tilde{G}$ qui contient $G'\tilde{x}$. On a l'\'egalit\'e $A_{\tilde{M}_{\natural}}=A_{I''}$. Fixons un L\'evi minimal  $M_{min}$ de $G$ contenu dans $M_{\natural}$ et un sous-groupe parabolique $P_{min}=M_{min}U_{min}\in {\cal P}(M_{min}) $. Pour simplifier, on peut supposer que $K$ est en bonne position   relativement \`a $M_{min}$.  On utilise le groupe $K$ pour d\'efinir les fonctions $H_{\tilde{Q}}$ sur $G(F)$, pour $\tilde{Q}\in {\cal L}^{\tilde{G}}$. Notons $\Delta$ l'ensemble des racines simples de $A_{M_{min}}$ dans $\mathfrak{u}_{min}$. Fixons $Y\in {\cal A}_{M_{min}}$. Pour tout $P'\in {\cal P}(M_{min})$, il y a un unique $w\in W^G$ tel que $wP_{min}w^{-1}=P'$. On pose $Y_{P'}=wY$. Pour $\tilde{Q}\in {\cal P}(\tilde{M}_{\natural})$, notons $Y_{\tilde{Q}}$ la projection orthogonale de $Y_{P'}$ sur ${\cal A}_{\tilde{M}_{\natural}}$, o\`u $P'$ est un \'el\'ement quelconque de ${\cal P}(M_{min})$ tel que $P'\subset Q$.
    Soit $g\in G(F)$.  On pose $Y(g)_{\tilde{Q}}=Y_{\tilde{Q}}-H_{\bar{\tilde{Q}}}(g)$, o\`u $\bar{\tilde{Q}}$ est le parabolique tordu oppos\'e \`a $\tilde{Q}$. La famille ${\cal Y}(g)=(Y(g)_{\tilde{Q}})_{\tilde{Q}\in {\cal P}(\tilde{M}_{\natural})}$ est $(\tilde{G},\tilde{M}_{\natural})$-orthogonale. Il existe une constante $c>0$ telle que, si 
      $$\inf\{\sigma(mg); m\in M_{\natural}(F)\}<c\,inf\{\alpha(Y ); \alpha\in\Delta\},$$
       on ait $\alpha(Y(g)_{\tilde{Q}})>0$ pour tout $\tilde{Q}=\tilde{M}_{\natural}U_{Q}\in {\cal P}(\tilde{M}_{\natural})$ et toute racine $\alpha$ de $A_{\tilde{M}_{\natural}}$ dans $\mathfrak{u}_{Q}$. Supposons cette condition v\'erifi\'ee. On note $\lambda\mapsto \sigma_{\tilde{M}_{\natural}}^{\tilde{G}}(\lambda,{\cal Y}(g))$ la fonction caract\'eristique dans ${\cal A}_{\tilde{M}_{\natural}}$ de l'enveloppe convexe des points de la famille ${\cal Y}(g)$. On pose
      $$v(g)=\nu(A_{I''})\int_{A_{I''}(F)}\sigma_{\tilde{M}_{\natural}}^{\tilde{G}}(H_{\tilde{M}_{\natural}}(a),{\cal Y}(g))da.$$
      
      On a
      
      (2) il existe $c>0$ et un sous-ensemble compact $\omega_{I}$ de $\mathfrak{i}(F)$ tels que les propri\'et\'es suivantes soient v\'erifi\'ees; soient $g\in G(F)$ et $X\in \mathfrak{i}(F)[>N^{-b}]$ tels que ${^g\tilde{f}}_{\tilde{x},\omega}^{\sharp}(X)\not=0$; alors $X\in \omega_{I}$ et il existe $t\in I(F)$ tel que $\sigma(tg)<c\,log(N)$.
      
      La preuve est la m\^eme que celle de [W2] 10.4(2). Cette propri\'et\'e assure que le membre de droite de l'\'egalit\'e ci-dessous est bien d\'efini.
      
      \ass{Proposition}{Il existe $c>0$ et un entier $N_{0}\geq1$ tels que, si $N\geq N_{0}$ et
      $$c\,log(N)\,<\,inf\{\alpha(Y); \alpha\in\Delta\},$$
      on ait l'\'egalit\'e
      $$\int_{I'(F)A_{I''}(F)\backslash G(F)}{^g\tilde{f}}_{\tilde{x},\omega}^{\sharp}(X)\kappa_{N,X''}(g)dg
=\int_{I'(F)A_{I''}(F)\backslash G(F)}{^g\tilde{f}}_{\tilde{x},\omega}^{\sharp}(X)v(g)dg$$
pour tout $X\in \mathfrak{i}(F)[>N^{-b}]\cap (\mathfrak{i}'(F)\times \mathfrak{i}''(F)^S)$.}

Preuve.  Pour simplifier la notation, fixons un ensemble $\omega_{I}$ comme en (2), posons 
$${\cal X}_{N}=\omega_{I}\cap  \mathfrak{i}(F)[>N^{-b}]\cap (\mathfrak{i}'(F)\times \mathfrak{i}''(F)^S)$$
et notons $\Omega_{N}$ l'ensemble des $g\in G(F)$ pour lesquels il existe $X\in {\cal X}_{N}$ tel que ${^g\tilde{f}}_{\tilde{x},\omega}^{\sharp}(X)\not=0$.
 Soit $Z\in {\cal A}_{M_{min}}$ tel que $\alpha(Z)>0$ pour tout $\alpha\in \Delta$ (on ne confond pas cet \'el\'ement avec le sous-espace $Z$ de $V$). En rempla\c{c}ant $Y$ par cet \'el\'ement, on construit une famille de points ${\cal Z}(g)$ pour tout $g\in G(F)$. On suppose
$$inf\{\alpha(Z); \alpha\in\Delta\}\geq c_{1}log(N),$$
ce qui assure que, pour  $g\in \Omega_{N}$, on a $Z(g)_{\tilde{Q}}\in {\cal A}_{\tilde{Q}}^+$ pour tout $\tilde{Q}\in {\cal P}(\tilde{M}_{\natural})$. Pour $\tilde{Q}=\tilde{L}U_{Q}\in {\cal P}(\tilde{M}_{\natural})$, notons $\lambda\mapsto \sigma_{\tilde{M}_{\natural}}^{\tilde{Q}}(\lambda,{\cal Z}(g))$ la fonction caract\'eristique dans ${\cal A}_{\tilde{M}_{\natural}}$ de la somme de ${\cal A}_{\tilde{L}}$ et de l'enveloppe convexe des $Z(g)_{\tilde{P}'}$ pour $\tilde{P}'\in {\cal P}(\tilde{M}_{min})$, $\tilde{P}'\subset \tilde{Q}$. Notons $\tau_{\tilde{Q}}$ la fonction caract\'eristique du sous-ensemble ${\cal A}_{\tilde{M}_{\natural}}^{\tilde{L}}\oplus {\cal A}_{\tilde{Q}}^+$ de ${\cal A}_{\tilde{M}_{\natural}}$. On a l'\'egalit\'e
$$\sum_{\tilde{Q}\in {\cal F}(\tilde{M}_{\natural})}\sigma_{\tilde{M}_{\natural}}^{\tilde{Q}}(\lambda,{\cal Z}(g))\tau_{\tilde{Q}}(\lambda-Z(g)_{\tilde{Q}})=1$$
pour tout $\lambda\in {\cal A}_{\tilde{M}_{\natural}}$. Pour $g\in  \Omega_{N}$ et $X\in {\cal X}_{N}$, on peut donc \'ecrire
$$v(g)=\nu(A_{I''})\sum_{\tilde{Q}\in {\cal F}(\tilde{M}_{\natural})}v(\tilde{Q},g),$$
$$\kappa_{N,X''}(g)=\nu(A_{I''})\sum_{\tilde{Q}\in {\cal F}(\tilde{M}_{\natural})}\kappa_{N,X''}(\tilde{Q},g),$$
o\`u
$$v(\tilde{Q},g)=\int_{A_{I''}(F)}\sigma_{\tilde{M}_{\natural}}^{\tilde{G}}(H_{\tilde{M}_{\natural}}(a),{\cal Y}(g))\sigma_{\tilde{M}_{\natural}}^{\tilde{Q}}(H_{\tilde{M}_{\natural}}(a),{\cal Z}(g))\tau_{\tilde{Q}}(H_{\tilde{M}_{\natural}}(a)-Z(g)_{\tilde{Q}})da,$$
$$\kappa_{N,X''}(\tilde{Q},g)=\int_{A_{I''}(F)}\kappa_{N}(\gamma_{X''}^{-1}ag)\sigma_{\tilde{M}_{\natural}}^{\tilde{Q}}(H_{\tilde{M}_{\natural}}(a),{\cal Z}(g))\tau_{\tilde{Q}}(H_{\tilde{M}_{\natural}}(a)-Z(g)_{\tilde{Q}})da.$$
Les fonctions $g\mapsto v(\tilde{Q},g)$ et $g\mapsto \kappa_{N,X''}(\tilde{Q},g)$ sont invariantes \`a gauche par $I'(F)A_{I''}(F)$. On peut donc d\'ecomposer le membre de gauche, resp. de droite, de l'\'egalit\'e de l'\'enonc\'e en
$$\nu(A_{I''})\sum_{\tilde{Q}\in {\cal F}(\tilde{M}_{\natural})}\Phi_{1}(\tilde{Q},X),$$
resp.
$$\nu(A_{I''})\sum_{\tilde{Q}\in {\cal F}(\tilde{M}_{\natural})}\Phi_{2}(\tilde{Q},X),$$
o\`u
$$\Phi_{1}(\tilde{Q},X)=\int_{I'(F)A_{I''}(F)\backslash G(F)}{^g\tilde{f}}_{\tilde{x},\omega}^{\sharp}(X)\kappa_{N,X''}(\tilde{Q},g)dg,$$
$$\Phi_{2}(\tilde{Q},X)=\int_{I'(F)A_{I''}(F)\backslash G(F)}{^g\tilde{f}}_{\tilde{x},\omega}^{\sharp}(X)v(\tilde{Q},g)dg.$$

On a:

(3) si $Z$ v\'erifie les in\'egalit\'es
$$sup\{\alpha(Z); \alpha\in\Delta\}\leq \left\lbrace\begin{array}{c}inf\{\alpha(Y); \alpha\in\Delta\},\\ log(N)^2\\ \end{array}\right.$$
on a $\Phi_{1}(\tilde{G},X)=\Phi_{2}(\tilde{G},X)$ pour tout $X\in {\cal X}_{N}$, pourvu que $N$ soit assez grand. 

En effet, soit $g\in  \Omega_{N}$. On v\'erifie comme en [W2] 10.4(8) que, sous ces hypoth\`eses, on a $\kappa_{N}(\gamma_{X''}^{-1}ag)=\sigma_{\tilde{M}_{\natural}}^{\tilde{G}}(H_{\tilde{M}_{\natural}}(a),{\cal Y}(g)) =1$ pour tout $a\in A_{I''}(F)$ tel que $\sigma_{\tilde{M}_{\natural}}^{\tilde{G}}(H_{\tilde{M}_{\natural}}(a),{\cal Z}(g))=1$. Donc $\kappa_{N,X''}(\tilde{G},g)=v(\tilde{G},g)$ et la conclusion.

Fixons $\tilde{Q}=\tilde{L}U_{Q}\in {\cal F}(\tilde{M}_{\natural})$, $\tilde{Q}\not=\tilde{G}$. On a:

(4) il existe $c>0$ tel que, si $c\,log(N)<inf \{\alpha(Z); \alpha\in\Delta\}$, on a $\Phi_{1}(\tilde{Q},X)=\Phi_{2}(\tilde{Q},X)=0$ pour tout $X\in {\cal X}_{N}$.

On \'ecrit
$$\Phi_{1}(\tilde{Q},X)=\int_{K_{min}}\int_{I'(F)A_{I''}(F)\backslash L(F)}\int_{U_{\bar{Q}}(F)}{^{\bar{u}lk}\tilde{f}}_{\tilde{x},\omega}^{\sharp}(X)\kappa_{N,X''}(\tilde{Q},\bar{u}lk)d\bar{u}\delta_{Q}(l)dl\,dk.$$
On v\'erifie que, si ${^{\bar{u}lk}\tilde{f}}_{\tilde{x},\omega}^{\sharp}(X)\not=0$ pour un $X\in {\cal X}_{N}$, on a une majoration $\sigma(\bar{u})<<log(N)$ et on peut repr\'esenter $l$ par un \'el\'ement tel que $\sigma(l)<<log(N)$. Le point est que, pour de tels \'el\'ements, on a l'\'egalit\'e 
$$\kappa_{N,X''}(\tilde{Q},\bar{u}lk)=\kappa_{N,X''}(\tilde{Q},lk).$$
Cela r\'esulte de l'assertion:

(5) soit $c>0$; il existe $c'>0 $ tel que, si $c'\,log(N)<inf \{\alpha(Z); \alpha\in\Delta\}$, les propri\'et\'es suivantes sont v\'erifi\'ees; soient $X\in {\cal X}_{N}$, $g\in G(F)$ et $\bar{u}\in U_{\bar{Q}}(F)$; supposons $\sigma(g),\sigma(\bar{u}),\sigma(\bar{u}g)<clog(N)$; alors on a l'\'egalit\'e $\kappa_{N,X''}(\tilde{Q},\bar{u}g)=\kappa_{N,X''}(\tilde{Q},g)$.

Reportons la preuve de cette assertion \`a plus tard. A l'int\'erieur de l'int\'egrale ci-dessus appara\^{\i}t donc l'int\'egrale
$$\int_{U_{\bar{Q}}(F)}{^{\bar{u}lk}\tilde{f}}_{\tilde{x},\omega}^{\sharp}(X)d\bar{u}.$$
Or celle-ci est nulle d'apr\`es 1.8(1). Donc $\Phi_{1}(\tilde{Q},X)=0$. On prouve de m\^eme que $\Phi_{2}(\tilde{Q},X)=0$, en utilisant une assertion similaire \`a (5) pour la fonction $v(\tilde{Q},g)$.

On a utilis\'e l'\'el\'ement auxiliaire $Z$. Il existe $c>0$ tel que, si
 $$c\,log(N)\,<\,inf\{\alpha(Y); \alpha\in\Delta\},$$
 on peut choisir cet \'el\'ement auxiliaire v\'erifiant les hypoth\`eses de (3) et (4). Ces relations entra\^{\i}nent la conclusion de l'\'enonc\'e.
 
 Il reste \`a prouver l'assertion  (5) et son analogue pour la fonction $v(\tilde{Q},g)$. La preuve de cette analogue est la m\^eme que celle de [W2] 10.5 et on ne la fera pas. Soient $c$, $X$, $g$, $\bar{u}$ v\'erifiant les hypoth\`eses de (5). Consid\'erons la formule int\'egrale qui d\'efinit $\kappa_{N,X''}(\tilde{Q},g)$. La fonction
 $$\lambda\mapsto \sigma_{\tilde{M}_{\natural}}^{\tilde{Q}}(\lambda,{\cal Z}(g))\tau_{\tilde{Q}}(\lambda-Z(g)_{\tilde{Q}})$$
 ne d\'epend  de $g$ que par l'interm\'ediaire des $H_{\bar{\tilde{P}}'}(g)$ pour $\tilde{P}'\in {\cal P}(\tilde{M}_{\natural})$ , $\tilde{P}'\subset \tilde{Q}$. Or ces termes sont insensibles au changement de $g$ en $\bar{u}g$. On en d\'eduit:
 $$(6) \qquad \kappa_{N,X''}(\tilde{Q},\bar{u}g)-\kappa_{N,X''}(\tilde{Q},g)=\int_{A_{I''}(F)}\sigma_{\tilde{M}_{\natural}}^{\tilde{Q}}(H_{\tilde{M}_{\natural}}(a),{\cal Z}(g))\tau_{\tilde{Q}}(H_{\tilde{M}_{\natural}}(a)-{\cal Z}(g)_{\tilde{Q}})$$
 $$(\kappa_{N}(\gamma_{X''}^{-1}a\bar{u}g)-\kappa_{N}(\gamma_{X''}^{-1}ag))da.$$
 Supposons prouv\'ee l'assertion suivante:
 
 (7) il existe $c''>0$, ne d\'ependant que de $c$ (et pas de $g$, $\bar{u}$, $X$), tel que, pour tout $a\in A_{I''}(F)$ tel que $\alpha(H_{\tilde{M}_{\natural}} (a))>c'' log(N)$ pour toute racine $\alpha$ de $A_{I''}$ dans $\mathfrak{u}_{Q}$, on a l'\'egalit\'e $\kappa_{N}(\gamma_{X''}^{-1}a\bar{u}g)=\kappa_{N}(\gamma_{X''}^{-1}ag)$.
 
  On peut trouver $c'>0$, ne d\'ependant que de $c$, tel que les deux conditions
  $$c'\,log(N)<inf \{\alpha(Z); \alpha\in\Delta\}$$
  et
  $$\sigma_{\tilde{M}_{\natural}}^{\tilde{Q}}(H_{\tilde{M}_{\natural}}(a),{\cal Z}(g))\tau_{\tilde{Q}}(H_{\tilde{M}_{\natural}}(a)-{\cal Z}(g)_{\tilde{Q}})=1$$
  entra\^{\i}nent $\alpha(H_{\tilde{M}_{\natural} }(a))>c'' log(N)$ pour toute racine $\alpha$ de $A_{I''}$ dans $\mathfrak{u}_{Q}$. Si l'on impose la premi\`ere condition \`a $Z$, la fonction que l'on int\`egre dans la formule (6) est donc nulle, d'o\`u la conclusion de (5). Cela nous ram\`ene \`a prouver (7). La preuve est essentiellement la m\^eme que celle du lemme 10.8 de [W2]. Signalons  seulement les  points \`a modifier. On voit comme dans cette r\'ef\'erence que l'on peut \'etendre les scalaires. Remarquons qu'en \'etendant les scalaires, on peut changer le tore $A_{I''}$ mais celui-ci ne peut que devenir plus gros et l'assertion (7) n'en devient que plus forte. Apr\`es extension des scalaires, on peut supposer que le groupe sp\'ecial orthogonal $G_{\tilde{x}}''$ de la restriction de $\tilde{x}$ \`a $V''$  est quasi-d\'eploy\'e et que $I''=A_{I''}$ en est un tore d\'eploy\'e maximal. On peut fixer un \'el\'ement $\tilde{P}_{\natural}\in {\cal P}(\tilde{M}_{\natural})$ et se limiter aux \'el\'ements $a\in A_{I''}(F)$ tels que $H_{\tilde{M}_{\natural}}(a)$ appartient \`a la chambre positive ferm\'ee pour ce parabolique tordu. Montrons que:
  
  (8) il existe $\delta\in G''_{\tilde{x}}(F)$ tel que $\delta\tilde{P}_{\natural}\delta^{-1}\subset \bar{\tilde{P}}$ et $A_{\tilde{M}}\subset \delta A_{I''}\delta^{-1}$.
  
  Le tore d\'eploy\'e $A_{\tilde{M}}$ est contenu dans $G''_{\tilde{x}}$ et $A_{I''}$ est un sous-tore maximal de ce groupe. On peut donc trouver $\delta\in G''_{\tilde{x}}(F)$ tel que  $A_{\tilde{M}}\subset \delta A_{I''}\delta^{-1}$. Alors $\delta^{-1}\bar{\tilde{P}}\delta\in {\cal F}(\tilde{M}_{\natural})$. Fixons $\tilde{P}'\in {\cal P}(\tilde{M}_{\natural})$ tel que $\tilde{P}'\subset \delta^{-1}\bar{\tilde{P}}\delta$. Le tore $A_{I''}$ \'etant un sous-tore maximal d'un groupe sp\'ecial orthogonal, il a une forme particuli\`ere: ses sous-espaces propres dans $V$ sont tous de dimension $1$, sauf \'eventuellement celui sur lequel $A_{I''}$ agit par l'identite. Le L\'evi $\tilde{M}_{\natural}$ \'etant le commutant de ce tore a lui-aussi une forme particuli\`ere: $M_{\natural}$ est un produit de $GL_{1}$ et \'eventuellement d'un unique bloc $GL_{k}$. Pour un tel L\'evi, deux \'el\'ements quelconques de ${\cal F}(\tilde{M}_{\natural})$ sont conjugu\'es par le groupe $Norm_{G(F)}(\tilde{M}_{\natural})$. Si $d$ est impair, l'application naturelle
  $$Norm_{G''_{\tilde{x}}(F)}(A_{I''})\to Norm_{G(F)}(\tilde{M}_{\natural})/M_{\natural}(F)$$
  est surjective. Donc $\tilde{P}_{\natural}$ et $\tilde{P}'$ sont conjugu\'es par un \'el\'ement de $Norm_{G''_{\tilde{x}}(F)}(A_{I''})$. Quitte \`a multiplier $\delta$ \`a droite par un tel \'el\'ement, on peut supposer $\tilde{P}'=\tilde{P}_{\natural}$ et la conclusion de (8) est v\'erifi\'ee. Supposons $d$ pair et introduisons le groupe orthogonal $O(V'')$ de la restriction \`a $V''$ de la forme $\tilde{x}$ ($G''_{\tilde{x}}$ en est sa composante neutre). L'application naturelle
  $$Norm_{O(V'')(F)}(A_{I''})\to Norm_{G(F)}(\tilde{M}_{\natural})/M_{\natural}(F)$$
  et on peut comme ci-dessus multiplier $\delta$ \`a droite par un \'el\'ement  de $Norm_{O(V'')(F)}(A_{I''})$ de sorte que $\tilde{P}'=\tilde{P}_{\natural}$. Si cet \'el\'ement appartient \`a $G''_{\tilde{x}}(F)$, on a fini. Sinon, on remarque que $A_{\tilde{M}}$ n'est pas maximal dans $\tilde{G}_{\tilde{x}}$, plus pr\'ecis\'ement, il agit trivialement sur un sous-espace non nul de $V''$: il fixe $z_{0}$. Il en r\'esulte qu'il existe un \'el\'ement $y\in Norm_{O(V'')(F)}(A_{I''})\cap \delta^{-1}M\delta$ tel que le d\'eterminant de $y$ agissant dans $V''$ soit $-1$. En multipliant $\delta$ \`a droite par un tel \'el\'ement, on obtient la conclusion de (8).
  
  Gr\^ace \`a (8), le m\^eme raisonnement qu'en [W2] 10.8 nous ram\`ene au cas o\`u $r=0$. Dans ce cas, $\kappa_{N}$ est par d\'efinition  le produit de deux fonctions $\kappa_{1,N}\kappa_{2,N}$. La premi\`ere est la fonction caract\'eristique de l'ensemble des $g\in G(F)$ tels que $g^{-1}z_{0}\in \mathfrak{p}_{F}^{-N}R$. La seconde est la fonction caract\'eristique de l'ensemble des $g\in G(F)$ tels que ${^tg}z_{0}^*\in \mathfrak{p}_{F}^{-N}R^{\vee}$. La d\'emonstration de [W2] 10.8 s'applique sans changement pour la premi\`ere fonction, c'est-\`a-dire que l'on d\'emontre l'analogue de (7) o\`u la conclusion est remplac\'ee par $\kappa_{1,N}(\gamma_{X''}^{-1}a\bar{u}g)=\kappa_{1,N}(\gamma_{X''}^{-1}ag)$.  Consid\'erons la seconde fonction. On a $\tilde{x}(z_{0})=\tilde{\zeta}(z_{0})=2\nu z_{0}^*$. La relation $\kappa_{N,2}(\gamma_{X''}^{-1}a\bar{u}g)=1$ \'equivaut donc \`a ${^tg}{^t\bar{u}}{^ta}{^t\gamma_{X''}^{-1}}\tilde{x}z_{0}\in 2\nu\mathfrak{p}_{F}^{-N}R^{\vee}$. Les produits sont ici les produits naturels et non pas les produits dans $\tilde{G}$. Faisons commuter $\tilde{x}$ aux \'el\'ements qui sont \`a sa gauche. La relation pr\'ec\'edente \'equivaut \`a $\tilde{x}g^{_{'}-1}\bar{u}^{_{'}-1}a^{_{'}-1}\gamma_{X''}'z_{0}\in 2\nu\mathfrak{p}_{F}^{-N}R^{\vee}$, o\`u $g'=\theta_{\tilde{x}}^{-1}(g)$, $\bar{u}'=\theta_{\tilde{x}}^{-1}(\bar{u})$, $a'=\theta_{\tilde{x}}^{-1}(a)$, $\gamma_{X''}'=\theta_{\tilde{x}}^{-1}(\gamma_{X''})$. Elle \'equivaut encore \`a $g^{_{'}-1}\bar{u}^{_{'}-1}a^{_{'}-1}\gamma_{X''}'z_{0}\in \mathfrak{p}_{F}^{-N}R'$, o\`u $R'=2\nu \tilde{x}^{-1}(R^{\vee})$. Mais $a'=a$ et $\gamma_{X''}'=\gamma_{X''}$ car ces \'el\'ements sont dans $G_{\tilde{x}}(F)$. Alors la condition pr\'ec\'edente est du m\^eme type que la condition $\kappa_{N,1}(\gamma_{X''}^{-1}a\bar{u}g)=1$: on a simplement remplac\'e $\bar{u}$ par $\bar{u}'$, $g$ par $g'$ et $R$ par $R'$. Ces objets v\'erifient les m\^emes hypoth\`eses que $\bar{u}$, $g$ et $R$. Alors la m\^eme d\'emonstration qu'en [W2] 10.8 s'applique et on d\'emontre l'analogue de (7) o\`u la conclusion est remplac\'ee par $\kappa_{2,N}(\gamma_{X''}^{-1}a\bar{u}g)=\kappa_{2,N}(\gamma_{X''}^{-1}ag)$. Evidemment, les deux analogues de (7) pour les fonctions $\kappa_{1,N}$ et $\kappa_{2,N}$ entra\^{\i}nent l'assertion (7) elle-m\^eme. Cela ach\`eve la preuve. $\square$
  
  \bigskip
  
  \subsection{Apparition des int\'egrales orbitales pond\'er\'ees}

  Le tore $I$ et la constante $b$ sont fix\'es comme dans le paragraphe pr\'ec\'edent.
  
  \ass{Lemme}{Il existe $N_{0}\geq1$ tel que, pour tout $N\geq N_{0}$ et tout $X\in \mathfrak{i}(F)[>N^{-b}]\cap (\mathfrak{i}'(F)\times \mathfrak{i}''(F)^S)$, on ait les \'egalit\'es
  $$\int_{I'(F) A_{I''}(F)\backslash G(F)}{^g\tilde{f}}_{\tilde{x},\omega}^{\sharp}(X)\kappa_{N,X''}(g)dg=0,$$
  si $A_{I'}\not=\{1\}$;
  $$\int_{I'(F) A_{I''}(F)\backslash G(F)}{^g\tilde{f}}_{\tilde{x},\omega}^{\sharp}(X)\kappa_{N,X''}(g)dg=\nu(I')\nu(A_{I''})\Theta^{J,\sharp}_{\tilde{f},\tilde{x},\omega}(X)$$
  si $A_{I'}=\{1\}$.}
  
  Preuve. Notons $v(g,Y)$ la fonction not\'ee $v(g)$ dans le paragraphe pr\'ec\'edent.
  La proposition 3.6 nous autorise \`a remplacer $\kappa_{N,X''}(g)$ par $v(g,Y)$ dans les int\'egrales de l'\'enonc\'e, pourvu que $Y$ v\'erifie une minoration
  $$inf\{\alpha(Y); \alpha\in\Delta\}>> log(N).$$
  Fixons $X$. Les int\'egrales sont  \`a support compact. On peut faire tendre $Y$ vers l'infini. Dans le cas non tordu, Arthur a calcul\'e en [A1] p.46 le comportement de $v(g,Y)$ quand $Y$ tend vers l'infini. Le m\^eme calcul vaut dans le cas tordu. Pour $Y$ dans un certain r\'eseau ${\cal R}\subset {\cal A}_{M_{min}}$, la fonction $Y\mapsto v(g,Y)$ est une somme de fonctions $Y\mapsto q_{\zeta}(Y)exp(\zeta(Y))$, o\`u $q_{\zeta}$ est un polyn\^ome et $\zeta\in Hom({\cal R},2\pi i{\mathbb Q}/2\pi i{\mathbb Z})$. De telles fonctions sont lin\'eairement ind\'ependantes. Puisque l'expression que l'on calcule est ind\'ependante de $Y$, on peut aussi bien remplacer $v(g,Y)$ par son "terme constant" $q_{0}(0)$. D'apr\`es [A1] p.92 (adapt\'e au cas tordu), on a
  $$q_{0}(0)=(-1)^{a_{\tilde{M}_{\natural}}}\sum_{\tilde{Q}\in {\cal F}(\tilde{M}_{\natural})} c'_{\tilde{Q}}v_{\tilde{M}_{\natural}}^{\tilde{Q}}(g).$$
  Les $c'_{\tilde{Q}}$ sont des constantes et on a $c'_{\tilde{G}}=1$. Cela conduit \`a l'\'egalit\'e
  $$ \int_{I'(F) A_{I''}(F)\backslash G(F)}{^g\tilde{f}}_{\tilde{x},\omega}^{\sharp}(X)\kappa_{N,X''}(g)dg=(-1)^{a_{\tilde{M}_{\natural}}}\sum_{\tilde{Q}\in {\cal F}(\tilde{M}_{\natural})}c'_{\tilde{Q}}\Phi(\tilde{Q}),$$ 
 o\`u
 $$\Phi(\tilde{Q})=\int_{I'(F) A_{I''}(F)\backslash G(F)}{^g\tilde{f}}_{\tilde{x},\omega}^{\sharp}(X)v_{\tilde{M}_{\natural}}^{\tilde{Q}}(g)dg.$$
 Pour $\tilde{Q}=\tilde{L}U_{Q}\not=\tilde{G}$, on d\'ecompose l'int\'egrale sur $I'(F) A_{I''}(F)\backslash G(F)$ en produit d'int\'egrales sur $I'(F)A_{I''}(F)\backslash L(F)$, $U_{Q}(F)$ et $K$. On voit appara\^{\i}tre une int\'egrale
 $$\int_{U_{Q}(F)}{^{ulk}\tilde{f}}_{\tilde{x},\omega}^{\sharp}(X)du,$$
 qui est nulle d'apr\`es 1.8(1). Il ne reste plus que le terme pour $\tilde{Q}=\tilde{G}$, qui vaut
 $$\Phi(\tilde{G})=mes(I(F)/I'(F)A_{I''}(F))D^{\tilde{G}_{\tilde{x}}}(X)^{-1/2} J^{\sharp}_{\tilde{M}_{\natural},\tilde{x},\omega}(X,\tilde{f}).$$
 Si $A_{I'}\not=\{1\}$, on a $A_{\tilde{M}_{\natural}}=A_{I''}\subsetneq A_{\tilde{G}_{\tilde{x}exp(X)}}=A_{I'}A_{I''}$, donc $ J^{\sharp}_{\tilde{M}_{\natural},\tilde{x},\omega}(X,\tilde{f})=0$ d'apr\`es 1.8(2). D'o\`u la premi\`ere assertion de l'\'enonc\'e. Supposons $A_{I'}=\{1\}$. Alors $\tilde{M}_{\natural}$ est le L\'evi tordu not\'e $\tilde{M}(X)$ en 1.8 et on obtient
 $$\Phi(\tilde{G})=(-1)^{a_{\tilde{M}_{\natural}}}\nu(I)mes(I(F)/I'(F)A_{I''}(F))\Theta^{J,\sharp}_{\tilde{f},\tilde{x},\omega}(X).$$
 On rappelle que l'on n'utilise pas les m\^emes mesures qu'en 1.8 (cf. 3.3), ce qui explique l'apparition de $\nu(I)$. Il reste \`a v\'erifier que
 $$\nu(I)mes(I(F)/I'(F)A_{I''}(F))=\nu(I')\nu(A_{I''})$$
 pour obtenir le (ii) de l'\'enonc\'e. $\square$

\bigskip

\subsection{Preuve du th\'eor\`eme 3.3}

Si ${\cal A}_{H'_{\tilde{x}}}\not=\{1\}$, on a $A_{I'}\not=\{1\}$ pour tous les tores $I$ intervenant dans la formule 3.5(2). En utilisant la relation 3.6(1) et le lemme 3.7, on voit que
$$lim_{N\to \infty}J_{\tilde{x},\omega,N}(\Theta,\tilde{f})=0.$$
On a aussi $J_{\tilde{x},\omega,geom}(\Theta,\tilde{f})=0$ d'apr\`es la d\'efinition de ce terme: il n'y a pas de sous-tore maximal de $H'_{\tilde{x}}$ qui soit anisotrope. D'o\`u la relation 3.5(1) dans ce cas.

Supposons ${\cal A}_{H'_{\tilde{x}}}=\{1\}$. Le m\^eme raisonnement nous d\'ebarrasse de tous les tores $I$ intervenant dans 3.5(2) tels que $I'$ n'est pas elliptique. Posons
$$J_{\tilde{x},\omega,\infty}(\Theta,\tilde{f})=\sum_{I=I'I''\in {\cal T}_{ell}(H_{\tilde{x}}')\times {\cal T}(G''_{\tilde{x}})}\nu(I')\vert W(G_{\tilde{x}},I)\vert ^{-1}$$
$$\int_{\mathfrak{i}'(F)\times \mathfrak{i}''(F)^S}\hat{j}_{S}(X')D^{G'_{\tilde{x}}}(X')D^{G''_{\tilde{x}}}(X'')^{1/2}\Theta^{J,\sharp}_{\tilde{f},\tilde{x},\omega}(X)dX.$$
Un calcul familier montre que cette expression est absolument convergente. Cette fois, la relation 3.6(1) et le lemme 3.7 montrent que
$$lim_{N\to \infty}J_{\tilde{x},\omega,N}(\Theta,\tilde{f})=J_{\tilde{x},\omega,\infty}(\Theta,\tilde{f}).$$
On a suppos\'e
$$\Theta_{\tilde{x},\omega}(X)=\hat{j}_{S}(X')\hat{j}^{H''_{\tilde{x}}}(S,X'').$$
Notons maintenant $\Theta''$ la fonction $X''\mapsto \hat{j}^{H''_{\tilde{x}}}(S,X'')$.
Un quasi-caract\`ere \`a support dans $\omega$ se d\'ecompose en combinaison lin\'eaire de produits d'un quasi-caract\`ere dans $\omega'$ et d'un quasi-caract\`ere dans $\omega''$. Ecrivons ainsi
$$\Theta^J_{\tilde{f},\tilde{x},\omega}(X)=\sum_{b\in B}\Theta'_{\tilde{f},b}(X')\Theta''_{\tilde{f},b}(X''),$$
o\`u $B$ est un ensemble fini d'indices. On a alors
$$J_{\tilde{x},\omega,\infty}(\Theta,\tilde{f})=\sum_{b\in B}J'_{b}J''_{b,\infty},$$
$$J_{\tilde{x},\omega,geom}(\Theta,\tilde{f})=\sum_{b\in B}J'_{b}J''_{b,geom},$$
o\`u
$$J'_{b}=\sum_{I'\in {\cal T}_{ell}(G'_{\tilde{x}})}\vert W(G'_{\tilde{x}},I')\vert ^{-1}\nu(I')\int_{\mathfrak{i}'(F)}\hat{j}_{S}(X')\Theta'_{\tilde{f},b}(X')D^{G'_{\tilde{x}}}(X')dX',$$
$$J''_{b,\infty}=\sum_{I''\in {\cal T}(G''_{\tilde{x}})}\vert W(G''_{\tilde{x}},I'')\vert ^{-1}\int_{\mathfrak{i}''(F)^S}\hat{\Theta}''_{\tilde{f},b}(X'')D^{G''_{\tilde{x}}}(X'')^{1/2}dX'',$$
$$J''_{b,geom}=\sum_{I''\in {\cal T}''}\nu(I'')\int_{\mathfrak{i}''(F)}c_{\Theta''}(X'')c_{\Theta''_{\tilde{f},b}}(X'')D^{H''_{\tilde{x}}}(X'')\Delta''(X)^rdX''.$$
D'apr\`es [W2] th\'eor\`eme 7.9 et lemme 11.2(ii), on a l'\'egalit\'e
$J''_{b,\infty}=J''_{b,geom}$ pour tout $b\in B$. D'o\`u l'\'egalit\'e $J_{\tilde{x},\omega,\infty}(\Theta,\tilde{f})=J_{\tilde{x},\omega,geom}(\Theta,\tilde{f})$, puis l'\'egalit\'e 3.5(1) qu'il fallait prouver. $\square$

 \bigskip
  
 \section{Majorations}
 
 \bigskip
 
 \subsection{Les r\'esultats}
 
 On \'enonce dans ce paragraphe les majorations qui seront prouv\'ees dans cette section. Les hypoth\`eses sont celles de 3.1.
 
 Pour un entier $N\geq1$ et un r\'eel $D$, posons
 $$I(N,D)=\int_{G(F)}\Xi^G(g)^2\kappa_{N}(g)\sigma(g)^Ddg.$$
 
 (1) Cette int\'egrale est convergente. Le r\'eel $D$ \'etant fix\'e, il existe un r\'eel $R$ tel que
 $$I(N,D)<<N^R$$
 pour tout entier $N\geq1$.
 
 Pour $u\in U(F)$ et $i=-r,...,r-1$, rappelons que l'on a not\'e $u_{i+1,i}$ la coordonn\'ee $<z_{i+1}^*,uz_{i}>$ de $u$. Pour un entier $c\geq1$, on note $U(F)_{c}$ le sous-ensemble des $u\in U(F)$ tels que $val_{F}(u_{i+1,i})\geq-c$ pour tout $i=-r,...,r-1$. C'est un sous-groupe de $U(F)$ conserv\'e par la conjugaison par $H(F)$. Pour un r\'eel $D$ et un \'el\'ement $m\in M(F)$, posons
 $$X(c,D,m)=\int_{U(F)_{c}}\Xi^G(um)\sigma(um)^Ddu.$$
 
 (2) Cette expression est convergente. Pour $D$ fix\'e, il existe un r\'eel $R$ tel que
 $$X(c,D,m)<<c^R\sigma(m)^R\delta_{P}(m)^{1/2}\Xi^M(m)$$
 pour tous $c\geq1$ et tout $m\in M(F)$.
 
 (3) Pour tout r\'eel $D$ et tout entier $c\geq1$, l'int\'egrale
 $$\int_{H(F)U(F)_{c}}\Xi^H(h)\Xi^G(hu)\sigma(hu)^Ddu\,dh$$
 est convergente.
 
 (4) Pour tout r\'eel $D$ et tout entier $c\geq1$, l'int\'egrale
 $$\int_{H(F)U(F)_{c}}\int_{H(F)U(F)_{c}}\Xi^G(hu)\Xi^H(h'h)\Xi^G(h'u')\sigma(hu)^D\sigma(h'u')^Ddu'\,dh'\,du\,dh$$
 est convergente.
 
 Soient $D$ et $C$ deux r\'eels, $c$, $c'$ et $N$ trois entiers. On suppose $C,c,c',N\geq1$. Pour $m\in M(F)$, $h\in H(F)$, $u,u'\in U(F)$, posons
 $$\phi(m,h,u,u';D)=\Xi^H(h)\Xi^G(u'm)\Xi^G(u^{-1}h^{-1}u'm)\kappa_{N}(m)\sigma(u')^D\sigma(u)^D\sigma(h)^D\sigma(m)^D\delta_{P}(m)^{-1}.$$
 Posons
 $$I(c,N,D)=\int_{M(F)}\int_{H(F)U(F)_{c}}\int_{U(F)}\phi(m,h,u,u';D)du'\,du\,dh\,dm,$$
 $$I(c,c',N,D)=\int_{M(F)}\int_{H(F)U(F)_{c}}\int_{U(F)-U(F)_{c'}}\phi(m,h,u,u';D)du'\,du\,dh\,dm,$$
$$I(c,c',N,C,D)=\int_{M(F)}\int_{H(F)U(F)_{c}}\int_{U(F)_{c'}}{\bf 1}_{\sigma\geq C}(hu)\phi(m,h,u,u';D)du'\,du\,dh\,dm.$$

(5) L'int\'egrale $I(c,N,D)$ est convergente. Les termes $c$ et $D$ \'etant fix\'es, il existe un r\'eel $R$ tel que
$$I(c,N,D)<<N^R$$
pour tout $N\geq1$.

(6) L'int\'egrale $I(c,c',N,D)$ est convergente. Les termes $c$ et $D$ \'etant fix\'es, pour tout r\'eel $R$, il existe $\alpha>0$ tel que
$$I(c,c',N,D)<<N^{-R}$$
pour tout $N\geq2$ et tout $c'\geq\alpha\,log(N)$.

(7) L'int\'egrale $I(c,c',N,C,D)$ est convergente. Les termes $c$ et $D$ \'etant fix\'es, pour tout r\'eel $R$, il existe $\alpha>0$ tel que
$$I(c,c',N,C,D)<<N^{-R}$$
pour tout $N\geq1$, tout $c'\geq1$ et tout $C\geq\alpha(log(N)+c')$. 

\bigskip

\subsection{Choix d'une base}

On fixe une base $(v_{i})_{i=r+1,...,d-r-1}$ de $W$.  On la compl\`ete en une base de $V$ en posant $v_{i}=z_{r+1-i}$ pour $i=1,...,r$, $v_{i}=z_{d-i-r}$ pour $i=d-r,...,d$. Dans cette section, le r\'eseau $R$ que l'on a fix\'e n'intervient que via les fonctions $\kappa_{N}$. Consid\'erons un autre r\'eseau $R'$, qui donne naissance \`a d'autres fonctions $\kappa'_{N}$. On v\'erifie qu'il existe $c\geq0$ tel que
$$\kappa'_{N}(g)\leq \kappa_{N+c}(g)\leq\kappa'_{N+2c}(g)$$
pour tout $g\in G(F)$ et tout $N\geq1$. Les majorations que l'on veut obtenir sont donc insensibles au changement de $\kappa_{N}$ en $\kappa'_{N}$. Cela nous autorise \`a supposer que $R$ est le r\'eseau de base $(v_{i})_{i=1,...,d}$ sur $\mathfrak{o}_{F}$. Donc $K=K_{d}$. On introduit le tore $A_{d}$ et le sous-groupe de Borel $B_{d}$. Notons $\mathfrak{S}_{d}$ le groupe des permutations de l'ensemble $\{1,...,d\}$. Pour tout $s\in \mathfrak{S}_{d}$, on note $A_{d}(F)_{s}^-$ le sous-ensemble des $a\in A_{d}(F)$ tels que
$$val_{F}(a_{s(1)})\geq val_{F}(a_{s(2)})\geq...\geq val_{F}(a_{s(d)}).$$
Pour $s$ \'egal \`a l'identit\'e, on pose simplement $A_{d}(F)^-=A_{d}(F)_{s}^-$.

\bigskip

\subsection{Comparaison de $\Xi^H$ et $\Xi^G$}

\ass{Lemme}{Supposons $r=0$. Il existe $\epsilon>0$ tel que $\Xi^G(h)<<exp(-\epsilon\sigma(h))\Xi^H(h)$ pour tout $h\in H(F)$.}

Preuve. Pour $r=0$, on a $H=GL_{d-1}$. En vertu de l'\'egalit\'e $H(F)=K_{d-1}A_{d-1}(F)^-K_{d-1}$, on peut supposer $h=a\in A_{d-1}(F)^-$. Si $val_{F}(a_{1})<0$, posons $k=1$. Sinon, soit $k$ le plus grand entier tel que $2\leq k\leq d$ tel que $val_{F}(a_{k-1})\geq0$. Introduisons l'\'el\'ement $b\in A_{d}(F)$ tel que $b_{i}=a_{i}$ pour $i=1,...,k-1$, $b_{k}=1$ et $b_{i}=a_{i-1}$ pour $i=k+1,...,d$. L'\'el\'ement $b$ appartient \`a $A_{d}(F)^-$ et est conjugu\'e \`a $a$. En vertu du lemme II.1.1 de [W4], il existe un entier $D$, ind\'ependant de $a$, tel que
$$\Xi^G(a)=\Xi^G(b)<<\delta_{B_{d}}(b)^{1/2}\sigma(a)^D.$$
On a aussi
$$\delta_{B_{d-1}}(a)^{1/2}<<\Xi^H(a).$$
On calcule
$$\delta_{B_{d}}(b)^{1/2}=\delta_{B_{d-1}}(a)^{1/2}(\prod_{i=1,...,k-1}\vert a_{i}\vert _{F}^{1/2})(\prod_{i=k,...,d-1}\vert a_{i}\vert _{F}^{-1/2}).$$
En vertu de la d\'efinition de $k$, cela \'equivaut \`a
$$\delta_{B_{d}}(b)^{1/2}=\delta_{B_{d-1}}(a)^{1/2}q^{-\Sigma(a)},$$
o\`u $q$ est le nombre d'\'el\'ements du corps r\'esiduel et
$$\Sigma(a)=\sum_{i=1,...,d-1}\vert val_{F}(a_{i})\vert .$$
On obtient
$$\Xi^G(a)<<\sigma(a)^Dq^{-\Sigma(a)}\Xi^H(a).$$
Mais $\Sigma(a)<<\sigma(a)$ et le lemme s'ensuit. $\square$

 \bigskip

\subsection{Majorations d'int\'egrales unipotentes}

 Pour simplifier les notations, on pose $B=B_{d}$. Pour tout $Q=M_{Q}U_{Q}\in {\cal F}(A_{d})$ et pour tout $g\in G(F)$, on fixe une d\'ecomposition $g=m_{Q}(g)u_{Q}(g)k_{Q}(g)$, avec $m_{Q}(g)\in M_{Q}(F)$, $u_{Q}(g)\in U_{Q}(F)$, $k_{Q}(g)\in K$. Dans le cas particulier o\`u $Q\in {\cal P}(A_{d})$, on note plut\^ot $a_{Q}(g)$ le terme $m_{Q}(g)$.

Supposons $r\geq1$. Notons $V_{r-1}$ le sous-espace de $V$ engendr\'e par $V_{0}$ et les vecteurs $z_{\pm j}$ pour $j=0,...,r-1$. Notons $P_{r}$ le sous-groupe parabolique de $G$ form\'e des \'el\'ements qui conservent le drapeau
$$Fz_{r}\subset  Fz_{r}\oplus V_{r-1}.$$
 On note $U_{r}$ son radical unipotent et $M_{r}$ sa composante de L\'evi qui contient $M$. Notons $U_{r,\natural}$ le sous-groupe des \'el\'ements $u\in U_{r}$ tels que $u_{r,r-1}=u_{1-r,-r}=0$. Pour deux r\'eels $b$ et $D$, posons
$$I_{r,\natural}(b,D)=\int_{U_{r,\natural}(F)}{\bf 1}_{\sigma\geq b}(u)\delta_{\bar{P}}(m_{\bar{P}}(u))^{1/2}\Xi^M(m_{\bar{P}}(u))\sigma(u)^Ddu.$$

\ass{Lemme}{Cette int\'egrale est convergente. Pour $D$ fix\'e, il existe $\epsilon>0$ tel que
$$I_{r,\natural}(b,D)<<exp(-\epsilon b)$$
pour tout $b\geq0$.}

Preuve. Supposons $r\geq2$. Remarquons que l'on peut aussi bien travailler avec l'int\'egrale
$$J_{r,\natural}(b,D)=\int_{U_{r,\natural}(F)}{\bf 1}_{\sigma\geq b}(u)\delta_{\bar{B}}(a_{\bar{B}}(u))^{1/2} \sigma(u)^Ddu.$$
En effet, d'apr\`es [W4] lemmes II.1.1 et II.3.2, il existe un r\'eel $D_{1}$ tel que
$$\delta_{\bar{P}}(m)^{1/2}\Xi^M(m)<<\delta_{\bar{B}}(a_{\bar{B}\cap M}(m))^{1/2}\sigma(m)^{D_{1}}$$
pour tout $m\in M(F)$. Pour $g\in G(F)$, on peut supposer $a_{\bar{B}}(g)=a_{\bar{B}\cap M}(m_{\bar{P}}(g))$ et on a $\sigma(m_{\bar{P}}(g))<<\sigma(g)$. Alors $I_{r,\natural}(b,D)<<J_{r,\natural}(b,D+D_{1})$.

 Notons $P'$ le sous-groupe parabolique de $G$ form\'e des \'el\'ements qui conservent la droite $Fz_{r}$. Notons $U'$ son radical unipotent et $M'$ sa composante de L\'evi contenant $M$. Posons $P''=M'\cap P_{r}$ et notons $U''=M'\cap P_{r}$ le radical unipotent de $P''$. On a $P''=M_{r}U''$. On pose $U'_{\natural}=U'\cap U_{r,\natural}$, $U''_{\natural}=U''\cap U_{r,\natural}$. On a $U_{r}=U'U''=U''U'$ et $U_{r,\natural}=U'_{\natural}U''_{\natural}=U''_{\natural}U'_{\natural}$. Posons
$$J'_{\natural}(b,D)=\int_{U'_{\natural}(F)}{\bf 1}_{\sigma\geq b}(u)\delta_{\bar{B}}(a_{\bar{ B}}(u))^{1/2}\sigma(u)^Ddu,$$
$$J''_{\natural}(b,D)=\int_{U''_{\natural}(F)}{\bf 1}_{\sigma\geq b}(u)\delta_{\bar{B}\cap M'}(a_{\bar{B}\cap M'}(u))^{1/2}\sigma(u)^Ddu.$$
On a

(1) ces int\'egrales sont convergentes; pour $D$ fix\'e, il existe $\epsilon_{D}>0$ tel que
$$J'_{\natural}(b,D)<< exp(-\epsilon_{D} b),\,\,J''_{\natural}(b,D)<< exp(-\epsilon_{D} b)$$
pour tout $b\geq0$.

L'assertion concernant $J'_{\natural}(b,D)$ est l'assertion (1) de [W3] 3.3. Consid\'erons $J''_{\natural}(b,D)$. On a $M'=GL_{1}\times GL_{d-1}$ et tout se passe dans le bloc $GL_{d-1}$. Appliquons dans ce bloc l'automorphisme $\theta_{d-1}$. Cela transforme le groupe $U''_{\natural}$ en l'analogue de $U'_{\natural}$ quand on remplace $d$ par $d-1$. Alors l'assertion r\'esulte de la m\^eme assertion (1) de [W3] 3.3, appliqu\'ee au groupe $GL_{d-1}$.

D'autre part, on v\'erifie facilement qu'il existe $\alpha>0$ tel que
$$(2) \qquad \delta_{\bar{B}}(a_{\bar{B}}(gg'))^{1/2}<<\delta_{\bar{B}}(a_{\bar{B}}(g))^{1/2}\delta_{\bar{B}}(a_{\bar{B}}(g'))^{1/2}exp(\alpha \sigma(g'))$$
pour tous $g,g'\in G(F)$.

Fixons un r\'eel auxiliaire $\beta>0$ que l'on pr\'ecisera plus tard. D\'ecomposons $J_{r,\natural}(b,D)$ en
$$J_{r,\natural}(b,D)=J_{\geq}(b,D)+J_{<}(b,D),$$
o\`u
$$J_{\geq}(b,D)=\int_{U''_{\natural}(F)}\int_{U'_{\natural}(F)}{\bf 1}_{\sigma\geq b}(u'u''){\bf 1}_{\sigma\geq \beta \sigma(u'')+b/2}(u')\delta_{\bar{B}}(a_{\bar{B}}(u'u''))^{1/2} \sigma(u'u'')^Ddu'\,du'',$$
$$J_{<}(b,D)=\int_{U''_{\natural}(F)}\int_{U'_{\natural}(F)}{\bf 1}_{\sigma\geq b}(u'u''){\bf 1}_{\sigma< \beta \sigma(u'')+b/2}(u')\delta_{\bar{B}}(a_{\bar{B}}(u'u''))^{1/2} \sigma(u'u'')^Ddu'\,du''.$$
Dans $J_{\geq}(b,D)$, on majore ${\bf 1}_{\sigma\geq b}(u'u'')$ par $1$ et on utilise (2). On obtient
$$J_{\geq}(b,D)<<\int_{U''_{\natural}(F)}exp(\alpha\sigma(u''))\delta_{\bar{B}}(a_{\bar{B}}(u''))^{1/2}\sigma(u'')^D$$
$$\int_{U'_{\natural}(F)}{\bf 1}_{\sigma\geq \beta \sigma(u'')+b/2}(u')\delta_{\bar{B}}(a_{\bar{B}}(u'))^{1/2}\sigma(u')^Ddu'\,du''.$$
L'int\'egrale int\'erieure est $J'_{\natural}(\beta \sigma(u'')+b/2,D)$. On d\'eduit de (1) la majoration
$$J_{\geq}(b,D)<<exp(-\epsilon_{D}b/2)\int_{U''_{\natural}(F)}exp((\alpha-\beta \epsilon_{D})\sigma(u''))\delta_{\bar{B}}(a_{\bar{B}}(u''))^{1/2}\sigma(u'')^Ddu''.$$
Remarquons que $\delta_{\bar{B}}(a_{\bar{B}}(u''))^{1/2}=\delta_{\bar{B}\cap M'}(a_{\bar{B}\cap M'}(u''))^{1/2}$. Introduisons le sous-groupe radiciel \'evident $U_{-r+1,-r}$ de $U$. On a $U''=U''_{\natural}U_{-r+1,-r}$. Soit $v\in U_{-r+1,-r}(F)\cap K$. On ne modifie pas l'int\'egrale  pr\'ec\'edente en rempla\c{c}ant $u''$ par $u''v$. On peut ensuite int\'egrer sur $v\in U_{-r+1,-r}(F)$. On obtient alors une int\'egrale sur un sous-ensemble ouvert de $U''(F)$, que l'on majore par l'int\'egrale sur tout le groupe. D'o\`u
$$J_{\geq}(b,D)<<exp(-\epsilon_{D}b/2)\int_{U''(F)}exp((\alpha-\beta \epsilon_{D})\sigma(u''))\delta_{\bar{B}\cap M'}(a_{\bar{B}\cap M'}(u''))^{1/2}\sigma(u'')^Ddu''.$$
On fixe maintenant $\beta$ tel que $\beta>0$ et  $\alpha-\beta \epsilon_{D}<0$. D'apr\`es le lemme II.4.2 de [W4], l'int\'egrale est convergente. On obtient
$$(3)\qquad J_{\geq}(b,D)<<exp(-\epsilon_{D}b/2).$$

Le groupe $U''$ normalise $U'$. Dans $J_{<}(b,D)$, effectuons le changement de variables $u'\mapsto u''u'u^{_{''}-1}$. Le terme ${\bf 1}_{\sigma<\beta\sigma(u'')+b/2}(u')$ devient ${\bf 1}_{\sigma<\beta\sigma(u'')+b/2}(u''u'u^{_{''}-1})$. Puisque $\sigma(u')\leq\sigma(u''u'u^{_{''}-1})+2\sigma(u'')$, ce terme est major\'e par ${\bf 1}_{\sigma<(\beta+2)\sigma(u'')+b/2}(u')$. Il y a aussi le terme ${\bf 1}_{\sigma\geq b}(u''u')$. Si les deux termes pr\'ec\'edents sont non nuls, on a
$$b\leq \sigma(u''u')\leq\sigma(u'')+\sigma(u')\leq (\beta+3)\sigma(u'')+b/2,$$
d'o\`u $\sigma(u'')\geq b/(2\beta+6)$ et ${\bf 1}_{\sigma\geq b/(2\beta+6)}(u'')=1$. On a donc
$$J_{<}(b,D)<<\int_{U''_{\natural}(F)}\int_{U'_{\natural}(F)}{\bf 1}_{\sigma\geq b/(2\beta+6)}(u''){\bf 1}_{\sigma<(\beta+2)\sigma(u'')+b/2}(u') \delta_{\bar{B}}(a_{\bar{B}}(u''u'))^{1/2}\sigma(u''u')^Ddu'\,du''.$$
Comme ci-dessus, on peut remplacer l'int\'egrale int\'erieure par une int\'egrale sur le groupe $U'(F)$ tout entier. On a $a_{\bar{B}}(u''u')=a_{\bar{B}\cap M'}(u'')a_{\bar{B}}(k_{\bar{B}\cap M'}(u'')u')$. Le groupe $U'(F)$ est normalis\'e par $M'(F)$. On effectue le changement de variables $u'\mapsto k_{\bar{B}\cap M'}(u'')^{-1}u'k_{\bar{B}\cap M'}(u'')$ et on obtient
$$J_{<}(b,D)<<\int_{U''_{\natural}(F)} {\bf 1}_{\sigma\geq b/(2\beta+6)}(u'')  \delta_{\bar{B}\cap M'}(a_{\bar{B}\cap M'}(u''))^{1/2}\sigma(u'')^D$$
$$\int_{U'(F)}{\bf 1}_{\sigma<(\beta+2)\sigma(u'')+b/2}(u')\delta_{\bar{B}}(a_{\bar{B}}(u'))^{1/2}\sigma(u')^Ddu'\,du''.$$
On peut majorer le terme ${\bf 1}_{\sigma<(\beta+2)\sigma(u'')+b/2}(u')$ par $\left((\beta+2)\sigma(u'')+b/2\right)^R\sigma(u')^{-R}$ pour n'importe quel r\'eel $R>0$. Remarquons que, quand ${\bf 1}_{\sigma\geq b/(2\beta+6)}(u'')=1$, ce terme est essentiellement major\'e par $\sigma(u'')^R\sigma(u')^{-R}$. D'o\`u
$$J_{<}(b,D)<<\int_{U''_{\natural}(F)} {\bf 1}_{\sigma\geq b/(2\beta+6)}(u'')  \delta_{\bar{B}\cap M'}(a_{\bar{B}\cap M'}(u''))^{1/2}\sigma(u'')^{D+R}$$
$$\int_{U'(F)} \delta_{\bar{B}}(a_{\bar{B}}(u'))^{1/2}\sigma(u')^{D-R}du'\,du''.$$
D'apr\`es le lemme II.4.2 de [W4], on peut fixer $R$ de sorte que l'int\'egrale int\'erieure soit convergente.  Ce nombre $R$ \'etant maintenant fix\'e, l'int\'egrale restante est $J''_{\natural}(b/(2\beta+6),D+R)$. Gr\^ace \`a (1), on obtient
$$J_{<}(b,D)<<exp(-\epsilon_{D+R}b/(2\beta+6)).$$
Jointe \`a (3), cette in\'egalit\'e entra\^{\i}ne celle de l'\'enonc\'e.

Supposons maintenant $r=1$. Dans ce cas, on a $P_{1}=P$ et $z_{0}=v_{d-1}$. Introduisons l'espace $V'$ engendr\'e par les vecteurs $v_{i}$ pour $i\in \{1,...,d\}-\{d-1\}$. Notons $G'\simeq GL_{d-1}$ son groupe d'automorphismes lin\'eaires. Posons $P'=G'\cap P$, $U'=G'\cap U$ et $M'=G'\cap M$. Le groupe $U_{1,\natural}$ est \'egal \`a $U'$ et $M'=GL(1)\times H\times GL(1)$. Soit $u\in U'(F)$.
On peut supposer $m_{\bar{P}}(u)=m_{\bar{P}'}(u)$. Ecrivons ce terme sous la forme
$m_{\bar{P}'}(u)=(a_{1},h,a_{d})$, avec $a_{1},a_{d}\in F^{\times}$ et $h\in H(F)$. On a
$$\delta_{\bar{P}}(m_{\bar{P}}(u))^{1/2}\Xi^M(m_{\bar{P}}(u))=\vert a_{1}\vert _{F}^{-(d-1)/2}\vert a_{d}\vert _{F}^{(d-1)/2}\Xi^{G_{0}}(h),$$
tandis que
$$\delta_{\bar{P}}(m_{\bar{P}'}(u))^{1/2}\Xi^{M'}(m_{\bar{P}'}(u))=\vert a_{1}\vert _{F}^{-(d-2)/2}\vert a_{d}\vert _{F}^{(d-2)/2}\Xi^{H}(h).$$
En vertu du lemme 4.3, on en d\'eduit
$$(4) \qquad \delta_{\bar{P}}(m_{\bar{P}}(u))^{1/2}\Xi^M(m_{\bar{P}}(u))<<\vert a_{1}\vert _{F}^{-1/2}\vert a_{d}\vert _{F}^{1/2}\delta_{\bar{P}}(m_{\bar{P}'}(u))^{1/2}\Xi^{M'}(m_{\bar{P}'}(u)).$$
On calcule $a_{1}$ et $a_{d}$ selon la m\'ethode habituelle. C'est-\`a-dire que l'on munit $V'$ et $\bigwedge^{d-2}V'$ de normes invariantes par $K$. On a
$$\vert a_{d}^{-1}\vert _{F}<<\vert a_{d}^{-1}v_{d}\vert $$
$$=\vert k_{\bar{P}'}(u)^{-1}u_{\bar{P}'}(u)^{-1}m_{\bar{P}'}(u)^{-1}v_{d}\vert =\vert u^{-1}v_{d}\vert .$$
Les coordonn\'ees de $u^{-1}v_{d}$ sont $1$ et les $(u^{-1})_{i,d}$, pour $i=1,...,d-2$. D'o\`u
$$(5) \qquad \vert a_{d}\vert _{F}^{-1}\geq sup(\{1\}\cup\{\vert (u^{-1})_{i,d}\vert _{F}; i=1,...,d-2\}).$$
On a
$$\vert det(h)^{-1}a_{d}^{-1}\vert _{F}<<\vert  m_{\bar{P}'}(u)^{-1}(v_{2}\wedge...\wedge v_{d-2}\wedge v_{d})\vert $$
$$=\vert k_{\bar{P}'}(u)^{-1}u_{\bar{P}'}(u)^{-1}m_{\bar{P}'}(u)^{-1}(v_{2}\wedge...\wedge v_{d-2}\wedge v_{d})\vert =\vert u^{-1}(v_{2}\wedge...\wedge v_{d-2}\wedge v_{d})\vert .$$
On peut supposer que $a_{1}a_{d}det(h)=det(u)=1$. Les coordonn\'ees $(u^{-1})_{1,j}$, pour $j=2,...,d-2$ sont aussi des coordonn\'ees de $u^{-1}(v_{2}\wedge...\wedge v_{d-2}\wedge v_{d})$: $(u^{-1})_{1,j}$ est la coordonn\'ee de $u^{-1}(v_{2}\wedge...\wedge v_{d-2}\wedge v_{d})$ sur $  v_{2}...\wedge v_{j-1}\wedge v_{1}\wedge v_{j+1}\wedge...\wedge v_{d-2}\wedge v_{d}$. La constante $1$ est aussi la coordonn\'ee de $u^{-1}(v_{2}\wedge...\wedge v_{d-2}\wedge v_{d})$ sur $v_{2}\wedge...\wedge v_{d-2}\wedge v_{d}$. D'o\`u
$$\vert a_{1}\vert _{F}\geq sup(\{1\}\cup\{\vert (u^{-1})_{1,j}\vert _{F}; j=2,...,d-2\}).$$
De cette in\'egalit\'e et de (5), on d\'eduit l'existence de $\alpha>0$ tel que
$$\vert a_{1}\vert _{F}\vert a_{d}\vert _{F}^{-1}\geq exp(\alpha\sigma(u)).$$
Alors (4) devient
$$\delta_{\bar{P}}(m_{\bar{P}}(u))^{1/2}\Xi^M(m_{\bar{P}}(u))<< exp(-\alpha\sigma(u)/2)\delta_{\bar{P}}(m_{\bar{P}'}(u))^{1/2}\Xi^{M'}(m_{\bar{P}'}(u)).$$
Donc
$$I_{1,\natural}(b,D)<<\int_{U'(F)}{\bf 1}_{\sigma\geq b}(u)exp(-\alpha\sigma(u)/2)\delta_{\bar{P}}(m_{\bar{P}'}(u))^{1/2}\Xi^{M'}(m_{\bar{P}'}(u))du.$$
Quand ${\bf 1}_{\sigma\geq b}(u)=1$, on peut majorer $exp(-\alpha\sigma(u)/2)$ par $exp(-\alpha b/4)exp(-\alpha\sigma(u)/4)$. D'apr\`es le lemme II.4.3 de [W4], l'int\'egrale est convergente et on obtient
$$I_{1,\natural}(b,D)<<exp(-\alpha b/4).$$
Cela ach\`eve la d\'emonstration. $\square$

\bigskip

\subsection{Majoration d'une int\'egrale de fonctions d'Harish-Chandra}

On suppose dans ce paragraphe $r=0$. Soit $D$ un r\'eel. Pour $h\in H(F)$ et $N\geq1$ un entier, posons
$$\chi(h,N,D)=\int_{G(F)}\Xi^G(hx)\Xi^G(x)\kappa_{N}(x)\sigma(x)^Ddx.$$

\ass{Lemme}{Cette int\'egrale est convergente. Le r\'eel $D$ \'etant fix\'e, il existe un r\'eel $R$ tel que
$$\chi(h,N,D)<<\Xi^G(h)N^R\sigma(h)^R$$
pour tout $h\in H(F)$ et tout entier $N\geq1$.}

Preuve. Comme en 4.3, on a $H=GL_{d-1}$. Ecrivons $h=k_{1}ak_{2}$, avec $k_{1},k_{2}\in K_{d-1}$ et $a\in A_{d-1}(F)$. On effectue le changement de variables $x\mapsto k_{2}^{-1}x$. Puisque $\Xi^G$ est biinvariante par $K_{d}$ et $\kappa_{N}$ est invariante \`a gauche par $H(F)$, on obtient l'\'egalit\'e
$\chi(h,N,D)=\chi(a,N,D)$. Cela nous ram\`ene au cas o\`u $h\in A_{d-1}(F)$, ce que l'on suppose d\'esormais.

On introduit le sous-groupe d'Iwahori standard $I$ de $G(F)$, c'est-\`a-dire le sous-ensemble des $k\in K_{d}$ tel que $val_{F}(k_{i,j})\geq1$ pour $i>j$. Notons $\Lambda$ le sous-ensemble des $a\in A_{d}(F)$ tels que toutes les coordonn\'ees $a_{i}$ sont des puissances de $\varpi_{F}$. On a l'\'egalit\'e
$$G(F)=\sqcup_{a\in \Lambda}IaK.$$
On a donc
$$\chi(h,N,D,a)=\sum_{a\in \Lambda}\chi(h,N,D,a),$$
o\`u
$$\chi(h,N,D,a)=\int_{IaK}\Xi^G(hx)\Xi^G(x)\kappa_{N}(x)\sigma(x)^Ddx.$$
Pour $s\in \mathfrak{S}_{d}$, posons $\Lambda_{s}^-=\Lambda\cap A_{d}(F)_{s}^-$, cf. 4.2. L'ensemble $\Lambda$ est r\'eunion des $\Lambda_{s}^-$. On peut donc fixer $s$ et se contenter de majorer
$$(1) \qquad \sum_{a\in \Lambda_{s}^-}\chi(h,N,D,a).$$
Quitte \`a r\'eindexer nos vecteurs de base, on peut supposer que $s$ est l'identit\'e. On note simplement $\Lambda^-$ l'ensemble correspondant.  Cette r\'eindexation change toutefois deux donn\'ees: le groupe $I$ n'est plus le m\^eme; le vecteur $z_{0}$ n'est plus \'egal \`a $v_{d}$, mais \`a un autre vecteur de base que nous notons $v_{k}$. D'apr\`es le lemme II.1.1 de [W4], il existe un r\'eel $R_{1}$ tel que  $\Xi^G(x)<<\delta_{B_{d}}(a)^{1/2}\sigma(a)^{R_{1}}$ pour tout $a\in \Lambda^-$ et tout $x\in IaK$. Pour $a\in \Lambda^-$, on a l'\'egalit\'e $IaK=(I\cap U_{d}(F))aK$ et on v\'erifie que l'on a une majoration $mes(IaK)<<\delta_{B_{d}}(a)^{-1}$.  On obtient
$$(2) \qquad \chi(h,N,D,a)<<\delta_{B_{d}}(a)^{-1/2}\sigma(a)^{R_{1}+D}Y(h,N,a),$$
o\`u
$$Y(h,N,a)=\int_{I\cap U_{d}(F)}\Xi^G(hua)\kappa_{N}(ua)du.$$
Posons 
$$X(h,N,a)=\int_{I\cap U_{d}(F)}\Xi^G(hua)du.$$
D\'efinissons quatre entiers $N_{1}(h)=N+1+sup(0,-val_{F}(h_{1}))$, $N_{d}(h)=N+1+sup(0,val_{F}(h_{d}))$, $b_{1}(h,N,a)=sup(0,val_{F}(a_{1})-N_{1}(h))$, $b_{d}(h,N,a)=sup(0,-val_{F}(a_{d})-N_{d}(h))$.
Montrons que

(3) il existe $\epsilon>0$ tel que
$$Y(h,N,a)<<exp(-\epsilon (b_{1}(h,N,a)+b_{d}(h,N,a)))X(h,N,a)$$
pour tout $h\in A_{d}(F)\cap H(F)$, tout $N\geq1$ et tout $a\in \Lambda^-$.

Supposons d'abord $k\not=1$ et $k\not=d$. Introduisons le sous-groupe $\Gamma$ de $U_{d}(F)$ form\'e des \'el\'ements $u(x,y,t)$ pour $x,y,t\in F$ tels que $u(x,y,t)v_{k}=v_{k}+xv_{1}$, $u(x,y)v_{d}=v_{d}+yv_{k}+tv_{1}$ et $u(x,y,t)$ fixe tout autre vecteur de base. Consid\'erons le sous-groupe $\Gamma_{h}$ de $\Gamma$ form\'e des \'el\'ements $u(x,y,t)$ tels que

- $val_{F}(x)\geq1$, $val_{F}(y)\geq1$, $val_{F}(t)\geq1$;

- $val_{F}(x)\geq -val_{F}(h_{1})+1$, $val_{F}(y)\geq val_{F}(h_{d})+1$, $val_{F}(t)\geq -val_{F}(h_{1})+val_{F}(h_{d})+1$.

Remarquons que les conditions portant sur $x$ et $y$  peuvent s'\'ecrire $val_{F}(x)\geq N_{1}(h)-N$  et $val_{F}(y)\geq N_{d}(h)-N$. Puisque $h\in H(F)$, on a $h_{k}=1$, d'o\`u l'\'egalit\'e $hu(x,y,t)h^{-1}=u(h_{1}x,h_{d}^{-1}y,h_{1}h_{d}^{-1}t)$. Les conditions ci-dessus assurent que $\gamma\in I\cap U_{d}(F)$ et $h\gamma h^{-1}\in I\cap U_{d}(F)$ pour tout $\gamma\in \Gamma_{h}$. Soit $\gamma\in \Gamma_{h}$. On ne modifie pas $Y(h,N,a)$ en rempla\c{c}ant $\Xi^G(hua)$ par $\Xi^G(h\gamma h^{-1} hua)$. On effectue ensuite le changement de variables $u\mapsto \gamma^{-1}u$. On peut enfin int\'egrer en $\gamma\in \Gamma_{h}$, \`a condition de diviser par $mes(\Gamma_{h})$.On obtient
$$(4) \qquad Y(h,N,a)=\int_{I\cap U_{d}(F)}\Xi^G(hua)\kappa_{N,h}(ua)du,$$
o\`u
$$\kappa_{N,h}(ua)=mes(\Gamma_{h})^{-1}\int_{\Gamma_{h}}\kappa_{N}(\gamma^{-1}ua)d\gamma.$$
Soient $\gamma=u(x,y,t)\in \Gamma_{h}$ et $u\in U_{d}(F)$. La condition $\kappa_{N}(\gamma^{-1}ua)=1$ \'equivaut \`a
$$a^{-1}u^{-1}\gamma v_{k}\in \mathfrak{p}_{F}^{-N}R\text{  et  }^ta{^tu}{^t\gamma}^{-1}v_{k}\in \mathfrak{p}_{F}^{-N}R,$$
ou encore
$$a_{1}^{-1}xv_{1}+a^{-1}u^{-1}v_{k}\in \mathfrak{p}_{F}^{-N}R\text{  et  }-a_{d}yv_{d}+{^ta}{^tu}v_{k}\in \mathfrak{p}_{F}^{-N}R.$$
 En notant $x_{u}$ la composante de $u^{-1}v_{k}$ sur $v_{1}$ et $y_{u}$ celle de $^tuv_{k}$ sur $v_{d}$, ces conditions entra\^{\i}nent 
 $$val_{F}(x+x_{u})\geq val_{F}(a_{1})-N \text{  et  }val_{F}(y-y_{u})\geq -val_{F}(a_{d})-N.$$
 Rappelons que l'on a aussi  $val_{F}(x)\geq N_{1}(h)-N$  et $val_{F}(y)\geq N_{d}(h)-N$. Les  relations pr\'ec\'edentes peuvent ne pas avoir de solution, auquel cas $\kappa_{N,h}(ua)=0$. Si elles en ont, elles d\'eterminent une classe modulo un sous-groupe $\Gamma'_{h}$ de $\Gamma_{h}$ et on a $\kappa_{N,h}(ua)=[\Gamma_{h}:\Gamma'_{h}]^{-1}$. L'indice de $\Gamma'_{h}$ dans $\Gamma_{h}$ se calcule ais\'ement. Les relations portant sur $x$ contribuent par $1$ si $val_{F}(a_{1})\leq N_{1}(h)$ et par $q^{N_{1}(h)-val_{F}(a_{1})}$ si $val_{F}(a_{1})\geq N_{1}(h)$. Autrement dit, elles contribuent par $q^{-b_{1}(h,N,a)}$. De m\^eme, les relations portant sur $y$ contribuent par $q^{-b_{d}(h,N,a)}$. Dans l'\'egalit\'e (4), on peut donc majorer $\kappa_{N,h}(ua)$ par $q^{-b_{1}(h,N,a)-b_{d}(h,N,a)}$ et la majoration  (3) s'en d\'eduit.
 
 Supposons maintenant $k=1$ (le cas $k=d$ est similaire). Pour $u\in U_{d}(F)$, la condition $\kappa_{N}(ua)=1$ entra\^{\i}ne $a^{-1}u^{-1}v_{1}\in \mathfrak{p}_{F}^{-N}R$, c'est-\`a-dire $a_{1}^{-1}v_{1}\in \mathfrak{p}_{F}^{-N}R$, ou encore $val_{F}(a_{1})\leq N$. Donc $b_{1}(h,N,a)=0$ et on peut oublier ce terme dans la relation (3). On introduit maintenant le sous-groupe $\Gamma$ de $U_{d}(F)$ form\'e des $u(y)$, pour $y\in F$, tels
 que $u(y)v_{d}=v_{d}+yv_{1}$ et $u(y)$ fixe tout autre vecteur de base. A l'aide de ce sous-groupe, le m\^eme raisonnement que ci-dessus conduit \`a la majoration (3).  
 
  Montrons que
 
 (5) il existe un r\'eel $R_{2}$ tel que
 $$X(h,N,a)<<\Xi^G(h)\delta_{B_{d}}(a)^{1/2}\sigma(h)^{R_{2}}\sigma(a)^{R_{2}}$$
 pour tout $h\in A_{d}(F)\cap H(F)$ et  tout $a\in \Lambda^-$.
 
 On peut fixer $s\in \mathfrak{S}_{d}$ tel que $h\in A_{d}(F)_{s}^-$. Introduisons le sous-groupe de Borel $B_{d,s}=A_{d}U_{d,s}$ de $G$ form\'e des \'el\'ements qui conservent le drapeau
 $$Fv_{s(1)}\subset Fv_{s(1)}\oplus Fv_{s(2)}\subset...\subset Fv_{s(1)}\oplus...\oplus Fv_{s(d)}.$$
 On a l'\'egalit\'e $I\cap U_{d}(F)=(I\cap U_{d}(F)\cap U_{d,s}(F))(I\cap U_{d}(F)\cap \bar{U}_{d,s}(F))$. Pour $u\in I\cap U_{d}(F)\cap U_{d,s}(F)$, on a $huh^{-1}\in I$. Le groupe $ I\cap U_{d}(F)\cap U_{d,s}(F)$ contribue donc trivialement \`a $X(h,N,a)$ et on obtient
 $$X(h,N,a)<<\int_{I\cap U_{d}(F)\cap \bar{U}_{d,s}(F)}\Xi^G(hua)du$$
 $$<<\delta_{0}(h)\int_{h(I\cap U_{d}(F)\cap \bar{U}_{d,s}(F))h^{-1}}\Xi^G(uha)du,$$
 o\`u $\delta_{0}(h)$ est la valeur absolue du d\'eterminant de $ad(h^{-1})$ agissant sur $\mathfrak{u}_{d}(F)\cap \bar{\mathfrak{u}}_{d,s}(F)$. Soit $u'\in I\cap U_{d}(F)\cap U_{d,s}(F)$. On ne modifie pas l'int\'egrale ci-dessus en rempla\c{c}ant $\Xi^G(uha)$ par $\Xi^G(u'uha)$. On peut ensuite int\'egrer en $u'$. En regroupant les deux int\'egrales, on obtient une int\'egrale sur un ouvert de $U_{d}(F)$. Sur cet ouvert, la variable d'int\'egration $u$ v\'erifie une majoration $\sigma(u)<<\sigma(h)$.  Pour tout r\'eel $R_{3}>0$, on peut donc glisser le terme $\sigma(h)^{R_{3}}\sigma(u)^{-R_{3}}$ dans l'int\'egrale et on obtient
 $$X(h,N,a)<<\delta_{0}(h)\sigma(h)^{R_{3}}\int_{U_{d}(F)}\sigma(u)^{-R_{3}}\Xi^G(uha)du.$$
 D'apr\`es [W4], proposition II.4.5, on peut choisir $R_{3}$ de sorte que cette int\'egrale soit convergente et essentiellement born\'ee par $\sigma(ha)^{R_{3}}\delta_{B_{d}}(ha)^{1/2}$. D'o\`u
 $$X(h,N,a)<<\delta_{0}(h)\delta_{B_{d}}(h)^{1/2}\delta_{B_{d}}(a)^{1/2}\sigma(h)^{2R_{3}}\sigma(a)^{R_{3}}.$$
 On v\'erifie que $\delta_{0}(h)\delta_{B_{d}}(h)^{1/2}=\delta_{B_{d,s}}(h)^{1/2}$.   D'apr\`es [W4] lemme II.1.1, ce terme est essentiellement born\'e par $\Xi^G(h)$. La majoration pr\'ec\'edente entra\^{\i}ne donc (5).

Gr\^ace \`a (2), (3) et (5), il existe $R_{5}$ tel que
$$\chi(h,N,D,a)<<exp(-\epsilon (b_{1}(h,N,a)+b_{d}(h,N,a)))\Xi^G(h)\sigma(a)^{R_{5}}\sigma(h)^{R_{5}}.$$
La somme (1) est born\'ee par
$$\sigma(h)^{R_{5}}\Xi^G(h)\sum_{a\in \Lambda^-}\sigma(a)^{R_{5}}exp(-\epsilon (b_{1}(h,N,a)+b_{d}(h,N,a))).$$
On v\'erifie qu'il existe $R_{6}$ tel que la s\'erie soit essentiellement born\'ee par $N_{1}(h)^{R_{6}}N_{d}(h)^{R_{6}}$. Ce terme est lui-m\^eme essentiellement born\'e par $N^{2R_{6}}\sigma(h)^{2R_{6}}$. Alors la majoration pr\'ec\'edente devient celle de l'\'enonc\'e. $\square$

\bigskip

\subsection{Preuve (sic!) des majorations de 4.1}

Dans la section 4 de [W3], on a d\'emontr\'e les analogues des majorations de 4.1 pour les groupes sp\'eciaux orthogonaux. On vient de d\'emontrer  en 4.3 ci-dessus l'analogue du lemme 4.9 de [W3], en 4.4 les analogues des lemmes 4.5 et 4.6 de [W3] et en 4.5 l'analogue du lemme 4.11 de [W3]. On l'a fait parce que les preuves pour $GL_{d}$ diff\`erent quelque peu de celles pour les groupes sp\'eciaux orthogonaux. En inspectant les preuves des autres paragraphes de [W3] (paragraphes 4.7, 4.8, 4.10 et 4.12 \`a 4.16), on constate qu'on peut les reprendre sans changement pour le groupe $GL_{d}$. On obtient ainsi une preuve des assertions de 4.1.

\bigskip

\section{Produits bilin\'eaires et facteurs $\epsilon$}

\bigskip

\subsection{Mod\`eles de Whittaker}

Fixons une base $(v_{i})_{i=1,...,d}$ de $V$ et identifions $G$ \`a $GL_{d}$. Soit $(\pi,E)$ une repr\'esentation admissible irr\'eductible de $G(F)$. On suppose qu'elle est temp\'er\'ee et qu'elle  admet un mod\`ele de Whittaker, que l'on construit comme en 2.2.

Pour $c\in {\mathbb Z}$, notons $\iota_{c}$ la fonction caract\'eristique du sous-ensemble de $A_{d}(F)$ form\'e des $a\in A_{d}(F)$ tels que $val_{F}(a_{i})\geq val_{F}(a_{i+1})-c$ pour tout $i=1,...,d-1$.  On sait qu'il existe un r\'eel $R$ et, pour tout $e\in E$, un entier $c$ de sorte que
$$(1) \qquad \vert W_{e}(a)\vert<< \iota_{c}(a)\delta_{B_{d}}(a)^{1/2}\sigma(a)^R$$
pour tout $a\in A_{d}(F)$.

 Pour $h=1,...,d$, on identifie $GL_{h}$ au sous-groupe de $G$ qui conserve le sous-espace de $V$ engendr\'e par $v_{1},...,v_{h}$ et fixe $v_{i}$ pour $i>h$. Pour $h,k,l\in {\mathbb N}$ tels que $h<k\leq l\leq d$, on note $\bar{U}_{h,k,l}$ le sous-groupe de $\bar{U}_{d}$ form\'e des \'el\'ements $\bar{u}\in \bar{U}_{d}$ tels que, pour $i,j=1,...,d$, $i\not=j$, $\bar{u}_{i,j}$ n'est non nul que si $(i,j)$ appartient au rectangle d\'efini par les in\'egalit\'es $k\leq i\leq l$, $1\leq j\leq h$.
  
\ass{Lemme}{On suppose $\pi$ temp\'er\'ee. Soit $h\in \{1,...,d-2\}$. Pour tout $e\in E$, il existe un sous-ensemble compact $\Omega\subset \bar{U}_{h,h+1,d-1}(F)$ tel que, pour tout $g\in GL_{h}(F)$, la fonction $\bar{u}\mapsto W_{e}(g\bar{u})$ sur $\bar{U}_{h,h+1,d-1}(F)$ soit \`a support dans $\Omega$.}

Preuve. La preuve  se trouve dans [JPSS] lemme 2.6. Rappelons-la. On remarque que $\bar{U}_{h,h+1,d-1}$ est le radical unipotent du sous-groupe parabolique de $GL_{d-1}$ qui est triangulaire inf\'erieur et de L\'evi $GL_{h}\times GL_{d-h-1}$. Supposons $W_{e}(g\bar{u})\not=0$. On \'ecrit $g\bar{u}=uak$, avec $u\in U_{d-1}(F)$, $a\in A_{d-1}(F)$, $k\in K_{d-1}$.   Gr\^ace \`a (1), on a une majoration $\vert a_{i}\vert _{F}\leq c$, o\`u $c$ est ind\'ependant de $g$. La base de $V$ d\'etermine une base du produit  $\bigwedge^{d-1-h}V$. On munit ce produit de la norme $\vert .\vert $ qui est le sup des valeurs absolues des coefficients dans la base en question. Ce produit est invariant par $K_{d}$. On a
$$\vert ^t(uak)(v_{h+1}\wedge...\wedge v_{d-1})\vert =\vert a_{h+1}...a_{d-1}\vert \leq c^{d-1-h}.$$
Mais un coefficient $\bar{u}_{i,j}$ est le coefficient de $ ^t(uak)(v_{h+1}\wedge...\wedge v_{d-1})={^t(g\bar{u})}(v_{h+1}\wedge...\wedge v_{d-1})$ sur le vecteur de base $v_{h+1}\wedge...v_{i-1}\wedge v_{j}\wedge v_{i+1}\wedge...\wedge v_{d-1}$. Donc les coefficients de $\bar{u}$ sont born\'es. $\square$

\bigskip

\subsection{Entrelacements}

Soient $(\pi,E_{\pi})\in Temp(G)$ et $(\rho,E_{\rho})\in Temp(H)$. On munit les espaces de ces repr\'esentations de produits scalaires invariants. Pour $c\in {\mathbb Z}$, on a d\'efini en 4.1 le sous-groupe $U(F)_{c}$  de  $U(F)$. Pour $e,e'\in E_{\pi}$ et $\epsilon,\epsilon'\in E_{\rho}$, posons
$${\cal L}_{\pi,\rho,c}(\epsilon'\otimes e',\epsilon\otimes e)=\int_{H(F)}\int_{U(F)_{c}}(\rho(h)\epsilon',\epsilon)(e',\pi(hu)e)\bar{\xi}(u)du\,dh.$$
D'apr\`es  4.1(3), cette int\'egrale est absolument convergente. 

\ass{Lemme}{Pour tous $\epsilon,\epsilon',e,e'$, il existe $c_{0}$ tel que ${\cal L}_{\pi,\rho,c}(\epsilon'\otimes e',\epsilon\otimes e)$ soit ind\'ependant de $c$ pour $c\geq c_{0}$.}

La preuve est similaire \`a celle du lemme 3.5 de [W3]. Rappelons-la. Pour un entier $c'\geq1$, notons $\omega(c')$ le sous-groupe des $a\in A(F)$ tels que $val_{F}(1-a_{i})\geq c'$ pour tout $i=\pm 1,...,\pm r$. Choisissons $c'$ tel que $e$ et $e'$ soient fix\'es par $\omega(c')$. Alors
$${\cal L}_{\pi,\rho,c}(\epsilon'\otimes e',\epsilon\otimes e)=mes(\omega(c'))^{-1}\int_{\omega(c')}\int_{H(F)}\int_{U(F)_{c}}(\rho(h)\epsilon',\epsilon)(\pi(a)e',\pi(hua)e)\bar{\xi}(u)du\,dh\,da.$$
Par le changement de variable $u\mapsto aua^{-1}$, on transforme cette expression en
$${\cal L}_{\pi,\rho,c}(\epsilon'\otimes e',\epsilon\otimes e)=mes(\omega(c'))^{-1} \int_{H(F)}\int_{U(F)_{c}}(\rho(h)\epsilon',\epsilon)(e',\pi(hu)e)\int_{\omega(c')}\bar{\xi}(aua^{-1})du\,dh\,da.$$
Mais il existe $c_{0}$, ne d\'ependant que de $c'$, tel que l'int\'egrale int\'erieure en $a$ soit nulle si $u\not\in U(F)_{c_{0}}$. D'o\`u le lemme. $\square$

On d\'efinit une forme sesquilin\'eaire ${\cal L}_{\pi,\rho}$ sur $E_{\rho}\otimes_{{\mathbb C}}E_{\pi}$ par
$${\cal L}_{\pi,\rho}(\epsilon'\otimes e',\epsilon\otimes e)=lim_{c\to \infty}{\cal L}_{\pi,\rho,c}(\epsilon'\otimes e',\epsilon\otimes e).$$

\bigskip

\subsection{Entrelacements et mod\`eles de Whittaker}

Fixons une base $(v_{i})_{i=1,...,m}$ de $W$ et compl\'etons-la en une base $(v_{i})_{i=1,...,d}$ de $V$ par $v_{m+i}=z_{r+1-i}$ pour $i=1,...,2r+1$.  Pour tout $h=1,...,d$, introduisons le sous-groupe parabolique standard $Q^h=L^hU^h$ tel que
$$L^h=GL_{h}\times GL_{1}\times...\times GL_{1}.$$
 Notons $\beta:GL_{r+1}\to G$ le plongement qui identifie $GL_{r+1}$ au sous-groupe des \'el\'ements de $G$ qui conservent le sous-espace de $V$ engendr\'e par $v_{m+1},...,v_{m+r+1}$ et fixent les autres vecteurs de base. On v\'erifie que notre groupe $U$ est produit des trois sous-groupes   $\beta(U_{r+1})$, $\bar{U}_{m,m+1,m+r}$  et $U^{m+r+1}$. Le caract\`ere $\xi$ de $U(F)$ est trivial sur $\bar{U}_{m,m+1,m+r}(F)$ et co\"{\i}ncide sur les deux autres facteurs avec le caract\`ere $\xi$ de $U_{d}(F)$.   Notons $U_{d}(F)_{c}$ le sous-groupe des $u\in U_{d}(F)$ tels que $val_{F}(u_{i,i+1})\geq-c$ pour tout $i=1,...,d-1$. Alors $U(F)_{c}$ est produit des trois sous-groupes   $\beta(U_{r+1}(F)_{c})$, $\bar{U}_{m,m+1,m+r}$  et $U^{m+r+1}(F)_{c}=U_{d}(F)_{c}\cap U^{m+r+1}(F)$. Soient $(\pi,E_{\pi})\in Temp(G)$ et $(\rho,E_{\rho})\in Temp(H)$. Fixons comme en 5.1 des mod\`eles de Whittaker de $\pi$ et $\rho$. Pour $e\in E_{\pi}$ et $\epsilon'\in E_{\rho}$, posons
$$L_{\pi,\rho}(\epsilon',e)=\int_{\bar{U}_{m,m+1,m+r}(F)}\int_{U_{m}(F)\backslash H(F)}\bar{W}_{\epsilon'}(h)W_{e}(h\bar{u})\vert det(h)\vert ^{-r}dh\,d\bar{u}.$$

\ass{Lemme}{(i) L'int\'egrale ci-dessus est absolument convergente.

(ii) Il existe $C>0$ tel que, pour tous $e,e'\in E_{\pi}$ et $\epsilon,\epsilon'\in E_{\rho}$, on ait l'\'egalit\'e
$${\cal L}_{\pi,\rho}(\epsilon'\otimes e',\epsilon\otimes e)=CL_{\pi,\rho}(\epsilon',e)\overline{L_{\pi,\rho}(\epsilon,e')}.$$}

Preuve. Consid\'erons $L(\epsilon',e)$. D'apr\`es le lemme 5.1, l'int\'egrale sur $\bar{U}_{m,m+1,m+r}(F)$ est \`a support compact et on n'a pas \`a s'en soucier. On d\'ecompose  $h\in U_{m}(F)\backslash H(F)$ en $h=ak$, avec $a\in A_{m}(F)$ et $k\in K_{m}$. La mesure devient $\delta_{B_{m}}(a)^{-1}da\,dk$. L'int\'egrale en $k$ ne nous importe pas. D'apr\`es  5.1(1), l'int\'egrale restante est born\'ee par
$$\int_{A_{m}(F)}\iota_{c'}(a)\delta_{B_{m}}(a)^{-1/2}\delta_{B_{d}}(a)^{1/2}\vert det(a)\vert _{F}^{-r}da,$$
pour un entier $c'$ convenable. On calcule
$$\delta_{B_{m}}(a)^{-1/2}\delta_{B_{d}}(a)^{1/2}\vert det(a)\vert _{F}^{-r}=\vert det(a)\vert _{F}^{1/2}$$
et on v\'erifie que l'int\'egrale ci-dessus est convergente. Cela prouve (i).

La premi\`ere \'etape pour prouver (ii) est de calculer
$$\int_{U(F)_{c}}(e',\pi(u)e)\bar{\xi}(u)du$$
pour un entier $c\geq1$ et deux \'el\'ements $e,e'\in E_{\pi}$. On a
$$(1) \qquad \int_{U(F)_{c}}(e',\pi(u)e)du=\int_{\bar{U}_{m,m+1,m+r}(F)}\int_{\beta(U_{r+1}(F)_{c})}$$
$$\int_{U^{m+r+1}(F)_{c}}(e',\pi(u'u\bar{u})e)\bar{\xi}(u'u)du'\,du\,d\bar{u}$$
et cette expression est absolument convergente.
Pour un entier $h=0,...,d-1$, notons $\omega_{[h+1,d-1]}(c)$ le sous-groupe des \'el\'ements $a\in A_{d}(F)$ tels que $a_{1}=...=a_{h}=1$, $a_{d}=1$ et $val_{F}(1-a_{i})\geq c$ pour tout $i=h+1,...,d-1$. Notons $c_{\psi}$ l'exposant du conducteur de $\psi$, c'est-\`a-dire que $\psi$ est trivial sur $\mathfrak{p}_{F}^{c_{\psi}}$ mais pas sur $\mathfrak{p}_{F}^{c_{\psi}-1}$. Posons
$$(2) \qquad I^h(e',e)=\int_{\omega_{[h+1,d-1]}(c+c_{\psi})}\int_{U_{h}(F)\backslash GL_{h}(F)}\bar{W}_{e'}(ag)W_{e}(ag)\vert det(g)\vert _{F}^{h+1-d}dg\,da.$$
On suppose $c\geq1$ et $c+c_{\psi}\geq 1$. D'apr\`es le lemme 3.6 de [W3], cette expression est absolument convergente et il existe $C_{r+1}>0$, ne d\'ependant que de $c$ et des mesures de Haar, tel que l'int\'egrale int\'erieure de l'expression (1) soit \'egale \`a $C_{r+1}I^{m+r}(e',\pi(u\bar{u})e)$. Pour $k=2,...,r+1$, notons $X^k$ le sous-groupe de $\beta(U_{r+1})$   form\'e des $u\in U_{d}$ tels que, pour $i,j=1,...,d$, $i\not=j$, on n'a $u_{i,j}\not=0$ que si $j=m+k$ et $i=m+1,...,m+k-1$. On a l'\'egalit\'e 
$\beta(U_{r+1})=X^{r+1}X^r...X^2$. Fixons un entier $N\geq c$. Pour $k=2,...,r+1$, notons ${\cal X}_{N}^k$ le sous-groupe des \'el\'ements $u\in X^k(F)$ tels que $val_{F}(u_{i,m+k})\geq-kN$ pour $i=m+1,...,m+k-2$ et $val_{F}(u_{m+k-1,m+k})\geq-c$. Notons $\bar{{\cal U}}_{N}$ le sous-groupe des \'el\'ements  $\bar{u}\in \bar{U}_{m,m+1,m+r}(F)$ tels que $val_{F}(\bar{u}_{i,j})\geq-(2r+1)N$ pour tous $i\not=j$. On a
$$\int_{U(F)_{c}}(e',\pi(u)e)\bar{\xi}(u)du=C_{r+1}lim_{N\to \infty}\int_{\bar{{\cal U}}_{N}}\int_{{\cal X}^2_{N}}...\int_{{\cal X}^{r+1}_{N}}$$
$$I^{m+r}(e',\pi(u_{r+1}...u_{2}\bar{u})e)\bar{\xi}(u_{r+1}...u_{2})du_{r+1}\,...du_{2}\,d\bar{u}.$$
Pour $k=1,...,r+1$, posons
$$J^k_{N}=\int_{\bar{{\cal U}}_{N}}\int_{{\cal X}^2_{N}}...\int_{{\cal X}^{k}_{N}}\int_{\bar{U}_{m,m+k,m+r}(F)}I^{m+k-1}(\pi(\bar{v})e',\pi(u_{k}...u_{2}\bar{u})e)\bar{\xi}(u_{k}...u_{2})d\bar{v}\,du_{k}\,...du_{2}\,d\bar{u}.$$
La formule pr\'ec\'edente se r\'ecrit
$$(3) \qquad \int_{U(F)_{c}}(e',\pi(u)e)\bar{\xi}(u)du=C_{r+1}lim_{N\to \infty}J^{r+1}_{N}.$$
Montrons que pour tout $k=2,...,r+1$, il existe $C'_{k}>0$ ne d\'ependant que de $c$ et des mesures de Haar de sorte que
$$(4) \qquad J^{k}_{N}=C'_{k}J^{k-1}_{N}$$
si $N$ est assez grand. Rempla\c{c}ons dans  $J_{N}$ le terme $I^{m+k-1}(\pi(\bar{v})e',\pi(u_{k}...u_{2}\bar{u})e)$ par son expression int\'egrale (2).  On obtient
$$(5) \qquad J^k_{N}=\int_{\bar{{\cal U}}_{N}}\int_{{\cal X}^2_{N}}...\int_{{\cal X}^{k}_{N}}\int_{\bar{U}_{m,m+k,m+r}(F)}\int_{\omega_{[m+k,d-1]}(c+c_{\psi})}\int_{U_{m+k-1}(F)\backslash GL_{m+k-1}(F)}$$
$$\bar{W}_{e'}(ag\bar{v})W_{e}(agu_{k}...u_{2}\bar{u})\vert det(g)\vert _{F}^{m+k-d}\bar{\xi}(u_{k}...u_{2})dg\,da\,d\bar{v}\,du_{k}\,...du_{2}\,d\bar{u}.$$
Cette expression est absolument convergente. En effet, les variables $u_{k}$,...,$u_{2}$, $\bar{u}$ restent dans des compacts. La variable $\bar{v}$ n'intervient que dans $\bar{W}_{e'}(ag\bar{v})$ et $g$ appartient \`a $GL_{m+k-1}(F)$. Si on ne tient pas compte de $a$, le lemme 5.1 nous dit que cette fonction est \`a support compact en $\bar{v}$. Le $a$ ne g\^ene pas: en le faisant commuter \`a $\bar{v}$, on obtient que $a\bar{v}a^{-1}$ reste dans un compact, ce qui \'equivaut \`a ce que $\bar{v}$ reste dans un compact, puisque $a$ lui-m\^eme reste dans un compact. Mais alors, la convergence absolue de notre expression r\'esulte de celle de  (2).

Notons $Y$ le sous-groupe des \'el\'ements de $GL_{m+k-1}$ form\'e des $y$  dont les seuls coefficients non nuls sont:

 - $y_{i,i}=1$ pour $i=1,...,m+k-2$;
 
  - les $y_{m+k-1,j}$ pour $j=1,...,m+k-1$.
  
   On note $dy$ le produit des mesures additives $dy_{m+k-1,j}$. On a la formule d'int\'egration
$$\int_{U_{m+k-1}(F)\backslash GL_{m+k-1}(F)}\varphi(g)dg=C'\int_{Y(F)}\int_{U_{m+k-2}(F)\backslash GL_{m+k-2}(F)}$$
$$\varphi(g'y)\vert det(g')\vert _{F}^{-1}dg'\vert y_{m+k-1,m+k-1}\vert ^{-1}dy$$
pour toute fonction int\'egrable $\varphi$ sur $U_{m+k-1}(F)\backslash GL_{m+k-1}(F)$, o\`u $C' $ est une constante positive. On utilise cette formule pour calculer l'int\'egrale en $g$  de la formule (5).On remplace donc $g$ par $g'y$, avec $g'\in GL_{m+k-2}(F)$ et $y\in Y(F)$.  On a l'\'egalit\'e $ag'yu_{k}=nag'y$, o\`u $n$ est un \'el\'ement de $U_{d}(F)$ qui n'a de coefficients non nuls, hormis les diagonaux, que sur la $(m+k)$-i\`eme colonne, et qui v\'erifie 
$$n_{m+k-1,m+k}=a_{m+k}^{-1}\sum_{i=m+1,...,m+k-1}y_{m+k-1,i}(u_{k})_{i,m+k}.$$
On a $W_{e}(ag'yu_{k}...u_{2}\bar{u})=\psi(n_{m+k-1,m+k})W_{e}(ag'yu_{k-1}...u_{2}\bar{u})$. L'int\'egrale en $u_{k}$ de la formule (5) est donc simplement
$$\int_{{\cal X}^k_{N}}\psi(-(u_{k})_{m+k-1,m+k}+a_{m+k}^{-1}\sum_{i=m+1,...,m+k-1}y_{m+k-1,i}(u_{k})_{i,m+k})du_{k}.$$
Notons $Y_{N}$ le sous-ensemble des $y\in Y(F)$ tels que $val_{F}(y_{m+k-1,i})\geq kN+c_{\psi}$ pour $i=m+1,...,m+k-2$ et $val_{F}(1-y_{m+k-1,m+k-1})\geq c+c_{\psi}$. Posons 
$$C(N)=mes(\mathfrak{p}_{F}^{-kN})^{k-2}mes(\mathfrak{p}_{F}^{-c}).$$
L'int\'egrale ci-dessus vaut $C(N) $ si $y\in Y_{N}$ et $0$ sinon. D\'ecomposons $Y_{N}$ en produit de trois groupes: son intersection avec $A_{d}(F)$, que l'on note $\omega_{[m+k-1,m+k-1]}(c+c_{\psi})$; le sous-groupe des $y$ unipotents tels que $y_{m+k-1,j}=0$ pour $j=m+1,...,m+k-2$, qui est \'egal \`a $\bar{U}_{m,m+k-1,m+k-1}(F)$; le sous-groupe des $y$ unipotents tels que $y_{m+k-1,j}=0$ pour $j=1,...,m$, que l'on note $Y'_{N}$. A ce point, on a obtenu l'\'egalit\'e
$$J^k_{N}=C'C(N)\int_{\bar{{\cal U}}_{N}}\int_{{\cal X}^2_{N}}...\int_{{\cal X}^{k-1}_{N}}\int_{\bar{U}_{m,m+k,m+r}(F)}\int_{\omega_{[m+k,d-1]}(c+c_{\psi})}\int_{U_{m+k-2}(F)\backslash GL_{m+k-2}(F)}$$
$$\int_{\omega_{[m+k-1,m+k-1]}(c+c_{\psi})}\int_{\bar{U}_{m,m+k-1,m+k-1}(F)}\int_{Y'_{N}}\bar{W}_{e'}(aga'\bar{u}'y'\bar{v})W_{e}(aga'\bar{u}'y'u_{k-1}...u_{2}\bar{u})$$
$$\vert det(g)\vert _{F}^{m+k-1-d}\bar{\xi}(u_{k-1}...u_{2})dy'\,d\bar{u}'\,da'\,dg\,da\,d\bar{v}\,du_{k-1}\,...du_{2}\,d\bar{u}.$$
  On a dej\`a dit que $\bar{v}$ restait dans un compact, qui est ind\'ependant de $N$. Alors l'\'el\'ement $\bar{v}^{-1}y'\bar{v}$ est proche de $1$ quand $N$ est grand, donc  $\bar{W}_{e'}(aga'\bar{u}'y'\bar{v})=\bar{W}_{e'}(aga'\bar{u}'\bar{v})$. Les termes $y'$, $u_{k-1}$,...,$u_{2}$ appartiennent \`a $\beta(GL_{r}(F))$, donc aussi $y''=(u_{k-1}...u_{2})^{-1}y'u_{k-1}...u_{2}$. De plus, les coefficients de $u_{k-1}...u_{2}$ sont de valuations $\geq -(k-1)N$ tandis que les coefficients non diagonaux de $y'$ sont de valuation $\geq kN+c_{\psi}$. Donc $y''$ est proche de $1$ quand $N$ est grand. Le groupe $\beta(GL_{r}(F))$ normalise $\bar{U}_{m,m+1,m+r}(F)$ et un \'el\'ement proche de $1$ de $\beta(GL_{r}(F))$ normalise $\bar{{\cal U}}_{N}$. Quitte \`a effectuer le changement de variables $\bar{u}\mapsto y^{_{''}-1}\bar{u}y''$, on remplace $W_{e}(aga'\bar{u}'y'u_{k-1}...u_{2}\bar{u})$ par $W_{e}(aga'\bar{u}'u_{k-1}...u_{2}\bar{u}y'')$, qui est \'egal \`a $W_{e}(aga'\bar{u}'u_{k-1}...u_{2}\bar{u})$. Alors $y'$ dispara\^{\i}t de notre formule. Posons $C''=C(N)mes(Y'_{N})$. Cette constante est ind\'ependante de $N$ et on a obtenu
  $$J^k_{N}=C'C''\int_{\bar{{\cal U}}_{N}}\int_{{\cal X}^2_{N}}...\int_{{\cal X}^{k-1}_{N}}\int_{\bar{U}_{m,m+k,m+r}(F)}\int_{\omega_{[m+k,d-1]}(c+c_{\psi})}\int_{U_{m+k-2}(F)\backslash GL_{m+k-2}(F)}$$
  $$\int_{\omega_{[m+k-1,m+k-1]}(c+c_{\psi})}\int_{\bar{U}_{m,m+k-1,m+k-1}(F)} \bar{W}_{e'}(aga'\bar{u}'\bar{v})W_{e}(aga'\bar{u}'u_{k-1}...u_{2}\bar{u})$$
  $$\vert det(g)\vert _{F}^{m+k-1-d}\bar{\xi}(u_{k-1}...u_{2})d\bar{u}'\,da'\,dg\,da\,d\bar{v}\,du_{k-1}\,...du_{2}\,d\bar{u}.$$
   De m\^eme que l'on a prouv\'e que $\bar{v}$ restait dans un compact, on montre maintenant que $\bar{u}'\bar{v}$ reste dans un compact, donc aussi $\bar{u}'$. Cet \'el\'ement appartient \`a $\bar{U}_{m,m+1,r}(F)$ et, comme ci-dessus, ce groupe est normalis\'e par $\beta(GL_{r}(F)$. Donc l'\'el\'ement $\bar{u}''=(u_{k-1}...u_{2})^{-1}\bar{u}'u_{k-1}....u_{2}$ appartient \`a $\bar{U}_{m,m+1,r}(F)$. Puisque $\bar{u}'$ reste dans un compact et que les coefficients de $u_{k-1}...u_{2}$ sont de valuation $\geq -(k-1)N$, on a $\bar{u}''\in \bar{{\cal U}}_{N}$. Quitte \`a effectuer le changement de variable $\bar{u}\mapsto \bar{u}^{_{''}-1}\bar{u}$, on remplace $W_{e}(aga'\bar{u}'u_{k-1}...u_{2}\bar{u})$ par $W_{e}(aga'u_{k-1}...u_{2}\bar{u})$. Le terme $\bar{u}'$ n'intervient plus que dans $\bar{W}_{e'}(aga'\bar{u}'\bar{v})$. Les int\'egrales en $\bar{u}'$ et $\bar{v}$ se combinent en une int\'egrale en $\bar{v}\in \bar{U}_{m,m+k-1,m+r}(F)$. Les \'el\'ements $a'$ et $g$ commutent. Ensuite, les int\'egrales en $a$ et $a'$ se combinent pour former une int\'egrale en $a\in \omega_{[m+k-1,d-1]}(c+c_{\psi})$. Cela conduit \`a l'\'egalit\'e
$$J^k_{N}=C'C''\int_{\bar{{\cal U}}_{N}}\int_{{\cal X}^2_{N}}...\int_{{\cal X}^{k-1}_{N}}\int_{\bar{U}_{m,m+k-1,m+r}(F)}\int_{\omega_{[m+k-1,d-1]}(c+c_{\psi})}\int_{U_{m+k-2}(F)\backslash GL_{m+k-2}(F)} $$
$$ \bar{W}_{e'}(ag\bar{v})W_{e}(agu_{k-1}...u_{2}\bar{u})\vert det(g)\vert _{F}^{m+k-1-d}\bar{\xi}(u_{k-1}...u_{2})dg\,da\,d\bar{v}\,du_{k-1}\,...du_{2}\,d\bar{u},$$
c'est-\`a-dire $J^k_{N}=C'C''J^{k-1}_{N}$. On a prouv\'e (4).

En utilisant (3) et (4), on obtient l'existence d'une constante $C_{1}>0$ telle que
$$ \int_{U(F)_{c}}(e',\pi(u)e)\bar{\xi}(u)du=C_{1}lim_{N\to \infty}J^{1}_{N}.$$
On a
 $$J^1_{N}=\int_{\bar{{\cal U}}_{N}}\int_{\bar{U}_{m,m+1,m+r}(F)}I^{m}(\pi(\bar{v})e',\pi(\bar{u})e) d\bar{v}\,d\bar{u}.$$ 
 Le m\^eme argument qui nous a dit que $\bar{v}$ restait dans un compact nous dit que $\bar{u}$ reste lui-aussi dans un compact. Pour $N$ grand, on peut remplacer $\bar{{\cal U}}_{N}$ par $\bar{U}_{m,m+1,m+r}(F)$ et on obtient
 $$\int_{U(F)_{c}}(e',\pi(u)e)\bar{\xi}(u)du=C_{1}\int_{\bar{U}_{m,m+1,m+r}(F)^2}\int_{\omega_{[m+1,d-1]}(c+c_{\psi})}\int_{U_{m}(F)\backslash GL_{m}(F)}$$
 $$\bar{W}_{e'}(ag\bar{v})W_{e}(ag\bar{u})\vert det(g)\vert _{F}^{m+1-d}dg\,da\, d\bar{v}\,d\bar{u}.$$
 Par les changements de variables $\bar{u}\mapsto a^{-1}\bar{u}a$ et $\bar{v}\mapsto a^{-1}\bar{v}a$, et en supposant que $c$ soit assez grand pour que $e$ et $e'$ soient invariants par  $\omega_{[m+1,d-1]}(c+c_{\psi})$, l'int\'egrale en $a$ dispara\^{\i}t. D'autre part, $m+1-d=-2r$.  On obtient
 $$(6) \qquad  \int_{U(F)_{c}}(e',\pi(u)e)\bar{\xi}(u)du=C_{1}\int_{\bar{U}_{m,m+1,m+r}(F)^2}\int_{U_{m}(F)\backslash GL_{m}(F)}$$
 $$\bar{W}_{e'}(g\bar{v})W_{e}(g\bar{u})\vert det(g)\vert _{F}^{-2r}dg\,d\bar{v}\,d\bar{u}.$$
 Rappelons que $H=GL_{m}$. Donc
 $${\cal L}_{\pi,\rho,c}(\epsilon'\otimes e',\epsilon\otimes e)=\int_{GL_{m}(F)}(\rho(h)\epsilon',\epsilon)\int_{U(F)_{c}}(e',\pi(uh)e)\bar{\xi}(u)du\,dh.$$
  On a
 $${\cal L}_{\pi,\rho,c}(\epsilon'\otimes e',\epsilon\otimes e)=lim_{N\to \infty}L_{N},$$
 o\`u
 $$L_{N}=\int_{GL_{m}(F)}(\rho(h)\epsilon',\epsilon)\int_{U(F)_{c}}(e',\pi(uh)e)\bar{\xi}(u){\bf 1}_{\sigma < N}(h)du\,dh.$$
 Fixons $N$. On peut remplacer l'int\'egrale int\'erieure par son expression (6). L'expression obtenue est absolument convergente. On effectue les changements de variables $\bar{u}\mapsto h\bar{u}h^{-1}$ puis $h\mapsto g^{-1}h$. Cela introduit un Jacobien $\vert det(g^{-1}h)\vert _{F}^{-r}$. On d\'ecompose ensuite $h$ en $u'a'k'$, avec $u'\in U_{m}(F)$, $a'\in A_{m}(F)$, $k'\in K_{m}$, et $g$ en $ak$, avec $a\in A_{m}(F)$ et $k\in K_{m}$. La mesure devient $\delta_{B_{m}}(aa')^{-1} du'\,da'\,dk'\,da\,dk$. On obtient
 $$L_{N}=C_{1}\int_{\bar{U}_{m,m+1,m+r}(F)^2}\int_{A_{m}(F)^2}\int_{K_{m}^2}\int_{U_{m}(F)}(\rho(u'a'k')\epsilon',\rho(ak)\epsilon)\bar{W}_{e'}(ak\bar{v})W_{e}(u'a'k'\bar{u})$$
 $$\delta_{B_{m}}(aa')^{-1}\vert det(aa')\vert_{F} ^{-r}{\bf 1}_{\sigma < N}(k^{-1}a^{-1}u'a'k')du'\,dk'\,dk\,da'\,da\,d\bar{v}\,d\bar{u}.$$
 On a $W_{e}(u'a'k'\bar{u})=\xi(u')W_{e}(a'k'\bar{u})$. Alors le m\^eme raisonnement qu'au lemme 5.2 montre que l'on peut remplacer l'int\'egrale sur $U_{m}(F)$ par une int\'egrale sur $U_{m}(F)_{c}$, pourvu que $c$ soit assez grand. Posons
 $$L=C_{1}\int_{\bar{U}_{m,m+1,m+r}(F)^2}\int_{A_{m}(F)^2}\int_{K_{m}^2}\int_{U_{m}(F)_{c}}(\rho(u'a'k')\epsilon',\rho(ak)\epsilon)\bar{W}_{e'}(ak\bar{v})W_{e}(a'k'\bar{u})$$
 $$\delta_{B_{m}}(aa')^{-1}\vert det(aa')\vert_{F} ^{-r} \xi(u')du'\,dk'\,dk\,da'\,da\,d\bar{v}\,d\bar{u}.$$
 On a
 
 (7) cette expression est absolument convergente.
 
 Consid\'erons l'int\'egrale
 $$\int_{U_{m}(F)_{c}}\vert (\rho(n
u'a'k')\epsilon',\rho(ak)\epsilon)\vert du'.$$
 On effectue le changement de variable $u'\mapsto au'a^{-1}$. Cela introduit un jacobien $\delta_{B_{m}}(a)$. La nouvelle variable parcourt un ensemble qui d\'epend de $a$, mais il existe $c_{0}>0$ tel que cet ensemble soit inclus dans $U_{m}(F)_{c+c_{0}\sigma(a)}$. L'int\'egrale ci-dessus est donc essentiellement born\'ee par
 $$\delta_{B_{m}}(a)\int_{U_{m}(F)_{c+c_{0}\sigma(a)}}\Xi^H(u'a^{-1}a')du'.$$
 D'apr\`es [W3] proposition 3.2, ceci est essentiellement major\'e par $\sigma(a)^R\sigma(a')^R\delta_{B_{m}}(aa')^{1/2}$ pour un certain r\'eel $R$. Revenons \`a l'expression $L$. Les variables $k$, $k'$, $\bar{v}$ et $\bar{u}$ restent dans des compacts, on peut les n\'egliger. Il reste \`a prouver que l'expression
 $$\int_{A_{m}(F)^2}\vert W_{e'}(a)W_{e}(a')\vert \delta_{B_{m}}(aa')^{-1/2}\vert det(aa')\vert _{F}^{-r}\sigma(a)^R\sigma(a')^Rda'\,da$$
 est convergente. Il suffit pour cela de reprendre le calcul de la preuve du (i) de l'\'enonc\'e. Cela prouve (7).
 
 Gr\^ace \`a (7), le th\'eor\`eme de convergence domin\'ee implique $lim_{N\to \infty}L_{N}=L$, d'o\`u
 $${\cal L}_{\pi,\rho,c}(\epsilon'\otimes e',\epsilon\otimes e)=L.$$
 En appliquant \`a $\rho$ le lemme 3.7 de [W3], on voit que, si $c$ est assez grand, il existe $C_{2}>0$ tel que
 $$\int_{U_{m}(F)_{c}}(\rho(u'a'k')\epsilon,\rho(ak)\epsilon)\xi(u')du'=C_{2}\bar{W}_{\epsilon'}(a'k')W_{\epsilon}(ak).$$
 Rempla\c{c}ons le membre de gauche par celui de droite dans l'expression de $L$. Rempla\c{c}ons ensuite $a'k'$ et $ak$ par $g',g\in U_{m}(F)\backslash GL_{m}(F)$. On reconna\^{\i}t alors
 $$L=C_{1}C_{2} L_{\pi,\rho}(\epsilon',e)\overline{L_{\pi,\rho}(\epsilon,e')},$$
 d'o\`u l'\'egalit\'e du (ii) de l'\'enonc\'e. On n'a pas \'et\'e pr\'ecis quant aux constantes, c'est-\`a-dire que l'on n'a pas v\'erifi\'e qu'elles ne d\'ependaient pas de l'entier auxiliaire $c$. Cela n'a pas d'importance. Puisque l'\'egalit\'e finale du (ii) de l'\'enonc\'e compare des termes ind\'ependants de $c$, ou bien tous les termes intervenant dans cette \'egalit\'e sont nuls, et on peut prendre $C$ quelconque, ou bien la constante $C$ ne peut  qu'\^etre ind\'ependante de $c$. $\square$
 
 \bigskip
 
 \subsection{Non-nullit\'e des entrelacements}
 
 \ass{Lemme}{Soient $\pi\in Temp(G)$ et $\rho\in Temp(H)$. Alors les formes sesquilin\'eaires $L_{\pi,\rho}$ et ${\cal L}_{\pi,\rho}$ sont non nulles.}
 
 Preuve. D'apr\`es le lemme pr\'ec\'edent, il suffit de prouver l'assertion relative \`a $L_{\pi,\rho}$. Pour $e\in E_{\pi}$, $\epsilon'\in E_{\rho}$ et $s\in {\mathbb C}$, posons
 $$L_{\pi,\rho}(\epsilon',e,s)= \int_{\bar{U}_{m,m+1,m+r}(F)}\int_{U_{m}(F)\backslash H(F)}\bar{W}_{\epsilon'}(h)W_{e}(h\bar{u})\vert det(h)\vert ^{s-r-1/2}dh\,d\bar{u}.$$
 Le m\^eme calcul que celui de la preuve du (i) du lemme pr\'ec\'edent montre que cette int\'egrale est absolument convergente pour $Re(s)>0$. Elle d\'efinit dans ce domaine une fonction holomorphe de $s$ et on a $L_{\pi,\rho}(\epsilon',e)=L_{\pi,\rho}(\epsilon',e,1/2)$. Mais la fonction $s\mapsto L_{\pi,\rho}(\epsilon',e,s)$ est celle \'etudi\'ee par Jacquet, Piatetski-Shapiro et Shalika. Ils montrent que pour tout $s$, on peut trouver $\epsilon'$ et $e$ tels que $L_{\pi,\rho}(\epsilon',e,s)\not=0$ ([JPSS] th\'eor\`eme 2.7(ii)). D'o\`u la conclusion. $\square$
 
\bigskip

\subsection{Entrelacements et groupes tordus}

 Soient $\tilde{\pi}\in Temp(\tilde{G})$ et $\tilde{\rho}\in Temp(\tilde{H})$.   On a introduit un nombre $\epsilon_{\nu}(\tilde{\pi},\tilde{\rho}) $ en 2.5.
 
 \ass{Proposition}{Soient $\epsilon,\epsilon'\in E_{\rho}$, $e,e'\in E_{\pi}$ et $\tilde{y}\in \tilde{H}(F)$. On a l'\'egalit\'e
 $${\cal L}_{\pi,\rho}(\epsilon'\otimes \tilde{\pi}(\tilde{y})e',\tilde{\rho}(\tilde{y})\epsilon\otimes e)= \epsilon_{\nu}(\tilde{\pi}^{\vee},\tilde{\rho}){\cal L}_{\pi,\rho}(\epsilon'\otimes e',\epsilon\otimes e).$$}
 
 Preuve. On v\'erifie que, pour tout $h\in H(F)$, on a l'\'egalit\'e
 $${\cal L}_{\pi,\rho}(\epsilon'\otimes \pi(h)e',\rho(h)\epsilon\otimes e)= {\cal L}_{\pi,\rho}(\epsilon'\otimes e',\epsilon\otimes e).$$
 Il suffit donc de prouver l'\'egalit\'e de la proposition pour un  \'el\'ement $\tilde{y}$ bien choisi. On fixe des bases comme en 5.3 et on choisit $\tilde{y}=\boldsymbol{\theta}_{m}$. Le lemme 5.3 nous ram\`ene \`a prouver l'\'egalit\'e
$$(1) \qquad L_{\pi,\rho}(\tilde{\rho}(\boldsymbol{\theta}_{m})\epsilon,\tilde{\pi}(\boldsymbol{\theta}_{m})e')=\overline{ \epsilon_{\nu}(\tilde{\pi}^{\vee},\tilde{\rho})}L_{\pi,\rho}(\epsilon,e').$$
 Une telle \'egalit\'e est ind\'ependante des normalisations de $\tilde{\pi}$ et $\tilde{\rho}$. On peut supposer $w(\tilde{\pi},\psi)=w(\tilde{\rho},\psi)=1$. Dans ce cas, on a $W_{\tilde{\rho}(\boldsymbol{\theta}_{m})\epsilon}(h)=W_{\epsilon}(\theta_{m}(h))$ pour tout $h\in H(F)$. Notons $\gamma$ l'\'el\'ement de $G(F)$ tel que $\boldsymbol{\theta}_{m}=\boldsymbol{\theta}_{d}\gamma$. On a de m\^eme
$$W_{\tilde{\pi}(\boldsymbol{\theta}_{m})e'}(g)=W_{\tilde{\pi}(\boldsymbol{\theta}_{d}\gamma)e'}(g)=W_{\pi(\gamma)e'}(\theta_{d}(g))=W_{e'}(\theta_{d}(g)\gamma)$$
pour tout $g\in G(F)$. Alors
$$L_{\pi,\rho}(\tilde{\rho}(\boldsymbol{\theta}_{m})\epsilon,\tilde{\pi}(\boldsymbol{\theta}_{m})e')=\int_{\bar{U}_{m,m+1,m+r}(F)}\int_{U_{m}(F)\backslash H(F)}\overline{W_{\epsilon}(\theta_{m}(h))}W_{e'}(\theta_{d}(h\bar{u})\gamma)\vert det(h)\vert ^{-r}dh\,d\bar{u}.$$
Notons $w_{m}\in GL_{m}(F)$ la matrice de permutation antidiagonale, $D_{m}\in GL_{m}(F)$ la matrice diagonale telle que $(D_{m})_{i,i}=(-1)^{m+1+i}$ et, pour $z\in F^{\times}$, notons $z_{m}\in GL_{m}(F)$ la matrice centrale de coefficients diagonaux $z$. Fixons $z\in F^{\times}$. On a $J_{m}=w_{m}D_{m}$, donc $\theta_{m}(h)=w_{m}D_{m}{^th}^{-1}J_{m}^{-1}$. De m\^eme $\theta_{d}(h\bar{u})\gamma=w_{d}D_{d}{^th}^{-1}{^t\bar{u}}^{-1}J_{d}^{-1}\gamma$. Effectuons le changement de variables $h\mapsto D_{m}z_{m}^{-1}h{^tJ_{m}^{-1}}$. On obtient
$$L_{\pi,\rho}(\tilde{\rho}(\boldsymbol{\theta}_{m})\epsilon,\tilde{\pi}(\boldsymbol{\theta}_{m})e')=\int_{\bar{U}_{m,m+1,m+r}(F)}\int_{U_{m}(F)\backslash H(F)}\overline{W_{\epsilon}(w_{m}z_{m}{^th}^{-1})}$$
$$W_{e'}(w_{d}D_{d}D_{m}z_{m}{^th}^{-1}J_{m}{^t\bar{u}}^{-1}J_{d}^{-1}\gamma)\vert det(h)\vert _{F}^{-r}\vert z\vert _{F}^{rm}dh\,d\bar{u}.$$
On a
$$\overline{W_{\epsilon}(w_{m}z_{m}{^th}^{-1})}=\overline{\omega_{\rho}(z)}\overline{W_{\epsilon}(w_{m}{^th}^{-1})}.$$
On v\'erifie que $D_{d}D_{m}z_{m}$ commute \`a $GL_{m}(F)$ et que $D_{d}D_{m}z_{m}J_{m}$ normalise $\bar{U}_{m,m+1,m+r}(F)$. Par le changement de variable $\bar{u}\mapsto {^t(D_{d}D_{m}z_{m}J_{m})}\bar{u} {^t(D_{d}D_{m}z_{m}J_{m})^{-1}}$, qui est de Jacobien $\vert z\vert _{F}^{-rm}$,
on obtient
$$L_{\pi,\rho}(\tilde{\rho}(\boldsymbol{\theta}_{m})\epsilon,\tilde{\pi}(\boldsymbol{\theta}_{m})e')=\overline{\omega_{\rho}(z)}\int_{\bar{U}_{m,m+1,m+r}(F)}\int_{U_{m}(F)\backslash H(F)}\overline{W_{\epsilon}(w_{m}{^th}^{-1})}$$
$$W_{e'}(w_{d}{^th}^{-1}{^t\bar{u}}^{-1}\gamma')\vert det(h)\vert _{F}^{-r}dh\,d\bar{u},$$
o\`u $\gamma'= D_{d}D_{m}z_{m}J_{m}J_{d}^{-1}\gamma=w_{m}z_{m}w_{d}\gamma$. Rappelons que, d'apr\`es notre construction du plongement de $\tilde{H}$ dans $\tilde{G}$, l'\'el\'ement $\boldsymbol{\theta}_{m}$, vu comme un \'el\'ement de $\tilde{G}(F)$, co\"{\i}ncide avec l'\'el\'ement $\boldsymbol{\theta}_{m}$ de d\'epart sur les vecteurs $v_{1},...,v_{m}$ et avec l'\'el\'ement $\tilde{\zeta}$ sur les vecteurs $v_{m+1},...,v_{d}$. C'est-\`a-dire que
$$\boldsymbol{\theta}_{m}(v_{i})=\left\lbrace\begin{array}{cc}(-1)^{i+[(m+1)/2]}v_{m+1-i}^*,&\text{ si }i=1,...,m,\\ (-1)^{i+m+1+r}2\nu v_{m+1+d-i}^*,&\text{ si }i=m+1,...,d.\\ \end{array}\right.$$
Puisque $\boldsymbol{\theta}_{m}=\boldsymbol{\theta}_{d}\gamma$, on calcule:
$$\gamma(v_{i})=\left\lbrace\begin{array}{cc}(-1)^{m+d+[(d+1)/2]+[(m+1)/2]}v_{d-m+i},& \text{ si }i=1,...,m,\\ (-1)^{r+1+[(d+1)/2]}2\nu v_{i-m},&\text{ si }i=m+1,...,d,\\ \end{array}\right.$$
puis
$$\gamma'(v_{i})=\left\lbrace\begin{array}{cc}(-1)^{m+d+[(d+1)/2]+[(m+1)/2]}zv_{i},& \text{ si }i=1,...,m,\\ (-1)^{r+1+[(d+1)/2]}2\nu v_{d+m+1-i},&\text{ si }i=m+1,...,d.\\ \end{array}\right.$$
  Posons $z'=(-1)^{r+1+[(d+1)/2]}2\nu$ et choisissons 
  $$z=(-1)^{m+d+[(d+1)/2]+[(m+1)/2]}z'=(-1)^{r+[(m+1)/2]}2\nu.$$
   Alors $\gamma'=z'_{d}w_{0}$ o\`u $w_{0}$ est la matrice de permutation obtenue en rempla\c{c}ant les constantes par $1$  dans la formule ci-dessus. On a
$$W_{e'}(w_{d}{^th}^{-1}{^t\bar{u}}^{-1}\gamma')=\omega_{\pi}(z')W_{e'}(w_{d}{^th}^{-1}{^t\bar{u}}^{-1}w_{0}),$$
et
$$L_{\pi,\rho}(\tilde{\rho}(\boldsymbol{\theta}_{m})\epsilon,\tilde{\pi}(\boldsymbol{\theta}_{m})e')=\overline{\omega_{\rho}(z)}\omega_{\pi}(z')\int_{\bar{U}_{m,m+1,m+r}(F)}\int_{U_{m}(F)\backslash H(F)}\overline{W_{\epsilon}(w_{m}{^th}^{-1})}$$
$$W_{e'}(w_{d}{^th}^{-1}{^t\bar{u}}^{-1}w_{0})\vert det(h)\vert _{F}^{-r}dh\,d\bar{u}.$$
D'apr\`es [JPSS] th\'eor\`eme 2.7(iii), l'int\'egrale ci-dessus est \'egale \`a $\omega_{\rho}(-1)^{d-1}\epsilon(1/2,\pi\times \bar{\rho},\psi)L_{\pi,\rho}(\epsilon,e')$. Pour obtenir l'\'egalit\'e (1), il reste \`a prouver l'\'egalit\'e
$$(2)\qquad \overline{\omega_{\rho}((-1)^{d-1}z)}\omega_{\pi}(z')\epsilon(1/2,\pi\times \bar{\rho},\psi)=\overline{ \epsilon_{\nu}(\tilde{\pi},\tilde{\rho})}.$$
On a suppos\'e $w(\tilde{\rho},\psi)=w(\tilde{\pi},\psi)=1$. D'apr\`es 2.2(1) et (2), on calcule $\overline{w(\tilde{\pi}^{\vee},\psi)}=w(\tilde{\pi},\psi^-)=\omega_{\pi}(-1)^{d-1}$. D'apr\`es 2.5(1),
$$\overline{\epsilon(1/2,\pi\times \rho,\psi)}=\epsilon(1/2,\bar{\pi}\times\bar{\rho},\psi^-)=\omega_{\pi}(-1)^m\omega_{\rho}(-1)^d\epsilon(1/2,\pi\times \rho,\psi).$$
Remarquons que, puisque $\rho$ est \`a la fois temp\'er\'ee, donc unitaire, et isomorphe \`a sa contragr\'ediente, on a $\rho\simeq \bar{\rho}\simeq \check{\rho}$. De m\^eme pour $\pi$.   De plus $\omega_{\rho}$ et $\omega_{\pi}$ prennent leurs valeurs dans $\{\pm 1\}$.  Le membre de droite de (2) vaut donc
$$\omega_{\pi}((-1)^{d-1+m+[m/2]}2\nu)\omega_{\rho}(-1)^{d+1+[d/2]}2\nu)\epsilon(1/2,\pi\times \bar{\rho},\psi).$$
Pour prouver (2) et la proposition, il reste \`a remarquer que $(-1)^{d-1}z=(-1)^{d+1+[d/2]}2\nu$ et $z'=(-1)^{d-1+m+[m/2]}2\nu$. $\square$

\bigskip

\subsection{Entrelacements et induction}

Soient $Q=LU_{Q}\in {\cal F}(A_{d})$ et $\tau\in Temp(L)$. Pour tout $\lambda\in i{\cal A}_{L,F}^*$, on d\'efinit $\tau_{\lambda}$ et la repr\'esentation induite $\pi_{\lambda}=Ind_{Q}^G(\tau_{\lambda})$. On la r\'ealise dans l'espace ${\cal K}_{Q,\tau}^G$ qui ne d\'epend pas de $\lambda$. Soit $\rho\in Temp(H)$. L'espace ${\cal K}_{Q,\tau}^G$ est muni d'un produit scalaire qui est invariant par $\pi_{\lambda}$ pour tout $\lambda$, et on utilise ce produit scalaire pour d\'efinir la forme sesquilin\'eaire ${\cal L}_{\pi_{\lambda},\rho}$.

\ass{Lemme}{(i) Soient $e,e'\in {\cal K}_{Q,\tau}^G$ et $\epsilon,\epsilon'\in E_{\rho}$. La fonction $\lambda\mapsto {\cal L}_{\pi_{\lambda},\rho}(\epsilon'\otimes e',\epsilon\otimes e)$ est $C^{\infty}$ sur $i{\cal A}_{L,F}^*$.

(ii) Il existe une famille finie $(\epsilon_{j})_{j=1,...,n}$ d'\'el\'ements de $E_{\rho}$, une famille finie $(e_{j})_{j=1,...,n}$ d'\'el\'ements de ${\cal K}_{Q,\tau}^G$ et une famille finie $(\varphi_{j})_{j=1,...,n}$ de fonctions $C^{\infty}$ sur $i{\cal A}_{L,F}^*$ de sorte que
$$\sum_{j=1,...,n}\varphi_{j}(\lambda){\cal L}_{\pi_{\lambda},\rho}(\epsilon_{j}\otimes e_{j},\epsilon_{j}\otimes e_{j})=1$$ pour tout $\lambda\in i{\cal A}_{L,F}^*$.}

La preuve est identique \`a celle du lemme 5.3 de [W3].

\bigskip

\section{La partie spectrale de la formule int\'egrale}

 \bigskip

 \subsection{Enonc\'e du th\'eor\`eme}
 
 On fixe une base $(v_{i})_{i=1,...,d}$ de $V$ et on utilise les notations introduites en 2.1. On prend pour L\'evi minimal $\tilde{M}_{min}=\tilde{A}_{d}$. Pour simplifier, on peut supposer que $(v_{i})_{i=1,...,d}$ est aussi une base du r\'eseau $R$ sur $\mathfrak{o}_{F}$. On choisit pour sous-groupe compact sp\'ecial de $G(F)$ le groupe $K=K_{d}$.
 
 Soient $\tilde{\rho}\in Temp(\tilde{H})$ et $\tilde{f}\in C_{c}^{\infty}(\tilde{G}(F))$ une fonction tr\`es cuspidale. Pour tout entier $N\geq1$, on a d\'efini $J_{N}(\Theta_{\tilde{\rho}},\tilde{f})$ en 3.3. Posons
 $$J_{spec}(\Theta_{\tilde{\rho}},\tilde{f})=\sum_{\tilde{L}\in {\cal L}^{\tilde{G}}}(-1)^{a_{\tilde{L}}}\vert W^{L}\vert \vert W^{G}\vert ^{-1}\sum_{{\cal O}\in \{\Pi_{ell}(\tilde{L})\}} [i{\cal A}_{{\cal O}}^{\vee}:i{\cal A}_{\tilde{L},F}^{\vee}]^{-1}$$
 $$2^{-s({\cal O})-a_{\tilde{L}'}}  \int_{i{\cal A}_{\tilde{L},F}^*} \epsilon_{\nu}(\tilde{\sigma}^{\vee},\tilde{\rho})J_{\tilde{L}}^{\tilde{G}}(\tilde{\sigma}_{\lambda},\tilde{f})d\lambda.$$
 On a fix\'e dans toute orbite ${\cal O}$ un point base que l'on a not\'e $\tilde{\sigma}$.
 
 \ass{Th\'eor\`eme}{On a l'\'egalit\'e
 $$lim_{N\to \infty}J_{N}(\Theta_{\rho},\tilde{f})=J_{spec}(\Theta_{\rho},\tilde{f}).$$}
 
 \bigskip
 
 \subsection{Utilisation de la formule de Plancherel}
 
 Fixons   $\tilde{y}\in \tilde{H}(F)$.  On d\'efinit une fonction $f$ sur $G(F)$ par $f(g)=\tilde{f}(g\tilde{y})$. Pour $g\in G(F)$, on a l'\'egalit\'e
 $$J(\Theta_{\tilde{\rho}},\tilde{f},g)=\int_{H(F)}\int_{U(F)}\Theta_{\tilde{\rho}}(h\tilde{y})\tilde{f}(g^{-1}hu\tilde{y}g)\xi(u)du\,dh,$$
 $$=\int_{H(F)}\int_{U(F)}\Theta_{\tilde{\rho}}(h\tilde{y})f(g^{-1}hu\theta_{\tilde{y}}(g))\xi(u)du\,dh.$$
 On exprime $f$ \`a l'aide de la formule de Plancherel, que l'on a reprise en [W3] 1.6. Pour tout $g\in G(F)$, on a l'\'egalit\'e
 $$f(g)=\sum_{L\in {\cal L}(A_{d})}\vert W^L\vert \vert W^G\vert ^{-1}\sum_{{\cal O}\in \{\Pi_{2}(L)\}}f_{{\cal O}}(g),$$
 o\`u
 $$f_{{\cal O}}(g)= [i{\cal A}_{{\cal O}}^{\vee}:i{\cal A}_{L,F}^{\vee}]^{-1}\int_{i{\cal A}_{L,F}^*}m(\tau_{\lambda})trace(Ind_{Q}^G(\tau_{\lambda},g^{-1})Ind_{Q}^G(\tau_{\lambda},f))d\lambda.$$
 Rappelons que $\Pi_{2}(L)$ est l'ensemble des repr\'esentations irr\'eductibles et de carr\'e int\'egrable de $L(F)$. L'ensemble $\{\Pi_{2}(L)\}$ est celui des orbites de l'action $\lambda\mapsto \tau_{\lambda}$ de $i{\cal A}_{L,F}^*$ dans $\Pi_{2}(L)$.   Pour tout ${\cal O}\in \{\Pi_{2}(L)\}$, on a fix\'e un point base $\tau$ et, pour tout $L$, on a fix\'e $Q\in {\cal P}(L)$. La fonction $\lambda\mapsto m(\tau_{\lambda})$ est la mesure de Plancherel.

  L'ensemble des orbites ${\cal O}$ pour lesquelles $f_{{\cal O}}\not=0$ est fini. On fixe un tel ensemble $\{\Pi_{2}(L)\}_{f}$. Fixons un sous-groupe ouvert compact $K_{f}$ de $K$ tel que $f$ soit biinvariante par $K_{f}$. Remarquons que toutes les fonctions $f_{{\cal O}}$ sont aussi biinvariantes par ce groupe. Pour tout $ g\in M(F)K$, fixons un sous-groupe ouvert compact $K'_{g}\subset H(F)$ tel que $gK'_{g}g^{-1}\subset K_{f}$ et $\theta_{\tilde{y}}(g)K'_{g}\theta_{\tilde{y}}(g^{-1})\subset K_{f}$. Fixons une base orthonorm\'ee ${\cal B}_{\rho}^{K'_{g}}$ du sous-espace des invariants $E_{\rho}^{K'_{g}}$. Alors, pour $g\in M(F)K$, on a l'\'egalit\'e
   $$J(\Theta_{\tilde{\rho}},\tilde{f},g)=\sum_{\epsilon\in {\cal B}_{\rho}^{K'_{g}}}\int_{H(F)}(\epsilon,\tilde{\rho}(h\tilde{y})\epsilon)\int_{U(F)}\sum_{L\in {\cal L}(A_{d})}\vert W^L\vert \vert W^G\vert ^{-1}$$
   $$\sum_{{\cal O}\in \{\Pi_{2}(L)\}_{f}}f_{{\cal O}}(g^{-1}hu\theta_{\tilde{y}}(g))\xi(u)du\,dh.$$
Pour $\epsilon\in E_{\rho}$, $L\in {\cal L}(A_{d})$, ${\cal O}\in \{\Pi_{2}(L)\}$ et $g\in G(F)$, posons
$$J_{L,{\cal O}}(\epsilon,f,g)=\int_{H(F)U(F)}(\epsilon,\tilde{\rho}(h\tilde{y})\epsilon)f_{{\cal O}}(g^{-1}hu\theta_{\tilde{y}}(g))\xi(u)du\,dh.$$
Comme en [W3] 6.2(2), on prouve que cette expression est absolument convergente. Pour $g\in M(F)K$, on a donc
$$(1) \qquad J(\Theta_{\tilde{\rho}},\tilde{f},g)=\sum_{\epsilon\in {\cal B}_{\rho}^{K'_{g}}}\sum_{L\in {\cal L}(A_{d})}\vert W^L\vert \vert W^G\vert ^{-1}\sum_{{\cal O}\in \{\Pi_{2}(L)\}_{f}}J_{L,{\cal O}}(\epsilon,f,g).$$
 
\bigskip

\subsection{Apparition des entrelacements}

Fixons $L\in {\cal L}(A_{d})$ et ${\cal O}\in \{\Pi_{2}(L)\}_{f}$. On a

(1) il existe $c_{0}\in {\mathbb N}$ tel que, pour tout $c\geq c_{0}$, tout $g\in M(F)K$ et tout $h\in H(F)$,
on ait l'\'egalit\'e
$$\int_{U(F)}f_{{\cal O}}(g^{-1}hu\theta_{\tilde{y}}(g))\xi(u)du=\int_{U(F)_{c}}f_{{\cal O}}(g^{-1}hu\theta_{\tilde{y}}(g))\xi(u)du.$$

 C'est la m\^eme preuve qu'au lemme 5.2, en utilisant le fait que $A$ commute \`a $M$.
 
 On fixe $\tau\in {\cal O}$ et $Q\in {\cal P}(L)$. On pose $\pi_{\lambda}=Ind_{Q}^G(\tau_{\lambda})$ pour tout $\lambda\in i{\cal A}_{L,F}^*$. On r\'ealise ces repr\'esentations dans l'espace ${\cal K}_{Q,\tau}^G$. On fixe une base orthonorm\'ee ${\cal B}_{{\cal O}}^{K_{f}}$ du sous-espace $({\cal K}_{Q,\tau}^G)^{K_{f}}$. Pour $g\in G(F)$, on a l'\'egalit\'e
 $$f_{{\cal O}}(g)=[i{\cal A}_{{\cal O}}^{\vee}:i{\cal A}_{L,F}^{\vee}]^{-1}$$
 $$\sum_{e\in {\cal B}_{{\cal O}}^{K_{f}}}\int_{i{\cal A}_{L,F}^*}m(\tau_{\lambda})(e,\pi_{\lambda}(g^{-1})\pi_{\lambda}(f)e)d\lambda.$$
 Pour $\epsilon\in E_{\rho}$, $g\in M(F)K$ et $c\geq c_{0}$, on a
 $$J_{L,{\cal O}}(\epsilon,f,g)=[i{\cal A}_{{\cal O}}^{\vee}:i{\cal A}_{L,F}^{\vee}]^{-1}\int_{H(F)U(F)_{c}}(\epsilon,\tilde{\rho}(h\tilde{y})\epsilon)$$
 $$\sum_{e\in {\cal B}_{{\cal O}}^{K_{f}}}\int_{i{\cal A}_{L,F}^*}m(\tau_{\lambda})(\pi_{\lambda}(\theta_{\tilde{y}}(g))e,\pi_{\lambda}((hu)^{-1}g)\pi_{\lambda}(f)e)\xi(u)d\lambda\,du\,dh.$$
 En changeant $h$ et $u$ en leurs inverses, on obtient
 $$J_{L,{\cal O}}(\epsilon,f,g)=[i{\cal A}_{{\cal O}}^{\vee}:i{\cal A}_{L,F}^{\vee}]^{-1}\int_{H(F)U(F)_{c}}(\rho(h)\epsilon,\tilde{\rho}(\tilde{y})\epsilon)$$
 $$\sum_{e\in {\cal B}_{{\cal O}}^{K_{f}}}\int_{i{\cal A}_{L,F}^*}m(\tau_{\lambda})(\pi_{\lambda}(\theta_{\tilde{y}}(g))e,\pi_{\lambda}(hug)\pi_{\lambda}(f)e)\bar{\xi}(u)d\lambda\,du\,dh.$$
  
   Pour $g\in M(F)K$ fix\'e, on a une majoration
 $$\vert (\pi_{\lambda}(\theta_{\tilde{y}}(g))e,\pi_{\lambda}(hug)\pi_{\lambda}(f)e)\vert <<\Xi^G(hu)$$
 pour tous $\lambda$, $h$ et $u$. Gr\^ace \`a 4.1(3), l'expression ci-dessus est absolument convergente. On peut permuter les int\'egrales:
 $$J_{L,{\cal O}}(\epsilon,f,g)=[i{\cal A}_{{\cal O}}^{\vee}:i{\cal A}_{L,F}^{\vee}]^{-1}\sum_{e\in {\cal B}_{{\cal O}}^{K_{f}}}\int_{i{\cal A}_{L,F}^*}m(\tau_{\lambda})$$
 $$\int_{H(F)U(F)_{c}}(\rho(h)\epsilon,\tilde{\rho}(\tilde{y})\epsilon) (\pi_{\lambda}(\theta_{\tilde{y}}(g))e,\pi_{\lambda}(hug)\pi_{\lambda}(f)e)\bar{\xi}(u)d\lambda\,du\,dh.$$
 L'int\'egrale int\'erieure est ${\cal L}_{\pi_{\lambda},\rho,c}(\epsilon\otimes \pi_{\lambda}(\theta_{\tilde{y}}(g))e,\tilde{\rho}(\tilde{y})\epsilon\otimes \pi_{\lambda}(g)\pi_{\lambda}(f)e)$. Quitte \`a accro\^{\i}tre $c_{0}$, c'est aussi  ${\cal L}_{\pi_{\lambda},\rho}(\epsilon\otimes \pi_{\lambda}(\theta_{\tilde{y}}(g))e,\tilde{\rho}(\tilde{y})\epsilon\otimes \pi_{\lambda}(g)\pi_{\lambda}(f)e)$ (l'utilisation directe du lemme 5.2 nous fournit pour chaque $\lambda$ un $c_{0}$ convenable; comme en (1) ci-dessus, on voit qu'en fait, on peut choisir $c_{0}$ ind\'ependant de $\lambda$). D'o\`u
 $$(2) \qquad J_{L,{\cal O}}(\epsilon,f,g)=[i{\cal A}_{{\cal O}}^{\vee}:i{\cal A}_{L,F}^{\vee}]^{-1}$$
 $$\sum_{e\in {\cal B}_{{\cal O}}^{K_{f}}}\int_{i{\cal A}_{L,F}^*}m(\tau_{\lambda}){\cal L}_{\pi_{\lambda},\rho}(\epsilon\otimes \pi_{\lambda}(\theta_{\tilde{y}}(g))e,\tilde{\rho}(\tilde{y})\epsilon\otimes \pi_{\lambda}(g)\pi_{\lambda}(f)e)d\lambda.$$

On fixe des familles $(\epsilon_{j})_{j=1,...,n}$, $(e_{j})_{j=1,...,n}$ et $(\varphi_{j})_{j=1,...,n}$ v\'erifiant 
 le lemme 5.6(ii).  Pour $e\in {\cal K}_{Q,\tau}^G$, $g\in M(F)K$ et $\lambda\in i{\cal A}_{L,F}^*$, posons
 $$X_{\lambda}(e,g)=\sum_{\epsilon\in {\cal B}_{\rho}^{K'_{g}}}{\cal L}_{\pi_{\lambda},\rho}(\epsilon\otimes  e',\underline{\epsilon}\otimes \underline{e}),$$
 o\`u $\underline{\epsilon}=\tilde{\rho}(\tilde{y})\epsilon$, $e'=\pi_{\lambda}(\theta_{\tilde{y}}(g))e$ et $\underline{e}=\pi_{\lambda}(g)\pi_{\lambda}(f)e$. Introduisons une forme sesquilin\'eaire $L_{\pi_{\lambda},\rho}$ comme en 5.3, soit $C_{\lambda}>0$ la constante telle que le (ii) du lemme 5.3 soit v\'erifi\'ee.  On a
 $${\cal L}_{\pi_{\lambda},\rho}(\epsilon\otimes  e',\underline{\epsilon}\otimes \underline{e})=C_{\lambda}\overline{L_{\pi_{\lambda},\rho}(\underline{\epsilon},e')}L_{\pi_{\lambda},\rho}(\epsilon,\underline{e}),$$
 tandis que la propri\'et\'e du lemme 5.6(ii) s'\'ecrit
 $$1=\sum_{j=1,...,n}C_{\lambda}\varphi_{j}(\lambda)\overline{L_{\pi_{\lambda},\rho}(\epsilon_{j},e_{j})}L_{\pi_{\lambda},\rho}(\epsilon_{j},e_{j}).$$
 Par multiplication, on obtient
 $$X_{\lambda}(e,g)=\sum_{j=1,...,n}\varphi_{j}(\lambda)X_{\lambda,j}(e,g),$$
 o\`u 
 $$X_{\lambda,j}(e,g)=\sum_{\epsilon\in {\cal B}_{\rho}^{K'_{g}}}C_{\lambda}^2\overline{L_{\pi_{\lambda},\rho}(\underline{\epsilon},e')}L_{\pi_{\lambda},\rho}(\epsilon,\underline{e})\overline{L_{\pi_{\lambda},\rho}(\epsilon_{j},e_{j})}L_{\pi_{\lambda},\rho}(\epsilon_{j},e_{j}).$$
 Fixons $j$. Le produit d'un $C_{\lambda}$ et des deux facteurs extr\^emes ci-dessus est \'egal \`a
 $${\cal L}_{\pi_{\lambda},\rho}(\epsilon_{j}\otimes e',\underline{\epsilon}\otimes e_{j}).$$
 Le produit d'un $C_{\lambda}$ et des deux facteurs restants est \'egal \`a
 $${\cal L}_{\pi_{\lambda},\rho}(\epsilon\otimes e_{j},\epsilon_{j}\otimes \underline{e}).$$
 Supposons $g\in M(F)K$. Par un argument d\'ej\`a utilis\'e plusieurs fois, on peut fixer $c_{0}$ ind\'ependant de $g$ et $\lambda$ de sorte que l'on puisse remplacer ${\cal L}_{\pi_{\lambda},\rho}$ par ${\cal L}_{\pi_{\lambda},\rho,c}$ dans ces expressions, pourvu que $c\geq c_{0}$. On obtient
 $$X_{\lambda,j}(e,g)=\sum_{\epsilon\in {\cal B}_{\rho}^{K'_{g}}}{\cal L}_{\pi_{\lambda},\rho,c}(\epsilon\otimes e_{j},\epsilon_{j}\otimes \underline{e})\int_{H(F) U(F)_{c}}(\rho(h)\epsilon_{j},\underline{\epsilon})(e',\pi_{\lambda}(hu)e_{j})\bar{\xi}(u)du\,dh,$$
 puis
  $$(3) \qquad X_{\lambda,j}(e,g) =\int_{H(F)U(F)_{c}}(\rho(h)\epsilon_{j}, \epsilon^{\sharp})(e',\pi_{\lambda}(hu)e_{j})\bar{\xi}(u)du\,dh,$$
  o\`u
  $$\epsilon^{\sharp}=\sum_{\epsilon\in {\cal B}_{\rho}^{K'_{g}}}{\cal L}_{\pi_{\lambda},\rho,c}(\epsilon\otimes e_{j},\epsilon_{j}\otimes \underline{e})\underline{\epsilon}.$$
  Ecrivons
  $${\cal L}_{\pi_{\lambda},\rho,c}(\epsilon\otimes e_{j},\epsilon_{j}\otimes \underline{e})=\int_{H(F)}(\rho(h')\epsilon,\epsilon_{j})\Lambda(h')dh',$$
  o\`u
  $$\Lambda(h')=\int_{U(F)_{c}}(e_{j},\pi_{\lambda}(h'u')\underline{e})\bar{\xi}(u')du'.$$
  Pour tout $\epsilon''\in E_{\rho}$, on a
  $$(\epsilon'',\epsilon^{\sharp})=\sum_{\epsilon\in {\cal B}_{\rho}^{K'_{g}}}{\cal L}_{\pi_{\lambda},\rho,c}(\epsilon\otimes e_{j},\epsilon_{j}\otimes \underline{e})(\epsilon'',\underline{\epsilon})=\int_{H(F)}\sum_{\epsilon\in {\cal B}_{\rho}^{K'_{g}}}(\epsilon,\rho(h^{_{'}-1})\epsilon_{j})(\epsilon'',\underline{\epsilon})\Lambda(h')dh'$$
  $$=\sum_{H(F)/K'_{g}}\sum_{\epsilon\in {\cal B}_{\rho}^{K'_{g}}}(\epsilon,\epsilon_{j}(h'))(\epsilon'',\underline{\epsilon}),$$
  o\`u
  $$\epsilon_{j}(h')=\int_{K'_{g}}\rho(k^{-1}h^{_{'}-1})\epsilon_{j}\Lambda(h'k)dk.$$
  La fonction $\Lambda$ est invariante \`a droite par $K'_{g}$ et $\epsilon_{j}(h')$ est invariant par $K'_{g}$. Donc
  $$\epsilon_{j}(h')=\sum_{\epsilon\in {\cal B}_{\rho}^{K'_{g}}}(\epsilon,\epsilon_{j}(h'))\epsilon,$$
  et
  $$\tilde{\rho}(\tilde{y})\epsilon_{j}(h'))=\sum_{\epsilon\in {\cal B}_{\rho}^{K'_{g}}}(\epsilon,\epsilon_{j}(h'))\underline{\epsilon}.$$
  On obtient alors
  $$(\epsilon'',\epsilon^{\sharp})=\sum_{H(F)/K'_{g}}(\epsilon'',\tilde{\rho}(\tilde{y})\epsilon_{j}(h'))=\int_{H(F)}(\epsilon'',\tilde{\rho}(\tilde{y}h^{_{'}-1})\epsilon_{j})\Lambda(h')dh'$$
  $$=\int_{H(F)U(F)_{c}}(\rho(\theta_{\tilde{y}}(h')\epsilon'',\tilde{\rho}(\tilde{y})\epsilon_{j})(e_{j},\pi_{\lambda}(h'u')\underline{e})\bar{\xi}(u')du'\,dh'.$$
  Reportons cette expression dans (3). On obtient
  $$(4) \qquad X_{\lambda,j}(e,g)=\int_{H(F)U(F)_{c}}\int_{H(F)U(F)_{c}}(\rho(\theta_{\tilde{y}}(h')h)\epsilon_{j},\tilde{\rho}(\tilde{y})\epsilon_{j})$$
  $$(e_{j},\pi_{\lambda}(h'u'g)\pi_{\lambda}(f)e)(\pi_{\lambda}(\theta_{\tilde{y}}(g))e,\pi_{\lambda}(hu)e_{j})\bar{\xi}(u)\bar{\xi}(u')du'\,dh'\,du\,dh.$$
  On a
  
  (5)  cette expression est absolument convergente.
  
  En effet, pour $g$ fix\'e, elle est major\'ee en valeur absolue par
  $$\int_{H(F)U(F)_{c}}\int_{H(F)U(F)_{c}}\Xi^G(h'u')\Xi^H(\theta_{\tilde{y}}(h')h)\Xi^G(hu)du'\,dh'\,du\,dh,$$
  qui est convergente d'apr\`es 4.1(4).
  
  Pour deux entiers $c,c'\in {\mathbb N}$ et pour $g\in M(F)K$, posons
  $$X_{\lambda,j,c,c'}(e,g)=\int_{H(F)U(F)_{c}}\int_{H(F)U(F)_{c'}}(\rho(h)\epsilon_{j},\tilde{\rho}(\tilde{y})\epsilon_{j})$$
  $$(e_{j},\pi_{\lambda}(h'u'g)\pi_{\lambda}(f)e))(\pi_{\lambda}(\theta_{\tilde{y}}(h'u'g))e,\pi_{\lambda}(hu)e_{j})\bar{\xi}(u)du'\,dh'\,du\,dh.$$
  Cette expression est absolument convergente comme la pr\'ec\'edente. On a
  
  (6) il existe $c_{0}$ ind\'ependant de $N$ et $\lambda$ tel que, pour $g\in M(F)K$,  $X_{\lambda,j,c,c'}(e,g)$ ne d\'epende pas de $c$ et $c'$, pourvu que $c\geq c_{0}$ et $c'\geq c_{0}$;  pour de tels $c,c'$, on a l'\'egalit\'e $X_{\lambda,j}(e,g)=X_{\lambda,j,c,c'}(e,g)$.
  
    Le proc\'ed\'e habituel montre qu'il existe $c_{0}$ tel que, pour $c\geq c_{0}$ et $g\in M(F)K$,  $X_{\lambda,j,c,c'}(e,g)$ ne d\'epend pas de $c$. Soient $c,c'\geq c_{0}$. Alors $X_{\lambda,j,c,c'}(e,g)=X_{\lambda,j,c',c'}(e,g)$. Par les changements de variables $h\mapsto \theta_{\tilde{y}}(h')h$, puis $u\mapsto h^{-1}\theta_{\tilde{y}}(u')hu$,   $X_{\lambda,j,c',c'}(e,g)$ est \'egal au membre de droite de (4) o\`u l'on a remplac\'e $c$ par $c'$. Mais celui-ci ne d\'epend pas de $c'$. Cela prouve (6).
 
 Pour $g\in M(F)K$, posons
 $$J_{L,{\cal O}}(\Theta_{\tilde{\rho}},\tilde{f},g)=[i{\cal A}_{{\cal O}}^{\vee}:i{\cal A}_{L,F}^{\vee}]^{-1}\sum_{e\in {\cal B}_{{\cal O}}^{K_{f}}}\sum_{j=1,...,n}\int_{i{\cal A}_{L,F}^*}m(\tau_{\lambda})\varphi_{j}(\lambda)X_{\lambda,j,c,c'}(e,g)d\lambda,$$
 o\`u $c,c'\geq c_{0}$, $c_{0}$ v\'erifiant (6). En revenant \`a la formule (2) et en rassemblant les calculs ci-dessus, on a montr\'e:
 
 (7) on a l'\'egalit\'e
 $$\sum_{\epsilon\in {\cal B}_{\rho}^{K_{N}}}J_{L,{\cal O}}(\epsilon,f,g)=J_{L,{\cal O}}(\Theta_{\tilde{\rho}},\tilde{f},g)$$
 pour tout $g\in M(F)K$.
 
 \bigskip
 
 \subsection{Une premi\`ere approximation}

D'apr\`es 6.2(1) et 6.3(7), on a l'\'egalit\'e
$$J(\Theta_{\tilde{\rho}},\tilde{f},g)=\sum_{L\in {\cal L}(A_{d})}\vert W^L\vert \vert W^G\vert ^{-1}\sum_{{\cal O}\in \{\Pi_{2}(L)\}_{f}}J_{L,{\cal O}}(\Theta_{\tilde{\rho}},\tilde{f},g)$$
pour tout $g\in M(F)K$. Soit $N\geq1$. Par d\'efinition, $J_{N}(\Theta_{\tilde{\rho}},\tilde{f})$ est l'int\'egrale de $J(\Theta_{\tilde{\rho}},\tilde{f},g)\kappa_{N}(g)$ sur $g\in H(F)U(F)\backslash G(F)$ ou, ce qui revient au m\^eme, l'int\'egrale de $J(\Theta_{\tilde{\rho}},\tilde{f},mk)\kappa_{N}(m)\delta_{P}(m)^{-1}$ sur $m\in H(F)\backslash M(F)$ et $k\in K$. C'est-\`a-dire
$$J_{N}(\Theta_{\tilde{\rho}},\tilde{f})=\int_{H(F)/ M(F)}\int_{K}\sum_{L\in {\cal L}(A_{d})}\vert W^L\vert \vert W^G\vert ^{-1}$$
$$\sum_{{\cal O}\in \{\Pi_{2}(L)\}_{f}}J_{L,{\cal O}}(\Theta_{\tilde{\rho}},\tilde{f},mk)\kappa_{N}(m)\delta_{P}(m)^{-1}dk\,dm.$$

Soient $L\in {\cal L}(A_{d})$ et ${\cal O}\in \{\Pi_{2}(L)\}_{f}$. Reprenons les notations du paragraphe pr\'ec\'edent. Soit $c_{0}$ v\'erifiant la relation (6) de ce paragraphe et soit $c\geq c_{0}$. Pour tout entier $C\in {\mathbb N}$, posons
$$J_{L,{\cal O},N,C}(\Theta_{\tilde{\rho}},\tilde{f})=[i{\cal A}_{{\cal O}}^{\vee}:i{\cal A}_{L,F}^{\vee}]^{-1}\sum_{e\in {\cal B}_{{\cal O}}^{K_{f}}}\sum_{j=1,...,n}\int_{i{\cal A}_{L,F}^*}m(\tau_{\lambda})\varphi_{j}(\lambda)\int_{H(F)U(F)_{c}}{\bf 1}_{\sigma< Clog(N)}(hu)$$
$$(\rho(h)\epsilon_{j},\tilde{\rho}(\tilde{y})\epsilon_{j})\bar{\xi}(u)\int_{G(F)}(e_{j},\pi_{\lambda}(g)\pi_{\lambda}(f)e))(\pi_{\lambda}(\theta_{\tilde{y}}(g))e,\pi_{\lambda}(hu)e_{j})\kappa_{N}(g)dg\,du\,dh\,d\lambda.$$

\ass{Lemme}{(i) Cette expression est absolument convergente.

(ii) Il existe $C$ tel que l'on ait la majoration
$$\vert J_{N}(\Theta_{\tilde{\rho}},\tilde{f})-\sum_{L\in {\cal L}(A_{d})}\vert W^L\vert \vert W^G\vert ^{-1}\sum_{{\cal O}\in \{\Pi_{2}(L)\}_{f}}J_{L,{\cal O},N,C}(\Theta_{\tilde{\rho}},\tilde{f})\vert << N^{-1}$$
pour tout entier $N\geq1$.}

La preuve est la m\^eme que celle du lemme 6.4 de [W3], en utilisant les majorations 4.1(5), (6) et (7).

On fixe $C$ v\'erifiant l'assertion (ii) ci-dessus. Dans les paragraphes suivants et jusqu'en 6.9, on fixe $L\in {\cal L}(A_{d})$ et ${\cal O}\in \{\Pi_{2}(L)\}_{f}$.

\bigskip

\subsection{Changement de fonction de troncature}

Notons $\Delta$ l'ensemble des racines simples de $A_{d}$ dans $\mathfrak{u}_{d}$. A toute racine $\alpha\in \Delta$ sont associ\'ees une coracine $\check{\alpha}$, un poids $\varpi_{\check{\alpha}}$ et un copoids $\varpi_{\alpha}$. Notons $\tilde{\Delta}$ l'ensemble des restrictions \`a $\tilde{A}_{d}$ des \'el\'ements de $\Delta$. Posons simplement $\theta=\theta_{d}$. L'ensemble $\tilde{\Delta}$ est en bijection avec l'ensemble des orbites dans $\Delta$ pour l'action de $\theta$. Notons $\alpha\mapsto (\alpha)$ cette bijection et notons $n(\alpha)$ le nombre d'\'el\'ements de l'orbite $(\alpha)$. On a alors l'\'egalit\'e 
$$\alpha=n(\alpha)^{-1}\sum_{\beta\in (\alpha)}\beta$$
dans ${\cal A}_{A_{d}}^*$. On utilise des \'egalit\'es analogues pour d\'efinir $\check{\alpha}$, $\varpi_{\check{\alpha}}$ et $\varpi_{\alpha}$. Par exemple
$$\varpi_{\check{\alpha}}=n(\alpha)^{-1}\sum_{\beta\in (\alpha)}\varpi_{\check{\beta}}.$$

Notons ${\cal A}_{A_{d}}^+$ le sous-ensemble des $\zeta\in {\cal A}_{A_{d}}$ tels que $\alpha(\zeta)\geq0$ pour tout $\alpha\in \Delta$. Fixons un r\'eel $\delta$ tel que $0<\delta<1$. Notons ${\cal D}$ l'ensemble des $Y\in {\cal A}_{A_{d}}^+\cap( {\cal A}_{\tilde{A}_{d},F}\otimes_{{\mathbb Z}}{\mathbb Q})$ tels que
$$inf\{\alpha(Y); \alpha\in \tilde{\Delta}\}\geq \delta sup\{\alpha(Y); \alpha\in \tilde{\Delta}\}.$$
Dans ce domaine, les fonctions $Y\mapsto inf\{\alpha(Y); \alpha\in \tilde{\Delta}\}$, $Y\mapsto sup\{\alpha(Y); \alpha\in \tilde{\Delta}\}$ et $Y\mapsto \vert Y\vert $ sont \'equivalentes.
Soit $Y\in {\cal D}$. Notons $\zeta\mapsto \varphi^+(\zeta,Y)$ la fonction caract\'eristique du sous-ensemble des $\zeta\in {\cal A}_{A_{d}}$ qui v\'erifient les deux conditions:

- $\alpha(\zeta)\geq 0$ pour tout $\alpha\in \Delta$;

- $\varpi_{\check{\alpha}}(Y-\zeta)\geq0$ pour tout $\alpha\in\tilde{\Delta}$.

Remarquons que ce sous-ensemble est compact modulo ${\cal A}_{G}$. Notons $g\mapsto \tilde{u}(g,Y)$ la fonction caract\'eristique de l'ensemble des $g\in G(F)$ pour lesquels il existe $k,k'\in K$ et $a\in A_{d}(F)$ de sorte que $g=kak'$ et $\varphi^+(H_{A_{d}}(a),Y)=1$.

 De m\^eme que l'on a d\'efini la fonction $f$, d\'efinissons une fonction $f'$ sur $G(F)$ par $f'(g)=\tilde{f}(g\boldsymbol{\theta}_{d})$. Fixons un sous-groupe ouvert compact $K_{f'}$ de $K$ tel que $f'$ soit biinvariante par $K_{f'}$, $K_{f'}$ soit distingu\'e dans $K$ et invariant par $\theta$. Fixons une base orthonorm\'ee ${\cal B}_{{\cal O}}^{K_{f'}}$ du sous-espace des invariants $({\cal K}_{Q,\tau}^{G})^{K_{f'}}$. Fixons $e',e''\in {\cal K}_{Q,\tau}^G$ et une fonction $\varphi$ sur $i{\cal A}_{L,F}^*$, que l'on suppose $C^{\infty}$. Pour $e\in {\cal K}_{Q,\tau}^G$, $g,g'\in G(F)$ et $\lambda\in i{\cal A}_{L,F}^*$, posons
$$\Phi(e,g,g',\lambda)=(\pi_{\lambda}(\theta(g))e,\pi_{\lambda}(g')e')(e'',\pi_{\lambda}(g)\pi_{\lambda}(f')e).$$
Posons
$$\Phi_{N}(g')=\sum_{e\in {\cal B}_{{\cal O}}^{K_{f'}}}\int_{G(F)}\int_{i{\cal A}_{L,F}^*}m(\tau_{\lambda})\varphi(\lambda) \Phi(e,g,g',\lambda)\kappa_{N}(g)\,d\lambda\,dg,$$
$$\Phi_{Y}(g')=\sum_{e\in {\cal B}_{{\cal O}}^{K_{f'}}}\int_{G(F)}\int_{i{\cal A}_{L,F}^*}m(\tau_{\lambda})\varphi(\lambda) \Phi(e,g,g',\lambda)\tilde{u}(g,Y) \,d\lambda\,dg.$$
 La premi\`ere expression est absolument convergente: cela r\'esulte de la convergence de
$$(1) \qquad \int_{G(F)}\kappa_{N}(g)\Xi^G(g)^2dg,$$
cf. 4.1(1). Montrons que la seconde est convergente dans l'ordre indiqu\'e. On peut l'\'ecrire
$$\Phi_{Y}(g')=\sum_{e\in {\cal B}_{{\cal O}}^{K_{f'}}}\int_{A_{G}(F)\backslash G(F)}\tilde{u}(g,Y)\int_{A_{G}(F)}\int_{i{\cal A}_{L,F}^*}m(\tau_{\lambda})\varphi(\lambda) \Phi(e,zg,g',\lambda)  \,d\lambda\,dz\,dg.$$ 
Parce que la fonction $g\mapsto \tilde{u}(g,Y)$ est \`a support compact modulo $A_{G}(F)$ et que la fonction \`a int\'egrer est localement constante, l'int\'egrale ext\'erieure est une somme finie. Il suffit de prouver la convergence dans l'ordre de la double int\'egrale int\'erieure. Notons $\omega$ la restriction \`a $A_{G}(F)$ du caract\`ere central de $\tau$. On a l'\'egalit\'e
$$\Phi(e,zg,g',\lambda)=\omega(z\theta(z)^{-1})e^{\lambda(H_{G}(z)-H_{G}(\theta(z)))}\Phi(e,g,g',\lambda)$$
pour tous $z$, $g$, $\lambda$.  Ainsi, l'int\'egrale int\'erieure en $\lambda$ d\'efinit une fonction de $z$ qui, au facteur $\omega(z\theta(z)^{-1})$ pr\`es, est une transform\'ee de Fourier partielle de la fonction que l'on a int\'egr\'ee, \'evalu\'ee au point $H_{G}(z)-H_{G}(\theta(z))$. Puisque la fonction de $\lambda$ est $C^{\infty}$, cette transform\'ee de Fourier est \`a d\'ecroissance rapide. On obtient une majoration
$$\vert \int_{i{\cal A}_{L,F}^*}m(\tau_{\lambda})\varphi(\lambda) \Phi(e,zg,g',\lambda)  \,d\lambda\vert <<(1+\vert H_{G}(z)-H_{G}(\theta(z))\vert) ^{-r}$$
pour tout r\'eel $r$. Mais la fonction $H\mapsto H-\theta(H)$ est injective sur ${\cal A}_{G}$. On a donc aussi une majoration
$$(1+\vert H_{G}(z)-H_{G}(\theta(z))\vert) ^{-r}<<\sigma(z)^{-r}.$$
Or l'int\'egrale
$$\int_{A_{G}(F)}\sigma(z)^{-r}\,dz$$
est convergente pour $r$ assez grand. L'int\'egrabilit\'e cherch\'ee en r\'esulte.

\ass{Proposition}{Soient $R$ et $\eta$ deux r\'eels, avec $0<\eta<1$. Il existe des r\'eels $c_{1},c_{2}>0$ tels que l'on ait la majoration
$$\vert \Phi_{N}(g')-\Phi_{Y}(g')\vert <<N^{-R}$$
pour tout $N\geq2$, tout $g'\in G(F)$ tel que $\sigma(g')\leq Clog(N)$ et tout $Y\in {\cal D}$ v\'erifiant les in\'egalit\'es $c_{1}N^{\eta}\leq  \vert Y\vert \leq c_{2 }N$.}

Preuve. Montrons d'abord:

(2) il existe $R_{1}\geq0$ tel que
$$\vert \Phi_{N}(g')\vert <<N^{R_{1}}\text{ et }\vert \Phi_{Y}(g')\vert <<N^{R_{1}}\vert Y\vert ^{R_{1}}$$
pour tout $N\geq1$, tout $Y\in {\cal D}$ et tout $g'\in G(F)$ tel que $\sigma(g')\leq Clog(N)$.

Pour $e\in {\cal K}_{Q,\tau}^G$, on a la majoration
$$\vert \Phi(e,g,g',\lambda)\vert <<\Xi^G(g^{_{'}-1}g)\Xi^G(g)$$
pour tous $\lambda,g,g'$. Gr\^ace \`a [W3] 3.3(5), il existe $R_{2}\geq0$ tel que
$$\Xi^G(g^{_{'}-1}g)<<exp(R_{2}\sigma(g'))\Xi^G(g).$$
D'apr\`es l'hypoth\`ese sur $g'$, on en d\'eduit
$$\vert \Phi(e,g,g',\lambda)\vert <<N^{CR_{2}}\Xi^G(g)^2.$$
D'apr\`es 4.1(1), il existe $R_{3}\geq0$ tel que l'int\'egrale (1) soit essentiellement major\'ee par $N^{R_{3}}$. On en d\'eduit l'assertion (2) pour $\Phi_{N}(g')$. 

Pour prouver l'assertion concernant la fonction $\Phi_{Y}(g')$, on reprend le raisonnement montrant la convergence de cette expression. On peut \'ecrire
$$\Phi_{Y}(g')=\sum_{e\in {\cal B}_{{\cal O}}^{K_{f'}}}\int_{A_{G}(F)\backslash G(F)}\tilde{u}(g,Y)\int_{A_{G}(F)}\int_{i{\cal A}_{L,F}^*/i{\cal A}_{G,F}^*}$$
$$\int_{i{\cal A}_{G,F}^*}m(\tau_{\lambda+\xi})\varphi(\lambda+\xi) \Phi(e,zg,g',\lambda+\xi)\,d\xi  \,d\lambda\,dz\,dg.$$
On peut identifier $i{\cal A}_{L,F}^*/i{\cal A}_{G,F}^*$ \`a $i{\cal A}_{L}^{G,*}/(i{\cal A}_{L}^{G,*}\cap(i{\cal A}_{L,F}^{\vee}+i{\cal A}_{G}^*))$. On peut identifier $A_{G}(F)\backslash G(F)$ \`a un ensemble de repr\'esentants sur lequel la fonction $H_{G}$ est born\'ee. On peut aussi d\'ecomposer $f'$ en $f'=\sum_{X\in {\cal X}}f'_{X}$, o\`u ${\cal X}$ est un sous-ensemble fini de ${\cal A}_{G,F}$ et $f'_{X}$ est une fonction dont le support est contenu dans l'ensemble des $x\in G(F)$ tels que $H_{G}(x)=X$. On a alors
$$\Phi(e,zg,g',\lambda+\xi)=\omega(z\theta(z)^{-1})\sum_{X\in {\cal X}}e^{\xi( B(z,g,g',X))} (\pi_{\lambda}(\theta(g))e,\pi_{\lambda}(g')e')(e'',\pi_{\lambda}(g)\pi_{\lambda}(f'_{X})e),$$
o\`u $B(z,g,g',X)=H_{G}(g')+H_{G}(g)-H_{G}(\theta(g))+X+H_{G}(z)-H_{G}(\theta(z))$. 
Par inversion de Fourier, on a une majoration
$$\vert \int_{i{\cal A}_{G,F}^*}m(\tau_{\lambda+\xi})\varphi(\lambda+\xi)e^{\xi( B(z,g,g',X))}\,d\xi\vert <<$$
$$(1+ \vert H_{G}(g')+H_{G}(g)-H_{G}(\theta(g))+X+H_{G}(z)-H_{G}(\theta(z))\vert )^{-r}$$
pour tout r\'eel $r$, et ce dernier terme est lui-m\^eme essentiellement major\'e par
$$(1+\vert H_{G}(g')\vert )^{r}(1+\vert H_{G}(z)-H_{G}(\theta(z))\vert )^{-r},$$
puis par $log(N)^r \sigma(z)^{-r}$. Alors $\Phi_{Y}(g')$ est major\'e par une somme finie d'expressions
$$log(N)^r\int_{A_{G}(F)\backslash G(F)}\tilde{u}(g,Y)\int_{A_{G}(F)}\int_{i{\cal A}_{L,F}^*/i{\cal A}_{G,F}^*}\sigma (z)^{-r}$$
$$\vert (\pi_{\lambda}(\theta(g))e,\pi_{\lambda}(g')e')(e'',\pi_{\lambda}(g)\pi_{\lambda}(f'_{X})e)\vert \,d\lambda\,dz\,dg.$$
Par les m\^emes calculs utilis\'es ci-dessus pour majorer $\Phi_{N}(g')$, ceci est essentiellement major\'e par
$$log(N)^rN^{R_{4}}\left(\int_{A_{G}(F)\backslash G(F)}\Xi^G(g)^2\tilde{u}(g,Y)\,dg\right)\left(\int_{A_{G}(F)}\sigma(z)^{-r}\,dz\right)$$
pour un r\'eel $R_{4}$ convenable. L'int\'egrale en $g$ est major\'ee par $\vert Y\vert ^{R_{5}}$ pour un r\'eel $R_{5}$ convenable. Il suffit de fixer $r$ tel que l'int\'egrale en $z$ soit convergente pour obtenir la majoration (2) pour $\Phi_{Y}(g')$.

On note $A_{d}(F)^+$ l'ensemble des $a\in A_{d}(F)$ tels que $\alpha(H_{A_{d}}(a))\geq 0$ pour tout $\alpha\in \Delta$. On d\'efinit  une fonction $D$ sur cet ensemble par
$$D(a)=mes(KaK)mes(K\cap A_{d}(F))^{-1}.$$
On a la formule d'int\'egration
$$\int_{G(F)}\phi(g)dg=\int_{K}\int_{K}\int_{A_{d}(F)^+}D(a)\phi(k_{1}ak_{2})da\,dk_{1}\,dk_{2}$$
pour toute fonction $\phi\in C_{c}^{\infty}(G(F))$. De $Y$ se d\'eduit comme en 3.6 une famille $(G,A_{d})$ orthogonale ${\cal Y}$. Pour tout $Q_{1}=L_{1}U_{1}\in {\cal F}(A_{d})$, on en d\'eduit des fonctions $\zeta\mapsto \sigma_{A_{d}}^{Q_{1}}(\zeta,{\cal Y})$ et $\zeta\mapsto \tau_{Q_{1}}(\zeta-Y_{Q_{1}})$ sur ${\cal A}_{A_{d}}$. Tous ces termes ont \'et\'e d\'efinis par Arthur. En fait, nous ne les utiliserons que dans le cas o\`u $B_{d}\subset Q_{1}$ et $\zeta\in {\cal A}_{A_{d}}^+$. Dans ce cas, on a
$$\sigma_{A_{d}}^{Q_{1}}(\zeta,{\cal Y})=1 \,\iff\, \varpi_{\alpha}(Y^{L_{1}}-\zeta^{L_{1}})\geq0\text{ pour tout }\alpha\in \Delta^{L_{1}},$$
$$\tau_{Q_{1}}(\zeta-Y_{Q_{1}})=1\,\iff\,\alpha(\zeta_{L_{1}}-Y_{L_{1}})>0\text{ pour tout }\alpha\in \Delta-\Delta^{L_{1}},$$
o\`u, pour tout $\zeta$, on note $\zeta^{L_{1}}$ et $\zeta_{L_{1}}$ les projections orthogonales de $\zeta$ sur ${\cal A}_{A_{d}}^{L_{1}}$, resp. ${\cal A}_{L_{1}}$, et $\Delta^{L_{1}}$ est l'ensemble des racines simples de $A_{d}$ dans $\mathfrak{u}_{d}\cap \mathfrak{l}_{1}$. On a l'\'egalit\'e
$$\sum_{Q_{1}\in {\cal F}(A_{d}),  B_{d}\subset Q_{1}}\sigma^{Q_{1}}_{A_{d}}(\zeta,{\cal Y})\tau_{Q_{1}}(\zeta-Y_{Q_{1}})=1$$
pour tout $\zeta\in{\cal A}_{A_{d}}^+$. Pour $Q_{1}\in {\cal F}(A_{d})$ contenant $B_{d}$, posons
$$\Phi_{N,Y,Q_{1}}(g')=\sum_{e\in {\cal B}_{{\cal O}}^{K_{f'}}}\int_{K}\int_{K}\int_{A_{d}(F)^+}\int_{i{\cal A}_{L,F}^*}m(\tau_{\lambda})\varphi(\lambda)$$
$$ \Phi(e,k_{1}ak_{2},g',\lambda)\kappa_{N}(k_{1}ak_{2})D(a)\sigma_{A_{d}}^{Q_{1}}(H_{A_{d}}(a),{\cal Y})\tau_{Q_{1}}(H_{A_{d}}(a)-Y_{Q_{1}})\,d\lambda\,da\,dk_{1}\,dk_{2}.$$
D\'efinissons de m\^eme $\Phi_{Y,Q_{1}}(g')$ en rempla\c{c}ant la fonction $\kappa_{N}(k_{1}ak_{2})$ par $\tilde{u}(k_{1}ak_{2},Y)$. Il r\'esulte de ce qui pr\'ec\`ede que l'on a
$$\Phi_{N}(g')=\sum_{Q_{1}\in {\cal F}(A_{d}),  B_{d}\subset Q_{1}}\Phi_{N,Y,Q_{1}}(g'),$$
$$\Phi_{Y}(g')=\sum_{Q_{1}\in {\cal F}(A_{d}),  B_{d}\subset Q_{1}}\Phi_{Y,Q_{1}}(g').$$
On va montrer:

(3) il existe $c_{2}>0$ tel que, si $\vert Y\vert \leq c_{2}N$, on a  la majoration
 $\vert \Phi_{N,Y,G}(g')-\Phi_{Y,G}(g')\vert <<N^{-R}$ pour tout  $N\geq1$ et tout $g'\in G(F)$ tel que $\sigma(g')\leq Clog(N)$;

(4) il existe $c_{1}, c_{2}>0$ tel que, si $c_{1}N^{\eta}\leq \vert Y\vert \leq c_{2}N$, on a les majorations
 $$\vert \Phi_{N,Y,Q_{1}}(g')\vert <<N^{-R}\text{ et }\vert \Phi_{Y,Q_{1}}(g')\vert <<N^{-R}$$
 pour tout $N\geq1$,  tout $g'$ comme ci-dessus et tout $Q_{1}\in {\cal F}(A_{d})$ tel que $B_{d}\subset Q_{1 }\subsetneq G$.

Cela d\'emontrera la proposition.

Prouvons (3).  Le support de la fonction $\zeta\mapsto \sigma_{A_{d}}^G(\zeta,{\cal Y})$ sur ${\cal A}_{A_{d}}^+$ est contenu dans celui de la fonction $\zeta\mapsto \varphi^+(\zeta,Y)$. On peut donc supprimer le terme $\tilde{u}(k_{1}ak_{2},Y)$ dans la d\'efinition de $\Phi_{Y,G}(g')$. La fonction $ \Phi_{N,Y,G}(g')-\Phi_{Y,G}(g')$ est donc d\'efinie par la m\^eme int\'egrale que $\Phi_{N,Y,G}(g')$, o\`u on remplace $\kappa_{N}(k_{1}ak_{2},Y)$ par $\kappa_{N}(k_{1}ak_{2},Y)-1$. On v\'erifie qu'il existe $c_{3}>0$ tel que le support de $\kappa_{N}$ contienne l'ensemble des $g\in G(F)$ tels que $\sigma(g)\leq c_{3}N$. On peut donc remplacer $\kappa_{N}(k_{1}ak_{2},Y)-1$ par la fonction caract\'eristique de l'ensemble des $a$ tels que $\sigma(a)> c_{3}N$. La condition $\sigma_{A_{d}}^G(H_{A_{d}}(a),{\cal Y})=1$ entra\^{\i}ne une majoration $\vert H_{A_{d}}^G(a)\vert <<\vert Y\vert $. Pour $c_{2}$ assez petit, la double condition ci-dessus entra\^{\i}ne $\vert H_{G}(a)\vert > c_{4}N$, pour un certain $c_{4}>0$. On peut alors reprendre la preuve de l'assertion (2) concernant $\Phi_{Y}(g')$. On obtient une majoration
$$\vert \Phi_{N,Y,G}(g')-\Phi_{Y,G}(g')\vert<<log(N)^rN^{R_{4}}\int_{A_{G}(F)\backslash  A_{d}(F)^+}\Xi^G(a)^2 \sigma_{A_{d}}(H_{A_{d}}(a),{\cal Y})D(a)\,da\,$$
$$\int_{z\in A_{G}(F); \vert H_{G}(z)\vert > c_{4}N-c_{5}} \sigma(z)^{-r}\,dz,$$
o\`u $c_{5}$ est un majorant de la fonction $H_{G}$ sur un ensemble de repr\'esentants du quotient $A_{G}(F)\backslash  A_{d}(F)$. La premi\`ere int\'egrale est major\'ee par
$$\int_{g\in A_{G}(F)\backslash G(F); \sigma(g)\leq c_{6}\vert Y\vert }\Xi(g)^2\,dg$$
pour $c_{6}>0$ convenable, donc par $\vert Y\vert ^{R_{6}}$ pour un r\'eel $R_{6}$ convenable, ou encore par $\vert N\vert ^{R_{6}}$. La seconde int\'egrale est essentiellement major\'ee par $\vert N\vert ^{-r}$. Le tout est essentiellement major\'e par
  $log(N)^rN^{R_{4}+R_{6}-r} $. En prenant $r$ assez grand, on obtient (3).

On fixe d\'esormais $Q_{1}\in {\cal F}(A_{d})$ tel que $B_{d}\subset Q_{1}\subsetneq G$. Introduisons un r\'eel $\epsilon$ tel que $0<\epsilon<1$. Comme ci-dessus, on a
$$\sum_{Q_{2}\in {\cal F}(A_{d}),  B_{d}\subset Q_{2}}\sigma^{Q_{2}}_{A_{d}}(\zeta,\epsilon{\cal Y})\tau_{Q_{2}}(\zeta-\epsilon Y_{Q_{2}})=1$$
pour tout $\zeta\in{\cal A}_{A_{d}}^+$. Supposons $\sigma_{A_{d}}^{Q_{1}}(\zeta,{\cal Y})\tau_{Q_{1}}(\zeta-Y_{Q_{1}})=1$. Pour $\alpha\in \Delta-\Delta^{L_{1}}$, on a $\alpha(\zeta-Y)>0$. Si $\sigma^{Q_{2}}_{A_{d}}(\zeta,\epsilon{\cal Y})=1$, on a $ \alpha(\zeta)<<\epsilon\vert Y\vert $ pour $\alpha\in \Delta^{L_{2}}$. Si $\epsilon$ est assez petit, ces in\'egalit\'es ne sont compatibles que si  $\Delta^{L_{2}}\subset \Delta^{L_{1}}$. C'est-\`a-dire que, pour un tel $\zeta$, la somme ci-dessus se limite aux $Q_{2}\subset Q_{1}$. Pour un tel $Q_{2}$ et pour $a\in A_{d}(F)$, posons
$$S_{Q_{1},Q_{2}}(a)=\sigma_{A_{d}}^{Q_{1}}(H_{A_{d}}(a),{\cal Y})\tau_{Q_{1}}(H_{A_{d}}(a)-Y_{Q_{1}})\sigma_{A_{d}}^{Q_{2}}(H_{A_{d}}(a),\epsilon{\cal Y})\tau_{Q_{2}}(H_{A_{d}}(a)-\epsilon Y_{Q_{2}}).$$
On peut d\'ecomposer
$$\Phi_{N,Y,Q_{1}}(g')=\sum_{Q_{2}\in {\cal F}(A_{d}); B_{d}\subset Q_{2}\subset Q_{1}}\Phi_{N,Y,Q_{1},Q_{2}}(g'),$$
o\`u
$$\Phi_{N,Y,Q_{1},Q_{2}}(g')=\sum_{e\in {\cal B}_{{\cal O}}^{K_{f'}}}\int_{K}\int_{K}\int_{A_{d}(F)^+}\int_{i{\cal A}_{L,F}^*}m(\tau_{\lambda})\varphi(\lambda)$$
$$ \Phi(e,k_{1}ak_{2},g',\lambda)\kappa_{N}(k_{1}ak_{2})D(a) S_{Q_{1},Q_{2}}(a)\,d\lambda\,da\,dk_{1}\,dk_{2}.$$
On peut d\'ecomposer de m\^eme $\Phi_{Y,Q_{1}}(g')$. Fixons $Q_{2}=L_{2}U_{2}$ avec $B_{d}\subset Q_{2}\subset Q_{1}$. On va montrer

(5) il existe $\epsilon_{0}>0$ et, pour $\epsilon<\epsilon_{0}$, il existe  $c_{1},c_{2}>0$ tels que, si $c_{1}N^{\eta}\leq  \vert Y\vert \leq c_{2}N$, on a les majorations
 $$\vert \Phi_{N,Y,Q_{1},Q_{2}}(g')\vert <<N^{-R}\text{ et }\vert \Phi_{Y,Q_{1},Q_{2}}(g')\vert <<N^{-R}$$

Pour simplifier la r\'edaction, on va fixer $c_{1}$ et $c_{2}$ et supposer $c_{1}N^{\eta}\leq\vert Y\vert \leq c_{2}log(N)$. On montrera que toutes les propri\'et\'es dont on a besoin sont v\'erifi\'ees si  $\epsilon$  est assez petit, $c_{1}$ est assez grand relativement \`a $\epsilon$ et $c_{2}$ est assez petit relativement \`a $\epsilon$.   

Soient $g'\in G(F)$ tel que $\sigma(g')\leq Clog(N)$, $k_{1},k_{2}\in K$ et $a\in A_{d}(F)^+$. On pose $\zeta=H_{A_{d}}(a)$ et on suppose $\sigma_{A_{d}}^{Q_{2}}(\zeta,\epsilon{\cal Y})\tau_{Q_{2}}(\zeta-\epsilon Y_{Q_{2}})=1$. Cette propri\'et\'e entra\^{\i}ne l'existence de $c_{3}>0$ tel que, pour toute racine $\alpha$ de $A_{d}$ dans $\mathfrak{u}_{2}$, on ait la minoration
$$\qquad c_{3}\inf\{\epsilon\beta(Y); \beta\in \Delta\}\leq \alpha(\zeta),$$
donc aussi, quitte \`a changer $c_{3}$,
$$(6) \qquad \epsilon c_{1}c_{3}N^{\eta}\leq \alpha(\zeta).$$
Ecrivons $\theta(k_{1}^{-1})g'=u'l'k'$, avec $u'\in \theta(U_{2})(F)$, $l'\in \theta(L_{2})(F)$ et $k'\in K$. Notons $A_{d}(F)^{L_{2},+}$ le sous-ensemble des $a'\in A_{d}(F)$ tels que $\alpha(H_{A_{d}}(a'))\geq0$ pour tout $\alpha\in \Delta^{L_{2}}$ et posons $K_{2}=K\cap L_{2}(F)$. On a $\theta(l')^{-1}a\in L_{2}(F)$ et on  peut \'ecrire $\theta(l')^{-1}a=k_{3}a'k_{4}$, avec $k_{3},k_{4}\in K_{2}$ et $a'\in A_{d}(F)^{L_{2},+}$.  On pose $\zeta'=H_{A_{d}}(a')$. On a:

(7) si $\epsilon c_{1}$ est assez grand, $k_{2}^{-1}a^{-1}\theta(u')^{-1}ak_{2}$ appartient \`a $K_{f'}$;

(8) pour tout $c>0$, $a'$ appartient \`a $A_{d}(F)^+$ et v\'erifie $\alpha(\zeta')> c N^{\eta/2}$ pour tout $\alpha\in \Delta-\Delta^{L_{2}}$ pourvu que $\epsilon c_{1}$ soit assez grand.

Ces deux propri\'et\'es r\'esultent ais\'ement de (6), cf. [W3] 6.6(7) et (8) pour plus de d\'etails. Soient $\lambda\in i{\cal A}_{L,F}^*$ et $e\in {\cal B}_{{\cal O}}^{K_{f'}}$. On a l'\'egalit\'e
$$\Phi(e,k_{1}ak_{2},g',\lambda)=(\pi_{\lambda}(\theta(\theta(l')^{-1}\theta(u')^{-1}ak_{2}))e,\pi_{\lambda}(k')e')(\pi_{\lambda}(k_{1}^{-1})e'',\pi_{\lambda}(ak_{2})\pi_{\lambda}(f')e).$$
Gr\^ace \`a (7), on peut supprimer $\theta(u')^{-1}$ dans cette expression, pourvu que $\epsilon c_{1}$ soit assez grand. On obtient
$$(9) \qquad \Phi(e,k_{1}ak_{2},g',\lambda)=(\pi_{\lambda}( \theta(a'k_{4}k_{2}))e,\pi_{\lambda}(\theta(k_{3})^{-1}k')e')(\pi_{\lambda}(k_{1}^{-1})e'',\pi_{\lambda}(ak_{2})\pi_{\lambda}(f')e).$$

Gr\^ace aux r\'esultats d'Harish-Chandra, on va approximer les produits scalaires par leurs termes constants faibles. Introduisons quelques notations. Posons
$$W(L_{2}\vert G\vert L)=\{s\in W^G; sLs^{-1}\subset L_{2}, B_{d}\cap L_{2}\subset sQs^{-1}\}.$$
Pour un \'el\'ement $s$ de cet ensemble, on note $\gamma(s):{\cal K}_{Q,\tau}^G\to {\cal K}_{sQs^{-1},\tau}^G$ l'op\'erateur d\'efini par $(\gamma(s)e_{0})(k)=e_{0}(s^{-1}k)$ pour tout $k\in K$. On pose $Q_{2,s}=(L_{2}\cap sQs^{-1})U_{2}$, $\underline{Q}_{2,s}=(L_{2}\cap sQs^{-1})\bar{U}_{2}$, o\`u $\bar{U}_{2}$ est le radical unipotent du parabolique $\bar{Q}_{2}$ oppos\'e \`a $Q_{2}$. Pour $\lambda\in i{\cal A}_{L,F}^*$, on d\'efinit les op\'erateurs
$$J_{Q_{2,s}\vert sQs^{-1}}((s\tau)_{s\lambda})\circ \gamma(s):{\cal K}_{Q,\tau}^G\to {\cal K}_{Q_{2,s},s\tau}$$
et
$$J_{\underline{Q}_{2,s}\vert sQs^{-1}}((s\tau)_{s\lambda})\circ \gamma(s):{\cal K}_{Q,\tau}^G\to {\cal K}_{\underline{Q}_{2,s},s\tau}.$$
Les termes $s\tau$ et $s\lambda$ ont la signification usuelle. Pour $e_{1}\in {\cal K}_{Q_{2,s},s\tau}^G$ et $e_{2}\in {\cal K}_{\underline{Q}_{2,s},s\tau}^G$, les restrictions de $e_{1}$ et $e_{2}$ \`a $K_{2}$ appartiennent au m\^eme espace ${\cal K}_{L_{2}\cap sQs^{-1},s\tau}^{L_{2}}$, qui est muni d'un produit scalaire naturel. On pose
$$(e_{1},e_{2})^{L_{2}}=\int_{K_{2}}(e_{1}(x),e_{2}(x))dx.$$
Posons
$$(\pi_{\lambda}(k_{1}^{-1})e'',\pi_{\lambda}(ak_{2})\pi_{\lambda}(f')e)_{Q_{2},\lambda}=$$
$$\sum_{s\in W(L_{2}\vert G\vert L)}(J_{Q_{2,s}\vert sQs^{-1}}((s\tau)_{s\lambda})\circ \gamma(s)\pi_{\lambda}(k_{1}^{-1})e'',J_{\underline{Q}_{2,s}\vert sQs^{-1}}((s\tau)_{s\lambda})\circ \gamma(s)\pi_{\lambda}(ak_{2})\pi_{\lambda}(f')e)^{L_{2}}.$$
Alors le lemme 6.5 de [W3] (qui est d\^u \`a Harish-Chandra) nous dit qu'il existe $R_{7}\geq0$ et $\mu>0$ tel que
$$\vert (\pi_{\lambda}(k_{1}^{-1})e'',\pi_{\lambda}(ak_{2})\pi_{\lambda}(f')e)-(\pi_{\lambda}(k_{1}^{-1})e'',\pi_{\lambda}(ak_{2})\pi_{\lambda}(f')e)_{Q_{2},\lambda}\vert <<$$
$$\delta_{Q_{2}}(a)^{-1/2}\Xi^{L_{2}}(a)\sigma(a)^{R_{7}}sup\{exp(-\mu\alpha(\zeta)); \alpha\in \Delta-\Delta^{L_{2}}\}.$$
Gr\^ace \`a (6),  cela entra\^{\i}ne que, pour tout entier $b\geq0$, on a
$$(10) \qquad \vert (\pi_{\lambda}(k_{1}^{-1})e'',\pi_{\lambda}(ak_{2})\pi_{\lambda}(f')e)-(\pi_{\lambda}(k_{1}^{-1})e'',\pi_{\lambda}(ak_{2})\pi_{\lambda}(f')e)_{Q_{2},\lambda}\vert <<$$
$$
\delta_{Q_{2}}(a)^{-1/2}\Xi^{L_{2}}(a)\sigma(a)^{R_{7}}N^{-b},$$
pourvu que $\epsilon c_{1}$ soit assez grand.

De m\^emes consid\'erations s'appliquent \`a l'autre produit de la formule (9). A cause du $\theta$ qui figure dans ce produit, on doit remplacer $Q_{2}$ par $\theta(Q_{2})$ dans les constructions. Pour $s\in W(\theta(L_{2})\vert G\vert L)$, on pose $\theta(Q)_{2,s}=(\theta(L_{2})\cap sQs^{-1})\theta(U_{2})$ et $\theta(\underline{Q})_{2,s}=(\theta(L_{2})\cap sQs^{-1})\theta(\bar{U}_{2})$ (les notations sont un peu ambigu\" es: $\theta(Q)_{2,s}$ n'est pas l'image par $\theta$ d'un hypoth\'etique parabolique $Q_{2,s}$).
On pose
$$(\pi_{\lambda}( \theta(a'k_{4}k_{2}))e,\pi_{\lambda}(\theta(k_{3})^{-1}k')e')_{\theta(Q_{2}),\lambda}=\sum_{s\in W(\theta(L_{2})\vert G\vert L)}$$
$$(J_{\theta(\underline{Q})_{2,s}\vert sQs^{-1}}((s\tau)_{s\lambda})\circ \gamma(s) \pi_{\lambda}(\theta(a'k_{4}k_{2}))e, J_{\theta(Q)_{2,s}\vert sQs^{-1}}((s\tau)_{s\lambda})\circ\gamma(s)\pi_{\lambda}(\theta(k_{3})^{-1}k')e')^{\theta(L_{2})}.$$
On a une majoration analogue \`a (10), o\`u, dans le deuxi\`eme membre, $Q_{2}$ et $a$ sont remplac\'es par $\theta(Q_{2})$ et $\theta(a')$ ou, ce qui revient au m\^eme, par $Q_{2}$ et $a'$. 

Posons
$$\Phi_{w}(e,k_{1}a k_{2},g',\lambda)=(\pi_{\lambda}( \theta(a'k_{4}k_{2})))e,\pi_{\lambda}(\theta(k_{3})^{-1}k')e')_{\theta(Q_{2}),\lambda}(\pi_{\lambda}(k_{1}^{-1})e'',\pi_{\lambda}(ak_{2})\pi_{\lambda}(f')e)_{Q_{ 2},\lambda}.$$
Notons $\Phi_{N,Y,Q_{1},Q_{2},w}(g')$ et $\Phi_{Y,Q_{1},Q_{2},w}(g')$ les termes obtenus en rempla\c{c}ant $\Phi(e,k_{1}a k_{2},g',\lambda)$ par $\Phi_{w}(e,k_{1}a k_{2},g',\lambda)$ dans les d\'efinitions de $\Phi_{N,Y,Q_{1},Q_{2}}(g')$ et $\Phi_{Y,Q_{1},Q_{2}}(g')$. D'apr\`es ce qui pr\'ec\`ede, $\Phi(e,k_{1}ak_{2},g',\lambda)-\Phi_{w}(e,k_{1}ak_{2},g',\lambda)$ est born\'e par la somme de
$$\delta_{Q_{2}}(a)^{-1/2}\Xi^{L_{2}}(a)\sigma(a)^{R_{7}}N^{-b}\vert (\pi_{\lambda}( \theta(a'k_{4}k_{2})))e,\pi_{\lambda}(\theta(k_{3})^{-1}k')e')\vert $$
et de
$$ \delta_{Q_{2}}(a')^{-1/2}\Xi^{L_{2}}(a')\sigma(a')^{R_{7}}N^{-b}\vert (\pi_{\lambda}(k_{1}^{-1})e'',\pi_{\lambda}(ak_{2})\pi_{\lambda}(f')e)_{Q_{2},\lambda}\vert .$$
Si on oublie les termes $N^{-b}$ un calcul analogue \`a celui utilis\'e dans la preuve de  (2) conduit \`a des majorations
$$\vert \Phi_{N,Y,Q_{1},Q_{2}}(g')-\Phi_{N,Y,Q_{1},Q_{2},w}(g')\vert <<N^{R_{8}},$$
$$\vert \Phi_{Y,Q_{1},Q_{2}}(g')-\Phi_{Y,Q_{1},Q_{2},w}(g')\vert <<N^{R_{8}}\vert Y\vert ^{R_{8}},$$
pour un certain r\'eel $R_{8}$. En r\'eintroduisant le terme $N^{-b}$, o\`u $b$ peut \^etre quelconque, et en se rappelant que l'on suppose $\vert Y\vert <<N$, les deux expressions ci-dessus sont major\'ees par $N^{-b}$ pourvu que $\epsilon c_{1}$ soit assez grand. Le bilan est que, pour d\'emontrer (5), on peut y remplacer $\Phi_{N,Y,Q_{1},Q_{2}}(g')$ et $\Phi_{Y,Q_{1},Q_{2}}(g')$ par $\Phi_{N,Y,Q_{1},Q_{2},w}(g')$ et $\Phi_{Y,Q_{1},Q_{2},w}(g')$.

On a
$$\Phi_{w}(e,k_{1}ak_{2},g',\lambda)=\sum_{s_{1}\in W(\theta(L_{2})\vert G\vert L), s_{2}\in W(L_{2}\vert G\vert L)}\Phi_{s_{1},s_{2}}(e,k_{1}ak_{2},g',\lambda),$$
o\`u
$$\Phi_{s_{1},s_{2}}(e,k_{1}ak_{2},g',\lambda)=$$
$$(J_{\theta(\underline{Q})_{2,s_{1}}\vert s_{1}Qs_{1}^{-1}}((s_{1}\tau)_{s_{1}\lambda})\circ \gamma(s_{1}) \pi_{\lambda}(\theta(a'k_{4}k_{2}))e, J_{\theta(Q)_{2,s_{1}}\vert s_{1}Qs_{1}^{-1}}((s_{1}\tau)_{s_{1}\lambda})\circ\gamma(s_{1})\pi_{\lambda}(\theta(k_{3})^{-1}k')e')^{\theta(L_{2})}$$
$$(J_{Q_{2,s_{2}}\vert s_{2}Qs_{2}^{-1}}((s_{2}\tau)_{s_{2}\lambda})\circ \gamma(s_{2})\pi_{\lambda}(k_{1}^{-1})e'',J_{\underline{Q}_{2,s_{2}}\vert s_{2}Qs_{2}^{-1}}((s_{2}\tau)_{s_{2}\lambda})\circ \gamma(s_{2})\pi_{\lambda}(ak_{2})\pi_{\lambda}(f')e)^{L_{2}}.$$
Au moins formellement, on peut d\'ecomposer de m\^eme $\Phi_{N,Y,Q_{1},Q_{2},w}(g')$ et $\Phi_{Y,Q_{1},Q_{2},w}(g')$ en somme de termes $\Phi_{N,Y,Q_{1},Q_{2},s_{1},s_{2}}(g')$ et $\Phi_{Y,Q_{1},Q_{2},s_{1},s_{2}}(g')$. Il y a un probl\`eme de convergence car les op\'erateurs d'entrelacement peuvent avoir des p\^oles. Ce probl\`eme dispara\^{\i}t gr\^ace \`a la multiplication par la mesure de Plancherel. En effet, soient $s_{1}\in W(\theta(L_{2})\vert G\vert L)$ et $s_{2}\in W(L_{2}\vert G\vert L)$. On a

(11) pour $e$ fix\'e, la fonction $m(\tau_{\lambda})\Phi_{s_{1},s_{2}}(e,k_{1}ak_{2},g',\lambda)$ est combinaison lin\'eaire de fonctions qui sont elles-m\^emes des produits $f_{1}(a',\lambda)f_{2}(a,\lambda)f_{3}(k_{1},k_{2},k_{3},k_{4},k')f_{4}(\lambda)$, o\`u

$$f_{1}(a',\lambda)=\delta_{Q_{2}}(a')^{-1/2}(Ind_{\theta(L_{2})\cap s_{1}Qs_{1}^{-1}}^{\theta(L_{2})}((s_{1}\tau)_{s_{1}\lambda},\theta(a'))e'_{1},e_{1}),$$
pour des \'el\'ements $e_{1},e'_{1}\in {\cal K}_{\theta(L_{2})\cap s_{1}Qs_{1}^{-1},s_{1}\tau}^{\theta(L_{2})}$; 

$$f_{2}(a,\lambda)=\delta_{Q_{2}}(a)^{1/2}(e'_{2},Ind_{L_{2}\cap s_{2}Qs_{2}^{-1}}^{L_{1}}((s_{2}\tau)_{s_{2}\lambda},a)e_{2}),$$
pour des \'el\'ements $e_{2},e'_{2}\in {\cal K}_{L_{2}\cap s_{2}Qs_{2}^{-1},s_{2}\tau}^{L_{2}}$;

$f_{3}$ est une fonction localement constante des variables $k_{1}$, $k_{2}$, $k_{3}$, $k_{4}$ et $k'$;

$f_{4}$ est une fonction $C^{\infty}$ de $\lambda$.

La preuve est la m\^eme qu'en [W3] 6.6(10). Le point central est le suivant. Soient $e_{1}\in {\cal K}_{Q,\tau}^G$ et $e'_{1}\in {\cal K}_{Q_{2,s_{1}},s_{1}\tau}^{G}$. Posons 
$$j_{Q_{2,s_{2}}}(\lambda)=(e'_{1},J_{Q_{2,s_{2}}\vert s_{2}Qs_{2}}((s_{2}\tau)_{s_{2}\lambda})\circ \gamma(s_{2})e_{1}).$$
D\'efinissons de m\^eme des fonctions $j_{\underline{Q}_{2,s_{2}}}(\lambda)$, $j_{\theta(Q)_{2,s_{1}}}(\lambda)$ et $j_{\theta(\underline{Q})_{2,s_{1}}}(\lambda)$. Alors la fonction
$$\lambda\mapsto m(\tau_{\lambda})\overline{j_{\theta(\underline{Q})_{2,s_{1}}}(\lambda)j_{Q_{2,s_{2}}}(\lambda)}j_{\underline{Q}_{2,s_{2}}}(\lambda)j_{\theta(Q)_{2,s_{1}}}(\lambda)$$
est $C^{\infty}$ sur $i{\cal A}_{L}^*$. Cela r\'esulte de la preuve du corollaire V.2.3 de [W4].

Gr\^ace \`a (11), les termes $\Phi_{N,Y,Q_{1},Q_{2},s_{1},s_{2}}(g')$ et $\Phi_{Y,Q_{1},Q_{2},s_{1},s_{2}}(g')$ sont bien d\'efinis. On peut fixer $s_{1}$ et $s_{2}$ et il nous suffit d'\'etablir la majoration
$$(12)\qquad \vert \Phi_{N,Y,Q_{1},Q_{2},s_{1},s_{2}}(g')\vert <<N^{-R}\text{ et }\vert \Phi_{Y,Q_{1},Q_{2},s_{1},s_{2}}(g')\vert <<N^{-R}.$$

Posons $s=s_{1}s_{2}^{-1}$.  Supposons d'abord

{\bf Hypoth\`ese.} Il n'existe pas de sous-groupe parabolique $Q'=L'U'$ de $G$ tel que $Q_{2}\subset Q'\subsetneq G$, $\theta(Q')=Q'$ et $s\in W^{L'}$.

Dans ce cas, on exprime $m(\tau_{\lambda})\Phi_{s_{1},s_{2}}(e,k_{1}ak_{2},g',\lambda)$ comme dans l'assertion (11). Quitte \`a multiplier $f_{4}$ par $\varphi$, $\Phi_{N,Y,Q_{1},Q_{2},s_{1},s_{2}}(g')$ s'\'ecrit comme combinaison lin\'eaire d'int\'egrales
$$\Psi_{N,Y}=\int_{K}\int_{K}\int_{A_{d}(F)^+}\int_{i{\cal A}_{L,F}^*}f_{1}(a',\lambda)f_{2}(a,\lambda)f_{3}(k_{1},k_{2},k_{3},k_{4},k')f_{4}(\lambda)$$
$$\kappa_{N}(k_{1}ak_{2})D(a)S_{Q_{1},Q_{2}}(a)d\lambda\,da\,dk_{1}\,dk_{2}.$$
 
Pour tout $x\in L_{2}$, choisissons des \'el\'ements $l(x)\in s_{2}L(F)s_{2}^{-1}$, $u(x)\in L_{2}(F)\cap s_{2}U_{Q}(F)s_{2}^{-1}$ et $k(x)\in K_{2}$ de sorte que $x=l(x)u(x)k(x)$. On a l'\'egalit\'e
$$f_{2}(a,\lambda)=\int_{K_{2}}f'_{2}(a,x)exp((s_{2}\lambda)(H_{L_{2}\cap s_{2}Qs_{2}^{-1}}(xa)))dx,$$
o\`u
$$f'_{2}(a,x)=\delta_{Q_{2}}(a)^{-1/2}\delta_{L_{2}\cap s_{2}Qs_{2}^{-1}}^{L_{2}}(l(xm))^{1/2}(e'_{2}(x),(s_{2}\tau)(l(xa))(e_{2}(k(xa))).$$
La fonction $f_{1}(a',\lambda)$ s'exprime de fa\c{c}on analogue: $L_{2}$ et $a$ doivent \^etre remplac\'es par $\theta(L_{2})$ et $\theta(a')$, et il y a un signe $-$ dans l'exponentielle car $a'$ intervient dans le premier terme du produit scalaire et non dans le second. Alors
$$\Psi_{N,Y}=\int_{K}\int_{K}\int_{A_{d}(F)^+}\int_{\theta(K_{2})}\int_{K_{2}}f'_{1}(a',y)f'_{2}(a,x)f_{3}(k_{1},k_{2},k_{3},k_{4},k')f_{5}(xa,y\theta(a'))$$
$$\kappa_{N}(k_{1}ak_{2})D(a)S_{Q_{1},Q_{2}}(a)dx\,dy\,da\,dk_{1}\,dk_{2},$$
 o\`u
$$f_{5}(xa,y\theta(a'))=\int_{i{\cal A}_{L,F}^*}f_{4}(\lambda)exp(-(s_{1}\lambda)(H_{\theta(L_{2})\cap s_{1}Qs_{1}^{-1}}(y\theta(a')))+(s_{2}\lambda)(H_{L_{2}\cap s_{2}Qs_{2}^{-1}}(xa)))d\lambda$$
$$=\int_{i{\cal A}_{L,F}^*}f_{4}(\lambda)exp(\lambda(\zeta(xa,y\theta(a'))))d\lambda,$$
o\`u
$$\zeta(xa,y\theta(a'))=s_{2}^{-1}(H_{L_{2}\cap s_{2}Qs_{2}^{-1}}(xa))-s_{1}^{-1}(H_{\theta(L_{2})\cap s_{1}Qs_{1}^{-1}}(y\theta(a'))).$$
Montrons que

(13)  si $\epsilon$ est assez petit et $ c_{1}$ est assez grand, on a
$$\vert \zeta(xa,y\theta(a'))\vert >> N^{\eta/2}+\vert H_{G}(a)\vert $$
pour tous $x$, $y$, $a$, $k_{1}$, $k_{2}$ intervenant dans $\Psi_{N,Y}$ avec $S_{Q_{1},Q_{2}}(a)\not=0$.

On pose  $\zeta=H_{A_{d}}(a)$, $\zeta'=H_{A_{d}}(a')$. On a $xa=aa^{-1}xa$, d'o\`u
$$H_{L_{2}\cap s_{2}Qs_{2}^{-1}}(xa)=H_{L_{2}\cap s_{2}Qs_{2}^{-1}}(a)+H_{L_{2}\cap s_{2}Qs_{2}^{-1}}(a^{-1}xa).$$
Parce que $x\in K_{2}\subset L_{2}(F)$, on a une majoration 
$$\sigma(a^{-1}xa)<<sup\{\alpha(\zeta); \alpha\in \Delta^{L_{2}}\}.$$  La condition $\sigma_{A_{d}}^{Q_{2}}(\zeta,\epsilon{\cal Y})=1$ entra\^{\i}ne $\alpha(a)<<\epsilon\vert Y\vert $ pour $\alpha\in \Delta^{L_{2}}$. D'o\`u $\sigma(a^{-1}xa)<< \epsilon\vert Y\vert $, puis $\vert H_{L_{2}\cap s_{2}Qs_{2}^{-1}}(a^{-1}xa)\vert <<\epsilon\vert Y\vert $. Pour la m\^eme raison, on a
$$\vert H_{L_{2}\cap s_{2}Qs_{2}^{-1}}(a)-H_{L_{2}}(a)\vert =\vert H_{L_{2}\cap s_{2}Qs_{2}^{-1}}(a)-\zeta_{L_{2}}\vert <<\epsilon\vert Y\vert ,$$
d'o\`u 
$$\vert H_{L_{2}\cap s_{2}Qs_{2}^{-1}}(xa)-\zeta_{L_{2}}\vert <<\epsilon\vert Y\vert .$$
Rappelons que $a'=k_{3}^{-1}\theta(l')^{-1}ak_{4}^{-1}$. Les termes $k_{3}$, $k_{4}$ appartiennent \`a $K_{2}$, $\theta(l')$ appartient \`a $L_{2}(F)$ et v\'erifie une majoration $\sigma(\theta(l'))<<Clog(N)$. On en d\'eduit une majoration $\alpha(\zeta')<<Clog(N)+\epsilon\vert Y\vert $ pour tout $\alpha\in \Delta^{L_{2}}$.  On d\'emontre alors comme ci-dessus
$$\vert H_{\theta(L_{2})\cap s_{1}Qs_{1}^{-1}}(y\theta(a'))-H_{\theta(L_{2})}(\theta(a'))\vert <<Clog(N)+\epsilon\vert Y\vert .$$
De la d\'efinition de $a'$ rappel\'ee ci-dessus r\'esulte que
$$H_{\theta(L_{2})}(\theta(a'))-H_{\theta(L_{2})}(\theta(a))\vert <<Clog(N).$$
Puisque $H_{\theta(L_{2})}(\theta(a))=\theta(\zeta_{L_{2}})$, on obtient
$$\vert H_{\theta(L_{2})\cap s_{1}Qs_{1}^{-1}}(y\theta(a'))-\theta(\zeta_{L_{2}})\vert <<Clog(N)+\epsilon\vert Y\vert .$$
De ces calculs r\'esulte la majoration
$$ (14) \qquad \vert \zeta(xa,y\theta(a'))-s_{2}^{-1}\zeta_{L_{2}}+s_{1}^{-1}\theta(\zeta_{L_{2}})\vert <<Clog(N)+\epsilon\vert Y\vert .$$
Pour tout sous-groupe $G'$ de $G$ stable par conjugaison par $A_{d}$, notons $\Sigma^{G'}$ l'ensemble des racines de $A_{d}$ dans $\mathfrak{g}'$. Montrons que

(15) l'un au moins des ensembles $s_{2}^{-1}(\Sigma^{U_{1}})\cap s_{1}^{-1}\theta(\Sigma^{\bar{Q}_{2}})$ ou $s_{2}^{-1}(\Sigma^{\bar{Q}_{2}})\cap s_{1}^{-1}\theta(\Sigma^{U_{1}})$ est non vide.

Supposons que le premier ensemble est vide. Alors $s_{2}^{-1}(\Sigma^{U_{1}})\subset s_{1}^{-1}\theta(\Sigma^{U_{2}})$, ou encore $sU_{1}s^{-1}\subset \theta(U_{2})$. Pour deux sous-groupes paraboliques $P'=M'U'$ et $P''=M''U''$ de $G$, les inclusions $U'\subset U''$ et $P''\subset P'$ sont \'equivalentes. Donc $\theta(Q_{2})\subset sQ_{1}s^{-1}$. Le sous-groupe parabolique $\theta(Q_{2})$ est standard, donc aussi $sQ_{1}s^{-1}$. Mais $Q_{1}$ est lui-aussi standard. Deux sous-groupes paraboliques standard n'\'etant conjugu\'es que s'ils sont \'egaux, on a $Q_{1}=sQ_{1}s^{-1}$, donc   $s\in W^{L_{1}}$.  L'inclusion $\theta(Q_{2})\subset sQ_{1}s^{-1}$ devient $\theta(Q_{2})\subset Q_{1}$. Supposons que le second ensemble de (15) soit vide. Un calcul similaire montre que $s\in W^{\theta(L_{1})}$. Posons $Q'=L'U'=Q_{1}\cap \theta(Q_{1})$. On a alors $Q_{2}\subset Q'$ et $s\in W^{L'}$. Cela contredit l'hypoth\`ese sur $s$. D'o\`u (15).

Supposons par exemple le premier ensemble ci-dessus non vide. Soit $\alpha$ un de ses \'el\'ements.
On a
$$\alpha(s_{2}^{-1}\zeta_{L_{2}})=(s_{2}\alpha)(\zeta)-(s_{2}\alpha)(\zeta^{L_{2}}),$$
o\`u $\zeta^{L_{2}}$ est la projection orthogonale de $\zeta$ sur ${\cal A}^{L_{2}}_{A_{d}}$. Parce que $s_{2}\alpha\in \Sigma^{U_{1}}$ et $\sigma_{A_{d}}^{Q_{1}}(\zeta,{\cal Y})\tau_{Q_{1}}(\zeta-Y_{Q_{1}})=1$, on a $ (s_{2}\alpha)(\zeta) >>\vert Y\vert $. Parce que $\sigma_{A_{d}}^{Q_{2}}(\zeta,\epsilon{\cal Y})=1$, on a $\vert (s_{2}\alpha)(\zeta^{L_{2}})\vert <<\epsilon \vert Y\vert $. Il existe donc $c_{4}>0$ tel que
$$\alpha(s_{2}^{-1}\zeta_{L_{2}})>>(1-c_{4}\epsilon)\vert Y\vert .$$
On a $\alpha(s_{1}^{-1}\theta(\zeta_{L_{2}}))=(\theta s_{1}\alpha)(\zeta_{L_{2}})$. L'\'el\'ement  $\theta s_{1}\alpha$ appartient \`a $\Sigma^{\bar{Q}_{2}}$. Sa projection sur ${\cal A}_{L_{2}}^*$ est une combinaison lin\'eaire \`a coefficients n\'egatifs ou nuls d'\'el\'ements de $\Delta-\Delta^{L_{2}}$. Donc
$$\alpha(s_{1}^{-1}\theta(\zeta_{L_{2}}))\leq 0.$$
D'o\`u
$$\alpha(s_{2}^{-1}\zeta_{L_{2}}-s_{1}^{-1}\theta(\zeta_{L_{2}}))>>(1-c_{4}\epsilon)\vert Y\vert ,$$
a fortiori
$$\vert s_{2}^{-1}\zeta_{L_{2}}^G-s_{1}^{-1}\theta(\zeta_{L_{2}}^G)\vert >>(1-c_{4}\epsilon)\vert Y\vert .$$
On a aussi
$$\vert  s_{2}^{-1}\zeta_{G}-s_{1}^{-1}\theta(\zeta_{G})\vert=\vert \zeta_{G}-\theta(\zeta_{G})\vert  >> \vert \zeta_{G}\vert $$
D'o\`u
$$\vert  s_{2}^{-1}\zeta_{L_{2}}-s_{1}^{-1}\theta(\zeta_{L_{2}})\vert >>(1-c_{4}\epsilon)\vert Y\vert +\vert \zeta_{G} \vert .$$
Un raisonnement analogue, avec la m\^eme conclusion, vaut si le deuxi\`eme ensemble de (15) est non vide. La minoration que l'on vient d'obtenir, jointe \`a (14) et \`a l'hypoth\`ese $\vert Y\vert \geq c_{1}N^{\eta}$, d\'emontre (13).

Le terme $f_{5}(xa,y\theta(a'))$ est la valeur au point $\zeta(xa,y\theta(a'))$ de la transform\'ee de Fourier d'une fonction $C^{\infty}$ sur $i{\cal A}_{L,F}^*$. Cette transform\'ee de Fourier  est \`a d\'ecroissance rapide. Gr\^ace \`a (13), pour tout entier $b$, on a une majoration
$$\vert f_{5}(xa,y\theta(a'))<<N^{-b}(1+\vert H_{G}(a)\vert )^{-b}$$
pourvu que $\epsilon$ soit assez petit et $c_{1}$ assez grand. D'o\`u
$$\Psi_{N,Y}<<N^{-b}\int_{K}\int_{K}\int_{A_{d}(F)^+}\int_{\theta(K_{2})}\int_{K_{2}}\vert f'_{1}(a',y)f'_{2}(a,x)f_{3}(k_{1},k_{2},k_{3},k_{4},k')\vert $$
$$(1+\vert H_{G}(a)\vert )^{-b}\kappa_{N}(k_{1}ak_{2})D(a)S_{Q_{1},Q_{2}}(a)dx\,dy\,da\,dk_{1}\,dk_{2}.$$
On n'a aucun mal \`a montrer qu'il existe $R_{9}$ tel que l'int\'egrale soit essentiellement born\'ee par $N^{R_{9}}$. Puisque $b$ est quelconque, on obtient la premi\`ere majoration de (12). La seconde s'obtient de la m\^eme fa\c{c}on. Le seul changement est que, dans l'int\'egrale ci-dessus, $\kappa_{N}(k_{1}ak_{2})$ est remplac\'e par $\tilde{u}(a,Y)$. L'int\'egrale est alors major\'ee par $N^{R_{9}}\vert Y\vert ^{R_{9}}$ et il faut utiliser la majoration $\vert Y\vert <<c_{2}N$ pour conclure.

Supposons maintenant que l'hypoth\`ese impos\'ee plus haut \`a $s$ ne soit pas v\'erifi\'ee. Dans ce cas, on va montrer que $\Phi_{N,Y,Q_{1},Q_{2},s_{1},s_{2}}(g')=\Phi_{Y,Q_{1},Q_{2},s_{1},s_{2}}(g')=0$.
 D'apr\`es la d\'efinition de ces termes, il suffit de prouver
 
 (16) pour $\lambda\in i{\cal A}_{L,F}^*$, $a\in A_{d}(F)^+$ et $k_{1},k_{2}\in K$, on a
 $$\sum_{e\in {\cal B}_{{\cal O}}^{K_{f'}}}\Phi_{s_{1},s_{2}}(e,k_{1}ak_{2},g',\lambda)=0.$$
 
 On note $X(\lambda)$ la somme ci-dessus. C'est la restriction \`a $i{\cal A}_{L,F}^*$ d'une fonction m\'eromorphe sur $({\cal A}_{L}\otimes_{{\mathbb R}}{\mathbb C})/i{\cal A}_{L,F}^{\vee}$. Il suffit  de d\'emontrer qu'elle est nulle en un point $\lambda$ g\'en\'eral. On peut donc supposer que tous les op\'erateurs d'entrelacement qui appara\^{\i}tront ci-dessous n'ont ni z\'ero ni p\^ole en $\lambda$. On a une \'egalit\'e
 $$\Phi_{s_{1},s_{2}}(e,k_{1}ak_{2},g',\lambda)=$$
 $$(J_{\theta(\underline{Q})_{2,s_{1}}\vert s_{1}Qs_{1}^{-1}}((s_{1}\tau)_{s_{1}\lambda})\circ\gamma(s_{1})\circ\pi_{\lambda}(\theta(a'k_{4}k_{2}))e,e_{1})^{\theta(L_{2})}$$
 $$(e_{2},J_{\underline{Q}_{2,s_{2}}\vert s_{2}Qs_{2}^{-1}}((s_{2}\tau)_{s_{2}\lambda})\circ\gamma(s_{2})\circ\pi_{\lambda}(ak_{2})\pi_{\lambda}(f')e)^{L_{2}},$$
 pour des \'el\'ements $e_{1}\in {\cal K}_{Q_{2,s_{1}},s_{1}\tau}^G$ et $e_{2}\in {\cal K}_{Q_{2,s_{2}},s_{2}\tau}^G$. On a $\pi_{\lambda}(k_{2})\pi_{\lambda}(f')=\pi_{\lambda}(f'')\pi_{\lambda}(\theta(k_{2}))$, o\`u $f''$ est d\'efini par $f''(g)=f(k_{2}^{-1}g\theta(k_{2}))$. Posons ${\cal B}_{\natural}^{K_{f'}}=\{\pi_{\lambda}(\theta(k_{2}))e; e\in {\cal B}_{{\cal O}}^{K_{f'}}\}$. C'est encore une base orthonorm\'ee de $({\cal K}_{Q,\tau}^G)^{K_{f'}}$ (rappelons que $K_{f'}$ est distingu\'e dans $K$) et on a
 $$X(\lambda)=\sum_{e\in {\cal B}_{\natural}^{K_{f'}}}(J_{\theta(\underline{Q})_{2,s_{1}}\vert s_{1}Qs_{1}^{-1}}((s_{1}\tau)_{s_{1}\lambda})\circ\gamma(s_{1})\circ\pi_{\lambda}(\theta(a'k_{4}))e,e_{1})^{\theta(L_{2})}$$
 $$(e_{2},J_{\underline{Q}_{2,s_{2}}\vert s_{2}Qs_{2}^{-1}}((s_{2}\tau)_{s_{2}\lambda})\circ\gamma(s_{2})\circ\pi_{\lambda}(a)\pi_{\lambda}(f')e)^{L_{2}}.$$
 Il existe une fonction $j_{1}(\lambda)$ qui est m\'eromorphe (au m\^eme sens que ci-dessus) et telle que
 $$J_{\theta(\underline{Q})_{2,s_{1}}\vert s_{1}Qs_{1}^{-1}}((s_{1}\tau)_{s_{1}\lambda})\circ\gamma(s_{1})=$$
 $$j_{1}(\lambda)J_{\theta(\underline{Q})_{2,s_{1}}\vert s\underline{Q}_{2,s_{2}}s^{-1}}((s_{1}\tau)_{s_{1}\lambda})\circ\gamma(s)\circ J_{\underline{Q}_{2,s_{2}}\vert s_{2}Qs_{2}^{-1}}((s_{2}\tau)_{s_{2}\lambda})\circ\gamma(s_{2}).$$
 L'ensemble 
 $$\{J_{\underline{Q}_{2,s_{2}}\vert s_{2}Qs_{2}^{-1}}((s_{2}\tau)_{s_{2}\lambda})\circ\gamma(s_{2})e; e\in {\cal B}_{\natural}^{K_{f'}}\}$$
 est une base de $({\cal K}_{\underline{Q}_{2,s_{2}},s_{2}\tau}^G)^{K_{f'}}$. Les propri\'et\'es d'adjonction et de composition des op\'erateurs d'entrelacement entra\^{\i}nent qu'elle est orthogonale et que tous ses \'el\'ements ont la m\^eme norme. Notons $j_{2}(\lambda)$ cette norme et divisons tout \'el\'ement de la base par $\sqrt{j_{2}(\lambda)}$. On obtient une base orthonorm\'ee de $({\cal K}_{\underline{Q}_{2,s_{2}},s_{2}\tau}^G)^{K_{f'}}$ que l'on note ${\cal B}_{\sharp}^{K_{f'}}$. On a l'\'egalit\'e
 $$X(\lambda)=j_{1}(\lambda)j_{2}(\lambda)\sum_{e\in {\cal B}_{\sharp}^{K_{f'}}}(J_{\theta(\underline{Q})_{2,s_{1}}\vert s\underline{Q}_{2,s_{2}}s^{-1}}((s_{1}\tau)_{s_{1}\lambda})\circ\gamma(s)\circ Ind_{\underline{Q}_{2,s_{2}}}^G((s_{2}\tau)_{s_{2}\lambda},\theta(a'k_{4}))e,e_{1})^{\theta(L_{2})}$$
 $$(e_{2},Ind_{\underline{Q}_{2,s_{2}}}^G((s_{2}\tau)_{s_{2}\lambda},a)Ind_{\underline{Q}_{2,s_{2}}}^G((s_{2 }\tau)_{s_{2}\lambda},f'')e)^{L_{2}}.$$ 
 Maintenant, le groupe $Q$ n'appara\^{\i}t plus. Pour simplifier les notations, on peut supposer $s_{2}=1$ et $Q=\underline{Q}_{2,s_{2}}$. Alors $s=s_{1}$, $Q\subset \bar{Q}_{2}$ et l'expression pr\'ec\'edente se simplifie en
 $$X(\lambda)=j_{1}(\lambda)j_{2}(\lambda)\sum_{e\in {\cal B}_{\sharp}^{K_{f'}}}(J_{\theta(\underline{Q})_{2,s}\vert sQs^{-1}}((s\tau)_{s\lambda})\circ\gamma(s)\circ\pi_{\lambda}(\theta(a'k_{4}))e,e_{1})^{\theta(L_{2})}$$
 $$(e_{2},\pi_{\lambda}(a)\pi_{\lambda}(f'')e)^{L_{2}}.$$
 Puisque l'hypoth\`ese indiqu\'ee plus haut sur $s$ n'est pas v\'erifi\'ee, on peut fixer un sous-groupe parabolique $Q'=L'U'$ de $G$ tel que $Q_{2}\subset Q'\subsetneq G$, $\theta(Q')=Q'$ et $s\in W^{L'}$.
 Les deux sous-groupes paraboliques $\theta(\underline{Q})_{2,s}$ et $sQs^{-1}$ sont inclus dans $\bar{Q}'$. Introduisons   la repr\'esentation $\pi'=Ind_{L'\cap Q}^{L'}(\tau_{\lambda})$, que l'on r\'ealise dans ${\cal K}_{L'\cap Q, \tau}^{L'}$.  On peut identifier ${\cal K}_{Q,\tau}^G$ \`a ${\cal K}_{\bar{Q}',\pi'}^G$. Modulo cette identification, on d\'efinit une forme sesquilin\'eaire $B$ sur ${\cal K}_{\bar{Q}',\pi'}^G\times {\cal K}_{\bar{Q}',\pi'}^G$ par
 $$B(e',e)=(J_{\theta(\underline{Q})_{2,s}\vert sQs^{-1}}((s\tau)_{s\lambda})\circ\gamma(s)\circ\pi_{\lambda}(\theta(a'k_{4}))e,e_{1})^{\theta(L_{2})}(e_{2},\pi_{\lambda}(a)\pi_{\lambda}(f'')e)^{L_{2}}.$$
 Alors 
 $$X(\lambda)=trace_{B}(\pi'(f'')).$$
 Pour prouver que cette expression est nulle, on va utiliser le lemme 1.11, ou plut\^ot sa variante "unitaire" \'evidente. Posons $\tilde{Q}'=Q'\boldsymbol{\theta}_{d}$, $\tilde{L}'=L'\boldsymbol{\theta}_{d}$. Puisque $\theta(Q')=Q'$, $\tilde{Q}'$ est bien un sous-groupe parabolique tordu de $\tilde{G}$, de L\'evi tordu $\tilde{L}'$. Il est diff\'erent de $\tilde{G}$. La fonction $f''$ co\"{\i}ncide avec $\tilde{f}''_{\boldsymbol{\theta}_{d}}$, o\`u $\tilde{f}''$ est d\'efinie par $\tilde{f}''(\tilde{x})=\tilde{f}(k_{2}^{-1}\tilde{x}k_{2})$. Cette fonction $\tilde{f}''$ est tr\`es cuspidale. Il reste \`a v\'erifier que $B$ v\'erifie l'hypoth\`ese $(H)_{\boldsymbol{\theta}_{d}}$ de 1.11. Il suffit pour cela de montrer que, pour $e,e'\in {\cal K}_{\bar{Q}',\pi'}^G$,   $B(e',e)$ ne d\'epend que de la restriction de $e$ et $e'$ \`a $K\cap \bar{Q}'(F)$.  Introduisons  
 l'espace   ${\cal K}_{L'\cap\theta(\underline{Q})_{2,s}}^{L'}$. On dispose de l'op\'erateur
 $$J^{L'}_{L'\cap \theta(\underline{Q})_{2,s}\vert L'\cap sQs^{-1}}((s\tau)_{s\lambda})\circ\gamma(s):{\cal K}_{L'\cap Q,\tau}^{L'}\to {\cal K}_{L'\cap \theta(\underline{Q})_{2,s},s\tau}^{L'}.$$
 On v\'erifie la formule
 $$B(e',e)=\delta_{Q'}(\theta(a'k_{4}))^{-1/2}\delta_{Q'}(a)^{-1/2}$$
 $$\int_{K_{2}}((J^{L'}_{L'\cap \theta(\underline{Q})_{2,s}\vert L'\cap sQs^{-1}}((s\tau)_{s\lambda})\circ\gamma(s)\circ\pi'(\theta(a'k_{4})(e'(1)))(x),e_{1}(x))dx$$
 $$\int_{K_{2}}(e_{2}(y),(\pi'(a)(e(1)))(y))dy.$$
 Donc $B(e',e)$ ne d\'epend que de $e(1)$ et $e'(1)$. On peut donc appliquer le lemme 1.11 (ou plut\^ot sa variante) et on conclut que $X(\lambda)=0$. Cela prouve (16) et la proposition. $\square$
 
 {\bf Erratum.} La preuve de la proposition 6.6 de [W3] contient une erreur: la minoration (11) de cette preuve n'entra\^{\i}ne pas la majoration $\vert f_{5}(xm',ym)\vert <<N^{-D}$ mais seulement  $\vert f_{5}(xm',ym)\vert <<log(N)^{-D}$. On doit remplacer (11) par une minoration $\vert \zeta(xm',ym)\vert >>N^{\epsilon}$ pour un $\epsilon>0$. Pour l'obtenir, il suffit de remplacer la condition $c_{1}log(N)<<\alpha(Y)$ pour tout $\alpha$ de l'\'enonc\'e de la proposition 6.6 par une condition $c_{1}N^{\eta}<<\alpha(Y)$ o\`u $\eta$ est un r\'eel fix\'e tel que $0<\eta<1$.

 \bigskip
 
 \subsection{Un r\'esultat extrait de la formule des traces locale tordue}
 
 Pour $i=1,...,4$, soit $e_{i}\in {\cal K}_{Q,\tau}^G$. Soient  $\varphi\in C_{c}^{\infty}(i{\cal A}_{L,F}^*)$ et $Y\in {\cal D}$. Posons 
  $$\Phi_{Y}((e_{i})_{i=1,...,4},\varphi)=\int_{G(F)}\int_{i{\cal A}_{L,F}^*}m(\tau_{\lambda})\varphi(\lambda) (\pi_{\lambda}(\theta(g))e_{1},e_{2})(e_{3},\pi_{\lambda}(g)e_{4})\tilde{u}(g,Y)dg\,d\lambda.$$
 
 Soit $\tilde{L}'\in {\cal L}^{\tilde{G}}$ tel que $L'$ contient $L$. Posons
 $$W^{\tilde{L}'}=\{\tilde{t}\in \tilde{L}'(F); \theta_{\tilde{t}}(A_{d})=A_{d}\}/A_{d}(F),$$
 $$W^{\tilde{L}'}(L)=\{\tilde{t}\in W^{\tilde{L}'}; \theta_{\tilde{t}}(L)=L\}/W^L.$$
 On repr\'esente tout \'el\'ement $\tilde{t}\in W^{\tilde{L}'}$ sous la forme $\tilde{t}=t\boldsymbol{\theta}_{d}$, o\`u $t\in K\cap Norm_{G(F)}(A_{d})$. Pour tout $\tilde{t}\in W^{\tilde{L}'}(L)$, l'automorphisme $\theta_{\tilde{t}}$ agit naturellement sur ${\cal A}_{L}$ et son ensemble de points fixes contient ${\cal A}_{\tilde{L}'}$. Notons $W^{\tilde{L}'}(L)_{reg}$ le sous-ensemble des $\tilde{t}\in W^{\tilde{L}'}(L)$ tels que cet ensemble de points fixes soit \'egal \`a ${\cal A}_{\tilde{L}'}$. Consid\'erons un tel \'el\'ement $\tilde{t}$. On d\'efinit la repr\'esentation $\theta_{\tilde{t}}(\tau) $ de $L(F)$ par $\theta_{\tilde{t}}(\tau)(l)=\tau(\theta_{\tilde{t}}^{-1}(l))$.  Notons $\Lambda_{{\cal O}}(\tilde{t})$ l'ensemble des $\lambda\in i{\cal A}_{L}^*$ tels que $\theta_{\tilde{t}}(\tau_{\lambda})\simeq \tau_{\lambda}$. Il est stable par translations par $i{\cal A}_{L,F}^{\vee}+i{\cal A}_{\tilde{L}'}$. Soit $\lambda$ un \'el\'ement de cet ensemble. On d\'efinit un signe $\epsilon_{\tau_{\lambda}}(\tilde{t})$ de la fa\c{c}on suivante. Introduisons l'ensemble $\Sigma^{L'}_{L}$ des racines de $A_{L}$ dans $\mathfrak{l}'$. On le munit du sous-ensemble "positif" form\'e des racines dans $\mathfrak{l}'\cap \mathfrak{q}$. A chaque racine, on peut associer une mesure de Plancherel $m_{\alpha}(\tau_{\lambda})$. Notons $n^{L'}(\tilde{t})$ le nombre des racines $\alpha\in \Sigma^{L'}_{L}$ telles que $\alpha>0$, $\theta_{\tilde{t}}(\alpha)<0$ et $m_{\alpha}(\tau_{\lambda})=0$. On pose $\epsilon_{\tau_{\lambda}}(\tilde{t})=(-1)^{n^{L'}(\tilde{t})}$. Introduisons l'op\'erateur d'entrelacement normalis\'e 
$$R_{Q\vert \theta_{\tilde{t}}(Q)}(\tau_{\lambda}): {\cal K}_{\theta_{\tilde{t}}(Q),\tau}^G\to {\cal K}_{Q,\tau}^G.$$
Ici, la normalisation choisie n'a pas d'importance, pouvu que cet op\'erateur pr\'eserve les produits scalaires. Fixons un automorphisme $A(\tau_{\lambda})$ de l'espace de $\tau$ tel que 
$$ \tau_{\lambda}(l)\circ A(\tau_{\lambda})=A(\tau_{\lambda})\circ\theta_{\tilde{t}}(\tau_{\lambda})(l) $$
pour  tout $l\in L(F)$.  D\'efinissons un op\'erateur
$$A(\tau_{\lambda},\tilde{t}):{\cal K}_{Q,\tau}^G\to {\cal K}_{\theta_{\tilde{t}}(Q),\tau}^G$$
par
$$ (A(\tau_{\lambda},\tilde{t})e)(k)=A(\tau_{\lambda})\circ e(\theta(t^{-1}k))$$
pour tous $e\in {\cal K}_{Q,\tau}^G$, $k\in K$. Pour $\tilde{Q}'=\tilde{L}'U'\in {\cal P}(\tilde{L}')$, posons $Q(Q')=(L'\cap Q)U'$. C'est un \'el\'ement de ${\cal P}(L)$. On d\'efinit les deux fonctions $\mu\mapsto j^{1,4}_{\tilde{Q}'}(\tilde{t},\lambda,\mu)$ et $\mu\mapsto j^{2,3}_{\tilde{Q}'}(\tilde{t},\lambda,\mu)$ sur $i{\cal A}_{\tilde{L}'}^*$ par
$$j^{1,4}_{\tilde{Q}'}(\tilde{t},\lambda,\mu)=(R_{Q\vert \theta_{\tilde{t}}(Q)}(\tau_{\lambda})A(\tau_{\lambda},\tilde{t})e_{1},J_{Q(\bar{Q}')\vert Q}(\tau_{\lambda})^{-1}J_{Q(\bar{Q}')\vert Q}(\tau_{\lambda+\mu}) e_{4}),$$
$$j^{2,3}_{\tilde{Q}'}(\tilde{t},\lambda,\mu)=(J_{Q(Q')\vert Q}(\tau_{\lambda})^{-1}J_{Q(Q')\vert Q}(\tau_{\lambda+\mu}) e_{2},R_{Q\vert \theta_{\tilde{t}}(Q)}(\tau_{\lambda})A(\tau_{\lambda},\tilde{t})e_{3}).$$
 Les familles $(j^{1,4}_{\tilde{Q}'}(\tilde{t},\lambda))_{\tilde{Q}'\in {\cal P}(\tilde{L}')}$ et $(j^{2,3}_{\tilde{Q}'}(\tilde{t},\lambda))_{\tilde{Q}'\in {\cal P}(\tilde{L}')}$ sont des $(\tilde{G},\tilde{L}')$-familles.  Notons $({\cal J}_{\tilde{Q}'}(\tilde{t},\lambda))_{\tilde{Q}'\in {\cal P}(\tilde{L}')}$ la famille produit:
$${\cal J}_{\tilde{Q}'}(\tilde{t},\lambda,\mu)=j^{1,4}_{\tilde{Q}'}(\tilde{t},\lambda,\mu)j^{2,3}_{\tilde{Q}'}(\tilde{t},\lambda,\mu).$$
Comme \`a toute $(\tilde{G},\tilde{L}')$-famille, on peut lui associer la fonction ${\cal J}_{\tilde{L}'}(\tilde{t},\lambda,\mu)$, puis le nombre ${\cal J}_{\tilde{L}'}(\tilde{t},\lambda)={\cal J}_{\tilde{L}'}(\tilde{t},\lambda,0)$. Admettons que ce nombre soit fonction $C^{\infty}$ de $\lambda$. Posons
$$\Phi((e_{i})_{i=1,...,4},\varphi)=\sum_{\tilde{L}'\in {\cal L}^{\tilde{G}}; L\subset L'}\sum_{\tilde{t}\in W^{\tilde{L}'}(L)_{reg}}\vert det(1-\theta_{\tilde{t}})_{\vert {\cal A}_{L}/{\cal A}_{\tilde{L}'}}\vert ^{-1}$$
$$\sum_{\lambda\in \Lambda_{{\cal O}}(\tilde{t})/(i{\cal A}^{\vee}_{L,F}+i{\cal A}_{\tilde{L}'}^*)}\epsilon_{\tau_{\lambda}}(\tilde{t})\int_{i{\cal A}_{\tilde{L}',F}^*}\varphi(\lambda+\chi){\cal J}_{\tilde{L}'}(\tilde{t},\lambda+\chi)d\chi.$$
Pour tout entier $k\in {\mathbb N}$, fixons une base ${\cal X}_{k}$ de l'espace des op\'erateurs diff\'erentiels sur $i{\cal A}_{L}^*$, \`a coefficients constants et d'ordre $\leq k$. Posons
$$\vert \varphi\vert _{k}=sup\{\vert X\varphi(\lambda)\vert ; \lambda\in i{\cal A}_{L}^*, X\in {\cal X}_{k}\}.$$
Introduisons l'espace $PolExp$ des fonctions $ \Phi$ sur ${\cal A}_{\tilde{A}_{d},F}\otimes_{{\mathbb Z}}{\mathbb Q}$ v\'erifiant la condition suivante. Soit ${\cal R}\subset {\cal A}_{\tilde{A}_{d},F}\otimes_{{\mathbb Z}}{\mathbb Q}$ un r\'eseau. Alors il existe un ensemble fini $\Xi_{{\cal R}}\subset i{\cal A}_{\tilde{A}_{d}}^*/i{\cal R}^{\vee}$ (o\`u ${\cal R}^{\vee}$ est l'ensemble des $\lambda\in {\cal A}_{\tilde{A}_{d}}^*$ tel que $\lambda(Y)\in 2\pi {\mathbb Z}$ pour tout $Y\in {\cal R}$) et, pour tout $\xi\in \Xi_{{\cal R}}$, il existe un polyn\^ome $p_{{\cal R},\xi}$ sur ${\cal A}_{\tilde{A}_{d}}$ de sorte que
$$\Phi(Y)=\sum_{\xi\in \Xi_{{\cal R}}}e^{\xi(Y)}p_{{\cal R},\xi}(Y)$$
pour tout $Y\in {\cal R}$.
Un tel d\'eveloppement est unique. On pose $c_{{\cal R},0}(\Phi)=p_{{\cal R},0}(0)$.

\ass{Proposition}{(i) La fonction $\lambda\mapsto {\cal J}_{\tilde{L}'}(\tilde{t},\lambda)$ est $C^{\infty}$.

(ii) Il existe une unique fonction $\Phi_{(e_{i})_{i=1,...,4},\varphi}\in PolExp$ de sorte que

(a) pour tout r\'eseau ${\cal R}\subset {\cal A}_{\tilde{A}_{d},F}\otimes_{{\mathbb Z}}{\mathbb Q}$, 
$$lim_{k\to \infty}p_{\frac{1}{k}{\cal R},0}(\Phi_{(e_{i})_{i=1,...,4},\varphi})=\Phi((e_{i})_{i=1,...,4},\varphi);$$

(b)  pour tout $R>0$, il existe un entier $k$ et    une constante $c>0$ d\'ependant de $R$ et des $e_{i}$ mais pas de $\varphi$, de sorte que
$$\vert \Phi_{Y}((e_{i})_{i=1,...,4},\varphi)-\Phi_{(e_{i})_{i=1,...,4},\varphi}(Y)\vert\leq c\vert \varphi\vert _{k}\vert Y\vert ^{-R}$$
 pour tout $Y\in {\cal D}$.}

Preuve.  L'existence d'un \'el\'ement de $\Phi_{(e_{i})_{i=1,...,4},\varphi}\in PolExp$ v\'erifiant la propri\'et\'e (ii)(b) est le lemme 3.19 de [W6], pr\'ecis\'e par la relation 3.19(1). Un tel \'el\'ement est forc\'ement unique. Le corollaire 3.24 de [W6] calcule $lim_{k\to \infty}p_{\frac{1}{k}{\cal R},0}(\Phi_{(e_{i})_{i=1,...,4},\varphi})$, mais exprime cette limite sous une forme diff\'erente de celle que l'on veut.  Dans la preuve de la proposition 3.24 de [W6], on montre qu'apr\`es sommation sur un certain ensemble de quadruplets $(e_{i})_{i=1,...,4}$, ces deux formes sont \'equivalentes.  En fait, la m\^eme preuve vaut pour un quadruplet fix\'e.  Cette preuve d\'emontre en m\^eme temps l'assertion (i). $\square$
  
 \bigskip
 
 \subsection{Utilisation de la formule des traces locale tordue}
  Revenons aux donn\'ees de 6.5. Pour $g'\in G(F)$ et $Y\in {\cal D}$, on a d\'efini $\Phi_{Y}(g')$ dans ce paragraphe.   Soit $\tilde{L}'\in {\cal L}^{\tilde{G}}$ tel que $L\subset L'$ et soient $\tilde{t}\in W^{\tilde{L}'}(L)_{reg}$ et $\lambda\in \Lambda_{{\cal O}}(\tilde{t})$. On d\'efinit les deux $(\tilde{G},\tilde{L}')$-familles suivantes. Pour $\tilde{Q}'\in {\cal P}(\tilde{L}')$ et $\mu\in i{\cal A}_{\tilde{L}'}^*$, on pose
 $$j_{\tilde{Q}'}(\tilde{t},\lambda,\mu) =\sum_{e\in {\cal B}_{{\cal O}}^{K_{f'}}}(R_{Q\vert \theta_{\tilde{t}}(Q)}(\tau_{\lambda})A(\tau_{\lambda},\tilde{t})e,J_{Q(\bar{Q}')\vert Q}(\tau_{\lambda})^{-1}J_{Q(\bar{Q}')\vert Q}(\tau_{\lambda+\mu}) \pi_{\lambda}(f')e),$$
 $$d_{\tilde{Q}'}(\tilde{t},\lambda,g',\mu)=(J_{Q(Q')\vert Q}(\tau_{\lambda})^{-1}J_{Q(Q')\vert Q}(\tau_{\lambda+\mu}) e'',R_{Q\vert \theta_{\tilde{t}}(Q)}(\tau_{\lambda})A(\tau_{\lambda},\tilde{t})\pi_{\lambda}(g')e').$$
 Posons
 $$(jd)_{\tilde{Q}'}(\tilde{t},\lambda,g',\mu)=j_{\tilde{Q}'}(\tilde{t},\lambda,\mu)d_{\tilde{Q}'}(\tilde{t},\lambda,g',\mu)c_{\tilde{Q}'}(\mu).$$
 La famille $((jd)_{\tilde{Q}'}(\tilde{t},\lambda,g'))_{\tilde{Q}'\in {\cal P}(\tilde{L}')}$ est une $(\tilde{G},\tilde{L}')$-famille   et on lui associe le nombre $(jd)_{\tilde{L}'}(\tilde{t},\lambda,g')$. Le (i) de la proposition du paragraphe pr\'ec\'edent entra\^{\i}ne que ce nombre est fonction $C^{\infty}$ de $\lambda$. On pose
 $$\Phi(g')=\sum_{\tilde{L}'\in {\cal L}^{\tilde{G}}; L\subset L'}\sum_{\tilde{t}\in W^{\tilde{L}'}(L)_{reg}}\vert det(1-\theta_{\tilde{t}})_{\vert {\cal A}_{L}/{\cal A}_{\tilde{L}'}}\vert^{-1} $$
$$\sum_{\lambda\in \Lambda_{{\cal O}}(\tilde{t})/(i{\cal A}^{\vee}_{L,F}+i{\cal A}_{\tilde{L}'}^*)}\epsilon_{\tau_{\lambda}}(\tilde{t})\int_{i{\cal A}_{\tilde{L}',F}^*}\varphi(\lambda+\chi)(jd)_{\tilde{L}'}(\tilde{t},\lambda+\chi,g')d\chi.$$

\ass{Proposition}{Pour tout r\'eel $R\geq1$,  il existe un entier $k\geq0$ tel que l'on ait une majoration
$$\vert \Phi_{Y}(g')-\Phi(g')\vert<< \sigma(g')^k\Xi^{G}(g')\vert Y\vert ^{-R }$$
pour tous $g'\in G(F)$ et $Y\in {\cal D}$.}

Preuve. Fixons un sous-groupe ouvert compact $K_{0}$ de $K$, conserv\'e par l'automorphisme $\theta$, et qui fixe $e''$. La fonction $g\mapsto \tilde{u}(g,Y)$ est biinvariante par $K$. Dans la d\'efinition de $\Phi_{Y}(g')$ donn\'ee en 6.5, on peut remplacer $g$ par $kg$ pour $k\in K_{0}$, puis int\'egrer en $k$, en divisant l'expression obtenue par $mes(K_{0})$. Cela remplace $\Phi_{Y}(g')$ par une expression analogue, o\`u $\pi_{\lambda}(g')e'$ est remplac\'e  par
$$mes(K_{0})^{-1}\int_{K_{0}}\pi_{\lambda}(kg')e'dk.$$
Fixons une base orthonorm\'ee ${\cal B}_{{\cal O}}^{K_{0}}$ de $({\cal K}_{Q,\tau}^G)^{K_{0}}$.  Le terme ci-dessus est \'egal \`a
$$\sum_{e_{0}\in {\cal B}_{{\cal O}}^{K_{0}}}(e_{0},\pi_{\lambda}(g')e')e_{0}.$$
De m\^eme, pour $e\in {\cal B}_{{\cal O}}^{K_{f'}}$, on peut remplacer $\pi_{\lambda}(f')e$ par son expression dans la base ${\cal B}_{{\cal O}}^{K_{f'}}$. On voit alors que
$$\Phi_{Y}(g')=\sum_{e,e_{4}\in {\cal B}_{{\cal O}}^{K_{f'}},e_{0}\in {\cal B}_{{\cal O}}^{K_{0}}}\Phi_{Y}(e,e_{0},e'',e_{4},\varphi_{e,e_{0},e_{4},g'}),$$
o\`u
$$\varphi_{e,e_{0},e_{4},g'}(\lambda)=\varphi(\lambda)(e_{4},\pi_{\lambda}(f')e)(e_{0},\pi_{\lambda}(g')e').$$
La proposition du paragraphe pr\'ec\'edent nous fournit une fonction 
$$\Phi_{g'}=\sum_{e,e_{4}\in {\cal B}_{{\cal O}}^{K_{f'}},e_{0}\in {\cal B}_{{\cal O}}^{K_{0}}}\Phi_{(e,e_{0},e'',e_{4}),\varphi_{e,e_{0},e_{4},g'}}.$$
Elle appartient \`a $PolExp$. Le (ii)(b) de la proposition nous dit que, pour tout r\'eel $R\geq1$, il existe un entier $k$ et une constante $c>0$ ind\'ependants de $g'$ telle que
$$\vert \Phi_{Y}(g')-\Phi_{g'}(Y)\vert \leq c\,sup_{e,e_{0},e_{4}}\vert \varphi_{e,e_{0},e_{4},g'}\vert _{k}\vert Y\vert ^{-R}$$
pour tout $Y\in {\cal D}$.
On montre que
$$sup_{e,e_{0},e_{4}}\vert \varphi_{e,e_{0},e_{4},g'}\vert _{k}<< \sigma(g')^k\Xi^G(g'),$$
cf. [W3] 6.7. La majoration pr\'ec\'edente devient
$$(1) \qquad \vert \Phi_{Y}(g')-\Phi_{g'}(Y)\vert \leq c\sigma(g')^k\Xi^G(g')\vert Y\vert ^{-R}.$$

 Montrons que cela entra\^{\i}ne que $\Phi_{g'}$ est constante. Pour cela, on peut fixer $g'$. Fixons $c_{1}$ et $c_{2}$ v\'erifiant les conditions de la proposition 6.5 pour $\eta=1/2$. Pour $Y\in {\cal D}$, notons $N_{Y}$ la partie enti\`ere de $2c_{2}^{-1}\vert Y\vert +1$. Soit $Y'\in {\cal D}$ tel que $\vert Y-Y'\vert \leq \vert Y\vert /2$. Si $\vert Y\vert $ est assez grand, les deux couples $(N_{Y},Y)$ et $(N_{Y},Y')$ v\'erifient les conditions de la proposition 6.5 (et on a aussi $\sigma(g')\leq Clog(N_{Y})$). Cette proposition appliqu\'ee aux deux couples entra\^{\i}ne
$$\vert \Phi_{Y}(g')-\Phi_{Y'}(g')\vert <<N_{Y}^{-R}<<\vert Y\vert ^{-R}.$$
Gr\^ace \`a (2), on en d\'eduit
$$\vert \Phi_{g'}(Y)-\Phi_{g'}(Y')\vert <<\vert Y\vert ^{-R}.$$
Or un \'el\'ement de $PolExp$ qui v\'erifie une telle majoration pour tout $Y\in {\cal D}$ tel que $\vert Y\vert $ soit assez grand et tout $Y'\in {\cal D}$ tel que $\vert Y-Y'\vert \leq \vert Y\vert /2$ ne peut qu'\^etre constant. Cela prouve l'assertion.

Puisque $\Phi_{g'}$ est constant, sa valeur n'est autre que $c_{{\cal R},0}(\Phi_{g'})$ pour tout r\'eseau ${\cal R}\subset {\cal A}_{\tilde{A}_{d},F}\otimes_{{\mathbb Z}}{\mathbb Q}$. L'assertion (ii)(a) de la proposition 6.6  implique alors que cette valeur est \'egale \`a
$$\sum_{e,e_{4}\in {\cal B}_{{\cal O}}^{K_{f'}},e_{0}\in {\cal B}_{{\cal O}}^{K_{0}}}\Phi(e,e_{0},e'',e_{4},\varphi_{e,e_{0},e_{4},g'}).$$
Mais on reconstitue ais\'ement cette somme: c'est $\Phi(g')$. Alors la majoration (1) devient celle de l'\'enonc\'e.
 $\square$
 
 \bigskip

\subsection{Simplification de $\Phi(g')$}

Soit $\tilde{L}'\in {\cal L}^{\tilde{G}}$ tel que $L\subset L'$.  Fixons $\tilde{S}'=\tilde{L}'U'\in {\cal P}(\tilde{L}')$.   Notons $\Lambda_{{\cal O},ell}^{\tilde{L}'}$ l'ensemble des $\lambda\in i{\cal A}_{L}^*$ tels que la repr\'esentation $Ind_{Q\cap L'}^{L'}(\tau_{\lambda})$ s'\'etende en une repr\'esentation elliptique de $\tilde{L}'(F)$. On a d\'ecrit ces repr\'esentations en 2.4. Pour tout \'el\'ement $\lambda$ de cet ensemble, fixons un prolongement unitaire  $ \tilde{\sigma}_{\lambda}$  de $Ind_{Q\cap L'}^{L'}(\tau_{\lambda})$ \`a $\tilde{L}'(F)$. A cette repr\'esentation est associ\'e un caract\`ere pond\'er\'e $J_{\tilde{L}'}^{\tilde{G}}(\tilde{\sigma}_{\lambda},.)$ sur $C_{c}^{\infty}(\tilde{G}(F))$, cf. 1.9. La repr\'esentation $\pi_{\lambda}=Ind_{Q}^G(\tau_{\lambda})$ s'identifie \`a l'induite de $S'(F)$ \`a $G(F)$ de la repr\'esentation pr\'ec\'edente (remarquons que, pour les groupes lin\'eaires, toutes ces induites sont irr\'eductibles). Du prolongement fix\'e $\tilde{\sigma}_{\lambda}$ se d\'eduit un prolongement $\tilde{\pi}_{\lambda}$ de $\pi_{\lambda}$ \`a $\tilde{G}(F)$.

 Remarquons que $\Lambda_{{\cal O},ell}^{\tilde{L}'}$ est invariant par translations par $i{\cal A}_{L,F}^{\vee}+i{\cal A}_{\tilde{L}'}^*$.  On peut supposer et on suppose que, pour $\chi\in i{\cal A}_{L,F}^{\vee}+i{\cal A}_{\tilde{L}'}^*$, $\tilde{\sigma}_{\lambda+\chi}=(\tilde{\sigma}_{\lambda})_{\chi}$. Rappelons que l'on a not\'e $a_{\tilde{L}'}$ la dimension de ${\cal A}_{\tilde{L}'}$ et que, pour tout $\lambda\in \Lambda_{{\cal O},ell}^{\tilde{L}'}$, on a d\'efini un entier $s(\tilde{\sigma}_{\lambda})$ en 2.4.
 
 \ass{Lemme}{Pour tout $g'\in G(F)$, on a l'\'egalit\'e
 $$\Phi(g')=\sum_{\tilde{L}'\in {\cal L}^{\tilde{G}}; L\subset L'} (-1)^{a_{\tilde{L}'}}\sum_{\lambda\in \Lambda^{\tilde{L}'}_{{\cal O},ell}/(i{\cal A}_{L,F}^{\vee}+i{\cal A}_{\tilde{L}'}^*)} 2^{-s(\tilde{\sigma}_{\lambda})-a_{\tilde{L}'}} $$
 $$\int_{i{\cal A}_{\tilde{L}',F}^*} (\tilde{\pi}_{\lambda+\chi}(\boldsymbol{\theta}_{d})e'',\pi_{\lambda+\chi}(g')e')J_{\tilde{L}'}^{\tilde{G}}(\tilde{\sigma}_{\lambda+\chi},\tilde{f})\varphi(\lambda+\chi)d\chi.$$}
 
 Preuve. Fixons $\tilde{L}'\in {\cal L}^{\tilde{G}}$ tel que $L\subset L'$, $\tilde{t}\in W^{\tilde{L}'}(L)_{reg}$ et $\lambda\in \Lambda_{{\cal O}}(\tilde{t})$. Consid\'erons le terme $(jd)_{\tilde{L}'}(\tilde{t},\lambda,g')$ du paragraphe pr\'ec\'edent. Rempla\c{c}ons dans les d\'efinitions des $(\tilde{G},\tilde{L}')$-familles les op\'erateurs d'entrelacement par des op\'erateurs normalis\'es. Cela remplace la famille produit $((jd)_{\tilde{Q}'}(\tilde{t},\lambda,g',\mu))_{\tilde{Q}'\in {\cal P}(\tilde{L}')}$ par un triple produit $((jdc)_{\tilde{Q}'}(\tilde{t},\lambda,g',\mu))_{\tilde{Q}'\in {\cal P}(\tilde{L}')}$, o\`u maintenant les familles $(j_{\tilde{Q}'}(\tilde{t},\lambda,\mu))_{\tilde{Q}'\in {\cal P}(\tilde{L}')}$ et $(d_{\tilde{Q}'}(\tilde{t},\lambda,g',\mu))_{\tilde{Q}'\in {\cal P}(\tilde{L}')}$ sont d\'efinies \`a l'aide d'op\'erateurs normalis\'es et $(c_{\tilde{Q}'}(\tilde{t},\lambda,g',\mu))_{\tilde{Q}'\in {\cal P}(\tilde{L}')}$ est form\'ee des facteurs de normalisation.
 On a donc maintenant
 $$j_{\tilde{Q}'}(\tilde{t},\lambda,\mu)=\sum_{e\in {\cal B}_{{\cal O}}^{K_{f'}}}(R_{Q\vert \theta_{\tilde{t}}(Q)}(\tau_{\lambda})A(\tau_{\lambda},\tilde{t})e,R_{Q(\bar{Q}')\vert Q}(\tau_{\lambda})^{-1}R_{Q(\bar{Q}')\vert Q}(\tau_{\lambda+\mu}) \pi_{\lambda}(f')e).$$
 Puisque $\theta_{\tilde{t}}(\tau_{\lambda})=\tau_{\lambda}$, la repr\'esentation $Ind_{L'\cap Q}^{L'}(\tau_{\lambda})$ s'\'etend en une repr\'esentation irr\'eductible de $\tilde{L}'(F)$. Celle-ci n'a pas de raison d'\^etre elliptique, mais on peut tout de m\^eme d\'efinir comme avant l'\'enonc\'e les repr\'esentations $\tilde{\sigma}_{\lambda}$ et $\tilde{\pi}_{\lambda}$.
 Consid\'erons l'op\'erateur 
 $$R_{Q\vert \theta_{\tilde{t}}(Q)}(\tau_{\lambda})A(\tau_{\lambda},\tilde{t}).$$
  En revenant \`a sa d\'efinition, on voit qu'il v\'erifie
 $$\pi_{\lambda}(g)R_{Q\vert \theta_{\tilde{t}}(Q)}(\tau_{\lambda})A(\tau_{\lambda},\tilde{t})=R_{Q\vert \theta_{\tilde{t}}(Q)}(\tau_{\lambda})A(\tau_{\lambda},\tilde{t})\pi_{\lambda}(\theta(g))$$
 pour tout $g\in G(F)$. Mais $\tilde{\pi}_{\lambda}(\boldsymbol{\theta}_{d})^{-1}$ v\'erifie la m\^eme propri\'et\'e. Puisque $\pi_{\lambda}$ est irr\'eductible, les deux op\'erateurs sont proportionnels. Puisqu'ils sont tous deux unitaires, il existe un nombre complexe $z$ de module $1$ tel que
 $$R_{Q\vert \theta_{\tilde{t}}(Q)}(\tau_{\lambda})A(\tau_{\lambda},\tilde{t})=z\tilde{\pi}_{\lambda}(\boldsymbol{\theta}_{d})^{-1}.$$
 L'ensemble $\{\tilde{\pi}_{\lambda}(\boldsymbol{\theta}_{d})e; e\in {\cal B}_{{\cal O}}^{K_{f'}}\}$ est encore une base orthonorm\'ee de $({\cal K}_{Q,\tau}^G)^{K_{f'}}$. Notons-la ${\cal B}_{\natural}^{K_{f'}}$. Par le changement de variables $e\mapsto \tilde{\pi}_{\lambda}(\boldsymbol{\theta}_{d})e$, on obtient
 $$j_{\tilde{Q}'}(\tilde{t},\lambda,\mu)=\sum_{e\in {\cal B}_{\natural}^{K_{f'}}}( ze,R_{Q(\bar{Q}')\vert Q}(\tau_{\lambda})^{-1}R_{Q(\bar{Q}')\vert Q}(\tau_{\lambda+\mu}) \pi_{\lambda}(f')\tilde{\pi}_{\lambda}(\boldsymbol{\theta}_{d})e).$$
 En se rappelant que $f'(g)=\tilde{f}(g\boldsymbol{\theta}_{d})$, on v\'erifie que
 $$ \pi_{\lambda}(f')\tilde{\pi}_{\lambda}(\boldsymbol{\theta}_{d})=\tilde{\pi}_{\lambda}(\tilde{f}).$$
 Alors
 $$j_{\tilde{Q}'}(\tilde{t},\lambda,\mu)=\bar{z}\,trace(R_{Q(\bar{Q}')\vert Q}(\tau_{\lambda})^{-1}R_{Q(\bar{Q}')\vert Q}(\tau_{\lambda+\mu}) \tilde{\pi}_{\lambda}(\tilde{f})).$$
 On peut encore remplacer le sous-groupe parabolique $Q$ par $Q(S')$. En effet,  les propri\'et\'es de composition des op\'erateurs d'entrelacement entra\^{\i}nent que
 $$j_{\tilde{Q}'}(\tilde{t},\lambda,\mu)=\bar{z}\,trace(R_{Q(\bar{Q}')\vert Q(S')}(\tau_{\lambda})^{-1}R_{Q(\bar{Q}')\vert Q(S')}(\tau_{\lambda+\mu}) r(\lambda,\mu)\tilde{\pi}_{\lambda}(\tilde{f}) ),$$
 o\`u on r\'ealise maintenant $\tilde{\pi}_{\lambda}$ dans ${\cal K}_{Q(S'),\tau}^G$ et o\`u
 $$r(\lambda,\mu)=R_{Q(S')\vert Q}(\tau_{\lambda+\mu})R_{Q(S')\vert Q}(\tau_{\lambda})^{-1}.$$
 Cet op\'erateur ne d\'epend pas de $Q'$ et une propri\'et\'e famili\`ere des $(\tilde{G},\tilde{L}')$-familles entra\^{\i}ne qu'il dispara\^{\i}t quand on calcule le nombre associ\'e $j_{\tilde{L}'}(\tilde{t},\lambda)$. Autrement dit, si on d\'efinit une $(\tilde{G},\tilde{L}')$-famille $(j'_{\tilde{Q}'}(\tilde{t},\lambda))_{\tilde{Q}'\in {\cal P}(\tilde{L}')}$ par
 $$j'_{\tilde{Q}'}(\tilde{t},\lambda,\mu)=\bar{z}\,trace(R_{Q(\bar{Q}')\vert Q(S')}(\tau_{\lambda})^{-1}R_{Q(\bar{Q}')\vert Q(S')}(\tau_{\lambda+\mu}) \tilde{\pi}_{\lambda}(\tilde{f})) ,$$
 on a l'\'egalit\'e $j_{\tilde{L}'}(\tilde{t},\lambda)=j'_{\tilde{L}'}(\tilde{t},\lambda)$. Mais on reconna\^{\i}t ce dernier terme en se reportant \`a la d\'efinition de 1.9. On obtient
 $$ j_{\tilde{L}'}(\tilde{t},\lambda)=\bar{z}J_{\tilde{L}'}^{\tilde{G}}(\tilde{\sigma}_{\lambda},\tilde{f}).$$
 Si $\tilde{Q}''=\tilde{L}''U''\in {\cal F}(\tilde{L}')$, le m\^eme calcul vaut pour les $(\tilde{L}'',\tilde{L}')$-familles d\'eduites, c'est-\`a-dire que l'on a
 $$j_{\tilde{L}'}^{\tilde{Q}''}(\tilde{t},\lambda)=\bar{z}J_{\tilde{L}'}^{\tilde{Q}''}(\tilde{\sigma}_{\lambda},\tilde{f}).$$
 On se rappelle que $\tilde{f}$ est tr\`es cuspidale. D'apr\`es  le lemme 1.13(i), le terme ci-dessus est nul si $\tilde{Q}''\not=\tilde{G}$. En appliquant les formules de descente habituelles, on obtient
 $$(jdc)_{\tilde{L}'}(\tilde{t},\lambda,g')=j_{\tilde{L}'}(\tilde{t},\lambda)d_{\tilde{L}'}^{\tilde{Q'}}(\tilde{t},\lambda,g')c_{\tilde{L}'}^{\tilde{Q}'}(\lambda),$$
 o\`u $\tilde{Q}'$ est un \'el\'ement quelconque de ${\cal P}(\tilde{L}')$. On v\'erifie que $c_{\tilde{L}'}^{\tilde{Q}'}(\lambda)=1$ et
 $$d_{\tilde{L}'}^{\tilde{Q}'}(\tilde{t},\lambda,g')=(e'',R_{Q\vert \theta_{\tilde{t}}(Q)}(\tau_{\lambda})A(\tau_{\lambda},\tilde{t})\pi_{\lambda}(g')e')$$
 $$=(e'',z\tilde{\pi}_{\lambda}(\boldsymbol{\theta}_{d})^{-1}\pi_{\lambda}(g')e')=z(\tilde{\pi}_{\lambda}(\boldsymbol{\theta}_{d})e'',\pi_{\lambda}(g')e').$$
 Les termes $z\bar{z}$ disparaissent puisque $z$ est de module $1$.Toujours d'apr\`es le lemme 1.13(i), le caract\`ere pond\'er\'e $J_{\tilde{L}'}^{\tilde{G}}(\tilde{\sigma}_{\lambda},\tilde{f})$ est nul si $\tilde{\sigma}_{\lambda}$ est une induite propre, ce qui \'equivaut \`a ce qu'elle ne soit pas elliptique. On a obtenu
 $$(jdc)_{\tilde{L}'}(\tilde{t},\lambda,g')=\left\lbrace\begin{array}{cc}(\tilde{\pi}_{\lambda}(\boldsymbol{\theta}_{d})e'',\pi_{\lambda}(g')e')J_{\tilde{L}'}^{\tilde{G}}(\tilde{\sigma}_{\lambda},\tilde{f}),&\text{ si }\lambda\in \Lambda_{{\cal O},ell}^{\tilde{L}'},\\ 0,&\text{ sinon.}\\ \end{array}\right.$$
 Il reste \`a remarquer que pour tout $\lambda\in \Lambda_{{\cal O},ell}^{\tilde{L}'}$, il existe exactement un $\tilde{t}\in W^{\tilde{L}'}(L)_{reg}$ tel que $\lambda\in \Lambda_{{\cal O}}(\tilde{t})$ et que, pour ce $\tilde{t}$, on a $\vert det(1-\theta_{\tilde{t}})_{\vert {\cal A}_{M}/{\cal A}_{\tilde{L}'}}\vert =2^{s(\tilde{\sigma}_{\lambda})+a_{\tilde{L}'}}$. Ces propri\'et\'es se v\'erifient ais\'ement sur la description de 2.4. Cela ach\`eve la preuve. $\square$
 
 \bigskip
 
 \subsection{Evaluation d'une limite}
 
 \ass{Lemme}{On a l'\'egalit\'e
 
 $$lim_{N\to \infty}J_{L,{\cal O},N,C}(\Theta_{\tilde{\rho}},\tilde{f})=[i{\cal A}_{{\cal O}}^{\vee}:i{\cal A}_{L,F}^{\vee}]^{-1}\sum_{\tilde{L}'\in {\cal L}^{\tilde{G}}; L\subset L'}(-1)^{a_{\tilde{L}'}} $$
  $$\sum_{\lambda\in \Lambda^{\tilde{L}'}_{{\cal O},ell}/((i{\cal A}_{L,F}^{\vee}+i{\cal A}_{\tilde{L}'})}2^{-s(\pi_{\lambda}^{L'})-a_{\tilde{L}'}}\epsilon_{\nu}((\tilde{\sigma}_{\lambda})^{\vee},\tilde{\rho}) \int_{i{\cal A}_{L',F}^*} J_{\tilde{L}'}^{\tilde{G}}(\tilde{\sigma}_{\lambda+\chi},\tilde{f})d\chi.$$}
 
 Preuve. Consid\'erons la d\'efinition de $J_{L,{\cal O},N,C}(\Theta_{\tilde{\rho}},\tilde{f})$ donn\'ee avant le lemme 6.4. Un calcul formel montre qu'elle ne d\'epend pas du choix de la base orthonorm\'ee ${\cal B}_{{\cal O}}^{K_{f}}$, ni de celui du groupe $K_{f}$ pourvu qu'il soit assez petit. En faisant entrer la sommation en $e$ dans la derni\`ere int\'egrale, on peut m\^eme remplacer ${\cal B}_{{\cal O}}^{K_{f}}$ par une base d\'ependant de $\lambda$. Soit $\gamma$ l'\'el\'ement de $G(F)$ tel que $\tilde{y}=\gamma\boldsymbol{\theta}_{d}$. On a les \'egalit\'es $f(g)=\tilde{f}(g\tilde{y})=\tilde{f}(g\gamma\boldsymbol{\theta}_{d})=f'(g\gamma)$ pour tout $g\in G(F)$. On peut donc supposer que $K_{f}=\gamma K_{f'}\gamma^{-1}$ et remplacer ${\cal B}_{{\cal O}}^{K_{f}}$ par
 $$\{\pi_{\lambda}(\gamma)e;e\in {\cal B}_{{\cal O}}^{K_{f'}}\}.$$
 On a $\pi_{\lambda}(f)\pi_{\lambda}(\gamma)e=\pi_{\lambda}(f')e$, et $\pi_{\lambda}(\theta_{\tilde{y}}(g))\pi_{\lambda}(\gamma)e=\pi_{\lambda}(\gamma)\pi_{\lambda}(\theta(g))e$. On obtient ainsi
 $$J_{L,{\cal O},N,C}(\Theta_{\tilde{\rho}},\tilde{f})= [i{\cal A}_{{\cal O}}^{\vee}:i{\cal A}_{L,F}^{\vee}]^{-1}\sum_{j=1,...,n}$$
 $$\int_{H(F)U(F)_{c}}{\bf 1}_{\sigma<Clog(N)}(hu)(\rho(h)\epsilon_{j},\tilde{\rho}(\tilde{y})\epsilon_{j})\bar{\xi}(u)\Phi_{N,j}(\gamma^{-1}hu)du\,dh,$$
 o\`u
 $$\Phi_{N,j}(\gamma^{-1}hu)=\sum_{e\in {\cal B}_{{\cal O}}^{K_{f'}}}\int_{i{\cal A}_{L,F}^*}m(\tau_{\lambda})\varphi_{j}(\lambda)$$
 $$\int_{G(F)}(e_{j},\pi_{\lambda}(f')e)(\pi_{\lambda}(\theta(g))e,\pi_{\lambda}(\gamma^{-1}hu)e_{j})\kappa_{N}(g)dg\,d\lambda.$$
 Ce terme $\Phi_{N,j}(\gamma^{-1}hu)$ est exactement le terme $\Phi_{N}(g')$ d\'efini en 6.5, associ\'e aux donn\'ees auxiliaires $e'=e''=e_{j}$, $\varphi=\varphi_{j}$, \'evalu\'e au point $g'=\gamma^{-1}hu$. Notons $\Phi_{Y,j}$ et $\Phi_{j}$ les fonctions associ\'ees \`a ces donn\'ees d\'efinies en 6.5 et 6.7. Remarquons que, si ${\bf 1}_{\sigma<Clog(N)}(hu)=1$, l'\'el\'ement $g'=\gamma^{-1}hu$ v\'erifie l'hypoth\`ese de la proposition 6.5, quitte \`a agrandir $C$, ce qui importe \'evidemment peu. Fixons $c_{1}$ et $c_{2}$ v\'erifiant les conditions de la proposition 6.5 pour $\eta=1/2$ et pour tout $j=1,...,n$. Si $N$ est assez grand, on peut choisir $Y\in {\cal D}$ tel que $c_{1}N^{1/2}\leq c_{2}N/2\leq \vert Y\vert \leq c_{2}N$. En appliquant successivement les propositions 6.5 et 6.7 pour cet $Y$, on obtient pour tout $R$ une majoration
 $$\vert \Phi_{N,j}(\gamma^{-1}hu)-\Phi_{j}(\gamma^{-1}hu)\vert<< (1+\sigma(hu)^{k(R)}\Xi^G(hu))N^{-R}$$
 pour tout $j$, tout $N$ assez grand et tous $h$, $u$ tels que ${\bf 1}_{\sigma<Clog(N)}(hu)=1$, o\`u $k(R)$ est un entier d\'ependant de $R$ comme la notation l'indique. On peut oublier le terme $\Xi^G(hu)$ qui est born\'e. Posons
 $$X_{N}=[i{\cal A}_{{\cal O}}^{\vee}:i{\cal A}_{L,F}^{\vee}]^{-1}\sum_{j=1,...,n}\int_{H(F)U(F)_{c}}{\bf 1}_{\sigma<Clog(N)}(hu)(\rho(h)\epsilon_{j},\tilde{\rho}(\tilde{y})\epsilon_{j})\bar{\xi}(u)\Phi_{j}(\gamma^{-1}hu)du\,dh.$$
 Alors
 $$\vert J_{L,{\cal O},N,C}(\Theta_{\tilde{\rho}},\tilde{f})-X_{N}\vert <<N^{-R}\int_{H(F)U(F)_{c}}{\bf 1}_{\sigma<Clog(N)}(hu)\Xi^H(h)\sigma(hu)^{k(R)}du\,dh$$
 $$<<log(N)^{k(R)}N^{-R}\int_{H(F)U(F)_{c}}{\bf 1}_{\sigma<Clog(N)}(hu)\Xi^H(h)du\,dh.$$
 On v\'erifie qu'il existe $R_{1}$ tel que l'int\'egrale soit born\'ee par $N^{R_{1}}$, cf. [W3] 4.3(1). Puisque $R$ est quelconque, on obtient que 
 $$lim_{N\to \infty}(J_{L,{\cal O},N,C}(\Theta_{\tilde{\rho}},\tilde{f})-X_{N})=0.$$
 Il reste \`a calculer $lim_{N\to \infty}X_{N}$. En utilisant le lemme 6.8, on a
 $$X_{N}=[i{\cal A}_{{\cal O}}^{\vee}:i{\cal A}_{L,F}^{\vee}]^{-1} \sum_{\tilde{L}'\in {\cal L}^{\tilde{G}}; L\subset L'} (-1)^{a_{\tilde{L}'}}\sum_{\lambda\in \Lambda^{\tilde{L}'}_{{\cal O},ell}/(i{\cal A}_{L,F}^{\vee}+i{\cal A}_{\tilde{L}'}^*)} 2^{-s(\tilde{\sigma}_{\lambda})-a_{\tilde{L}'}}$$
 $$ \int_{i{\cal A}_{\tilde{L}',F}^*} J_{\tilde{L}'}^{\tilde{G}}(\tilde{\sigma}_{\lambda+\chi},\tilde{f})
 \sum_{j=1,...,n}\varphi_{j}(\lambda+\chi)$$
 $$\int_{H(F)U(F)_{c}}{\bf 1}_{\sigma<Clog(N)}(hu)(\rho(h)\epsilon_{j},\tilde{\rho}(\tilde{y})\epsilon_{j})(\tilde{\pi}_{\lambda+\chi}(\boldsymbol{\theta}_{d})e_{j},\pi_{\lambda+\chi}(\gamma^{-1}hu)e_{j}) \bar{\xi}(u)du\,dh\,d\chi.$$
 Les permutations d'int\'egrales sont justifi\'ees par la convergence absolue de cette expression, qui provient elle-m\^eme de la convergence de
 $$\int_{H(F)U(F)_{c}}\Xi^H(h)\Xi^G(hu)du\,dh,$$
 cf. 4.1(3). Pour la m\^eme raison, on peut appliquer le th\'eor\`eme de convergence domin\'ee pour calculer
 $$(1) \qquad lim_{N\to \infty}X_{N}=[i{\cal A}_{{\cal O}}^{\vee}:i{\cal A}_{L,F}^{\vee}]^{-1} \sum_{\tilde{L}'\in {\cal L}^{\tilde{G}}; L\subset L'} (-1)^{a_{\tilde{L}'}}\sum_{\lambda\in \Lambda^{\tilde{L}'}_{{\cal O},ell}/(i{\cal A}_{L,F}^{\vee}+i{\cal A}_{\tilde{L}'}^*)} 2^{-s(\tilde{\sigma}_{\lambda})-a_{\tilde{L}'}}$$
 $$ \int_{i{\cal A}_{\tilde{L}',F}^*} J_{\tilde{L}'}^{\tilde{G}}(\tilde{\sigma}_{\lambda+\chi},\tilde{f})\sum_{j=1,...,n}\varphi_{j}(\lambda+\chi)$$
 $$\int_{H(F)U(F)_{c}} (\rho(h)\epsilon_{j},\tilde{\rho}(\tilde{y})\epsilon_{j})(\tilde{\pi}_{\lambda+\chi}(\boldsymbol{\theta}_{d})e_{j},\pi_{\lambda+\chi}(\gamma^{-1}hu)e_{j}) \bar{\xi}(u)du\,dh\,d\chi.$$
 On a l'\'egalit\'e
  $$(\tilde{\pi}_{\lambda+\chi}(\boldsymbol{\theta}_{d})e_{j},\pi_{\lambda+\chi}(\gamma^{-1}hu)e_{j})  =(\tilde{\pi}_{\lambda+\chi}(\tilde{y})e_{j},\pi_{\lambda+\chi}(hu)e_{j}) .$$
 On reconna\^{\i}t alors la derni\`ere int\'egrale en $h$, $u$: elle vaut 
 $${\cal L}_{\pi_{\lambda+\chi},\rho}(\epsilon_{j}\otimes \tilde{\pi}_{\lambda+\chi}(\tilde{y})e_{j},\tilde{\rho}(\tilde{y})\epsilon_{j}\otimes e_{j}).$$
 D'apr\`es la proposition 5.5, c'est aussi
 $$ \epsilon_{\nu}((\tilde{\pi}_{\lambda+\chi})^{\vee},\tilde{\rho}){\cal L}_{\pi_{\lambda+\chi},\rho}(\epsilon_{j}\otimes e_{j},\epsilon_{j}\otimes e_{j}).$$
 Par d\'efinition des familles $(\epsilon_{j})_{j=1,...,n}$, $(e_{j})_{j=1,...,n}$ et $(\varphi_{j})_{j=1,...,n}$, on a l'\'egalit\'e
 $$\sum_{j=1,...,n}\varphi_{j}(\lambda+\chi){\cal L}_{\pi_{\lambda+\chi},\rho}(\epsilon_{j}\otimes e_{j},\epsilon_{j}\otimes e_{j})=1.$$
 D'apr\`es 2.5(2), on a $\epsilon_{\nu}((\tilde{\pi}_{\lambda+\chi})^{\vee},\tilde{\rho})=\epsilon_{\nu}((\tilde{\sigma}_{\lambda})^{\vee},\tilde{\rho})$.   Mais alors le membre de droite de la formule (1) devient celui de l'\'egalit\'e de l'\'enonc\'e. $\square$

 \bigskip
 
 \subsection{Preuve du th\'eor\`eme 6.1}
 
 Les lemmes 6.4 et 6.9 prouvent que $J_{N}(\Theta_{\tilde{\rho}},\tilde{f})$ a une limite quand $N$ tend vers l'infini et ils calculent cette limite. En intervertissant les sommations en $L$ et $\tilde{L}'$, on obtient
 $$(1) \qquad lim_{N\to\infty}J_{N}(\Theta_{\tilde{\rho}},\tilde{f})=\sum_{\tilde{L}'\in {\cal L}^{\tilde{G}}}(-1)^{a_{\tilde{L}'}}\vert W^G\vert ^{-1}X(\tilde{L}'),$$
 o\`u
 $$X(\tilde{L}')=\sum_{L\in {\cal L}^{L'}(A_{d})}\vert W^L\vert \sum_{{\cal O}\in \{\Pi_{2}(L)\}_{f}}[i{\cal A}_{{\cal O}}^{\vee}:i{\cal A}_{L,F}^{\vee}]^{-1}$$
 $$\sum_{\lambda\in \Lambda^{\tilde{L}'}_{{\cal O},ell}/((i{\cal A}_{L,F}^{\vee}+i{\cal A}_{\tilde{L}'}^*)}2^{-s(\sigma_{\lambda})-a_{\tilde{L}'}} \epsilon_{\nu}((\tilde{\sigma}_{\lambda})^{\vee},\tilde{\rho})\int_{i{\cal A}_{L',F}^*}J_{\tilde{L}'}^{\tilde{G}}(\tilde{\sigma}_{\lambda+\chi},\tilde{f})d\chi.$$
Fixons $\tilde{L}'$. Dans la formule ci-dessus, on peut remplacer la sommation sur ${\cal O}\in \{\Pi_{2}(L)\}_{f}$ par une sommation sur tout $\{\Pi_{2}(L)\}$: pour une orbite dans le compl\'ementaire, les fonctions $J_{\tilde{L}'}^{\tilde{G}}(\tilde{\sigma}_{\lambda+\chi},\tilde{f})$ sont nulles. Consid\'erons l'ensemble ${\cal Z}$ des triplets $(L,{\cal O},\lambda)$ intervenant dans  $X(\tilde{L}')$. A un tel triplet, associons l'orbite de $\tilde{\sigma}_{\lambda}$ sous l'action de $i{\cal A}_{\tilde{L}'}^*$. On obtient une application $\iota:{\cal Z}\to \{\Pi_{ell}(\tilde{L}')\}$. Cette application est surjective. D'o\`u
$$X(\tilde{L}')=\sum_{{\cal O}'\in \{\Pi_{ell}(\tilde{L}')\}}x({\cal O}')\int_{i{\cal A}_{L',F}^*}J_{\tilde{L}'}^{\tilde{G}}(\tilde{\sigma}'_{\chi},\tilde{f})d\chi,$$
o\`u, dans chaque orbite ${\cal O}'$, on a fix\'e un point base $\tilde{\sigma}'$, et o\`u on a pos\'e
$$(2) \qquad x({\cal O}')=\sum_{z=(L,{\cal O},\lambda)\in {\cal Z}, \iota(z)={\cal O}'}\vert W^L\vert [i{\cal A}_{{\cal O}}^{\vee}:i{\cal A}_{L,F}^{\vee}]^{-1}2^{-s(\tilde{\sigma}_{\lambda})-a_{\tilde{L}'}}\epsilon_{\nu}((\tilde{\sigma}_{\lambda})^{\vee},\tilde{\rho}).$$
Fixons ${\cal O}'$ et notons ${\cal Z}'$ l'ensemble des triplets $(L,{\cal O},\lambda)$ v\'erifiant les conditions suivantes. Les termes $L$ et ${\cal O}$ sont comme ci-dessus. L'\'el\'ement $\lambda$ appartient \`a $\Lambda^{\tilde{L}'}_{{\cal O},ell}/i{\cal A}_{L,F}^{\vee}$. On impose que $\tilde{\sigma}_{\lambda}\simeq \tilde{\sigma}'$. Les projections de $\Lambda^{\tilde{L}'}_{{\cal O},ell}/i{\cal A}_{L,F}^{\vee}$ sur $\Lambda^{\tilde{L}'}_{{\cal O},ell}/(i{\cal A}_{L,F}^{\vee}+i{\cal A}_{\tilde{L}'}^*)$ induisent une surjection de ${\cal Z}$ sur la fibre de $\iota$ au-dessus de ${\cal O}'$.  Deux \'el\'ements $(L,{\cal O},\lambda)$ et $(L_{1},{\cal O}_{1},\lambda_{1})$ de ${\cal Z}'$ ont m\^eme image dans ${\cal Z}$ si et seulement si $L_{1}=L$, ${\cal O}_{1}={\cal O}$ et il existe $\chi\in i{\cal A}_{\tilde{L}'}^*$ tel que $\lambda_{1}=\lambda+\chi$. La condition $\tilde{\sigma}_{\lambda}=\tilde{\sigma}_{\lambda_{1}}=\tilde{\sigma}'$ impose que $\chi\in i{\cal A}_{{\cal O}'}^{\vee}$. Autrement dit, ce sont des orbites pour les actions du groupe $i{\cal A}_{{\cal O}'}^{\vee}$ sur les espaces $i{\cal A}_{L}^*/i{\cal A}_{L,F}^{\vee}$. Puisque $i{\cal A}_{\tilde{L}'}^*\cap i{\cal A}_{L,F}^{\vee}=i{\cal A}_{\tilde{L}',F}^{\vee}$, toutes les fibres ont pour nombre d'\'el\'ements $[ i{\cal A}_{{\cal O}'}^{\vee}:i{\cal A}_{\tilde{L}',F}^{\vee}]$. Dans (2), on peut remplacer la sommation sur les  $z\in {\cal Z}$ tel que $\iota(z)={\cal O}'$ par une sommation sur les $z\in {\cal Z}'$, \`a condition de diviser par $[ i{\cal A}_{{\cal O}'}^{\vee}:i{\cal A}_{\tilde{L}',F}^{\vee}]$. Fixons $(L,{\cal O},\lambda)\in {\cal Z}'$ et soit $(L_{1},{\cal O}_{1},\lambda_{1})$ un autre \'el\'ement. D'apr\`es Harish-Chandra, la condition $\tilde{\sigma}_{\lambda}\simeq \tilde{\sigma}'\simeq \tilde{\sigma}_{1,\lambda_{1}}$ \'equivaut \`a l'existence de $w\in W^{L'}$ tel que $L_{1}=wLw^{-1}$, ${\cal O}_{1}=w{\cal O}$ et $\tau_{1,\lambda_{1}}\simeq (w\tau)_{w\lambda}$. Il en r\'esulte d'abord que les termes associ\'es \`a $(L_{1},{\cal O}_{1},\lambda_{1})$ qui apparaissent dans (2) sont \'egaux \`a ceux associ\'es \`a $(L,{\cal O},\lambda)$. Donc
$$(3) \qquad x({\cal O}')=\vert {\cal Z}'\vert [ i{\cal A}_{{\cal O}'}^{\vee}:i{\cal A}_{\tilde{L}',F}^{\vee}]^{-1}\vert W^L\vert [i{\cal A}_{{\cal O}}^{\vee}:i{\cal A}_{L,F}^{\vee}]^{-1}2^{-s(\tilde{\sigma}_{\lambda})-a_{\tilde{L}'}}\epsilon_{\nu}((\tilde{\sigma}_{\lambda})^{\vee},\tilde{\rho}).$$
D'autre part, pour $w\in W^{L'}$, soit $\Lambda(w)$ l'ensemble des $\lambda_{1}\in \Lambda^{\tilde{L}'}_{w{\cal O},ell}/i{\cal A}_{wLw^{-1},F}^{\vee}$ tels que $\tau_{1,\lambda_{1}}\simeq (w\tau)_{w\lambda}$ (o\`u $\tau_{1}$ est le point base fix\'e dans ${\cal O}_{1}=w{\cal O}$). Posons ${\cal Z}''=\{(w,\lambda_{1}); w\in W; \lambda_{1}\in \Lambda(w)\}$. L'application
$$\begin{array}{ccc}{\cal Z}''&\to &{\cal Z}'\\ (w,\lambda_{1})&\mapsto &(wLw^{-1},w{\cal O},\lambda_{1})\\ \end{array}$$
est surjective. Ses fibres ont toutes le m\^eme nombre d'\'el\'ements, \`a savoir le nombre de $w\in W^{L'}$ tels que $wLw^{-1}=L$ et $w\tau_{\lambda}\simeq \tau_{\lambda}$. Il r\'esulte de la description de 2.4 que cet ensemble se r\'eduit \`a $W^L$. D'autre part, les ensembles $\Lambda(w)$ ont tous le m\^eme nombre d'\'el\'ements, \`a savoir $[i{\cal A}_{{\cal O}}^{\vee}:i{\cal A}_{L,F}^{\vee}]$. On obtient
$$\vert {\cal Z}''\vert =\vert W^{L'}\vert  [i{\cal A}_{{\cal O}}^{\vee}:i{\cal A}_{L,F}^{\vee}],$$
 $$\vert {\cal Z}'\vert =\vert {\cal Z}''\vert \vert W^L\vert ^{-1} =\vert W^{L'}\vert  [i{\cal A}_{{\cal O}}^{\vee}:i{\cal A}_{L,F}^{\vee}]\vert W^L\vert ^{-1}.$$ 
 On reporte cette valeur dans l'expression (3). Cela  calcule explicitement $X(\tilde{L}')$.  On reporte cette valeur dans l'expression (1) et on obtient le th\'eor\`eme 6.1. $\square$
 
 \bigskip
 
 \section{Une formule int\'egrale calculant une valeur d'un facteur $\epsilon$}
 
 \bigskip
 
 \subsection{Enonc\'e du th\'eor\`eme}
 
 On consid\`ere la situation de 3.1. Soient $\tilde{\pi}\in Temp(\tilde{G})$ et $\tilde{\rho}\in Temp(\tilde{H})$. On reprend les d\'efinitions de 3.2, en particulier on introduit  l'ensemble ${\cal T}$. Soit $\tilde{T}\in {\cal T}$. Aux caract\`eres $\Theta_{\tilde{\pi}}$ et $\Theta_{\tilde{\rho}}$ sont associ\'ees des fonctions $c_{\Theta_{\tilde{\pi}}}$ et $c_{\Theta_{\tilde{\rho}}}$ sur $\tilde{T}(F)_{/\theta}$. On les note simplement $c_{\tilde{\pi}}$ et $c_{\tilde{\rho}}$. On pose
 $$\epsilon_{geom,\nu}(\tilde{\pi},\tilde{\rho})=\sum_{\tilde{T}\in {\cal T}}\vert W(H,\tilde{T})\vert ^{-1}\int_{\tilde{T}(F)_{/\theta}}c_{\tilde{\pi}}(\tilde{t})c_{\tilde{\rho}}(\tilde{t})D^{\tilde{H}}(\tilde{t})\Delta_{r}(\tilde{t})d\tilde{t}.$$
 On rappelle que l'on a d\'efini un nombre $\epsilon_{\nu}(\tilde{\pi},\tilde{\rho})$ en 2.5.
 
 \ass{Th\'eor\`eme}{Pour tous $\tilde{\pi}\in Temp(\tilde{G})$ et $\tilde{\rho}\in Temp(\tilde{H})$, on a l'\'egalit\'e 
 $$\epsilon_{geom,\nu}(\tilde{\pi},\tilde{\rho})= \epsilon_{\nu}(\tilde{\pi},\tilde{\rho}).$$}
 
 \bigskip
 
 \subsection{Terme g\'eom\'etrique et induction}
 
 Soient $ \Gamma$ un quasi-caract\`ere sur $\tilde{G}(F)$ et $\Theta$ un quasi-caract\`ere sur $\tilde{H}(F)$. Posons
 $$\Xi_{\nu}(\Gamma,\Theta)=\Xi_{-\nu}(\Theta,\Gamma)=\sum_{\tilde{T}\in {\cal T}}\vert W(H,\tilde{T})\vert ^{-1} \int_{\tilde{T}(F)_{/\theta}}c_{\Gamma}(\tilde{t})c_{\Theta}(\tilde{t})D^{\tilde{H}}(\tilde{t})\Delta_{r}(\tilde{t})d\tilde{t}.$$
 Consid\'erons un L\'evi tordu $\tilde{L}$ de $\tilde{G}$. On reprend les notations de 2.3, en particulier, on introduit une d\'ecomposition 
   $$V=V_{u}\oplus...\oplus V_{1}\oplus V_{0}\oplus V_{-1}\oplus...\oplus V_{-u}$$
telle que $\tilde{L}$ soit l'ensemble des \'el\'ements  $\tilde{x}\in \tilde{G}$ tels que $\tilde{x}(V_{j})=V_{-j}^*$ pour tout $j=-u,...,u$. On a
  $$L=GL_{d_{u}}\times...\times GL_{d_{1}}\times GL_{d_{0}}\times GL_{d_{1}}\times ...\times GL_{d_{u}}.$$
  On introduit   le groupe tordu $\tilde{G}_{0}$ analogue de $\tilde{G}$ quand on remplace $d$ par $d_{0}$. 
  
  {\bf Remarque}: les termes $V_{0}$ et $G_{0}$ n'ont plus la m\^eme signification qu'en 3.1.
  
  Soient $\Gamma_{0}$ un quasi-caract\`ere sur $\tilde{G}_{0}(F)$ et, pour tout $j=1,...,u$, soit $\Gamma_{j}$ un quasi-caract\`ere sur $GL_{d_{j}}(F)$. Pour $j=1,...,u$, le quasi-caract\`ere $\Gamma_{j}$ admet un d\'eveloppement au voisinage de l'\'el\'ement neutre, avec des coefficients $c_{\Gamma_{j},{\cal O}}(1)$, o\`u ${\cal O}$ parcourt les orbites nilpotentes de $\mathfrak{gl}_{j}(F)$. Il y a une unique orbite nilpotente r\'eguli\`ere, notons-la ${\cal O}_{reg}$. On pose
 $$\Xi(\Gamma_{j})=c_{\Gamma_{j},{\cal O}_{reg}}(1).$$
 
 Si $d_{0}=dim(V_{0})>m=dim(W)$, posons  $r_{0}=(d_{0}-m-1)/2$ et notons $Z_{0}$, resp. $Z_{0}^*$ le sous-espace de $Z$, resp. $Z^*$ engendr\'e par  les vecteurs $z_{i}$, resp. $z_{i}^*$, pour $i= -r_{0},...,r_{0}$. Quitte \`a conjuguer $\tilde{L}$, on peut supposer que $V_{0}=W\oplus Z_{0}$. On note $\tilde{\zeta}_{0}\in Isom(Z_{0},Z_{0}^*)$ la restriction de $\tilde{\zeta}$, cf. 3.1, et on d\'efinit le plongement $\iota_{0}:\tilde{H}\to \tilde{G}_{0}$ comme en 3.1, en rempla\c{c}ant $\tilde{\zeta}$ par $\tilde{\zeta}_{0}$. En appliquant les constructions pr\'ec\'edentes au couple $(\tilde{H},\tilde{G}_{0})$, on d\'efinit le nombre $\Xi_{\nu}(\Gamma_{0},\Theta)$.
 
 Si $d_{0}<m$, posons $r_{0}=(m-d_{0}-1)/2$. Quitte \`a conjuguer $\tilde{L}$, on peut supposer $V_{0}\subset W$. On fixe un suppl\'ementaire $Z_{0}$ de $V_{0}$ dans $W$ et une base $(z_{0,i})_{i=-r_{0},...,r_{0}}$ de $Z_{0}$. On identifie $W^*$ \`a $V_{0}^*\oplus Z_{0}^*$. On note $(z_{0,i}^*)_{i=-r_{0},...,r_{0}}$ la base de $Z_{0}^*$ duale de $(z_{0,i})_{i=-r_{0},...,r_{0}}$. On d\'efinit $\tilde{\zeta}_{0}\in Isom(Z_{0},Z_{0}^*)$ par $\tilde{\zeta}(z_{0,i})=(-1)^{i+1}2\nu z_{0,i}^*$. Remarquons que le signe n'est pas le m\^eme qu'en 3.1: on a remplac\'e $\nu$ par $-\nu$. En utilisant ces donn\'ees, on d\'efinit un plongement $\iota_{0}:\tilde{G}_{0}\to \tilde{H}$. En appliquant les constructions pr\'ec\'edentes au couple $(\tilde{G}_{0},\tilde{H})$, on d\'efinit le nombre $\Xi_{\nu}(\Gamma_{0},\Theta)$ (le changement de signe de $\nu$ compense l'inversion entre les r\^oles des deux groupes tordus).

  Consid\'erons la fonction $\Gamma^{\tilde{L}}$ sur $\tilde{L}(F)$ d\'efinie par
 $$\Gamma^{\tilde{L}}(\tilde{x})=\Gamma_{0}(\tilde{x}_{0})\prod_{j=1,...,t}\Gamma_{j}({^t\tilde{x}}_{-j}^{-1}\tilde{x}_{j}).$$
 On v\'erifie que c'est un quasi-caract\`ere. On peut l'induire en un quasi-caract\`ere $\Gamma=Ind_{L}^G(\Gamma^{\tilde{L}})$, cf. 1.12, puis on d\'efinit comme ci-dessus $\Xi_{\nu}(\Gamma,\Theta)$.
 
 \ass{Lemme}{Sous ces hypoth\`eses, on a l'\'egalit\'e
 $$\Xi_{\nu}(\Gamma,\Theta)=\Xi_{\nu}(\Gamma_{0},\Theta)\prod_{j=1,...,u}\Xi(\Gamma_{j}).$$}
 
 Preuve. Supposons d'abord $d_{0}>m$.   En plus de l'hypoth\`ese d\'ej\`a faite sur $V_{0}$, on peut supposer que, pour tout $j=1,...,u$, il y a un sous-ensemble $I_{j}\subset \{r_{0}+1,...,r\}$ tel que  $V_{j}$, resp. $V_{-j}$, soit engendr\'e par les vecteurs $z_{i}$ pour $i\in D_{j}$, resp. les $z_{-i}$ pour $i\in D_{j}$.    Alors  $\tilde{H}\subset \tilde{L}$. Consid\'erons les formules qui d\'efinissent $\Xi_{\nu}(\Gamma,\Theta)$ et $\Xi_{\nu}(\Gamma_{0},\Theta)$. On voit que les ensembles ${\cal T}$ qui y apparaissent sont les m\^emes. Fixons $\tilde{T}\in{\cal T}$. Les fonctions $c_{\Theta}$ sont aussi les m\^emes, ainsi que les fonctions $D^{\tilde{H}}$. La fonction $\Delta_{r}$ de la premi\`ere formule est chang\'ee en $\Delta_{r_{0}}$. En introduisant la d\'ecomposition $W=W'\oplus W''$ attach\'ee \`a $\tilde{T}$, on a
 $$\Delta_{r}(\tilde{t})=\vert 2\vert _{F}^{a}\vert det((1-t)_{\vert W'})\vert _{F}^{r-r_{0}}\Delta_{r_{0}}(\tilde{t})$$
 pour tout $\tilde{t}\in \tilde{T}(F)$, o\`u $t={^t\tilde{t}}^{-1}\tilde{t}$ et
 $$a=r^2+r-r_{0}^2-r_{0}+(r-r_{0})dim(W'').$$
 Pour d\'emontrer l'\'egalit\'e cherch\'ee, il suffit de prouver que, pour un \'el\'ement g\'en\'eral $\tilde{t}\in \tilde{T}(F)$, on a l'\'egalit\'e
 $$(1) \qquad c_{\Gamma}(\tilde{t})\vert 2\vert _{F}^{a}\vert det((1-t)_{\vert W'})\vert _{F}^{r-r_{0}}=c_{\Gamma_{0}}(\tilde{t})\prod_{j=1,...,t}\Xi(\Gamma_{j}).$$
    A $\tilde{T}$ est associ\'ee une d\'ecomposition $V=W'\oplus V''$ et $V''$ est muni d'une forme quadratique $\tilde{\zeta}_{G,T}$. De m\^eme, on a $V_{0}=W'\oplus V''_{0}$, o\`u $V''_{0}=V''\cap V_{0}$ et $V''_{0}$ est muni de la forme quadratique $\tilde{\zeta}_{G_{0},T}$, \'egale \`a la restriction de $\tilde{\zeta}_{G,T}$. Soit $\tilde{t}\in \tilde{T}(F)$ en position g\'en\'erale. Alors $G_{\tilde{t}}=T_{\theta}\times G''_{\tilde{t}}$, o\`u $G''_{\tilde{t}}=SO(V'')$, et $G_{0,\tilde{t}}=T_{\theta}\times G''_{0,\tilde{t}}$, o\`u $G''_{0,\tilde{t}}=SO(V''_{0})$. On a $c_{\Theta}(\tilde{t})=c_{\Theta,{\cal O}}(\tilde{t})$, o\`u ${\cal O}$ est une certaine orbite nilpotente r\'eguli\`ere de $\mathfrak{g}_{\tilde{t}}''(F)$.  Ce coefficient est calcul\'e par le lemme 1.12. Dans ce qui suit, on utilise les notations de ce lemme.Montrons que 
 
 (2)  l'ensemble ${\cal X}^{\tilde{L}}(\tilde{t})$ n'a qu'un \'el\'ement, que l'on peut supposer \^etre $\tilde{t}$ lui-m\^eme.
 
   Soit $g\in G(F)$ tel que $g\tilde{t}g^{-1}\in \tilde{L}(F)$. Alors $g^{-1}A_{\tilde{L}}g$ est inclus dans $G_{\tilde{t}}$. C'est un tore d\'eploy\'e. Par d\'efinition de ${\cal T}$, $T_{\theta}$ ne contient pas de tore d\'eploy\'e non trivial, donc $g^{-1}A_{\tilde{L}}g\subset G''_{\tilde{t}}$. Les espaces $V_{j}$ sont d\'etermin\'es par le tore $A_{\tilde{L}}$: ce sont des espaces propres pour l'action de ce tore dans $V$. L'inclusion pr\'ec\'edente entra\^{\i}ne que $g^{-1}V_{j}\subset V''$ pour tout $j=\pm 1,...,\pm u$. On a mieux. Parce que $g\tilde{t}g^{-1}$ appartient \`a $\tilde{L}(F)$, les espaces $V_{j}$ sont isotropes pour la forme bilin\'eaire $g\tilde{t}g^{-1}$ et deux espaces $V_{j}$ et $V_{-j}$ sont en dualit\'e pour cette forme.  Il en est donc de m\^eme pour les espaces $g^{-1}V_{j}$, relativement \`a la forme $\tilde{\zeta}_{G,T}$. Mais les espaces $V_{j}$ eux-m\^emes v\'erifient les m\^emes conditions, et, bien s\^ur, $dim(V_{j})=dim(g^{-1}V_{j})$. Sauf dans le cas exceptionnel expliqu\'e ci-dessous, deux telles familles de sous-espaces isotropes, dont les dimensions se correspondent, se d\'eduisent l'une de l'autre par l'action d'un \'el\'ement de $G''_{\tilde{t}}(F)$. Le cas exceptionnel est celui o\`u $dim(V'') $ est paire et $V''$ est \'egal \`a la somme des sous-espaces isotropes. Il y a dans ce cas deux classes de conjugaison de telles familles. Mais ce cas ne se produit pas car l'in\'egalit\'e $d_{0}>m$ entra\^{\i}ne que $V''_{0}\not=\{0\}$ et $V''$ ne peut pas \^etre \'egal \`a la somme des $V_{j}$ pour $j=\pm 1,...,\pm u$. Donc, quitte \`a multiplier $g$ \`a droite par un \'el\'ement de $G''_{\tilde{t}}(F)$, on peut supposer $gV_{j}=V_{j}$ pour tout $j=\pm 1,...,\pm u$. Les \'egalit\'es $g\tilde{t}g^{-1}(V_{j})=V_{-j}^*=\tilde{t}(V_{j})$ entra\^{\i}nent ${^tg}V_{j}^*=V_{j}^*$. Par dualit\'e, cela implique que $gV_{0}$ est annul\'e par $V_{j}^*$ pour tout $j=\pm 1,...,\pm u$, donc $gV_{0}=V_{0}$. Mais alors $g\in L(F)$ et $g\tilde{t}g^{-1}$ appartient \`a la classe de conjugaison par $L(F)$ de $\tilde{t}$. Cela prouve (2).
   
   L'ensemble $\Gamma_{\tilde{t}}$ du lemme 1.12 est celui des $g\in G(F)$ tels que $g\tilde{t}g^{-1}=\tilde{t}$, autrement dit $\Gamma_{\tilde{t}}=Z_{G}(\tilde{t})(F) =T(F)^{\theta}\times O(V'')(F)=T_{\theta}(F)\times O(V'')(F)$. Remarquons que, dans le lemme 1.12, l'\'el\'ement $g$ n'intervient que par l'interm\'ediaire de l'orbite $g{\cal O}$. Or les orbites nilpotentes r\'eguli\`eres dans l'alg\`ebre de Lie d'un groupe sp\'ecial orthogonal sont invariantes par l'action du groupe orthogonal tout entier.  On peut donc remplacer la somme en $g$ par la valeur en $g=1$ multipli\'ee par le nombre d'\'el\'ements de $\Gamma_{\tilde{t}}/G_{\tilde{t}}(F)$, c'est-\`a-dire par $[O(V'')(F):SO(V'')(F)] $.  On calcule 
   $$(3) \qquad Z_{L}(\tilde{t})(F) \simeq T(F)^{\theta}\times O(V''_{0})\times \prod_{j=1,...,u}GL_{d_{j}}(F)=\simeq T_{\theta}(F)\times O(V''_{0})\times \prod_{j=1,...,u}GL_{d_{j}}(F).$$
   Les deux premiers termes sont des sous-groupes de $GL_{d_{0}}(F)$.  Pour $j=1,...,u$, un \'el\'ement $g_{j}\in GL_{d_{j}}(F)$ agit de fa\c{c}on naturelle dans $V_{j}$ et par $g_{j}^{\sharp}$ dans $V_{-j}$, o\`u $g_{j}^{\sharp}$ est tel que le couple $(g_{j},g_{j}^{\sharp})$ appartienne au groupe sp\'ecial orthogonal de $V_{j}\oplus V_{-j}$.On en d\'eduit l'\'egalit\'e
   $$ [O(V'')(F):SO(V'')(F)] =[Z_{L}(\tilde{t})(F):L_{\tilde{t}}(F)] ,$$
   et la formule du lemme 1.12 devient simplement
   $$(4) \qquad c_{\Gamma,{\cal O}}(\tilde{t})=D^{\tilde{G}}(\tilde{t})^{-1/2}D^{\tilde{L}}(\tilde{t})^{1/2}c_{\Gamma^{\tilde{L}},{\cal O}^L}(\tilde{t}),$$
   o\`u ${\cal O}^L$ est l'unique orbite nilpotente r\'eguli\`ere de $\mathfrak{l}_{\tilde{t}}(F)$ telle que $[{\cal O}:{\cal O}^L]=1$.  Le calcul du produit des deux premiers facteurs est similaire au calcul du terme $d(\tilde{x})$ que l'on a fait en 3.4. On le laisse au lecteur. Le r\'esultat est
   $$(5) \qquad D^{\tilde{G}}(\tilde{t})^{-1/2}D^{\tilde{L}}(\tilde{t})^{1/2}=\vert 2\vert _{F}^{-a+b}\vert det((1-t)_{\vert W'})\vert _{F}^{r_{0}-r},$$
   o\`u $b=\sum_{j=1,...,u}d_{j}(d_{j}-1)/2$. Il reste \`a calculer $c_{\Gamma^{\tilde{L}},{\cal O}^L}(\tilde{t})$. Rappelons que, pour ${\cal O}'\in Nil(\mathfrak{l}_{\tilde{t}})$, le coefficient $c_{\Gamma^{\tilde{L}},{\cal O}'}(\tilde{t})$ est d\'etermin\'e par le d\'eveloppement en combinaison lin\'eaire de transform\'ees de Fourier d'int\'egrales orbitales nilpotentes  de la fonction $X\mapsto \Gamma^{\tilde{L}}(\tilde{t}exp(X))$ au voisinage de $0$ dans $\mathfrak{l}_{\tilde{t}}(F)$. En vertu de (2), $X$ se d\'ecompose en $X_{0}+\sum_{j=1,...,u}X_{j}$, avec $X_{0}\in \mathfrak{g}_{0,\tilde{t}}(F)$ et, pour tout $j=1,...,u$, $X_{j}\in \mathfrak{gl}_{d_{j}}(F)$. De m\^eme, une orbite ${\cal O}'$ se d\'ecompose en ${\cal O}_{0}+\sum_{j=1,...,u}{\cal O}_{j}$. Si on pose $\tilde{x}=\tilde{t}exp(X)$, on a $\tilde{x}_{0}=\tilde{t}exp(X_{0})$ tandis que ${^t\tilde{x}}_{-j}\tilde{x}_{j}=exp(2X_{j})$ pour $j=1,...,u$. Donc 
   $$\Gamma^{\tilde{L}}(\tilde{t}exp(X))=\Gamma_{0}(\tilde{t}exp(X_{0}))\prod_{j=1,...,u}\Gamma_{j}(exp(2X_{j})).$$
   On en d\'eduit que, pour une orbite ${\cal O}'={\cal O}_{0}+\sum_{j=1,...,u}{\cal O}_{j}$, on a l'\'egalit\'e
   $$c_{\Gamma^{\tilde{L}},{\cal O}'}(\tilde{t})=c_{\Gamma_{0},{\cal O}_{0}}(\tilde{t})\prod_{j=1,...,u}\gamma({\cal O}_{j})c_{\Gamma_{j},{\cal O}_{j}}(1),$$
   o\`u $\gamma({\cal O}_{j})$ est le nombre complexe tel que $\hat{j}^{GL_{d_{j}}}({\cal O}_{j},2X)=\gamma({\cal O}_{j})\hat{j}^{GL_{d_{j}}}({\cal O}_{j},X)$.  On v\'erifie sur la d\'efinition de ${\cal O}$ donn\'ee en 3.2 que ${\cal O}^L$ est la somme de l'orbite ${\cal O}_{0}$ de $\mathfrak{g}''_{0,\tilde{t}}(F)$ qui intervient dans la d\'efinition de $c_{\Gamma_{0}}(\tilde{t})$ et, pour tout $j=1,...,u$, de l'unique orbite r\'eguli\`ere de $\mathfrak{gl}_{d_{j}}(F)$, notons-la ${\cal O}_{j,reg}$. Pour cette orbite, $\gamma({\cal O}_{j,reg})=\vert 2\vert _{F}^{-d_{j}(d_{j}-1)/2}$, cf. par exemple [W2] 2.6(1). Alors
   $$(6) \qquad c_{\Gamma^{\tilde{L}},{\cal O}^L}(\tilde{t})=c_{\Gamma_{0},{\cal O}_{0}}(\tilde{t})\vert 2\vert _{F}^{-b}\prod_{j=1,...,u}\Xi(\Gamma_{j}).$$
   Les formules (4), (5) et (6) entra\^{\i}nent  l'\'egalit\'e (1) cherch\'ee.
   
   Passons au cas o\`u $d_{0}<m$. Consid\'erons l'espace $Z''=Z\oplus Z_{0}$ muni de la forme quadratique $\tilde{\zeta}''=\tilde{\zeta}\oplus \tilde{\zeta}_{0}$. Il est de dimension $2r''$, o\`u $r''=r+r_{0}+1$. Le noyau anisotrope du sous-espace $Z$ est la forme $x\mapsto 2\nu x^2$ tandis que le noyau anisotrope du sous-espace $Z_{0}$ est $x\mapsto -2\nu x^2$. L'espace total $Z\oplus Z_{0}$ est donc hyperbolique. On peut en  fixer une base  $(z''_{i})_{i=\pm 1,...\pm r''}$ hyperbolique, c'est-\`a-dire telle que $\tilde{\zeta}''z''_{i}=z^{_{''}*}_{-i}$. Quitte \`a conjuguer $\tilde{L}$, on peut supposer que, pour tout $j=1,...,u$, il existe un sous-ensemble $I_{j}\subset \{1,...,r''\}$ tel que $V_{j}$ soit engendr\'e par les $z''_{i}$ pour $i\in I_{j}$ tandis que $V_{-j}$ est engendr\'e par les $z''_{-i}$ pour $i\in I_{j}$. Alors l'image du compos\'e des deux plongements $\tilde{G}_{0}\subset \tilde{H}\subset \tilde{G}$ est incluse dans $\tilde{L}$.

 Les ensembles de tores tordus qui interviennent dans les d\'efinitions de $\Xi_{\nu}(\Gamma,\Theta)$ et $\Xi_{\nu}(\Gamma_{0},\Theta)$ ne sont plus les m\^emes, notons-les ${\cal T}$ et ${\cal T}_{0}$.  Le premier est un ensemble de repr\'esentants des classes de conjugaison par $H(F)$ dans un ensemble $\underline{\cal T}$. Le second est un ensemble de repr\'esentants des classes de conjugaison par $G_{0}(F)$ dans un ensemble $\underline{{\cal T}}_{0}$. On a plong\'e $\tilde{G}_{0}$ dans $\tilde{H}$. Via ce plongement, on a
   
   (7) $\underline{\cal T}_{0}\subset \underline{\cal T}$.
   
   Soit $\tilde{T}\in \underline{\cal T}_{0}$. A $\tilde{T}$ sont associ\'ees des d\'ecompositions $V_{0}=V'\oplus V''_{0}$, $W=V'\oplus W''$, $V=V'\oplus V''$.  Les espaces $V''_{0}$, $W''$ et $V''$ sont munis de formes quadratiques. Pour montrer que $\tilde{T}$ appartient \`a $\underline{\cal T}$, la seule condition non \'evidente est de v\'erifier que les groupes sp\'eciaux orthogonaux de $W''$ et $V''$ sont quasi-d\'eploy\'es. Celui de $W''$ l'est par d\'efinition de $\underline{\cal T}_{0}$. Celui de $V''_{0}$ aussi. Mais $V''=Z''\oplus V''_{0}$.  Comme on l'a d\'ej\`a dit, la forme quadratique sur $Z''$ est hyperbolique.  Donc  le groupe sp\'ecial orthogonal de $V''_{0}$ est quasi-d\'eploy\'e  si et seulement si celui de $V''$ l'est. D'o\`u (7).
   
   Montrons que
   
   (8) deux \'el\'ements de $\underline{\cal T}_{0}$ sont conjugu\'es par un \'el\'ement de $G_{0}(F)$ si et seulement s'ils le sont par un \'el\'ement de $H(F)$.
   
   Soient $\tilde{T}_{1},\tilde{T}_{2}\in \underline{\cal T}_{0}$ et $h\in H(F)$ tel que $h\tilde{T}_{1}h^{-1}=\tilde{T}_{2}$. On introduit les d\'ecompositions d'espaces relatives aux deux tores tordus, que l'on affecte d'indices $1$ et $2$. On a forc\'ement $hW''_{1}=W''_{2}$, plus pr\'ecis\'ement $h$ induit une isom\'etrie entre les espaces quadratiques $W''_{1}$ et $W''_{2}$. Mais $Z_{0}$ est contenu (comme espace quadratique) dans ces deux espaces. D'apr\`es le th\'eor\`eme de Witt, il existe un \'el\'ement $\gamma\in O(W''_{2})(F)$ tel que $\gamma h$ fixe tout \'el\'ement de $Z_{0}$. On prolonge $\gamma$ en le faisant agir par l'identit\'e sur $V'_{0,2}$ et on remplace $h$ par $\gamma h$. Alors $h$ agit trivialement sur $Z_{0}$.  Il envoie forc\'ement l'orthogonal $V''_{0,1}$ de $Z_{0}$ dans $W''_{1}$ sur son analogue $V''_{0,2}$.  L'\'egalit\'e $h\tilde{T}_{1}h^{-1}=\tilde{T}_{2}$ l'oblige aussi \`a envoyer $V'_{1}$ sur $V'_{2}$. Donc il conserve $V_{0}=V'_{1}  \oplus V''_{0,1}=V'_{2}\oplus V''_{0,2}$. Alors $h\in G_{0}(F)$, ce qui prouve (8).
   
   D'apr\`es (7) et (8), on peut supposer que ${\cal T}_{0}\subset {\cal T}$. Montrons que
   
   (9) si $\tilde{T}\in {\cal T}-{\cal T}_{0}$, la fonction $c_{\Gamma}$ est nulle en un point g\'en\'eral de $\tilde{T}(F)$. 
   
  Soit $\tilde{t}\in \tilde{T}(F)$ en position g\'en\'erale, supposons $c_{\Gamma}(\tilde{t})\not=0$.   Le terme $c_{\Gamma}(\tilde{t})$ est calcul\'e par le lemme 1.12 qui entra\^{\i}ne qu'il existe $g\in G(F)$ tel que $g\tilde{t}g^{-1}\in \tilde{L}(F)$. Fixons un tel $g$ et notons $W=W'\oplus W''$ et $V=W'\oplus V''$ les d\'ecompositions attach\'ees \`a $\tilde{T}$. On a $g^{-1}A_{\tilde{L}}g\subset G_{\tilde{t}}=T_{\theta}\times G''_{\tilde{t}}$, avec les notations de la premi\`ere partie de la preuve. Puisque $T_{\theta}$ est anisotrope, on a $g^{-1}A_{\tilde{L}}g\subset G''_{\tilde{t}}$. Comme  pr\'ec\'edemment, cela entra\^{\i}ne que, pour $j=\pm 1,...,\pm u$, les espaces $g^{-1}V_{j}$ sont inclus dans $V''$, que ce sont des sous-espaces isotropes  et que $g^{-1}V_{j}$ est en dualit\'e avec $V_{-j}$. L'espace quadratique $V''$ contient donc un sous-espace hyperbolique de dimension $2r''$. Puisqu'il contient $Z$, l'orthogonal $W''$ de $Z$ dans $V''$ contient n\'ecessairement un sous-espace isomorphe \`a $Z_{0}$ (muni de $\tilde{\zeta}_{0}$). Fixons un tel sous-espace $W''_{0}$ et notons $W''_{1}$ son orthogonal dans $W''$. Fixons un \'el\'ement $h\in H(F)$ tel que $h(W''_{0})=Z_{0}$, $h(W'\oplus W''_{1})=V_{0}$ et $h$ induise une isom\'etrie de $W''_{0}$ muni de $\tilde{\zeta}_{H,T}$ sur $Z_{0}$ muni de $\tilde{\zeta}_{0}$.  Quitte \`a remplacer $\tilde{T}$ par $h\tilde{T}h^{-1}$, on est ramen\'e au cas o\`u $W'\subset V_{0}$, $Z_{0}\subset W''$ et $\tilde{\zeta}_{H,T}$ a m\^eme restriction \`a $Z_{0}$ que $\tilde{\zeta}_{0}$. Alors $\tilde{T}\subset \tilde{G}_{0}$. Un raisonnement analogue \`a celui de la preuve de (7) montre que $\tilde{T}\in \underline{\cal T}_{0}$. En revenant \`a notre tore de d\'epart, on a montr\'e que $\tilde{T}$ \'etait conjugu\'e \`a un \'el\'ement de $\underline{\cal T}_{0}$ par un \'el\'ement de $H(F)$. 
 Cela contredit l'hypoth\`ese que $\tilde{T}\in {\cal T}-{\cal T}_{0}$. D'o\`u (9).
 
Par d\'efinition, $\Xi_{\nu}(\Gamma,\Theta)$ est une somme index\'ee par ${\cal T}$. Pour tout $\tilde{T}\in {\cal T}$ , notons $\Xi_{\nu,\tilde{T}}(\Gamma,\Theta)$ le terme correspondant de cette somme. Pour tout $\tilde{T}\in {\cal T}_{0}$, notons de m\^eme $\Xi_{\nu,\tilde{T}}(\Gamma_{0},\Theta)$ la contribution de $\tilde{T}$ \`a $\Xi_{\nu}(\Gamma_{0},\Theta)$. Gr\^ace \`a (9), pour d\'emontrer l'\'egalit\'e de l'\'enonc\'e, il suffit de prouver  l'\'egalit\'e

(10) $\Xi_{\nu,\tilde{T}}(\Gamma,\Theta)=\Xi_{\nu,\tilde{T}}(\Gamma_{0},\Theta)\prod_{j=1,...,u}\Xi(\Gamma_{j})$

\noindent
pour tout $\tilde{T}\in {\cal T}_{0}$. Fixons $\tilde{T}\in {\cal T}_{0}$ et introduisons les d\'ecompositions habituelles $V_{0}=V'\oplus V''_{0}$, $W=V'\oplus W''$. Nous allons comparer les diff\'erents termes qui interviennent dans les d\'efinitions de $\Xi_{\nu,\tilde{T}}(\Gamma,\Theta)$ et $\Xi_{\nu,\tilde{T}}(\Gamma_{0},\Theta)$. Posons
$$C=\left\lbrace\begin{array}{cc}1 ,&\text{ si }V'\not=V_{0},\\ 2 ,&\text{ si }V'=V_{0}.\\ \end{array}\right.$$
Montrons que
$$(11) \qquad \vert W(H,\tilde{T})\vert = C\vert W(G_{0},\tilde{T})\vert .$$ 
 Introduisons le groupe lin\'eaire $G'$ de l'espace $V'$ et l'espace tordu $\tilde{G}'=Isom(V',V^{_{'}*})$. On le plonge dans $\tilde{G_{0}}$ en prolongeant tout \'el\'ement de $\tilde{G}'$ par l'isomorphisme $\tilde{\zeta}_{G_{0},T}:V''_{0}\to V^{_{''}*}_{0}$. On  a $\tilde{T}\subset \tilde{G}'$ et les \'egalit\'es
 $$Norm_{H}(\tilde{T})=Norm_{G'}(\tilde{T})\times O(W''),\,\,Norm_{G_{0}}(\tilde{T})=Norm_{G'}(\tilde{T})\times O(V''_{0}).$$
 D'o\`u
 $$\vert W(H,\tilde{T})\vert =\vert Norm_{G'}(\tilde{T})(F)/T(F)\vert \vert O(W'')(F)/SO(W'')(F)\vert ,$$
  $$\vert W(G_{0},\tilde{T})\vert =\vert Norm_{G'}(\tilde{T})(F)/T(F)\vert \vert O(V''_{0})(F)/SO(V_{0}'')(F)\vert .$$
  L'espace $W''$ n'est pas nul, donc $\vert O(W'')(F)/SO(W'')(F)\vert =2$. Si $V''_{0}\not=\{0\}$, on a aussi $\vert O(V''_{0})(F)/SO(V''_{0})(F)\vert =2$. Si $V''_{0}=\{0\}$, $\vert O(V''_{0})(F)/SO(V''_{0})(F)\vert =1$. D'o\`u (11).
  
   Consid\'erons un \'el\'ement $\tilde{t}\in \tilde{T}(F)$, en position g\'en\'erale. Les fonctions $c_{\Theta}(\tilde{t})$ sont  les m\^emes dans les deux formules. Quand on passe d'une formule \`a l'autre, $\nu$ est chang\'e en $-\nu$, mais le r\^ole des groupes est lui-aussi chang\'e: $H$ est le plus petit groupe pour l'une et le plus grand pour l'autre. Quand on se rappelle les d\'efinitions, on voit que ces deux changements se compensent. Il intervient des fonctions $\Delta_{r}(\tilde{t})$ et $\Delta_{r_{0}}(\tilde{t})$. On a l'\'egalit\'e
  $$\Delta_{r}(\tilde{t})=\vert 2\vert _{F}^{a'}\vert det((1-t)_{\vert V'})\vert _{F}^{r-r_{0}}\Delta_{r_{0}}(\tilde{t}),$$
  o\`u $a'=r^2+r-r_{0}^2-r_{0}+rdim(W'')-r_{0}dim(V''_{0})$. Un calcul similaire \`a celui du terme $d(\tilde{x})$ en 3.4 conduit aux \'egalit\'es
  $$D^{\tilde{H}}(\tilde{t})=D^{\tilde{G}'}(\tilde{t})\vert 2\vert _{F}^{dim(W'')(dim(W'')+1)/2}\vert det((1-t)_{\vert V'})\vert _{F}^{dim(W'')},$$
  $$D^{\tilde{G}_{0}}(\tilde{t})=D^{\tilde{G}'}(\tilde{t})\vert 2\vert _{F}^{dim(V_{0}'')(dim(V_{0}'')+1)/2}\vert det((1-t)_{\vert V'})\vert _{F}^{dim(V_{0}'')}.$$
  Puisque $dim(W'')=dim(V''_{0})+2r_{0}+1$, on obtient
  $$D^{\tilde{H}}(\tilde{t})\Delta_{r}(\tilde{t})=\vert 2\vert _{F}^{b'}\vert det((1-t)_{\vert V'})\vert _{F}^{r+r_{0}+1}D^{\tilde{G}_{0}}(\tilde{t})\Delta_{r_{0}}(\tilde{t}),$$
  o\`u $b'= (r+r_{0}+1)^2+(r+r_{0}+1)dim(V''_{0})$.  Compte tenu de cette \'egalit\'e et de (11), il suffit pour prouver (10) d'\'etablir l'\'egalit\'e
  
  (12) $c_{\Gamma}(\tilde{t})=C \vert 2\vert _{F}^{-b'}\vert det((1-t)_{\vert V'})\vert _{F}^{-(r+r_{0}+1)}c_{\Gamma_{0}}(\tilde{t})$.
  
  On utilise le lemme 1.12. Montrons que 
  
  (13) l'ensemble ${\cal X}^{\tilde{L}}(\tilde{t})$ est r\'eduit \`a un \'el\'ement. 
  
  Soit $g\in G(F)$ tel que $g\tilde{t}g^{-1}\in  \tilde{L}(F)$. Posons $\tilde{t}_{\natural}=g\tilde{t}g^{-1}$. On peut certainement conjuguer $\tilde{t}_{\natural}$ par un \'el\'ement de $L(F)$ de sorte que, pour $j=1,....,u$, l'application $\tilde{t}_{\natural,j}:V_{j}\to V_{-j}^*$ soit \'egale \`a $\tilde{\zeta}_{G,T,j}$. Supposons qu'il en soit ainsi, posons $\tilde{T}_{\natural}=g\tilde{T}g^{-1}$, soit $V=V'_{\natural}\oplus V''_{\natural}$ la d\'ecomposition associ\'ee \`a $\tilde{T}_{\natural}$. Par un raisonnement d\'ej\`a fait plusieurs fois, l'inclusion $A_{\tilde{L}}\subset G_{\tilde{t}_{\natural}}$ entra\^{\i}ne que, pour $j=\pm 1,...\pm u$, les espaces  $V_{j}$ sont inclus dans $V''_{\natural}$, qu'ils sont isotropes et que $V_{j}$ est en dualit\'e avec $V_{-j}$ pour la forme $\tilde{\zeta}_{G,T_{\natural}}$.  Cela entra\^{\i}ne d'abord que $V'_{\natural}\subset V_{0}$ (il doit \^etre annul\'e par $V_{j}^*$ pour tout $j=\pm 1,...,\pm u$). Cela entra\^{\i}ne aussi que, pour $j=1,...,u$, ${^t\tilde{t}}_{\natural,-j}^{-1}\tilde{t}_{\natural,j}=1={^t\tilde{\zeta}}_{G,T,-j}^{-1}\tilde{\zeta}_{G,T,j}$. D'apr\`es l'hypoth\`eses faite sur $\tilde{t}_{\natural,j}$, cela entra\^{\i}ne que $\tilde{t}_{\natural,-j}=\tilde{\zeta}_{G,T,-j}$. Donc les formes $\tilde{\zeta}_{G,T_{\natural}}$ et $\tilde{\zeta}_{G,T}$ co\"{\i}ncident sur $Z+Z_{0}=\sum_{j=\pm1,...,\pm u}V_{j}$. En utilisant le th\'eor\`eme de Witt, on peut trouver un \'el\'ement $\gamma\in O(V''_{\natural})(F)$ tel que $\gamma g$ conserve $Z+Z_{0}$ et y agisse par l'identit\'e. On prolonge $\gamma$ par l'identit\'e de $V'_{\natural}$ et on pose $g_{\natural}=\gamma g$. On a encore $g_{\natural}\tilde{t}g_{\natural}^{-1}=\tilde{t}_{\natural}$. Cela entra\^{\i}ne que $g$ envoie $V'$ dans $V'_{\natural}$ et l'orthogonal $V''_{0}$ de $Z+Z_{0}$ dans $V''$ dans l'orthogonal $V''_{\natural,0}$ de $Z+Z_{0}$ dans $V''_{\natural}$. Puisque $V'\oplus V''_{0}=V_{0}=V'_{\natural}\oplus V''_{\natural,0}$, $g_{\natural}$ conserve $V_{0}$. On sait aussi que $g_{\natural}$ agit par l'identit\'e sur $Z\oplus Z_{0}=\sum_{j=\pm 1,...,\pm u}V_{j}$. Mais alors $g_{\natural}\in L(F)$ et $\tilde{t}_{\natural}$ est conjugu\'e \`a $\tilde{t}$ par un \'el\'ement de $L(F)$. D'o\`u (13). 
  
  On peut donc supposer ${\cal X}^{\tilde{L}}(\tilde{t})=\{\tilde{t}\}$. On a $\Gamma_{\tilde{t}}=Z_{G}(\tilde{t})(F)=T(F)^{\theta}\times O(V'')(F)$. Comme dans la premi\`ere partie de la preuve, on peut remplacer la somme en $g$ du lemme 1.12 par sa valeur en $g=1$, multipli\'ee par $\vert \Gamma_{\tilde{t}}/G_{\tilde{t}}(F)\vert $, c'est-\`a-dire par
  $$\vert T(F)^{\theta}/T_{\theta}(F)\vert \vert O(V'')(F)/SO(V'')(F)\vert=\vert O(V'')(F)/SO(V'')(F)\vert .$$
  Certainement $V''$ est non nul, donc  ce terme vaut $2$. Le groupe $Z_{L}(\tilde{t})(F)$ est d\'ecrit par (3). Donc
 $$[Z_{L}(\tilde{t})(F):L_{\tilde{t}}(F)]^{-1}=\vert T(F)^{\theta}/T_{\theta}(F)\vert^{-1} \vert O(V_{0}'')(F)/SO(V_{0}'')(F)\vert^{-1}=\vert O(V_{0}'')(F)/SO(V_{0}'')(F)\vert^{-1} .$$
 Ce terme vaut $2$ si $V_{0}''\not=\{0\}$, $1$ sinon. On en d\'eduit que le produit des deux facteurs pr\'ec\'edents vaut $C$ et la formule du lemme 1.12 devient simplement
 $$(14) \qquad c_{\Gamma}(\tilde{t})=c_{\Gamma,{\cal O}}(\tilde{t})=CD^{\tilde{G}}(\tilde{t})^{-1/2}D^{\tilde{L}}(\tilde{t})^{1/2}c_{\Gamma^{\tilde{L}},{\cal O}^L}(\tilde{t}),$$
 avec les m\^emes notations qu'en (4). On a encore la formule (6). L'orbite ${\cal O}_{0}$ est la bonne: ici encore, quand on passe du couple $(\tilde{G},\tilde{H})$ au couple $(\tilde{G}_{0},\tilde{H})$, on change $\nu$ en $-\nu$ et on \'echange les r\^oles des deux groupes, et ces deux changements se compensent. Enfin, on calcule le produit $D^{\tilde{G}}(\tilde{t})^{-1/2}D^{\tilde{L}}(\tilde{t})^{1/2}$ en imitant le calcul de 3.4 et on obtient
 $$D^{\tilde{G}}(\tilde{t})^{-1/2}D^{\tilde{L}}(\tilde{t})^{1/2}=\vert 2\vert _{F}^{-b'+b}\vert det((1-t)_{\vert V'})\vert _{F}^{-(r+r_{0}+1)}.$$
 Alors la formule (14) entra\^{\i}ne (12), ce qui ach\`eve la preuve. $\square$
 
 \bigskip
 
 \subsection{D\'ebut de la preuve: le cas des repr\'esentations induites}
 
 Nous d\'emontrons le th\'eor\`eme 7.1 par r\'ecurrence sur $dim(G)+dim(H)$. Soient $\tilde{\pi}\in Temp(\tilde{G})$ et $\tilde{\rho}\in Temp(\tilde{H})$.  On suppose dans ce paragraphe que $\tilde{\pi}$ n'est pas elliptique. Il existe donc un sous-groupe parabolique tordu $\tilde{Q}=\tilde{L}U_{Q}$ de $\tilde{G}$, avec $\tilde{Q}\not=\tilde{G}$, et une repr\'esentation $\tilde{\sigma}\in Temp(\tilde{L})$ telle que $\tilde{\sigma}$ soit elliptique et $\tilde{\pi}=Ind_{Q}^G(\tilde{\sigma})$. On \'ecrit $\tilde{L}$ et $\sigma $ comme en 2.3.   On r\'ealise $\tilde{\sigma}$ comme dans ce paragraphe: il suffit pour cela de choisir le prolongement $\tilde{\sigma}_{0}$ de $\sigma_{0}$ \`a $\tilde{G}_{0}(F)$ de sorte  que $w(\tilde{\sigma}_{0},\psi)=w(\tilde{\pi},\psi)$. Notons $\Gamma_{0}$ le caract\`ere $\Theta_{\tilde{\sigma}_{0}}$ et, pour $j=1,...,u$, notons $\Gamma_{j}$ le caract\`ere $\Theta_{\sigma_{j}}$ multipli\'e par $\omega_{\sigma_{j}}((-1)^{d+1})$. D'apr\`es le lemme 2.3, le caract\`ere $\Theta_{\tilde{\sigma}}$ co\"{\i}ncide avec le quasi-caract\`ere $\Gamma^{\tilde{L}}$ du paragraphe pr\'ec\'edent. Le caract\`ere $\Theta_{\tilde{\pi}}$ co\"{\i}ncide avec $\Gamma=Ind_{L}^G(\Gamma^{\tilde{L}})$. En appliquant le lemme 7.2, on a
 $$ \epsilon_{geom,\nu}(\tilde{\pi},\tilde{\rho})=\Xi_{\nu}(\Gamma,\Theta_{\tilde{\rho}})=\Xi_{\nu}(\Gamma_{0},\Theta_{\tilde{\rho}})\prod_{j=1,...,u}\Xi(\Gamma_{j})=\epsilon_{geom,\nu}(\tilde{\sigma}_{0},\tilde{\rho})\prod_{j=1,...,u}\Xi(\Gamma_{j}).$$
 
 Soit $j\in \{1,...,u\}$. On a $\Xi(\Gamma_{j})=\omega_{\sigma_{j}}((-1)^{d+1})\Xi(\Theta_{\sigma_{j}})$. Le dernier terme est le coefficient de l'orbite r\'eguli\`ere dans le d\'eveloppement du caract\`ere $\Theta_{\sigma_{j}}$ au voisinage de $1$. D'apr\`es un r\'esultat de Rodier ([R], th\'eor\`eme p.161 et remarque 2, p.162), ce coefficient est $1$ si $\sigma_{j}$ admet un mod\`ele de Whittaker, $0$ sinon. Mais $\sigma_{j}$ est temp\'er\'ee donc admet un mod\`ele de Whittaker. Donc $\Xi(\Theta_{\sigma_{j}})=1$ et $\Xi(\Gamma_{j})=\omega_{\sigma_{j}}((-1)^{d+1})$. 

 On a l'\'egalit\'e 
$$\epsilon_{geom,\nu}(\tilde{\sigma}_{0},\tilde{\rho})=\epsilon_{\nu}(\tilde{\sigma}_{0},\tilde{\rho}).$$
En effet, si $d_{0}>m$, cela r\'esulte de l'hypoth\`ese de r\'ecurrence. Si $d_{0}<m$, on  utilise les \'egalit\'es $\epsilon_{geom,\nu}(\tilde{\sigma}_{0},\tilde{\rho})=\epsilon_{geom,-\nu}(\tilde{\rho},\tilde{\sigma}_{0})$ et $\epsilon_{\nu}(\tilde{\sigma}_{0},\tilde{\rho})=\epsilon_{-\nu}(\tilde{\rho},\tilde{\sigma}_{0})$, plus l'hypoth\`ese de r\'ecurrence pour le couple $(\tilde{\rho},\tilde{\sigma}_{0})$.

En rassemblant ces \'egalit\'es, on obtient
$$\epsilon_{geom,\nu}(\tilde{\pi},\tilde{\rho})= \epsilon_{\nu}(\tilde{\sigma}_{0},\tilde{\rho}) \prod_{j=1,...,u}\omega_{\sigma_{j}}((-1)^{d+1}).$$
 Compte tenu du fait que $d+1$ est de m\^eme parit\'e que $m$, la relation 2.5(2) nous dit que le membre de droite est \'egal \`a $ \epsilon_{\nu}(\tilde{\pi},\tilde{\rho})$.  L'\'egalit\'e ci-dessus devient celle que l'on cherchait \`a \'etablir.

\bigskip

\subsection{Le cas elliptique}

Il reste \`a \'etablir le th\'eor\`eme dans le cas o\`u $\tilde{\pi}$ est elliptique, ce que l'on suppose d\'esormais. La th\'eorie des pseudo-coefficients est valable pour le groupe tordu $\tilde{G}$. C'est une cons\'equence de la formule des traces locale tordue, que nous avons admise. Mais on peut aussi le v\'erifier gr\^ace aux constructions de Schneider et Stuhler, cf. [W5] corollaire 2.2. On peut donc choisir une fonction $\tilde{f}\in C_{c}^{\infty}(\tilde{G}(F))$ qui est cuspidale et telle que:

(1) pour $\tilde{\pi}'\in \Pi_{ell}(\tilde{G})$, $\tilde{\pi}'\not\simeq \tilde{\pi}^{\vee}$, $\Theta_{\tilde{\pi}'}(\tilde{f})=0$;

(2) $\Theta_{\tilde{\pi}^{\vee}}(\tilde{f})=1$.

Gr\^ace au lemme 1.13(iii), on peut supposer que $\tilde{f}$ est tr\`es cuspidale. Pour tout entier $N\geq1$, on d\'efinit $J_{N}(\Theta_{\tilde{\rho}},\tilde{f})$ comme en 3.3. On  comparant les th\'eor\`emes 3.3 et 6.1, on obtient l'\'egalit\'e
$$(3) \qquad J_{geom}(\Theta_{\tilde{\rho}},\tilde{f})=J_{spec}(\Theta_{\tilde{\rho}},\tilde{f}).$$

En utilisant la notation de 7.2 et la d\'efinition des fonctions $c_{\tilde{f}}$, on a l'\'egalit\'e
$$J_{geom}(\Theta_{\tilde{\rho}},\tilde{f})=\Xi_{\nu}(\Theta^J_{\tilde{f}},\Theta_{\tilde{\rho}}).$$
 Le lemme 1.13(ii) nous dit que
$$\Theta^J_{\tilde{f}}=\sum_{\tilde{L}\in {\cal L}^{\tilde{G}}}\vert W^L\vert \vert W^G\vert ^{-1}(-1)^{a_{\tilde{L}}}Ind_{L}^G(\Theta^L_{\phi_{\tilde{L}}(\tilde{f})}).$$
Pour tout $\tilde{L}$, la fonction $\phi_{\tilde{L}}(\tilde{f})$ appartient \`a l'espace ${\cal H}_{ac}(\tilde{L}(F))$ introduit en 1.10 et doit \^etre interpr\'et\'ee comme une somme
$$\sum_{\zeta\in {\cal A}_{\tilde{L},F}}{\bf 1}_{H_{\tilde{L}}=\zeta}\phi_{\tilde{L}}(\tilde{f}),$$
Chaque terme est une fonction cuspidale \`a support compact. D'apr\`es sa d\'efinition, la distribution $\Theta^L_{\phi_{\tilde{L}}(\tilde{f})}$ est calcul\'ee par la relation 1.10(3).  On a
$$\Theta^L_{\phi_{\tilde{L}}(\tilde{f})}=\sum_{{\cal O}\in \{\Pi_{ell}(\tilde{L})\}}c({\cal O})\Theta_{{\cal O}},$$
o\`u $\Theta_{{\cal O}}$ est le quasi-caract\`ere d\'efini par
$$\Theta_{{\cal O}}(\tilde{x})=\int_{i{\cal A}_{\tilde{L},F}^*}\Theta_{\tilde{\sigma}_{\lambda}}(\tilde{x})\Theta_{(\tilde{\sigma}_{\lambda})^{\vee}}(\phi_{\tilde{L}}(\tilde{f}){\bf 1}_{H_{\tilde{L}}=H_{\tilde{L}}(\tilde{x})})d\lambda.$$
On a not\'e ici $\tilde{\sigma}$ le point base de ${\cal O}$.  D'o\`u
$$J_{geom}(\Theta_{\tilde{\rho}},\tilde{f})=\sum_{\tilde{L}\in {\cal L}^{\tilde{G}}}\vert W^L\vert \vert W^G\vert ^{-1}(-1)^{a_{\tilde{L}}}\sum_{{\cal O}\in \{\Pi_{ell}(\tilde{L})\}}c({\cal O})I_{geom}(\tilde{L},{\cal O}),$$
o\`u
$$I_{geom}(\tilde{L},{\cal O})=\Xi_{\nu}( Ind_{L}^G( \Theta_{{\cal O}}),\Theta_{\tilde{\rho}}).$$

Consid\'erons la d\'efinition de $J_{spec}(\Theta_{\tilde{\rho}},\tilde{f})$ donn\'ee en 6.1. On peut y \'echanger les r\^oles des repr\'esentations et de leurs contragr\'edientes.  D'apr\`es 2.4, on peut  remplacer le produit $[i{\cal A}_{{\cal O}}^{\vee}:i{\cal A}_{\tilde{L},F}^{\vee}]^{-1}2^{-s({\cal O})-a_{\tilde{L}}}$ par $c({\cal O})$. Enfin, par d\'efinition de la fonction $\phi_{\tilde{L}}(\tilde{f})$, l'int\'egrale
$$\int_{i{\cal A}_{\tilde{L},F}^*}J_{\tilde{L}}^{\tilde{G}}((\tilde{\sigma}_{\lambda})^{\vee},\tilde{f})d\lambda$$
est \'egale \`a $mes(i{\cal A}_{\tilde{L},F}^*)\Theta_{\tilde{\sigma}^{\vee}}(\tilde{f}_{\tilde{L}})$,
o\`u $\tilde{f}_{\tilde{L}}=\phi_{\tilde{L}}(\tilde{f}){\bf 1}_{H_{\tilde{L}}=0}$. On obtient
$$J_{spec}(\Theta_{\tilde{\rho}},\tilde{f})=\sum_{\tilde{L}\in {\cal L}^{\tilde{G}}}\vert W^L\vert \vert W^G\vert ^{-1}(-1)^{a_{\tilde{L}}}\sum_{{\cal O}\in \{\Pi_{ell}(\tilde{L})\}}c({\cal O})I_{spec}(\tilde{L},{\cal O}),$$
o\`u
$$I_{spec}(\tilde{L},{\cal O})= \epsilon_{\nu}(\tilde{\sigma},\tilde{\rho})mes(i{\cal A}_{\tilde{L},F}^*)\Theta_{\tilde{\sigma}^{\vee}}(\tilde{f}_{\tilde{L}}).$$

Consid\'erons d'abord le cas $\tilde{L}=\tilde{G}$. Alors $\phi_{\tilde{G}}(\tilde{f})=\tilde{f}_{\tilde{G}}=\tilde{f}$ et les orbites n'ont qu'un \'el\'ement. Les propri\'et\'es (1) et (2) entra\^{\i}nent que $I_{geom}(\tilde{G},{\cal O})=0=I_{spec}(\tilde{G},{\cal O})$ si ${\cal O}\not=\{\tilde{\pi}\}$. Pour ${\cal O}=\{\tilde{\pi}\}$, on a
$$I_{geom}(\tilde{G},{\cal O})=\epsilon_{geom,\nu}(\tilde{\pi},\tilde{\rho})\text{ et }I_{spec}(\tilde{G},{\cal O})= \epsilon_{\nu}(\tilde{\pi},\tilde{\rho}).$$
Alors l'\'egalit\'e (3) se r\'ecrit
$$\epsilon_{geom,\nu}(\tilde{\pi},\tilde{\rho})- \epsilon_{\nu}(\tilde{\pi},\tilde{\rho})=\sum_{\tilde{L}\in {\cal L}^{\tilde{G}}, \tilde{L}\not=\tilde{G}}\vert W^L\vert \vert W^G\vert ^{-1}(-1)^{a_{\tilde{L}}}$$
$$\sum_{{\cal O}\in \{\Pi_{ell}(\tilde{L})\}}c({\cal O})(I_{spec}(\tilde{L},{\cal O})-I_{geom}(\tilde{L},{\cal O})).$$
Pour d\'emontrer le th\'eor\`eme 7.1, il suffit de fixer $\tilde{L}\not=\tilde{G}$ et ${\cal O}\in \{\Pi_{ell}(\tilde{L})\}$ et de prouver l'\'egalit\'e
$$(5) \qquad I_{geom}(\tilde{L},{\cal O})=I_{spec}(\tilde{L},{\cal O}).$$
Fixons donc de tels $\tilde{L}$ et ${\cal O}$. On peut d\'ecomposer 
$$\Theta_{{\cal O}}=\sum_{\zeta\in {\cal A}_{\tilde{L},F}}\Theta_{{\cal O},\zeta},$$
o\`u, pour $\tilde{x}\in \tilde{L}(F)$,
$$\Theta_{{\cal O},\zeta}(\tilde{x})={\bf 1}_{H_{\tilde{L}}=\zeta}(\tilde{x})\Theta_{\tilde{\sigma}}(\tilde{x})\int_{i{\cal A}_{\tilde{L},F}^*}exp(\lambda(\zeta))\Theta_{(\tilde{\sigma}_{\lambda})^{\vee}}(\phi_{\tilde{L}}(\tilde{f}{\bf 1}_{H_{\tilde{L}}=\zeta}))d\lambda$$
$$={\bf 1}_{H_{\tilde{L}}=\zeta}(\tilde{x})\Theta_{\tilde{\sigma}}(\tilde{x})mes(i{\cal A}_{\tilde{L},F}^*)\Theta_{\tilde{\sigma}^{\vee}}(\phi_{\tilde{L}}(\tilde{f}){\bf 1}_{H_{\tilde{L}}=\zeta}).$$
Si l'on regroupe les $\zeta$ selon leur classe de conjugaison par $W^G$, la somme des quasi-caract\`eres induits $Ind_{L}^G(\Theta_{{\cal O},\zeta})$  devient une somme de quasi-caract\`eres \`a supports disjoints. On en d\'eduit ais\'ement que
$$I_{geom}(\tilde{L},{\cal O})=\sum_{\zeta\in {\cal A}_{\tilde{L},F}}\Xi_{\nu}(Ind_{L}^G(\Theta_{{\cal O},\zeta}),\Theta_{\tilde{\rho}}).$$
On \'ecrit $\tilde{L}$ et $\tilde{\sigma}$ comme dans le paragraphe pr\'ec\'edent. Le m\^eme raisonnement que dans ce paragraphe , c'est-\`a-dire essentiellement le lemme 7.2, nous permet de calculer l'expression ci-dessus. On a  $\Xi_{\nu}(Ind_{L}^G(\Theta_{{\cal O},\zeta}),\Theta_{\tilde{\rho}})=0$ si $\zeta\not=0$: les termes $\Xi(\Gamma_{j})$ de 7.2 sont nuls. Si $\zeta=0$,  $\phi_{\tilde{L}}(\tilde{f}){\bf 1}_{H_{\tilde{L}}=0}=\tilde{f}_{\tilde{L}}$ et on obtient
$$I_{geom}(\tilde{L},{\cal O})=\Xi_{\nu}(Ind_{L}^G(\Theta_{{\cal O},0}),\Theta_{\tilde{\rho}})$$
$$=mes(i{\cal A}_{\tilde{L},F}^*)\epsilon_{geom,\nu}(\tilde{\sigma}_{0},\tilde{\rho})\big(\prod_{j=1,...,u}\omega_{\sigma_{j}}((-1)^{d+1})\big)\Theta_{\tilde{\sigma}^{\vee}}(\tilde{f}_{\tilde{L}}).$$
Toujours comme dans le paragraphe pr\'ec\'edent, l'hypoth\`ese de r\'ecurrence nous dit que ceci est \'egal \`a
$$mes(i{\cal A}_{\tilde{L},F}^*) \epsilon_{\nu}(\tilde{\sigma},\tilde{\rho})\Theta_{\tilde{\sigma}^{\vee}}(\tilde{f}_{\tilde{L}}),$$
c'est-\`a-dire \`a $I_{spec}(\tilde{L},{\cal O})$. Cela prouve (5) et le th\'eor\`eme. $\square$

\bigskip

{\bf Bibliographie}

\bigskip

[AGRS] A. Aizenbud, D. Gourevitch, S. Rallis, G. Schiffmann:{ \it Multiplicity one theorems}, pr\'epublication 2007

[A1] J.Arthur: {\it A local trace formula}, Publ. Math. IHES 73 (1991), p.5-96 

[A2] ----------: {\it The trace formula in invariant form}, Annals of Math. 114 (1981), p.1-74

[A3] ----------: {\it Intertwining operators and residues I. Weighted characters}, J. Funct. Analysis 84 (1989), p.19-84

[A4] -----------: {\it On elliptic tempered characters}, Acta Math. 171 (1993), p.73-138

[C] L. Clozel: {\it Characters of non connected reductive $p$-adic groups}, Can. J. of Math. 39 (1987), p.149-167

[JPSS] H. Jacquet, I. Piatetskii-Shapiro, J. Shalika: {\it Rankin-Selberg convolutions}, American J. of Math.105 (1983), p.367-464

[R] F. Rodier: {\it Mod\`ele de Whittaker et caract\`eres de repr\'esentations}, in Non commutative harmonic analysis, J. Carmona, J. Dixmier, M. Vergne \'ed., Springer LN 466 (1981), p.151-171

[W1] J.-L. Waldspurger: {\it A propos du lemme fondamental pond\'er\'e tordu}, Math. Ann. 343 (2009), p.103-174

[W2] -----------------------: {\it Une formule int\'egrale reli\'ee \`a la conjecture locale de Gross-Prasad}, pr\'epublication (2009)

[W3] -----------------------: {\it Une formule int\'egrale reli\'ee \`a la conjecture locale de Gross-Prasad, $2^{\grave{e}me}$ partie: extension aux repr\'esentations temp\'er\'ees}, pr\'epublication (2009)

[W4] ------------------------: {\it La formule de Plancherel pour les groupes $p$-adiques d'apr\`es Harish-Chandra}, J. of the Inst. Math. Jussieu 2 (2003), p.235-333

[W5] ------------------------: {\it Le groupe $GL_{N}$ tordu, sur un corps $p$-adique, $1^{\grave{e}re}$ partie}, Duke Math. J. 137 (2007), p.185-234

 [W6] -----------------------: {\it La formule des traces locale tordue}, pr\'epublication 2012

\bigskip

CNRS-Institut de Math\'ematiques de Jussieu

2 place Jussieu

75005 Paris

e-mail: waldspur@math.jussieu.fr

\end{document}